\documentclass[preprint,3p,times]{elsarticle}
\usepackage{example}

\begin{document}

\begin{frontmatter}

\title{Annealed adaptive importance sampling method in PINNs for solving high dimensional partial differential equations}
\author[1]{Zhengqi Zhang}
\ead{zhengqizhang@zhejianglab.com}
\author[1]{Jing Li\corref{cor1}}
\ead{lijing@zhejianglab.com}
\author[1]{Bin Liu\corref{cor1}}
\ead{bins@ieee.org}
\cortext[cor1]{Corresponding author}
\affiliation[1]{organization={Zhejiang Lab}, city={Hangzhou}, postcode={311121},country={China}}
\begin{abstract}
    {Physics-informed neural networks (PINNs) have emerged as powerful tools for solving a wide range of partial differential equations (PDEs). However, despite their user-friendly interface and broad applicability, PINNs encounter challenges in accurately resolving PDEs, especially when dealing with singular cases that may lead to unsatisfactory local minima. To address these challenges and improve solution accuracy, we propose an innovative approach called Annealed Adaptive Importance Sampling (AAIS) for computing the discretized PDE residuals of the cost functions, inspired by the Expectation Maximization algorithm used in finite mixtures to mimic target density. Our objective is to approximate discretized PDE residuals by strategically sampling additional points in regions with elevated residuals, thus enhancing the effectiveness and accuracy of PINNs. Implemented together with a straightforward resampling strategy within PINNs, our AAIS algorithm demonstrates significant improvements in efficiency across a range of tested PDEs, even with limited training datasets. Moreover, our proposed AAIS-PINN method shows promising capabilities in solving high-dimensional singular PDEs. The adaptive sampling framework introduced here can be integrated into various PINN frameworks.}
  \end{abstract}

\end{frontmatter}

\section{INTRODUCTION}
\label{sec:intro}
Partial differential equations (PDEs) are formulated to model physical phenomena. Over the past few decades, traditional numerical solvers for PDEs have undergone significant advancements. With robust mathematical guarantees and efficient implementation, classical solvers have flourished, offering high accuracy and numerical stability. Nonetheless, these solvers often come with computational expenses and encounter significant challenges when applied to complex systems. In recent years, with the rapid development of computing resources and machine learning algorithms, physics-informed neural networks(PINNs) introduced in \cite{RaissiParisGE:2019:JCP} have garnered significant attention for their utility in a wide array of realistic simulations, such as fluid dynamics \cite{CaiMaoWangYinGE:2021:fluid,RaissiYazdaniGE:2020:Science,NSFnets:2021}, optics \cite{Chen2020:optics}, finance \cite{Bai2022:Finance1}, see more details in  reviews \cite{Cuomo2022:PINN:review,RaissiYazdaniGE:2020:Science} and reference therein. Applying an automatic differentiation mechanism, PINNs could encode PDEs into loss functions with soft boundary and initial conditions. Compared to classical numerical solvers, PINNs are preferred with easy coding algorithms for both forward and inverse problems, meshless structures, and their ability to circumvent the curse of dimensionality. Therefore, PINNs are widely applied in various PDE-related problems, including uncertainty quantification \cite{Yang2019:UQ}, stochastic differential equations \cite{Guo2022:SDE}, fractional differential equations \cite{Pang2019:fPINN,Guo2022:MCfPINN} and more. However, despite wide success of PINNs, it would fail in solving complicated PDEs satisfactorily reflecting on "failure modes" \cite{WangYuParis:2022:NTK,WangParis:2021:SISC,WangParis:2023:ICML,Krishnapriyan2021:NIPs:FailMode}. To illustrate specifically, the non-convex optimization inherent in PINNs often leads the neural network to converge to trivial solutions, which represent local solutions to the PDE.
{For simple PDEs (single-scale, single-mode), conventional PINNs \cite{RaissiParisGE:2019:JCP} can quickly obtain satisfactory solutions. However, for multiscale PDEs, conventional PINNs often perform unsatisfactorily as the low-frequency local solution diverges from the exact solution \cite{WangYuParis:2022:NTK}. Consequently, efficient implementations of the PINNs method have emerged in recent years. For example, loss re-weighting methods \cite{McClenny2023:JCP,WangParis:2021:SISC,Wang2022:causality} and adaptive sampling strategies \cite{GaoYanZhou2023:FIPINN:I,GaoTangYanZhou2023:FIPINN:II,LuLu:2023:RAD,Lu2021:deepxde,Tang2023:DAS,Jiao2023:GAS}  have similar ideas in finding trade-offs between loss and probability distribution on weight or sampling. Additionally, some other studies aim to enhance the representation capacity of PINNs by employing new neural network architectures \cite{Jagtap2020:JCP:adaptiveactf,Liu2020:CICP,HuJinZhou:2023}. For a comprehensive review and additional references, please refer to \cite{WangParis:2023:arxiv:PINN}.}

Adaptive sampling strategies have become a cornerstone in the training of PINNs, offering a robust approach to enhance learning efficiency. The underlying concept is to strategically sample additional points around regions with higher loss values, thereby directing the training process to prioritize these areas, which may lead to a better result. In \cite{Subramanian2023:Movepoint}, authors demonstrate that resampling of collocation points is necessary and provide a resampling scheme based on gradient of loss. In \cite{Lu2021:deepxde}, a residual-based adaptive refinement method(RAR) is provided based on the large residual of loss function. In \cite{GaoYanZhou2023:FIPINN:I,GaoTangYanZhou2023:FIPINN:II}, authors propose a novel sampling method called failure-informed PINNs(FIPINN) to model the failure function via loss. Also, DAS-PINN\cite{Tang2023:DAS} uses KR-net {based on Knothe-Rosenblatt rearrangement} to approximate proposed distribution while training. These methods have shown to achieve better results than standard PINNs. However, existing strategies encounter challenges with multi-modal loss landscapes. Techniques such as domain decomposition and Gaussian Mixture methods have been suggested to address this, as cited in \cite{GaoTangYanZhou2023:FIPINN:II} and \cite{Jiao2023:GAS}. The Residual-based Adaptive Distribution method (RAD), discussed in \cite{LuLu:2023:RAD}, demonstrates satisfactory performance in low-dimensional PDEs using a simple Monte-Carlo approach. Yet, the scalability of these methods to high-dimensional cases remains an open question, with the Monte-Carlo method potentially failing in such scenarios, and there is a paucity of research addressing this issue.

{Recently, our research has been captivated by the Annealed Adaptive Importance Sampling (AAIS) methodology as articulated in \cite{Liu2014}. This approach has demonstrated its efficacy in the context of sampling from complex, multi-modal target distributions, with the resultant samples subsequently utilized to approximate marginal likelihoods. Recognizing the substantial benefits of AAIS in addressing high-dimensional, multi-modal distributions, we have identified a strategic opportunity to integrate this algorithm into the framework of adaptive sampling methods, particularly pertinent to the high-dimensional PDEs.} 

{Therefore, in this paper, we propose a novel adaptive sampling method called AAIS-PINNs to augment the efficiency of PINNs model training. The AAIS-PINNs framework is based  on the Expectation Maximization(EM) algorithm \cite{Liu2014,Cappe2008,McLachlan2008:EM} for multi-modal distributions, and we aim to employ finite mixtures to emulate the distributions derived from PDE residuals. To rigorously evaluate the performance of the EM algorithm, we have adopted the Effective Sample Size (ESS) metric. This choice enables a precise quantification of the approximation quality between the target density and the proposed distributions, please refer to \cite{Elvira2022:Ess:rethink,Martino2017:SP:Ess,Kong1994:Ess} for more information about ESS. }

{The AAIS-PINNs necessitates fewer nodes for exploration compared to RAD, offering enhanced efficiency, particularly in high-dimensional problems}. In our experiments in Section \ref{sec:Experiments}, we observe that the RAD method struggles with nine-dimensional Poisson problems, since uniformly sampling in high dimensions is very inefficient. Moreover, based on the NTK theory\cite{WangYuParis:2022:NTK,Lau2024:pinnacle}, where authors propose that the PINNs would learn the low frequency parts of solutions firstly, we observe that the adaptive sampling methods can increase the frequency of the loss function or the absolute error, revealing that the low frequency parts of the solution are well-resolved. In Section \ref{sec:Experiments}, the frequency-increasing phenomenon in the residuals appear in most cases compared to PINNs using uniformly sampling method(\textit{Uni}).

 We summarize the main contributions of this paper as follows:
\begin{itemize}
 \item \textbf{Proposed Adaptive Sampling Methodology:} {We introduce AAIS-PINNs, encompassing both Gaussian mixture (\textit{AAIS-g}) and Student's t-distribution (\textit{AAIS-t}) approaches. This novel methodology represents an adaptive sampling technique designed to approximate complex target distributions within PINNs. Our approach significantly enhances the accuracy of capturing multi-modal residuals in high-dimensional PDEs, a task not previously addressed in the literature. The AAIS-PINNs have demonstrated robust performance across a spectrum of high-dimensional PDEs, showcasing their potential in various PDE-related machine learning applications. The adaptability of the AAIS algorithm also makes it readily integrable into a wide array of existing and future machine learning algorithms for PDEs.}
\item \textbf{Robust Resampling Framework for PINNs:} {We introduce a straightforward yet robust resampling framework specifically tailored for PINNs. This framework is designed to maintain a controlled size of training datasets for consistency, while also strategically incorporating adaptive points to mitigate the risk of local minima. Our approach, when applied to a variety of forward PDE problems, consistently demonstrates a frequency-increasing phenomenon. This consistent observation across different PDE problems underscores the efficiency and effectiveness of our proposed algorithms, with particular significance in high-dimensional scenarios.}
\end{itemize}

This paper is organized as follows: in Section \ref{sec:PINNs} we introduce the basic knowledge of PINNs. In Section \ref{sec:AAIS} the idea of AAIS algorithm and the definition of effective sample size are presented. In Section \ref{sec:AAIS-PINNs} we propose a simple yet efficient resampling framework of PINNs, adaptable to various sampling methods, including uniform sampling and our proposed AAIS methods. Finally, some numerical experiments are conducted in Section \ref{sec:Experiments} followed by the conclusion Section \ref{sec:conclusion}.
\section{Physics informed neural networks}
\label{sec:PINNs}
In this section, we briefly introduce PINNs based on the formula of \cite{RaissiParisGE:2019:JCP}. Given a $d$-dimensional domain $\Omega$, and the boundary $\partial\Omega$, consider the following problem:
\begin{align}
    \label{eqn:PDE}
    \begin{aligned}
        \mathcal{N}(\x;u(\x))&=0,~~\x\in\Omega,\\
        \mathcal{B}(\x;u(\x))&=0,~~\x\in\partial\Omega,
    \end{aligned}
\end{align}
where $\mathcal{N}$ is a differential operator, $\mathcal{B}$ is the boundary operator, $u(\x)$ is the unknown solution to the partial differential equation \eqref{eqn:PDE}. Denote $u(\x;\theta)$ a neural network representation of $u(\x)$ with parameter set $\theta$, which is the combination of all trainable parameters of fully connected neural networks.

The parameter set $\theta$ is obtained by the following minimization problem of loss function $\L(u(\x;\theta))$:
\begin{equation}
    \label{eqn:loss}
    \theta =\arg\min_{\theta\in\Theta}\L(u(\cdot;\theta))=\arg\min_{\theta\in\Theta}\L_{in}(u(\x;\theta))+\L_b(u(\x;\theta))
\end{equation}
where $\Theta$ is the parameter space. The most common choice of loss function is weighted $L^2$ norm on the entire domain $\Omega$ with measure $\omega(\x)$ i.e.,
\begin{align*}
    \L_{in}(\theta)=\|\N(\x;u(\x;\theta))\|^2_{L^2(\Omega, \omega)},~~    \L_b(\theta)=\|\B(\x;u(\x;\theta))\|^2_{L^2(\partial\Omega,\omega)}.
\end{align*}
However, in real world we could not compute $\L_{in}$ and $\L_b$ accurately but sampling points on the domain would lead to a discretized loss function
\begin{equation}
    \label{eqn:DiscreteLoss}
    \hat\L(u(;\theta)) = \hat\L_{in}(u(\x^{in};\theta))+\hat\L_b(u(\x^b;\theta)),
\end{equation}
where the discretized loss would be weighted MSE loss, i.e.
\begin{equation}
    \label{eqn:EachDiscreteloss}
    \begin{aligned}
    \hat\L_{in}(u(\x^{in};\theta))&=\frac{1}{N_{in}}\sum_{i=1}^{N_{in}}\omega^{in}_i|\N(\x^{in}_i;u(\x_i^{in};\theta))|^2,\\
    \hat\L_b(u(\x^b;\theta))&=\frac{1}{N_b}\sum_{i=1}^{N_b}\omega^b_i|\B(\x^b_i;u(\x_i^b;\theta))|^2,
    \end{aligned}
\end{equation}
where $\{\x^{in}_i\}_{i=1}^{N_{in}}$, $\{\x^{b}_i\}_{i=1}^{N_b}$ are sampling points in the domain and on the boundary respectively, with $\{\omega^{in}_i\}_{i=1}^{N_{in}}$,  $\{\omega^{b}_i\}_{i=1}^{N_{b}}$ the discrete weight corresponding to collocation points. PINNs would encode the PDE into loss function simply via automatic differentiation mechanism.

The aforementioned loss reweighting strategies and adaptive sampling point strategies try to minimize the statistical errors generated from collocation points, on the other hand some methods improving neural network structures would strengthen the representation capability of PINNs. To be more specifically, based on the formula in \cite{Tang2023:DAS,Jiao2023:GAS}, we let $\E$ be the expectation then the total error introduced from PINNs could be split into two parts
\begin{equation}
    \label{eqn:ErrSplit}
    \E(\|u(\cdot;\hat\theta-u(\cdot)\|_{L^2(\Omega,\omega)}))\le \E(\|u(\cdot;\hat\theta)-u(\cdot;\theta)\|_{L^2(\Omega,\omega)})+\|u(\cdot;\theta)-u(\cdot)\|_{L^2(\Omega,\omega)}
\end{equation}
where
\begin{equation*}
    \hat\theta=\arg\min_{\theta\in\Theta}\min\hat\L(u(\cdot;\theta)),~~\theta=\arg\min_{\theta\in\Theta}\min\L(u(\cdot;\theta)).
\end{equation*}
We could see that the first part of \eqref{eqn:ErrSplit} is due to the statistical error from discretizing the loss function with Monte Carlo methods and the second part is from the approximation capabilities of neural network on the parameter space $\Theta$. In this work, we only consider how to decrease the error in the PDE loss $\hat\L(u(\x^{in};\theta))$ from first part of \eqref{eqn:ErrSplit} and assume the boundary is well-approximated.
\section{Annealed adaptive importance sampling}
\label{sec:AAIS}
In this section we briefly introduce the annealed adaptive importance sampling(AAIS) method based on EM algorithm \cite{Cappe2008,Liu2014}.

Given a target density $\Q$ and the designed mixture
$$
q(\cdot|\psi):=\sum_{m=1}^M p_m(\cdot|\xi_m)=\sum_{m=1}^M\al_mf_m(\cdot|\xi_m)
$$
where $\psi\triangleq\{M,\{\al_m,\xi_m\}_{m=1}^M\}$ is the set of all tunable model parameters, $p_i$, $i=1,2,..., M$ is each component of $q$, $M$ is the number of modes of $q$, $\al_m$, $\xi_m$ is the mixture weights and parameter corresponding to each component $f_m$, the weights satisfy $\al_m>0$, $\sum_{m=1}^M\al_m=1$. We desired to approximate $\Q$ with proposed $q$, from \cite[Section 2.1]{Cappe2008}, it is equivalent to consider the following maximizing likelihood estimation
\begin{equation}
    \label{eqn:MaxEntropy}
    \psi=\arg\max\int\log\left(\sum_{m=1}^M\al_mf_m(x|\xi_m)\right)\Q(x)\d x,
\end{equation}
which implies the application of EM algorithm for finite mixtures.

Moreover, to evaluate the performance of approximation the classical idea is to calculate the log-likelihood \eqref{eqn:MaxEntropy}, which may be inefficient when calculating integration. Here we define the effective sample size(ESS) from importance sampling\cite{Kong1994:Ess}: let $\{X_i\}_{i=1}^N$ be the point set sampled from distribution $q$, then ESS can be defined by
\begin{equation}
    \label{eqn:ESS}
    \ESS(\Q;q)=1-\frac{Var[\tilde w_i]}{Var[\tilde w_i]+(\E[\tilde w_i])^2} = \frac{(\sum_{i=1}^N\tilde w_i)^2}{\sum_{i=1}^N\tilde w_i^2},
\end{equation}
where $\tilde w_i = \Q(X_i)/q(X_i,\psi)$. We could easily see that $\ESS(\Q;q)\in(0,1)$ which is normalized. And when variation of importance weights $\tilde w_i$ become smaller, ESS would be larger, implying better approximation of $q$. When $Var[\tilde w_i]$ vanishes, that means all IS weights are equal, the proposal is an ideal approximation of $\Q$ up to a constant, i.e.
$$
    q(x|\psi)=\frac{\Q(x)}{\int \Q(x)\d x}.
$$
ESS has been widely used in the last decade as a measurement of importance sampling methods for its simplicity, see more discussion of effective sample size in \cite{Kong1994:Ess,Martino2017:SP:Ess,Elvira2022:Ess:rethink}. Here we adopt the definition of \eqref{eqn:ESS} as the same in \cite{Kong1994:Ess}.

The EM algorithm for multivariate distributions is introduced in \cite{McLachlan2008:EM,Cappe2008}, where the authors consider fixed number of components and adjust the component weights and parameters. In \cite{Liu2014}, a novel approach is proposed, which involves gradually adding new components to the designed proposal under iteration, along with an annealed strategy. In the following we will introduce the idea of annealed importance sampling strategy(AAIS).

\subsection{EM operation} Firstly, we introduce the EM algorithm for fix-size($M$ fix) finite mixtures. Here we provide some details for both Gaussian mixtures and Student's t-mixtures.

At iteration $t$, $t\in\mathbb{N}$, the proposal $q^t$ at the current iteration has parameters $M$ and $\psi^t$ including $\al_m^t$ and $\xi_m^t$, $m=1,2,..., M$. Let $\{X_i^t\}_{i=1}^N$ be sampled points from $q^t$, and define the posterior probabilities describing the role each component plays in generating each sampling point,
\begin{equation*}
    \rho^t_m(X_i^t) = \frac{\al_m^tf_m(X_i^t;\xi^t_m)}{q^t(X_i^t;\psi^t)}.
\end{equation*}
Following the formula of \cite{Cappe2008}, in the next iteration $t+1$, the new parameter set $\psi^{t+1}$ of proposal $q^{t+1}$ shall be updated by EM steps
\begin{equation}
    \label{eqn:updateEM}
    \begin{aligned}
        &\text{{\bf E} step}:~~\al_m^{t+1} = \sum_{i=1}^N w_i^t \rho_m^t(X_i^t),\\
        &\text{{\bf M} step}:~~\xi_m^{t+1} = \arg\max_{\xi}\sum_{i=1}^N w_i^t\rho_m^t(X_i^t)\log(f_m(X_i^t;\xi))
    \end{aligned}
\end{equation}
where $w_i^t$ is the normalized importance weights
\begin{equation*}
    w_i^t = \frac{\tilde w_i^t}{\sum_{i=1}^N \tilde w_i^t},~~\tilde w_i^t = \frac{\Q(X_i^t)}{q^t(X_i^t)},~~i=1,2,..., N.
\end{equation*}

\paragraph{The Gaussian Component Case} In the following we assume that the proposal $q$ is a mixture with multiple Gaussians, i.e. $f_m(\cdot|\xi_m)$ is a Gaussian distribution. That means the parameter $\xi_m$ of each $f_m$ includes two parts: the mean $\mu_m$ and the covariance $\Sigma_m$. Then by the formula of \cite[Section 2.2]{McLachlan2008:EM}, \cite[section 3]{Cappe2008} the EM algorithm for updating Gaussian mixtures could be written as
\begin{equation}
    \label{EM:Gaussian}
    \begin{aligned}
        \al_m^{t+1} &= \sum_{i=1}^N w_i^t\rho_m^t(X_i^t),\\
        \mu_m^{t+1} &= \sum_{i=1}^N w_i^t\rho_m^t(X_i^t)X_i^t/\al_m^{t+1},\\
        \Sigma_m^{t+1} &= \sum_{i=1}^N w_i^t\rho_m^t(X_i^t)(X_i^t-\mu_m^{t+1})(X_i^t-\mu_m^{t+1})^T/\al_m^{t+1}.
    \end{aligned}
\end{equation}
\paragraph{The Student's t Component Case} Similar as Gaussian mixtures, we select the components to be Student's t-distributions, as they offer greater efficiency in importance sampling owing to their heavy tail property. The parameter $\xi_m$ would include mean $\mu_m$ and covariance $\Sigma_m$ with fixed degree of freedom $v$. Then the EM algorithm for the Student's t-mixture would be(\cite[section 3.4]{Liu2014}, \cite[section 4]{Cappe2008})
\begin{equation}
    \label{EM:T}
    \begin{aligned}
        \al_m^{t+1} &= \sum_{i=1}^N w_i^t\rho_m^t(X_i^t),\\
        \mu_m^{t+1} &= \sum_{i=1}^N \frac{w_i^t\rho_m^t(X_i^t)\delta(X_i^t,\mu_m^t;\Sigma_m^t)\mu_m^tX_i^t}{\sum_{i=1}^Nw_i^t\rho_m(X_i^t)\delta(X_i^t,\mu_m^t;\Sigma_m^t)},\\
        \Sigma_m^{t+1} &= \sum_{i=1}^N w_i^t\rho_m^t(X_i^t)\delta(X_i^t,\mu_m^t;\Sigma_m^t)(X_i^t-\mu_m^{t+1})(X_i^t-\mu_m^{t+1})^T/\al_m^{t+1},
    \end{aligned}
\end{equation}
where
\begin{equation*}
    \delta(X,\mu;\Sigma)= \frac{v+d}{v+(X-\mu)^T(\Sigma)^{-1}(X-\mu)},
\end{equation*}
see more details of $\delta(X,\mu;\Sigma)$ in \cite[section 2.6]{McLachlan2008:EM}.

\subsection{AAIS algorithm}
In the following we will discuss a simplified version of AAIS algorithm from \cite{Liu2014} for Gaussian and Student's t distribution, i.e. the parameter $\xi_m$ of each component of $q$ only includes the mean $\mu_m$ and covariance $\Sigma_m$. Firstly, we must build intermediate target density functions, $\Q_1,\Q_2,..., \Q_I$ during iteration, where
\begin{equation}
    \label{eqn:target}
    \Q_i(\cdot)\triangleq q(\cdot|\psi)^{1-\lambda_k}\Q^{\lambda_k}(\cdot),~~k=1,2,..., I,
\end{equation}
$0=\lambda_1\le ...\le \lambda_I=1$ is the temperature ladder, $I$ is given. If the temperature ladder is appropriate, the target density $\Q$ would be approximated smoothly from initial guess. Then we need the following operation in the AAIS algorithm.
\paragraph{Initial operation} From the beginning, we must set one initial proposal $q^0$ from target function $\Q$. Firstly, we uniformly sample $N_{A}$ points  $\{X_i\}_{i=1}^{N_A}$ in the domain then find the point $X_s$ maximizing target $\Q(x)$. Then let the mean of the single component proposal $q^0$ be $X_s$, given an initial covariance $\Sigma_0$, we execute EM algorithm with $n$ points sampled from $q^0$ until the $\ESS(\Q;q^0)$ would be larger than a given threshold $T_a$ or reach the maximum cycle limit $C_u$. The \textbf{Initial} algorithm is stated in Algorithm \ref{alg:init}.
\begin{algorithm}[htbp]
    \caption{Initial algorithm for AAIS}\label{alg:init}
     \KwIn{Target function $\Q$, number of searching points $N_S$, $N_A$, $n$, threshold  $T_a$, cycle limit $C_u$}
     Uniformly sample $N_S$ points $\{X_i\}_{i=1}^{N_S}$ in the domain\;
     Find the maximizing target function point $X_s$\;
     Construct a proposal $q^0$ with mean $X_s$ and covariance $\Sigma_0$\;
     \For{$\ell = 0$ to $C_u$}
     {
        Sample $n$ points from $q^0$\;
        Calculate importance weights $\{\tilde w_i\}_{i=1}^n$ and $\ESS(\Q;q^0)$\;
        Apply \textbf{EM} algorithm for $q^0$ once\;
        \If{$\ESS(\Q;q^0)\ge T_a$}
        {\textbf{Break}}
     }
    \KwOut{The initial guess proposal $q^0$.}
\end{algorithm}
\paragraph{Delete operation} At iteration $q^t$, if some components of it can be neglected, i.e., the corresponding weights are less than a given threshold $T_d$, we delete them from $q$ and rearrange the survival components that makes sum of their weights equal to 1.
\paragraph{Update operation} When the $\ESS(\Q;q^t)$ is not large enough, we could add new component based on importance sampling weights. Firstly, sample $N_A$ points $\{X_i^t\}_{i=1}^{N_A}$ based on $q^t$, calculate the importance weights $\{\tilde w_i^t\}_{i=1}^{N_A}$ based on target $\Q$. Find the largest weighted sample $X_s^t$, given an initial covariance $\Sigma_0$ and mean $X_s^t$, construct a new proposal $p$. Apply EM algorithm for $p$ according to target function $\Q$ until the $\ESS(\Q;p)$ is satisfactory or reach the maximum number of cycles. That means $p$ is ready to combine with $q^t$, then obtain $n$ samples from $p$, check if we need to merge $p$ with $q^t$. If no merge operation is applied, we update the proposal $q^{t+1}$ with $q^{t+1} = \sigma q^t +(1-\sigma)p$, where $\sigma$ is the update weight. The details are given in Algorithm \ref{alg:update}
\begin{algorithm}[htbp]
    \caption{Update algorithm for AAIS}\label{alg:update}
     \KwIn{The number of searching points $n$, the new proposal $p^*$, the designed proposal $q^t$, the threshold of merge $T_m$, the update weight $\sigma$.}
     Sample $n$ points from $p^*$;
     Calculate weighted ESS value $\al_m^t\ESS(p_m^t;p^*)$ for each component $p_m^t$ of $q^t$\;
     \eIf{any of $\al_m^t\ESS(p_m^t;p^*)$ is bigger than $T_m$}
     {
        Find the largest component $m^*$ maximizing $\al_{m^*}^t\ESS(p_{m^*}^t;p^*)$\;
        Obtain $q^{t+1}$ by updating $q^t$ with the $m^*$-th component changed into $p_{m^*}=\al_{m^*}^tf(\cdot|\mu_{m^*}^t+\mu^*,\Sigma_{m^*}^t+\Sigma^*)$, where $\mu^*$, $\Sigma^*$ are parameters of $p^*$.
     }
     {
        Obtain $q^{t+1} = \sigma q^t+(1-\sigma)p^*$\;
     }
    \KwOut{The updated proposal $q^{t+1}$.}
\end{algorithm}

The AAIS algorithms given by Algorithm \ref{alg:AAIS} includes all above operations and could be applied using both Gaussian mixtures(AAIS-g) and Student's t-mixtures(AAIS-t). The EM algorithms for different selection of distributions follow \eqref{EM:Gaussian} and \eqref{EM:T} respectively. And notice that if AAIS-t is used we must specify the degree of freedom $v$ as one more input.
\begin{algorithm}[htbp]
    \caption{AAIS algorithms (AAIS-g and AAIS-t)}\label{alg:AAIS}
     \KwIn{Target function $\Q$, number of searching points $N_S$, $N_A$,  $n$, threshold of delete, merge, update $T_d$, $T_m$, $T_a$, updated weight $\sigma$, cycle limit for update and initial guess $C_u$, anneal ladder $\{\lambda_k\}_{k=1}^I$, threshold ESS ladder $\{\eta_k\}_{k=1}^I$, iteration ladder $\{C_k\}_{k=1}^I$, degree of freedom $v$ if using AAIS-t.}
    Let $t=0$\;
    To obtain $q^0$, execute \textbf{Initial} operation(Algorithm \ref{alg:init}) with $\Q$, $N_S$, $N_A$, $n$, $T_a$ and $C_u$ \;
    \For{$k=1$ to $I$}
    {
        Define the target function $\Q_k=(q^t)^{1-\lambda_k}\Q^{\lambda_k}$\;
        Sample $N_A$ points $\{X_i^t\}_{i=1}^{N_A}$ from $q^t$\;
        Calculate importance weights $\{w_i^t\}_{i=1}^{N_A}$ and $\ESS(\Q_k;q^t)$\;

        \For{$j=0$ to $C_k$}
        {   \If{$\ESS(\Q_k;q^t)>=\eta_k$}
                {\textbf{Break}}
            Find the largest importance weight point $X_s^t$\;
            Set the candidate $p$ with mean $X_s^t$ and given covariance $\Sigma_0$\;
            \For{$\l=0$ to $C_u$}
                {
                    Sample $n$ points $\{X_i\}_{i=1}^n$ from $p$\;
                    Calculate importance weights $\{w_i\}_{i=1}^n$ and $\ESS(\Q_k;p)$\;
                    Execute \textbf{EM} algorithm for $p$ once\;
                    \If{$\ESS(\Q_k;p)>=T_a$}
                        {\textbf{Break}}

                }
            Execute \textbf{Update} operation(Algorithm \ref{alg:update}) with $p$, $n$, $T_m$, $\sigma$ and $q^t$, obtain $q^{t+1}$\;
            Execute \textbf{EM} algorithm for $q^{t+1}$ twice\;
            Execute \textbf{Delete} operation for $q^{t+1}$\;
            Sample $N_A$ points $\{X_i^{t+1}\}_{i=1}^{N_A}$ from $q^{t+1}$\;
            Calculate importance weights $\{w_i^{t+1}\}_{i=1}^{N_A}$ and $\ESS(\Q_k;q^{t+1})$\;
            $t\leftarrow t+1$\;
        }
    }
    Find the maximum ESS proposal $q^*:=\arg\max_{q^t}\ESS(\Q_I;q^t)$\;
    \KwOut{The proposal $q^*$.}
\end{algorithm}
\begin{remark}
    For Algorithm \ref{alg:AAIS}, a detailed statement about how to choose annealed ladder is discussed in \cite[Section 3.6]{Liu2014}. Compared to the merge algorithm in \cite[Section 3.2]{Liu2014}, we design the merging in update operation based solely on the ESS, rather than mutual information, for simplicity. For the choice of prior setting parameter, we always choose $T_d$ as 1 percent of the number of components of proposal, $T_a=0.15$, $T_m=0.85$, $N_A=\lfloor 0.1N_S \rfloor$, $n=\lfloor 0.1N_A\rfloor$ and $\sigma=0.5$. The initial covariance matrix $\Sigma_0$ is set to be a diagonal matrix with identical diagonal entries. For dimensions smaller than 2, the diagonal entries are set to $100\cdot n^{-2}$, while for dimensions equal to or greater than 2, the diagonal entries are set to $0.1$. The selection of these prior set parameters is informed by our experiences and the optimal settings for each parameter may require further investigation in future research.
\end{remark}
\section{Resampling framework of PINNs}
\label{sec:AAIS-PINNs}
In this section we propose a simple but efficient resampling PINNs framework, which could be applied with adaptive sampling strategies. The idea is to generate new sampling points during training and combine them with the training datasets controlling the size fixed. Hence, drawing upon the insights gained from RAD results \cite{LuLu:2023:RAD} and the concept of sampling additional points around regions with higher loss values, we opt for a density function proportional to the residual function $\Q(\x), \x\in\Omega$ as
\begin{equation}
    \label{eqn:PDERes}
    \Q(\x) = |\N(\x;u(\x;\theta))|^2,~~\x\in\Omega,
\end{equation}
with $\N(\cdot;u(\cdot,\theta))$ the PDE loss.


To sum up, the Algorithm for training PINNs with resampling is showed in Algorithm \ref{alg:PINNSamplingResample}. Here $\A$ could be uniformly sampling method, RAD method from \cite{LuLu:2023:RAD}, or our AAIS methods from Section \ref{sec:AAIS}. Note that at each iteration we keep the size of training datasets as the same by selecting nodes in the domain from the last training datasets and newly generated points in the domain and resampling $S^j_{b}$, which is  motivated by \cite{Subramanian2023:Movepoint} that remaining collocation points unchanged during training is suboptimal that would lead to a local behavior of PINNs.
\begin{algorithm}[htbp]
    \caption{PINNs using resample strategy}\label{alg:PINNSamplingResample}
     \KwIn{number of points of initial sampling $N_{in}$ in the domain and $N_{b}$ on the boundary, number of points of from  sampling method $N_\D$, maximum iteration $M$, sampling method $\A$.
     }
     Let $j=0$, uniformly sampling $N_{in}$ points in the domain $\Omega$, $N_{b}$ points on the boundary $\partial\Omega$, denoted with $\S^0_{in}$, $\S^0_{b}$ respectively\;
     Use $\S^0_{in}$ and $\S^0_{b}$ to pre-train PINNs\;
     \For{$t=1$ to $M$}
     {  Define the target function $\Q(\cdot)=|\N(\cdot;u(\cdot;\theta))|^2$\;
        Generate sampling set $\D$ according to one sampling method $\A$ and $\Q(\cdot)$\;
        Uniformly choose ($N_{in}-N_{D}$) points from $S^t_{in}$, combine them with $\D$ to generate $\S^{t+1}_{in}$\;
        Uniformly resample $N_b$ points on the boundary to obtain $\S^{t+1}_{b}$\;
        $t\leftarrow t+1$\;
        Use $\S^t_{in}$ and $\S^t_{b}$ to train PINNs\;

     }
    \KwOut{The neural network solution $u(\cdot;\theta)$.}
\end{algorithm}

The selection of parameters, the hyperparameter settings and the AAIS algorithm settings are problem dependent given in section \ref{sec:Experiments}.
\begin{remark}
    Algorithm \ref{alg:PINNSamplingResample} is a framework of training PINNs with resampling training datasets. The resampling idea is widely used in many works like \cite{GaoTangYanZhou2023:FIPINN:II,LuLu:2023:RAD} to fix the sample size in order to show the efficiency of adaptive sampling. We hope no complicated training tricks showed here, and the algorithm could be applied in other PINNs improvements like \cite{WangParis:2023:arxiv:PINN,Brecht2023:Hard}.
\end{remark}
\section{Experiments}
\label{sec:Experiments}
In this section in order to perform the efficiency of sampling, we present numerical results of three different sampling methods:
\begin{itemize}
    \item The uniformly sampling method\cite{RaissiParisGE:2019:JCP}. This method uses uniform sampling to update training dataset.
    \item The RAD method \cite{LuLu:2023:RAD}. This method uses a simple numerical integration to calculate $p(\x) = \Q(\x)/\int\Q(\x)\d \x$, that is, uniformly sample $N_S$ points $\{\x_i\}_{i=1}^{N_S}$ in the domain and calculate $p(\x_i) = \Q(\x_i)/(\sum_{j=1}^{N_S} \Q(\x_j)/{N_S})$ then updates training dataset from $\{\x_i\}_{i=1}^{N_S}$ based on the probability $p(\x_i)$.
    \item The AAIS algorithm presented in Section \ref{sec:AAIS}.
\end{itemize}
 In the following we always use \textit{Uni} referring to the uniformly sampling method, \textit{RAD} referring to the residual-based adaptive distribution method, \textit{AAIS-t} and \textit{AAIS-g} referring to the AAIS algorithm with Student's t-mixtures and Gaussian mixtures. According to the PINNs sampling framework in Section \ref{sec:PINNs}, we firstly pre-train the neural network with a small dataset uniformly sampled in the domain, then we will combine our training dataset $\S^j_{in}$ with newly generated dataset $\D$ with resample strategy.

For the all tested numerical methods, unless otherwise specified, the parameters keep the same in the following:  the annealed ladder $\{\lambda_k\}_{k=1}^3(I=3)$ defined in Section \ref{sec:AAIS} is always set as $[0.7, 0.9, 1.0]$, $C_u=10$, $\{C_k\}_{k=1}^3$ is set $[100, 100, 100]$, the ESS break threshold $\{\eta_k\}_{k=1}^3$ is set $[0.9, 0.88, 0.85]$, the degree of freedom $v$ of \textit{AAIS-t} is set to be 3. These settings for AAIS is a trade-off of high ESS values and computational costs according to our experiences. We always choose the collocation point weight in \eqref{eqn:EachDiscreteloss} to be 1, i.e., $\omega_i^{in}=\omega_i^{b}=1$.

The fully connected neural network with 7 hidden layers and 20 neurons in each layer is used to model the numerical solution of PDEs. We select \textit{tanh} as our activation function, Adam and lbfgs optimizers are applied to optimize loss function with learning rate 0.0001 and 0.3. To measure the accuracy of sampling strategies, we use the $L^2$ relative error and the $L^\infty$ error
$$
    e_r(u(\cdot;\theta)) = \frac{\sqrt{\sum_{i=1}^N (u(\x_i;\theta)-u^*(\x_i))^2}}{\sqrt{\sum_{i=1}^N (u^*(\x_i))^2}},~~e_\infty(u(\cdot;\theta)) = \max_{i}|u(\x_i;\theta)-u^*(\x_i)|
$$
where $\x_i$, $i=1,...,N$ belongs to the test datasets with size $N$.
\subsection{Summary of results and findings}
Firstly, we summarize our numerical results and findings for different problems.
We investigate the error decay during training, analyze the behavior of the PDE target density $\Q(\x)$ defined in \eqref{eqn:PDERes}, and explore the limits of adaptive sampling methods.

For PDEs exhibiting single and multiple high singularities, such as the Poisson problem with single and multiple peaks, we observe that adaptive sampling methods outperform the \textit{Uni} approach, underscoring the significance of adaptive residual-based sampling. Moreover, in low dimensions, \textit{RAD} outperforms the other three sampling methods, owing to its ability to closely mimic the target residual
$\Q$ with a sufficient number of search points (approximately 100k). Additionally, the \textit{AAIS-t} algorithm demonstrates superior performance compared to \textit{AAIS-g}, attributed to the heavy tail property of Student's t-distribution.

However, in high dimensions or when there are limitations on the size of searching points $N_S$, \textit{AAIS-t} may outperform the \textit{RAD} since simple Monte-Carlo integration may fail. This represents a significant advantage of our proposed AAIS-PINN algorithm, which could prove invaluable in solving high-dimensional multi-peak PDEs with constrained computing resources.

Furthermore, we conduct additional tests on various PDEs in \ref{sec:Appendix}. Our findings reveal that for PDEs such as Burgers' equation and Allen-Cahn equation (adjusted weight), adaptive sampling methods yield superior results. However, for complex PDEs like the KdV equation, we observe no significant difference between adaptive sampling methods and conventional PINNs. {This calls for further investigation to elucidate the reasons behind this phenomenon and to potentially refine the adaptive sampling strategies for broader efficacy.}

\subsection{Two-dimensional Poisson problems}
In this part, we focus on two-dimensional Poisson problems with low regularities, the solution of which has multi-peaks. This problem is always considered as a test problem for adaptive sampling efficiency of PINNs \cite{Tang2023:DAS,GaoYanZhou2023:FIPINN:I,Jiao2023:GAS}.
\begin{equation}
    \label{pde:Poisson1Peak}
    \begin{aligned}
        -\Delta u(x, y) &= f (x,y),~~ (x,y)\in\Omega,\\
        u(x,y) &= g(x,y),~~(x,y)\in\partial\Omega,
    \end{aligned}
\end{equation}
with $\Omega = (-1,1)^2$ and the source term $f$ and boundary term $g$ are defined by the exact solution
\begin{equation*}
    u^*(x,y) = \sum_{i=1}^c\exp\left[-1000((x-x_i)^2+(y-y_i)^2)\right]
\end{equation*}
the centers $\{(x_i, y_i)\}_{i=1}^c$ are case-dependent.
\subsubsection{One peak}
\begin{figure}[htbp]
    \centering
    \includegraphics[scale=0.25]{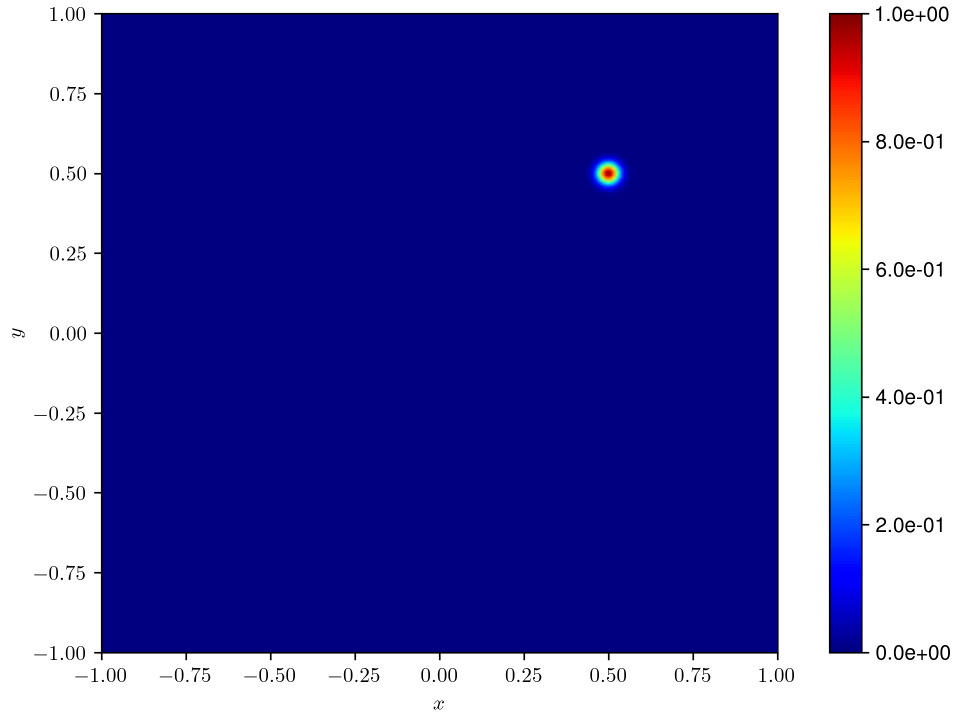}
    \caption{Exact solution for Poisson problems with one peak in \eqref{pde:Poisson1Peak}.}
    \label{fig:PS1Pexact}
\end{figure}

 We consider the one peak case $(c=1)$ and let $(x_1, y_1) = (0.5,0.5)$, the exact solution is plotted in Figure \ref{fig:PS1Pexact}. The sparsity of high singular solutions often leads PINNs to become trapped in local minima, as demonstrated in the following analyses.
 We sample total 2000 points in the domain $\Omega$ with 500 points updated in each iteration, i.e. $\#\S_{in}^t=1500$ and $\#\D=500$, and sample 500 points on the boundary $\partial\Omega$.

\begin{figure}[htbp]
    \centering
    \begin{subfigure}{.5\textwidth}
        \centering
        \includegraphics[scale=0.3]{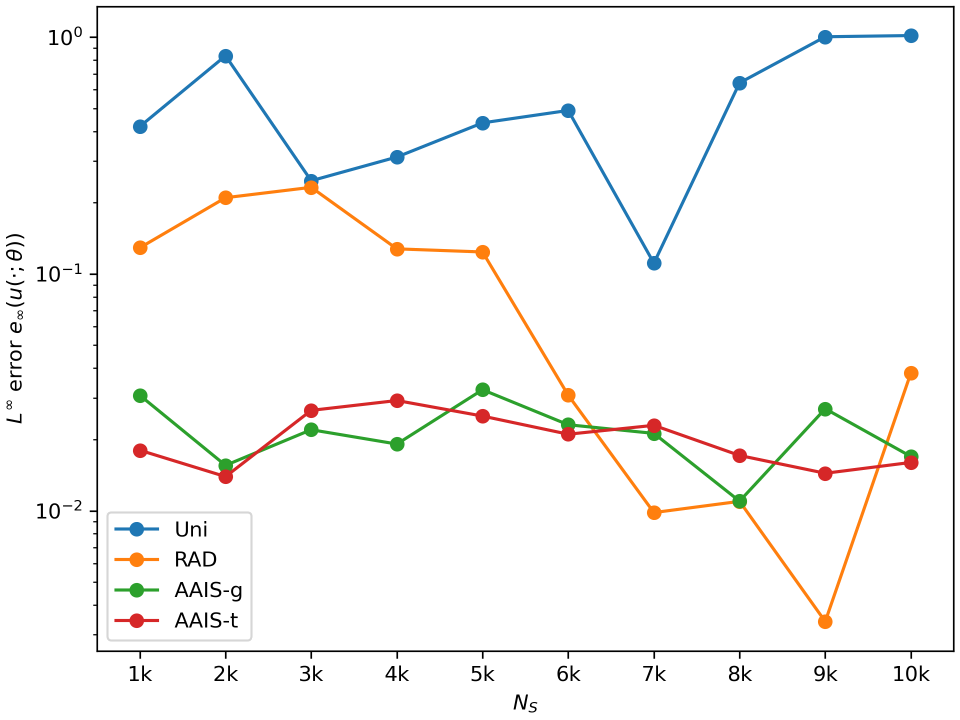}
        \end{subfigure}%
        \begin{subfigure}{.5\textwidth}
        \centering
        \includegraphics[scale=0.3]{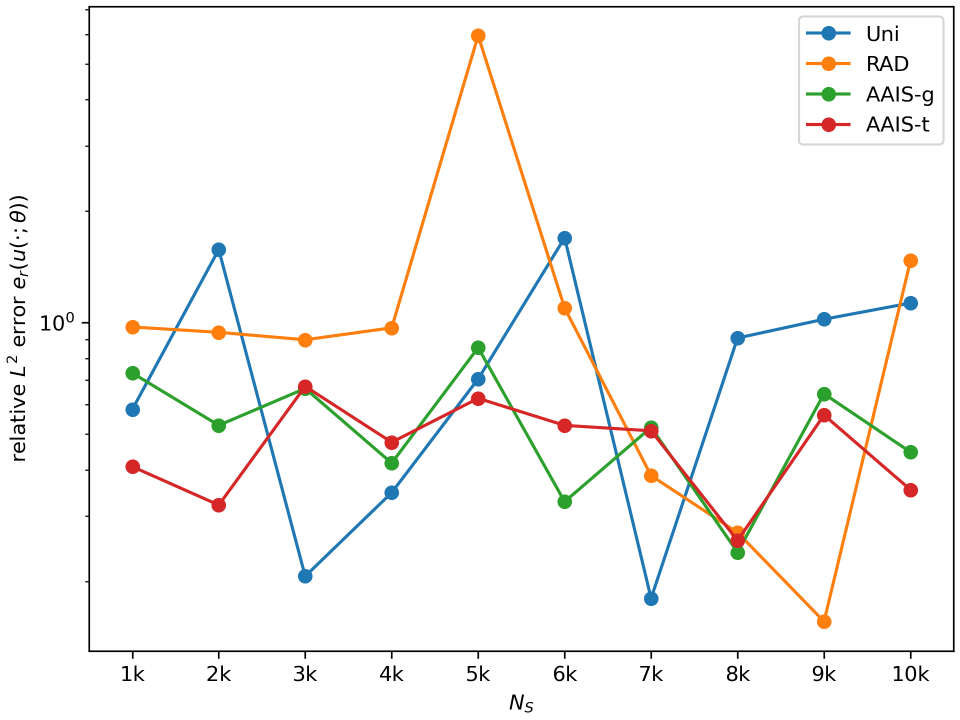}
        \end{subfigure}%
        \caption{$L^2$ relative error and $L^\infty$ error during training for one peak Poisson equation with training schedule of lbfgs 1000 epochs each versus different searching points $N_S$. Left: the $L^\infty$ error $e_\infty(u(\cdot;\theta))$. Right: the relative $L^2$ error $e_r(u(\cdot;\theta))$.
        }
    \label{fig:PS1PErr_eff}
\end{figure}

Firstly, we investigate the efficiency of these four algorithms. We set Adam 500 epochs and lbfgs 1000 epochs for pre-training and adaptive training with max iteration $M=10$. Here we let $N_S$ for \textit{RAD}, \textit{AAIS-g} and \textit{AAIS-t} to be $1000, 2000, ... , 10000$, and $N_A=N_S$ for AAIS algorithms in Algorithm \ref{alg:AAIS}.  The errors are showed in Figure \ref{fig:PS1PErr_eff}. We could see that with fewer epochs, all algorithms would behave similarly viewing the $L^2$ relative error but differently seeing the $L^\infty$ error, implying that the PINNs would arrive at local minimums, but the adaptive algorithms would jump from them. It could be more clearly seen in the following cases with more epochs. Moreover, it could be seen that the \textit{RAD} method behave similarly as \textit{Uni} method since small $N_S$ leads to a bad adaptive sampling. However, our proposed AAIS algorithm would not be influenced by the small $N_S$, showing the efficiency of importance sampling. The ineffective adaptive sampling would be the main reason why \textit{RAD} fails in high-dimensional Poisson problems with multi-peaks, which is showed in the followings.

\begin{figure}[htbp]
    \centering
    \begin{subfigure}{.33\textwidth}
        \centering
        \includegraphics[scale=0.2]{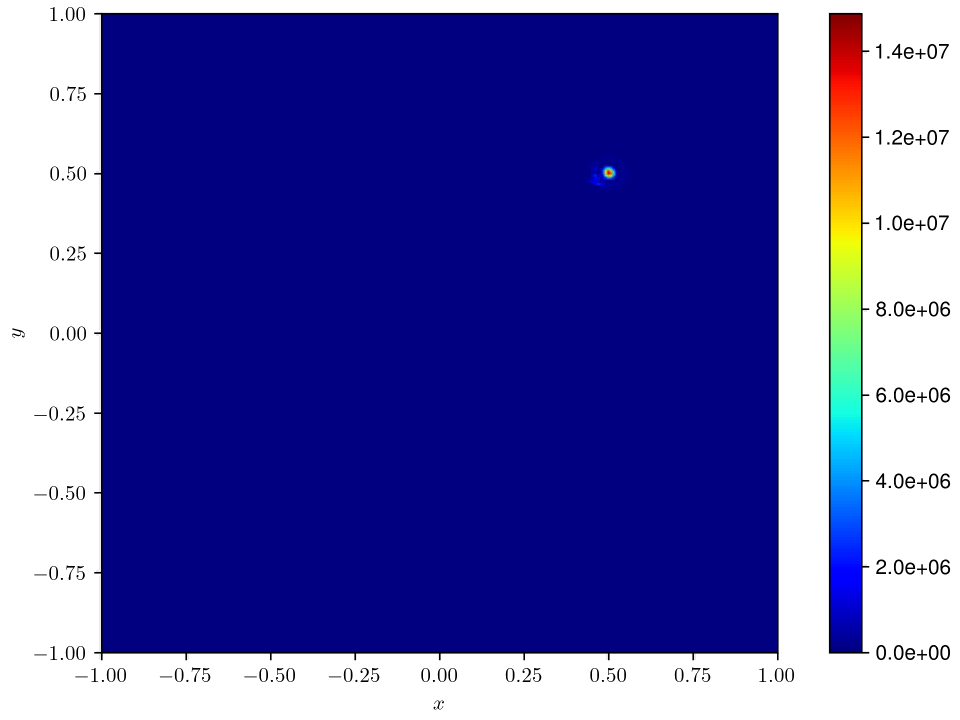}
        \caption{Residual $\Q$ after pre-training.}
    \end{subfigure}%
    \begin{subfigure}{.33\textwidth}
        \centering
        \includegraphics[scale=0.2]{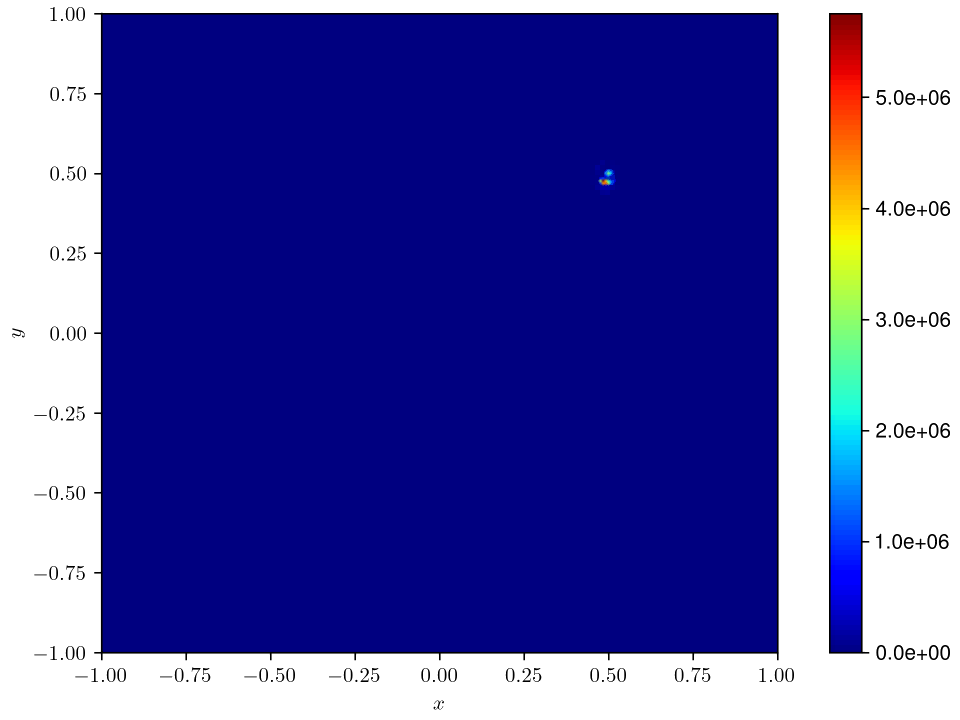}
        \caption{Residual $\Q$ after 4-th training.}
    \end{subfigure}%
    \begin{subfigure}{.33\textwidth}
        \centering
        \includegraphics[scale=0.2]{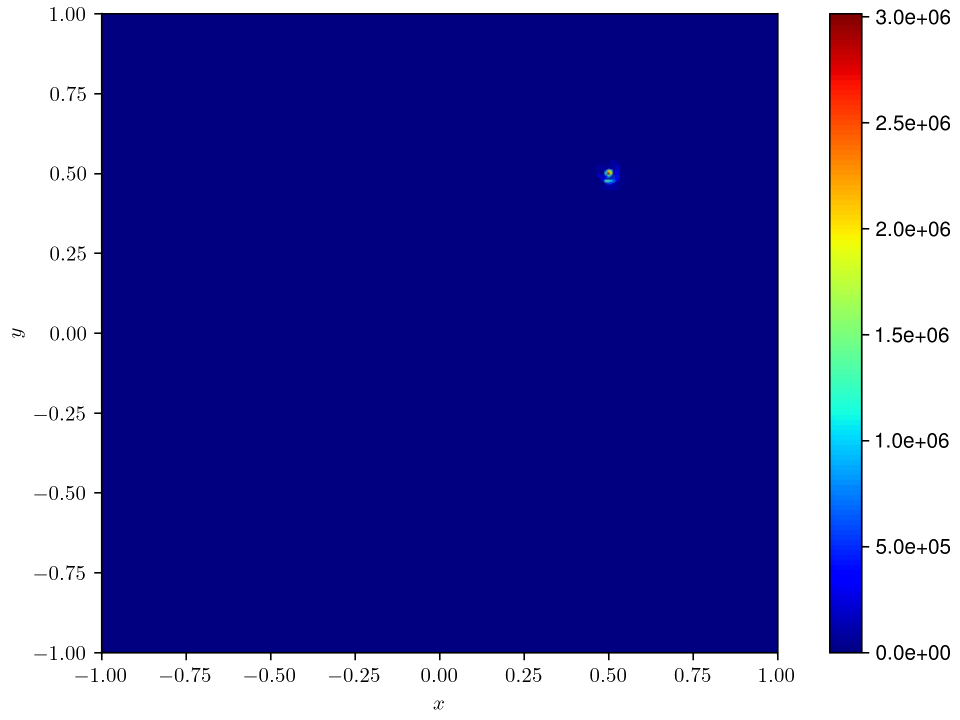}
        \caption{Residual $\Q$ after 9-th training.}
    \end{subfigure}%
    \caption{Residual $\Q$ for \textit{Uni} methods for training schedule of lbfgs 1000 epochs where $N_S=7000$.}
    \label{fig:PS1PlossUni7k}
\end{figure}

\begin{figure}[htbp]
    \centering
    \begin{subfigure}{.25\textwidth}
        \centering
        \includegraphics[height=0.75\textwidth,width=1.0\textwidth]{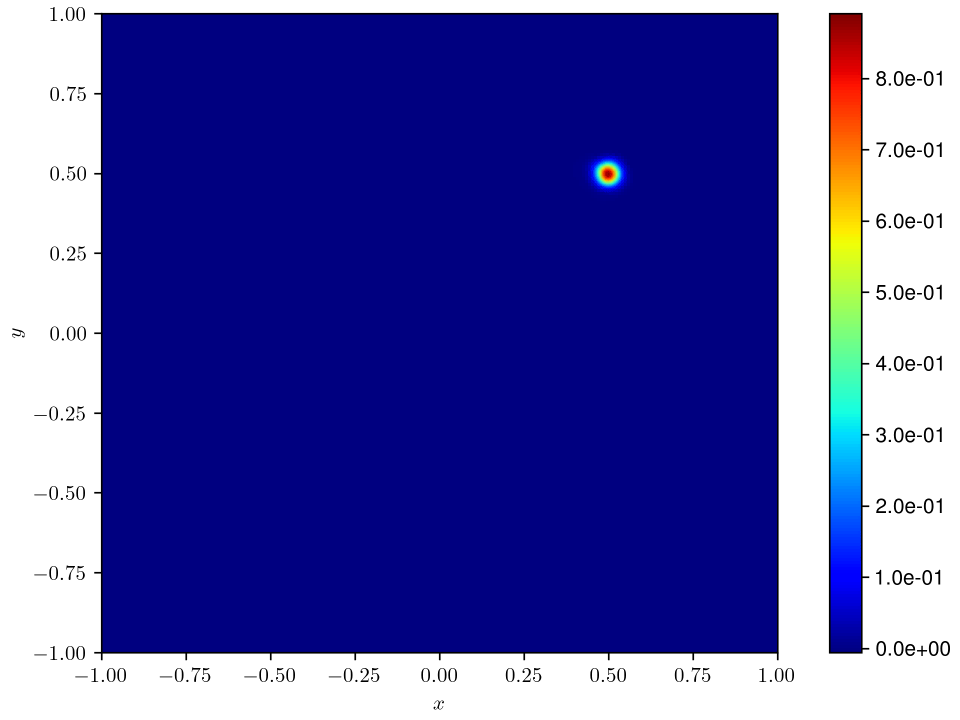}
    \end{subfigure}%
    \begin{subfigure}{.25\textwidth}
        \centering
        \includegraphics[height=0.75\textwidth,width=1.0\textwidth]{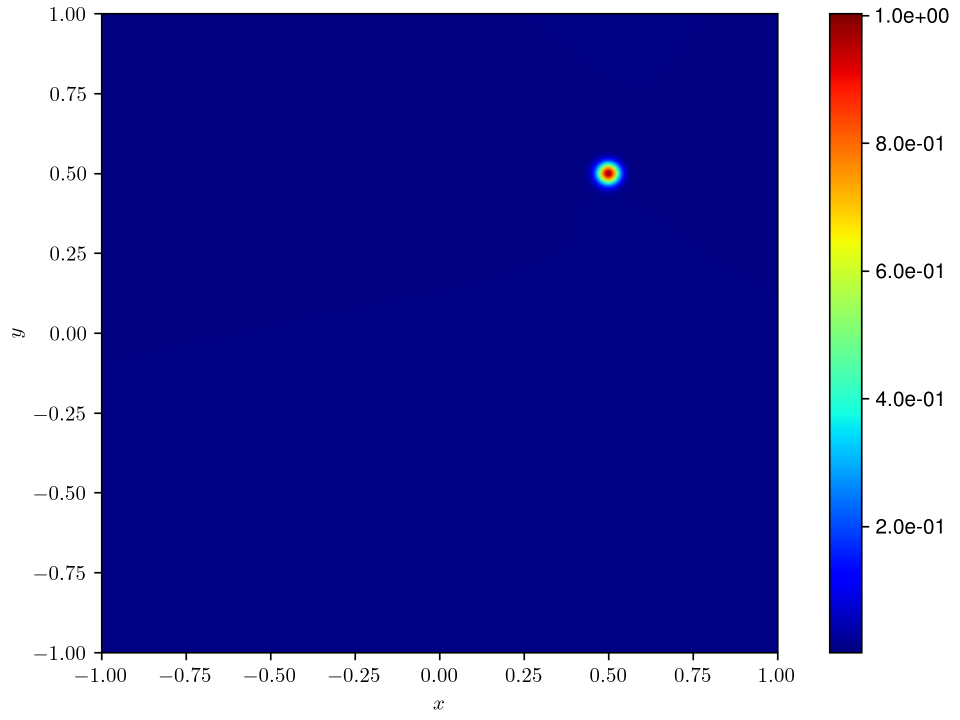}
    \end{subfigure}%
    \begin{subfigure}{.25\textwidth}
        \centering
        \includegraphics[height=0.75\textwidth,width=1.0\textwidth]{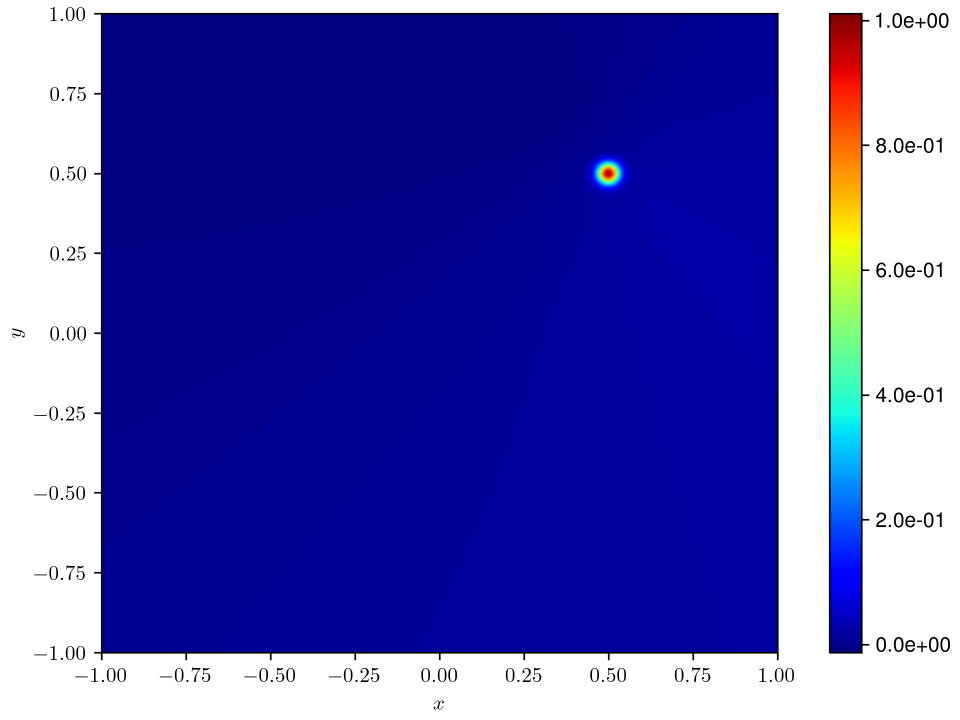}
    \end{subfigure}%
    \begin{subfigure}{.25\textwidth}
        \centering
        \includegraphics[height=0.75\textwidth,width=1.0\textwidth]{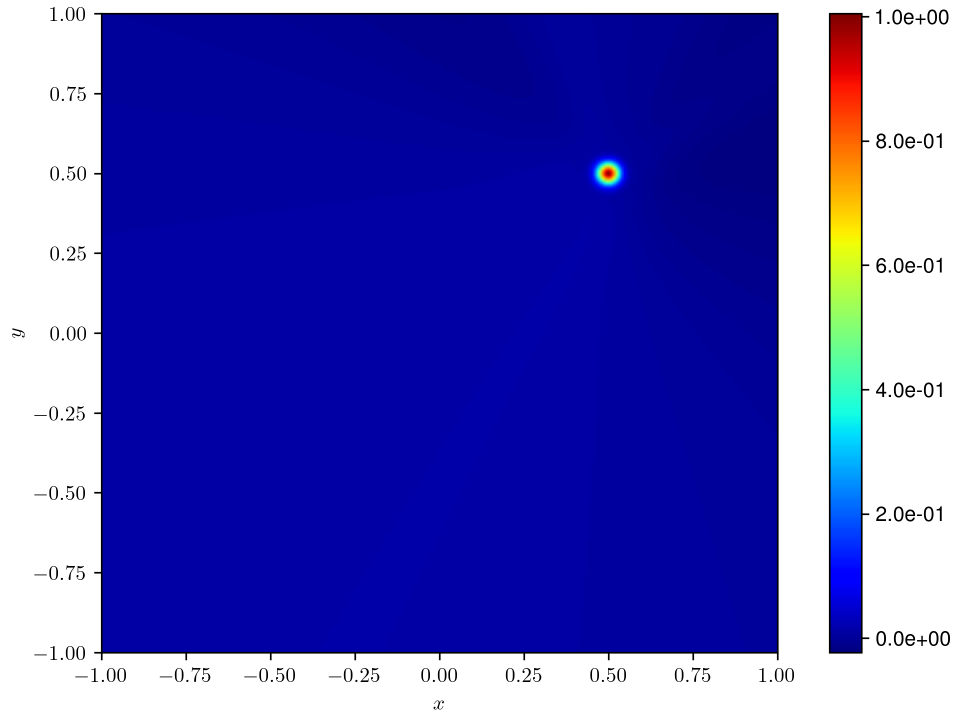}
    \end{subfigure}%
    \newline
    \raggedleft
    \begin{subfigure}{.25\textwidth}
        \centering
        \includegraphics[height=0.75\textwidth,width=1.0\textwidth]{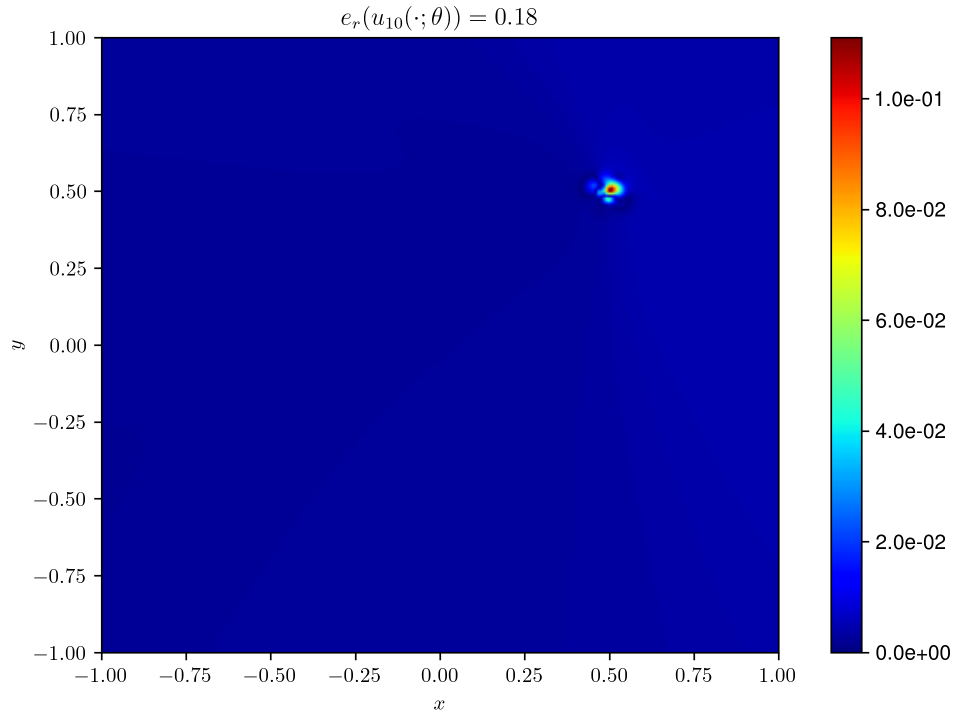}
        \caption{\textit{Uni}}
    \end{subfigure}%
    \begin{subfigure}{.25\textwidth}
        \centering
        \includegraphics[height=0.75\textwidth,width=1.0\textwidth]{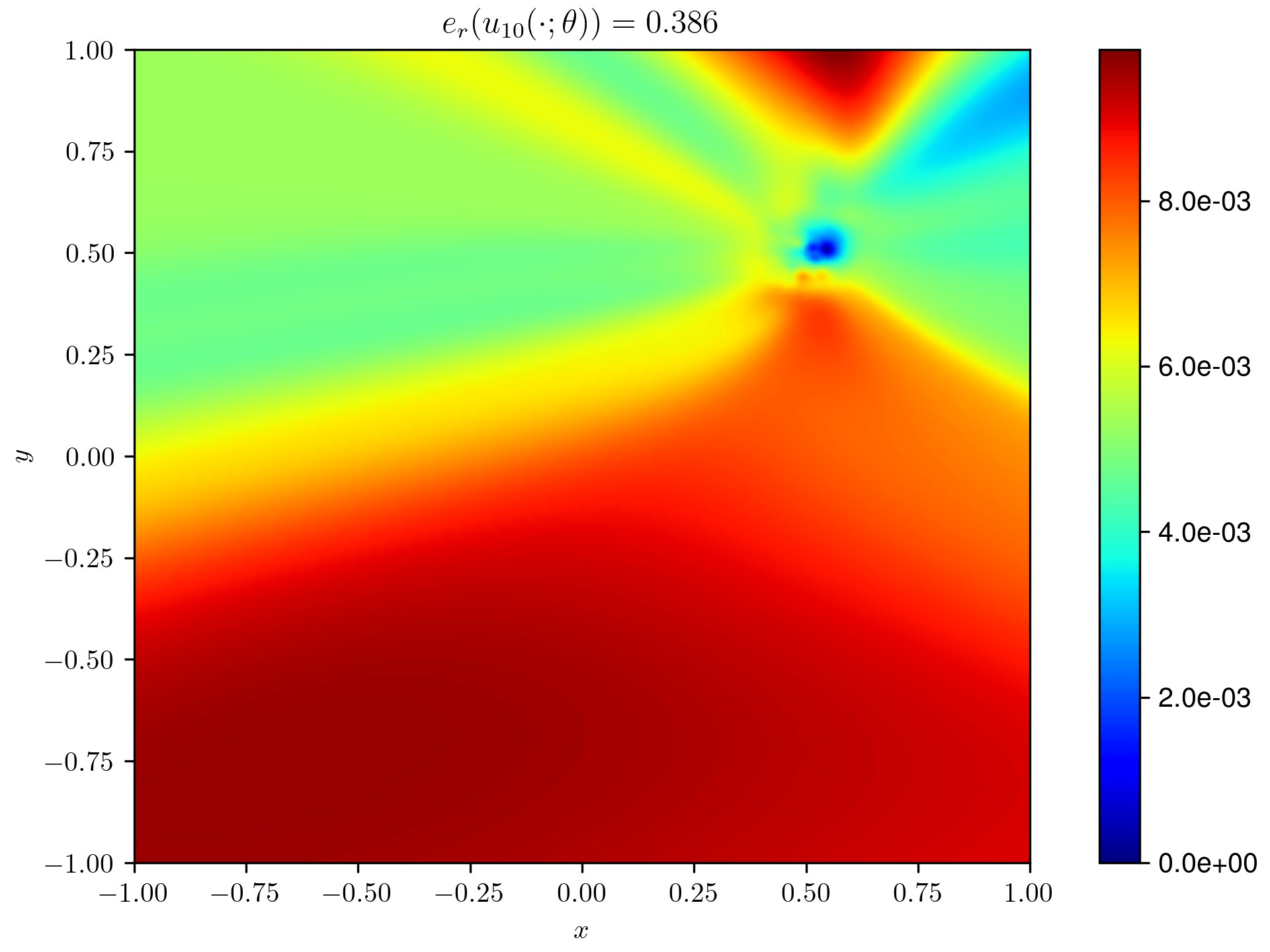}
        \caption{\textit{RAD}}
    \end{subfigure}%
    \begin{subfigure}{.25\textwidth}
        \centering
        \includegraphics[height=0.75\textwidth,width=1.0\textwidth]{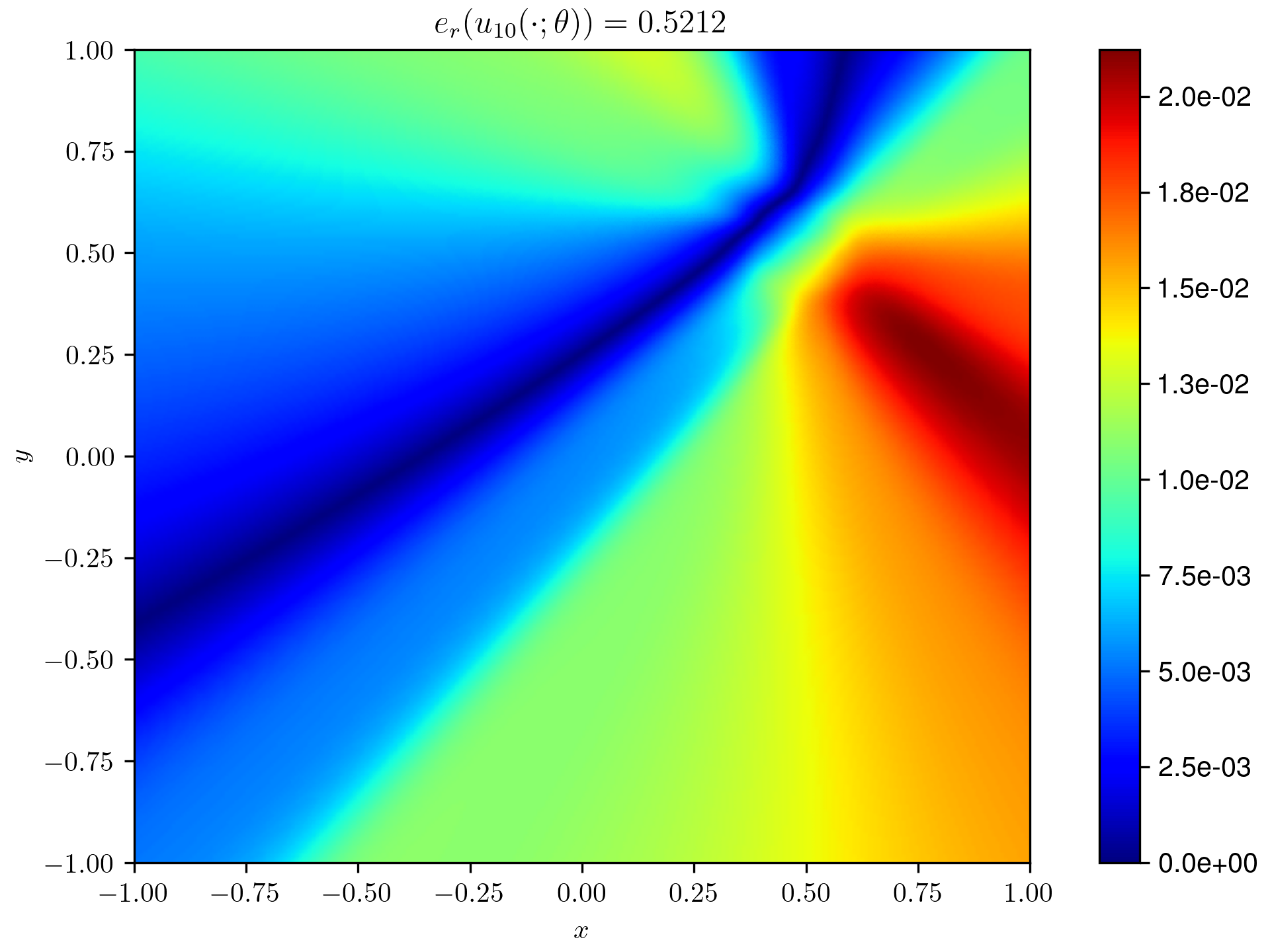}
        \caption{\textit{AAIS-g}}
    \end{subfigure}%
    \begin{subfigure}{.25\textwidth}
        \centering
        \includegraphics[height=0.75\textwidth,width=1.0\textwidth]{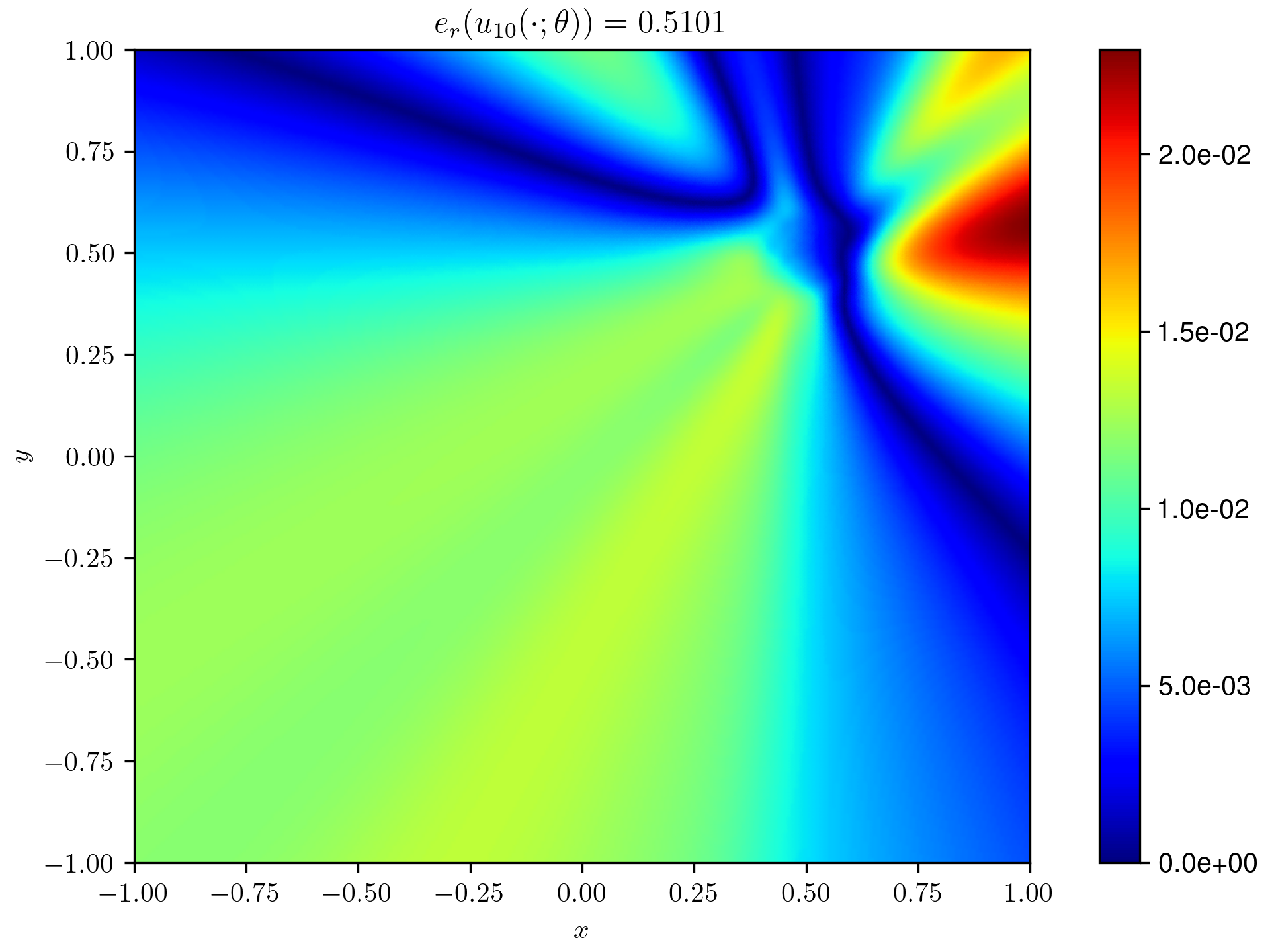}
        \caption{\textit{AAIS-t}}
    \end{subfigure}%
    \caption{Profiles of absolute error and neural network solutions for one peak Poisson equation with training schedule of 1000 epochs for lbfgs. First row: numerical solutions. Second row: absolute error. For \textit{Uni}, the solution is generated after 5th iteration. For other three sampling methods, the solutions are from the 10th iteration.}
    \label{fig:PS1PAbs7k}
\end{figure}

Consider the case $N_S=7000$, the \textit{Uni} method would be overfitted during the last few iterations due to the less training at singularity, which could be seen in the Figure \ref{fig:PS1PlossUni7k} that  the residual $\Q(\x)=|\N(\x;u(\x;\theta))|^2$ would concentrate on the singularity in the training, implying the unsatisfactory solution behavior presented in Figure \ref{fig:PS1PAbs7k}(a). Even it is the minimum $L^2$ error from the 10 cases of $N_S$, it still shows a huge mismatch at the singularity.

\begin{figure}[htbp]
    \centering
    \begin{subfigure}{.33\textwidth}
        \centering
        \includegraphics[height=0.75\textwidth,width=1.0\textwidth]{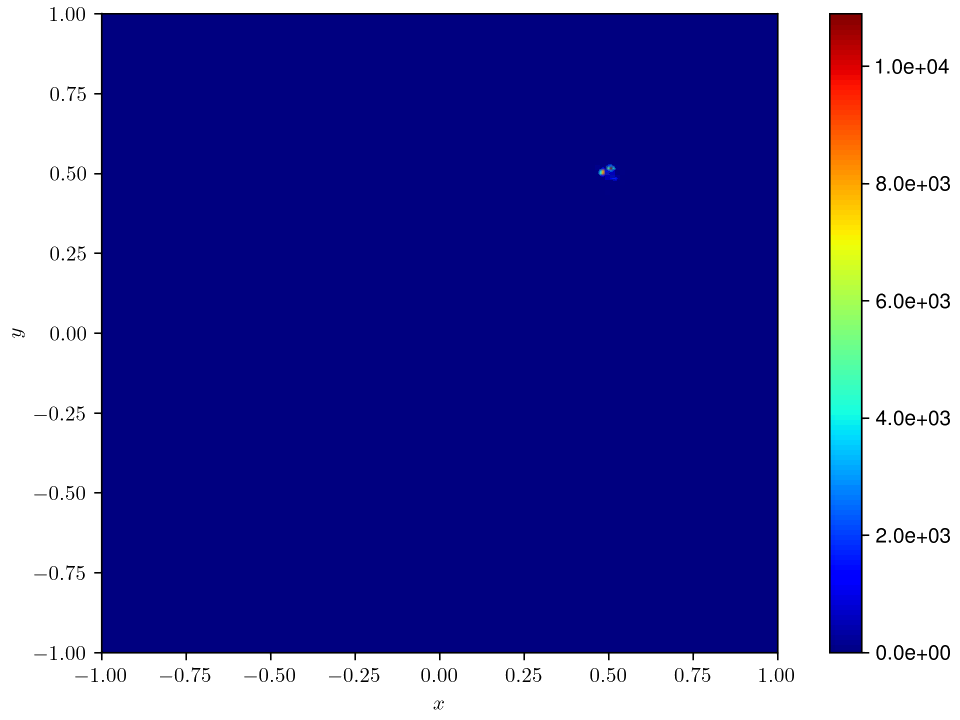}
    \end{subfigure}%
    \begin{subfigure}{.33\textwidth}
        \centering
        \includegraphics[height=0.75\textwidth,width=1.0\textwidth]{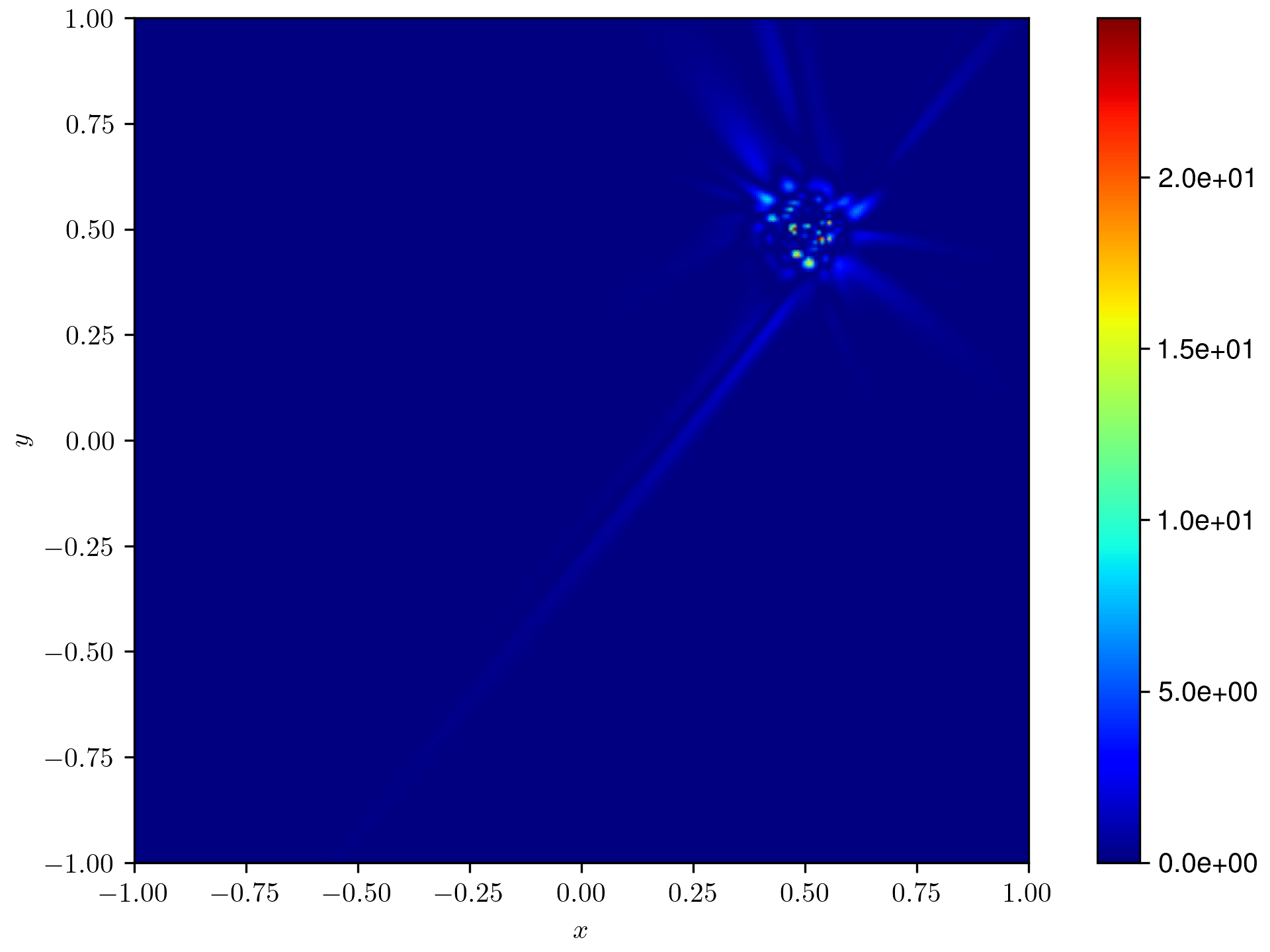}
    \end{subfigure}%
    \begin{subfigure}{.33\textwidth}
        \centering
        \includegraphics[height=0.75\textwidth,width=1.0\textwidth]{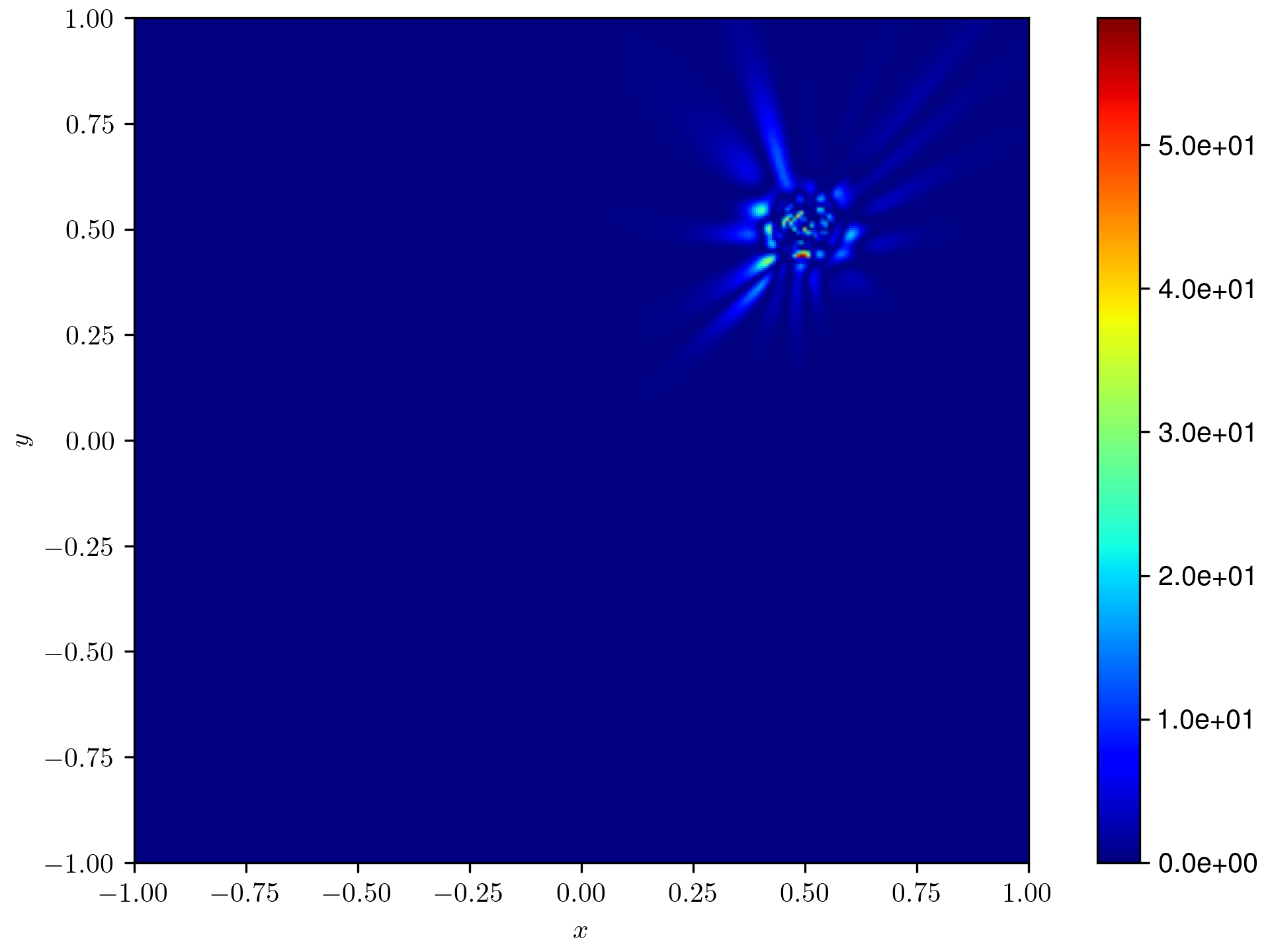}
    \end{subfigure}%
    \newline
    \raggedleft
    \begin{subfigure}{.33\textwidth}
        \centering
        \includegraphics[height=0.75\textwidth,width=1.0\textwidth]{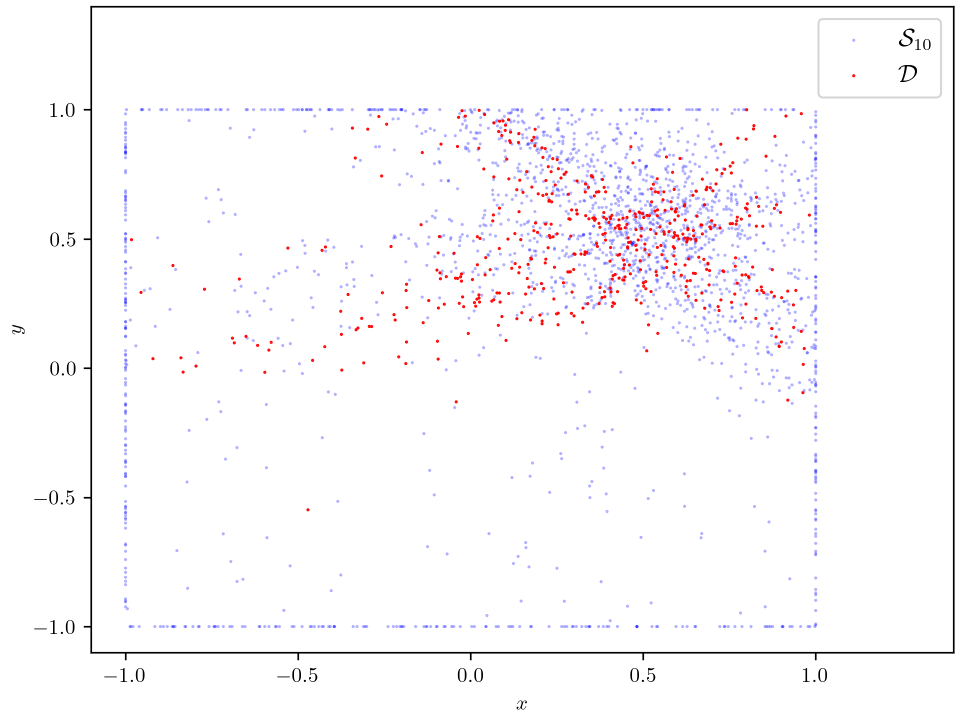}
    \end{subfigure}%
    \begin{subfigure}{.33\textwidth}
        \centering
        \includegraphics[height=0.75\textwidth,width=1.0\textwidth]{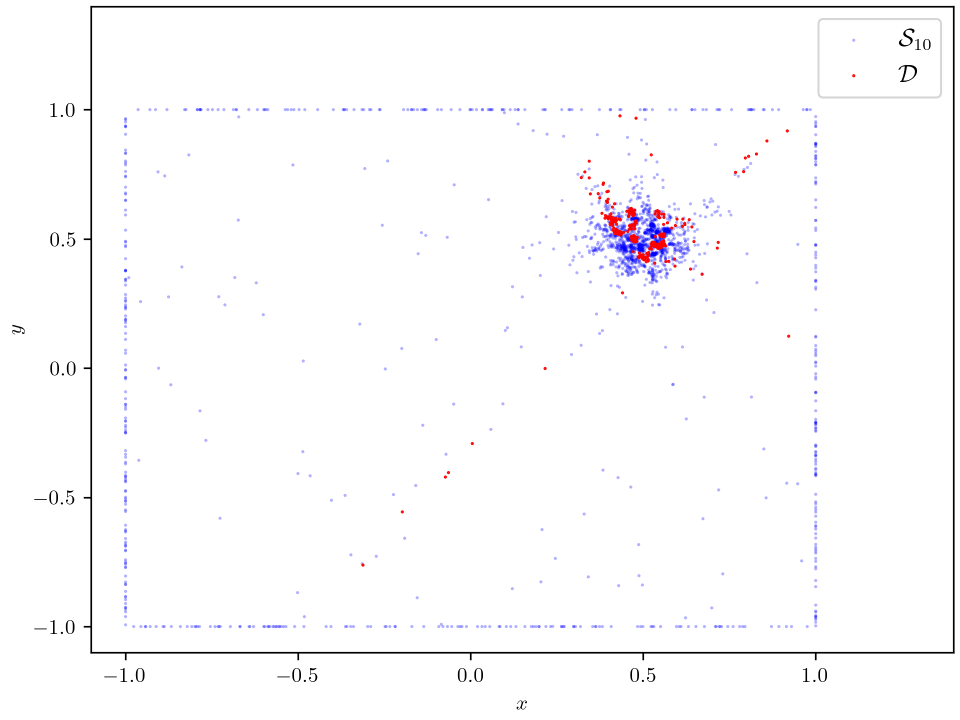}
    \end{subfigure}%
    \begin{subfigure}{.33\textwidth}
        \centering
        \includegraphics[height=0.75\textwidth,width=1.0\textwidth]{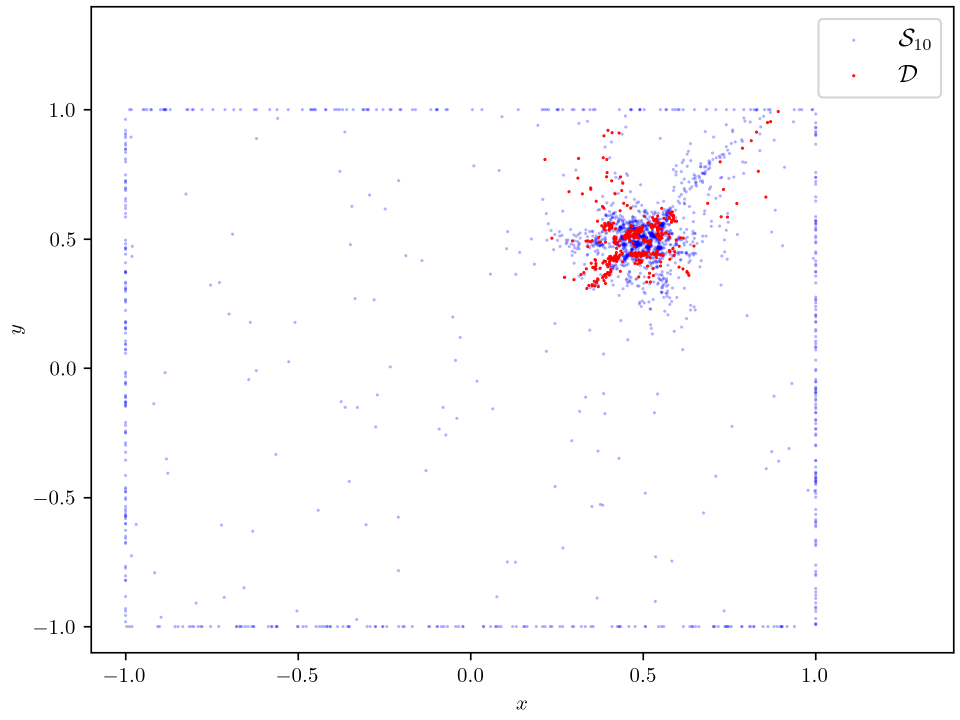}
    \end{subfigure}%
    \caption{Residual $\Q$ and training datasets  for \textit{RAD}, \textit{AAIS-g} and \textit{AAIS-t} for training schedule of lbfgs 1000 epochs after 9-th training. Left column: \textit{RAD}. Middle column: \textit{AAIS-g}. Right column: \textit{AAIS-t}. $\S_j$ means the nodes sampled from the dataset used in $(j-1)$-th iteration and $\D$ is the adaptive sampling nodes from the residual $\Q$ correspondingly. }
    \label{fig:PS1PlossIS7k}
\end{figure}

However, for adaptive sampling methods the results would be very different. They obtain a similar result  compared to \cite{GaoYanZhou2023:FIPINN:I}, \cite{Jiao2023:GAS} with fewer points and epochs. The primary reason for this phenomenon is the heightened focus on the singularity, as illustrated in the residual and node plots in Figure \ref{fig:PS1PlossIS7k}, where the nodes progressively converge around the singularity with the resampling strategy. For \textit{RAD} method, due to the small $N_S$, the loss function still concentrates around the singularity but nodes cluster around it. Moreover, for \textit{AAIS-t} and \textit{AAIS-g} algorithms, it is noticed that the frequency of residual increases during training, revealing that the PINNs are well-trained around the singularity and low-frequency of the solution is firstly and well learned according to the NTK theory. The absolute error and the neural network solution in Figure \ref{fig:PS1PAbs7k} also support the fact that compared to the numerical results from \textit{Uni} method, the adaptive sampling methods generate a better solution that the singularity would hide from the absolute error and the frequency of the absolute error increases.

\begin{figure}[htbp]
    \centering
    \begin{subfigure}{.5\textwidth}
        \centering
        \includegraphics[scale=0.3]{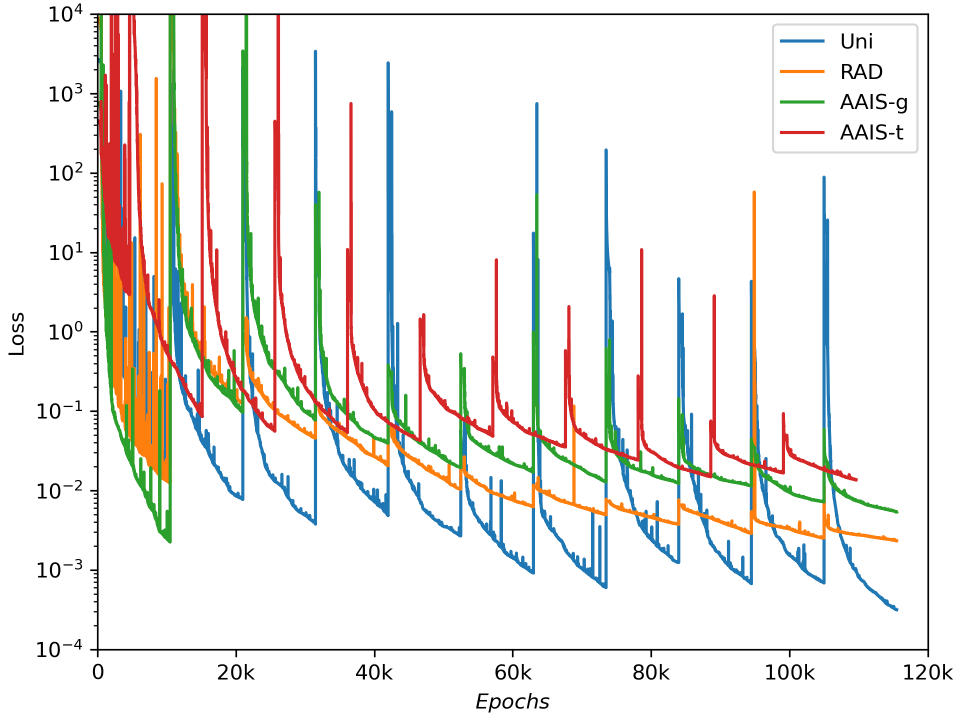}
        \end{subfigure}%
        \begin{subfigure}{.5\textwidth}
        \centering
        \includegraphics[scale=0.3]{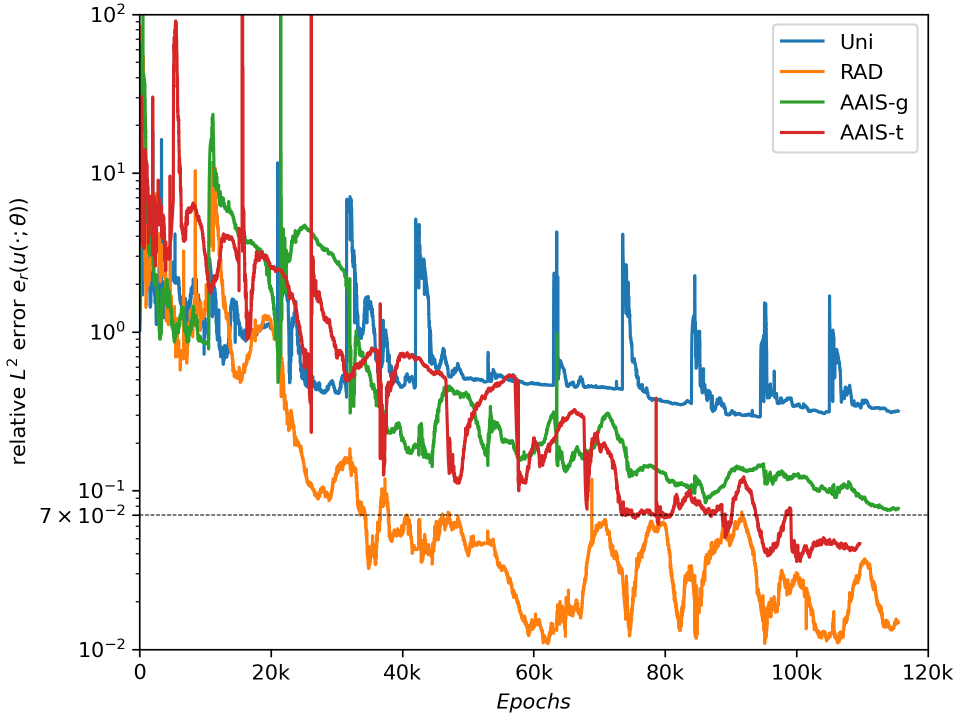}
        \end{subfigure}%
        \caption{Loss and relative errors during training for one peak Poisson equation with training schedule of lbfgs 10000 epochs each. Left: the loss function. Right: the relative $L^2$ error $e_r(u(\cdot;\theta))$.
        }
    \label{fig:PS1PErr10000e}
\end{figure}

\begin{figure}[htbp]
    \centering
    \begin{subfigure}{.25\textwidth}
        \centering
        \includegraphics[height=0.75\textwidth,width=1.0\textwidth]{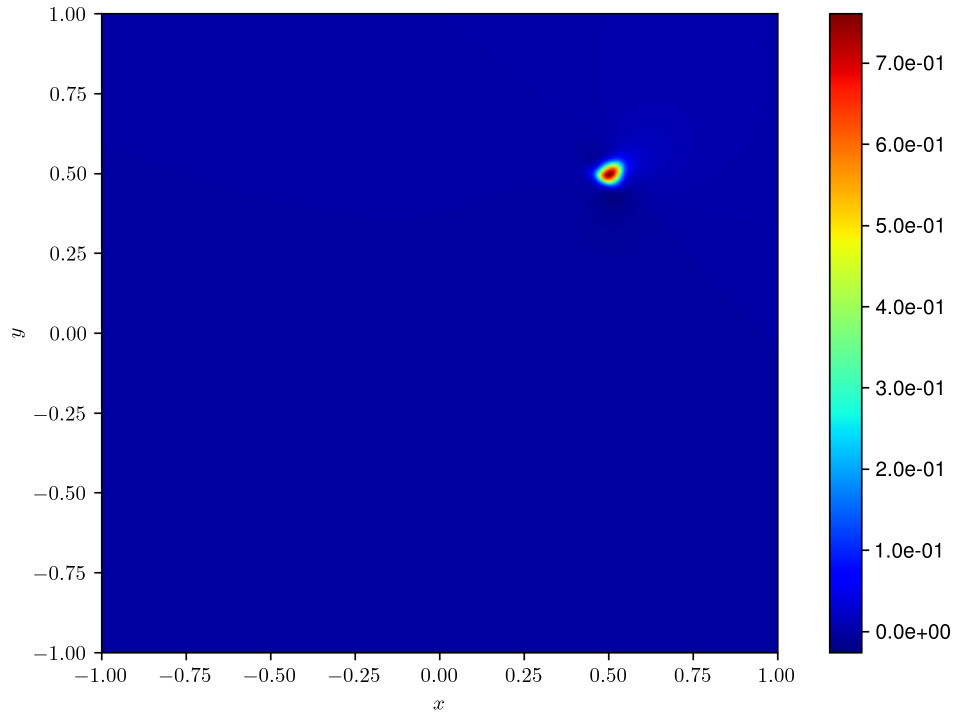}
    \end{subfigure}%
    \begin{subfigure}{.25\textwidth}
        \centering
        \includegraphics[height=0.75\textwidth,width=1.0\textwidth]{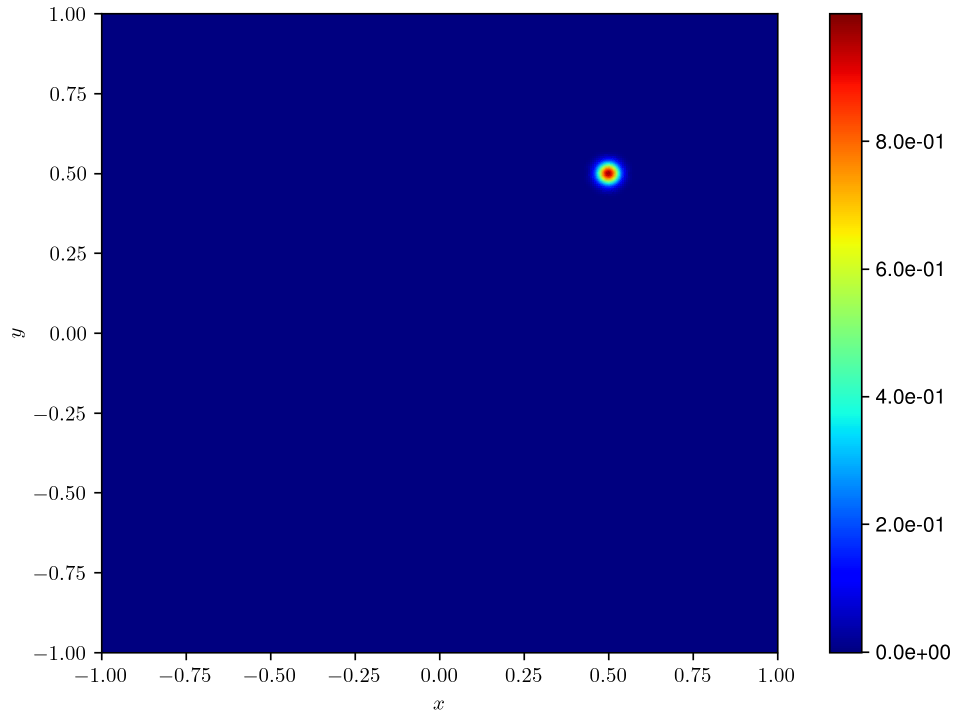}
    \end{subfigure}%
    \begin{subfigure}{.25\textwidth}
        \centering
        \includegraphics[height=0.75\textwidth,width=1.0\textwidth]{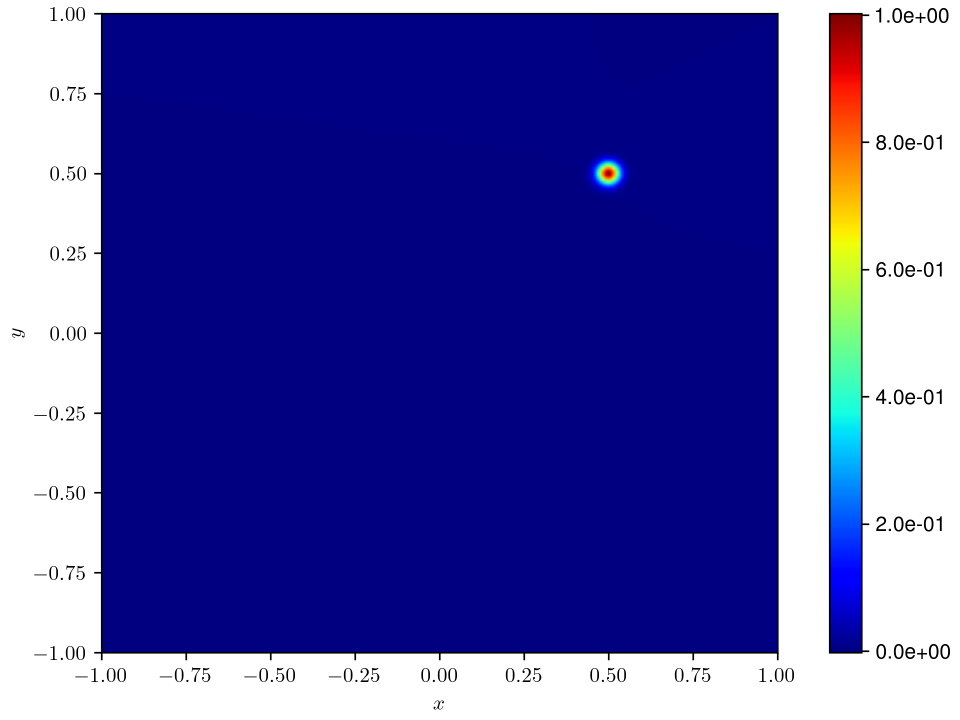}
    \end{subfigure}%
    \begin{subfigure}{.25\textwidth}
        \centering
        \includegraphics[height=0.75\textwidth,width=1.0\textwidth]{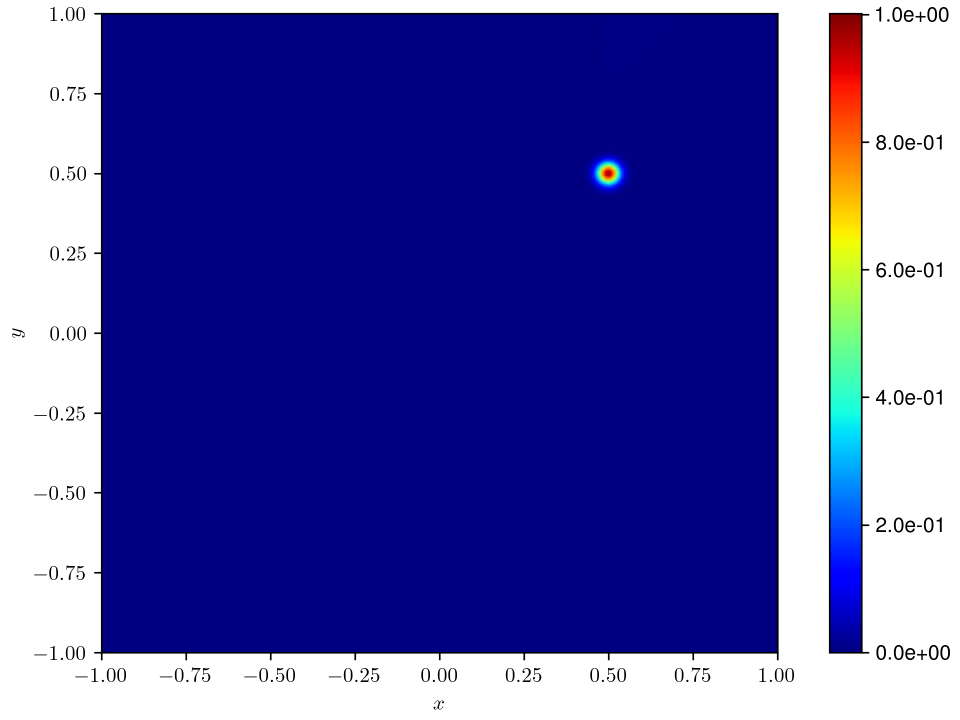}
    \end{subfigure}%
    \newline
    \raggedleft
    \begin{subfigure}{.25\textwidth}
        \centering
        \includegraphics[height=0.75\textwidth,width=1.0\textwidth]{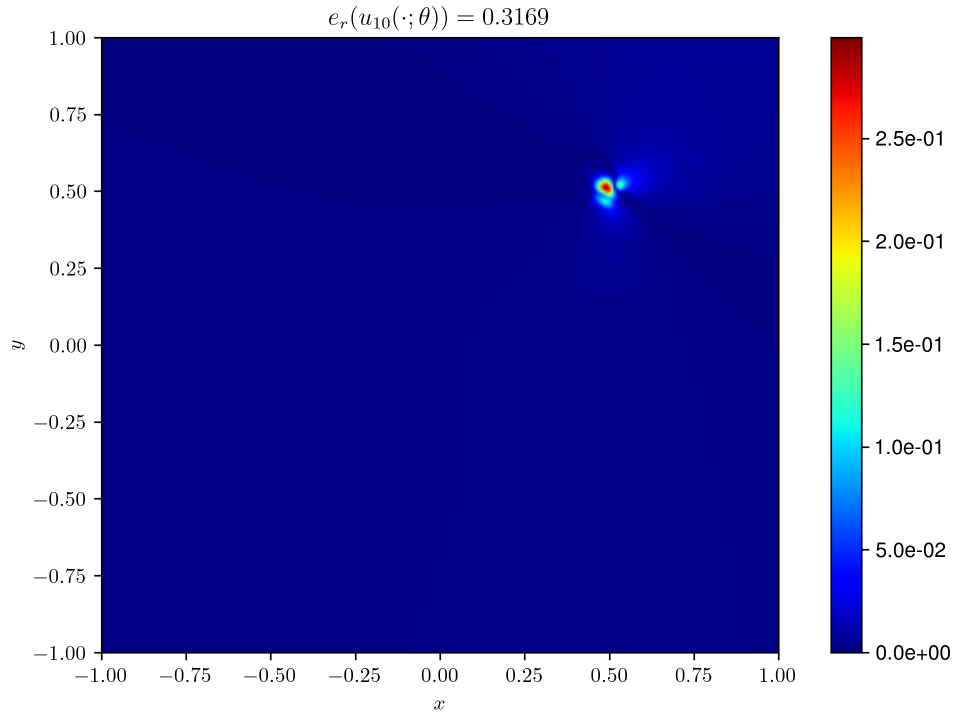}
        \caption{\textit{Uni}}
    \end{subfigure}%
    \begin{subfigure}{.25\textwidth}
        \centering
        \includegraphics[height=0.75\textwidth,width=1.0\textwidth]{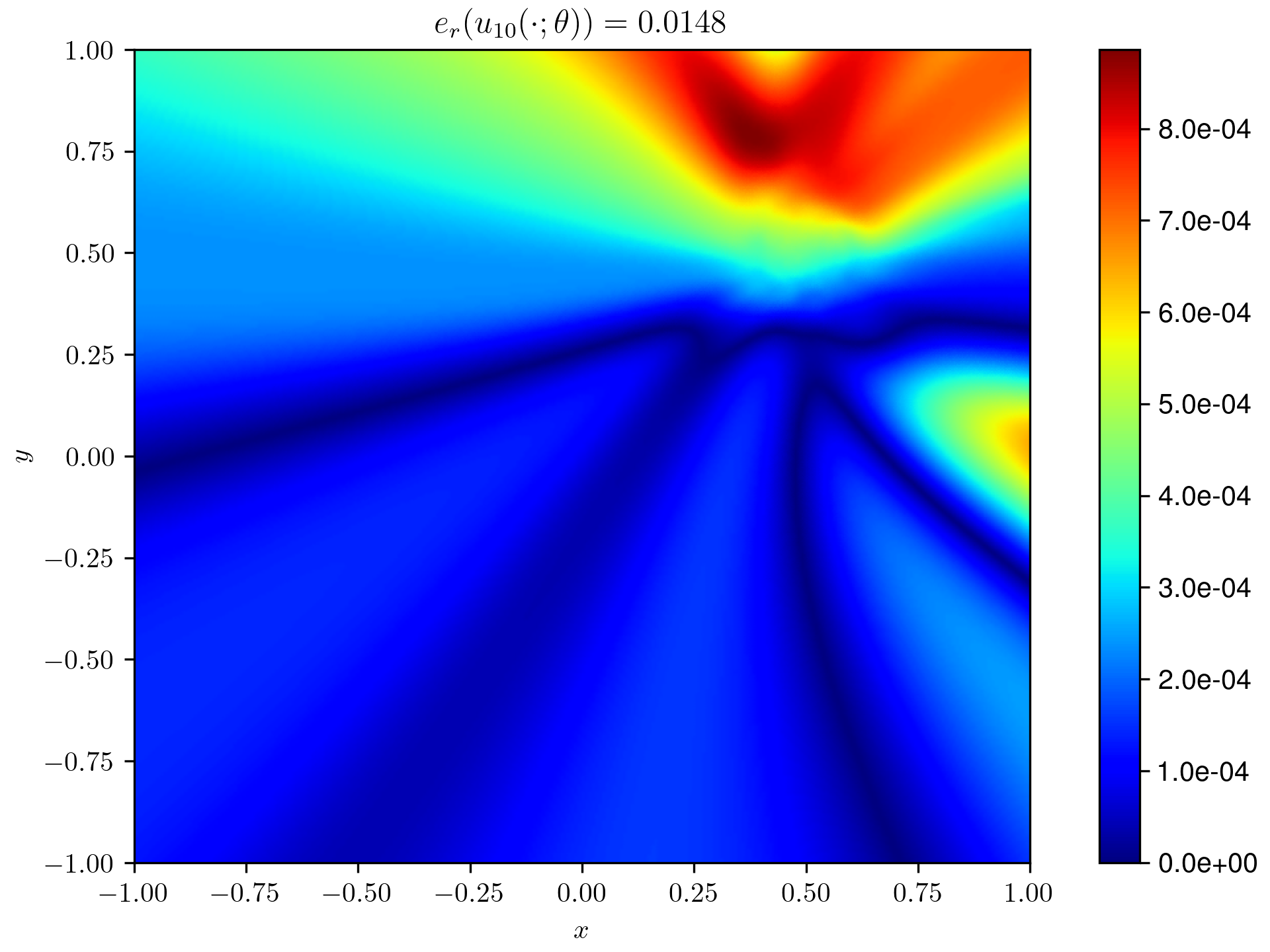}
        \caption{\textit{RAD}}
    \end{subfigure}%
    \begin{subfigure}{.25\textwidth}
        \centering
        \includegraphics[height=0.75\textwidth,width=1.0\textwidth]{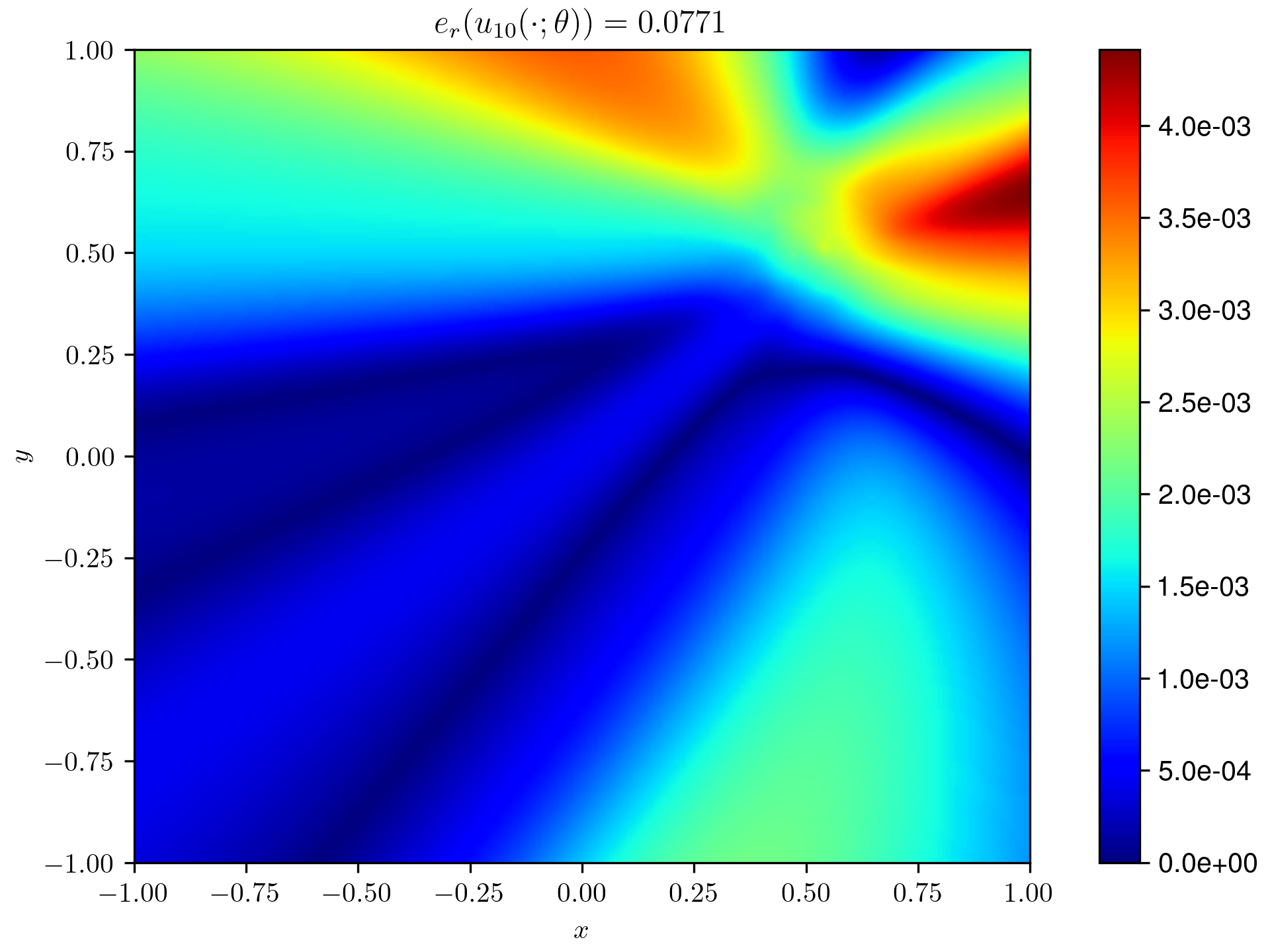}
        \caption{\textit{AAIS-g}}
    \end{subfigure}%
    \begin{subfigure}{.25\textwidth}
        \centering
        \includegraphics[height=0.75\textwidth,width=1.0\textwidth]{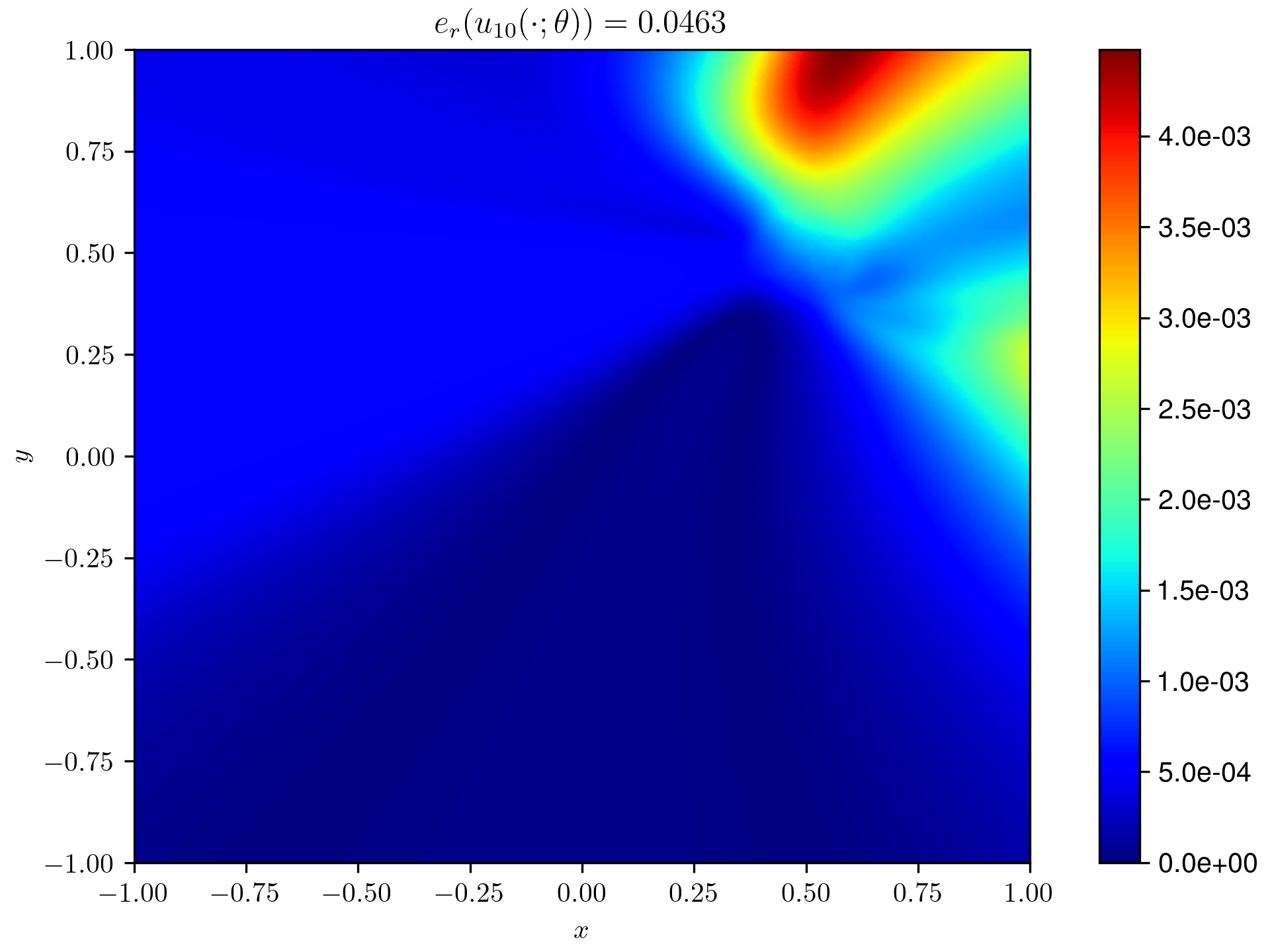}
        \caption{\textit{AAIS-t}}
    \end{subfigure}%
    \caption{Profiles of absolute error and neural network solutions after 10th training for one peak Poisson equation with training schedule of 10000 epochs for lbfgs. First row: numerical solutions. Second row: absolute error.}
    \label{fig:PS1PAbs10000e}
\end{figure}

In the following we keep the parameters fixed of above setting and let the maximum epochs of lbfgs optimizer at pre-train and each iteration to be 10000. And $N_S=100000$ for \textit{RAD} and $N_S=10N_A=60000$ for AAIS algorithms. With more training epochs the relative error of adaptive sampling methods behave better than the above schedule of 1000 epochs, see Figure \ref{fig:PS1PErr10000e}. We could see that the \textit{RAD} method would arrive at relative error $1\%$ which is significant smaller than the previous adaptive sampling work \cite{GaoYanZhou2023:FIPINN:I,Jiao2023:GAS}, and \textit{AAIS-g, AAIS-t} also behave better than the mentioned work. The profiles of solutions are listed in Figure \ref{fig:PS1PAbs10000e}.

\subsubsection{Nine peaks}
\begin{figure}[htbp]
    \centering
    \includegraphics[scale=0.25]{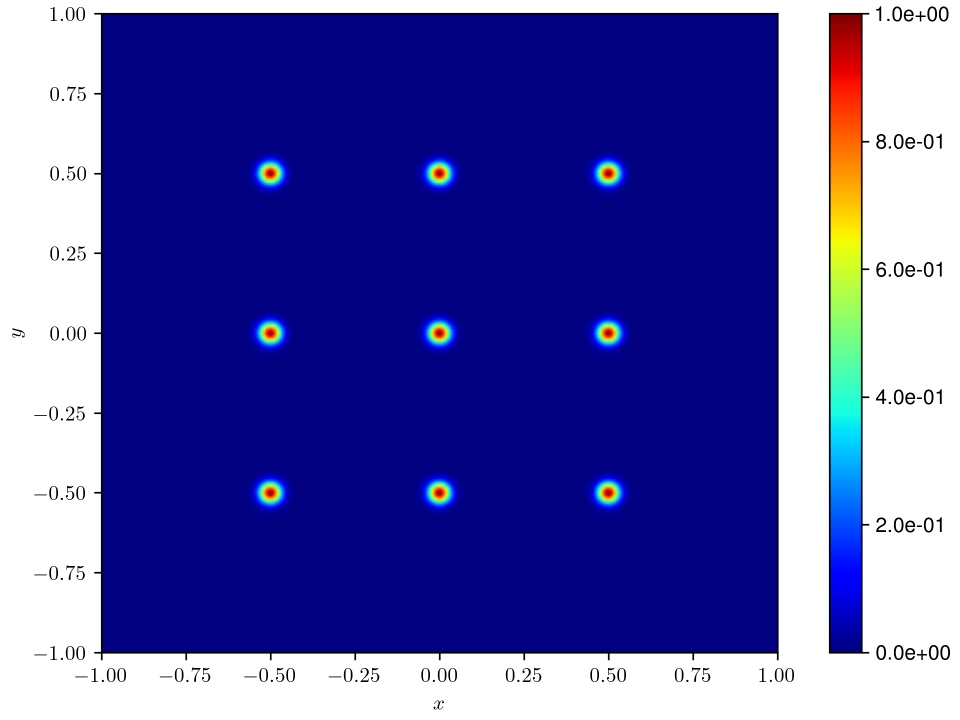}
    \caption{Exact solution for Poisson problems with nine peaks in \eqref{pde:Poisson1Peak}.}
    \label{fig:PS9Pexact}
\end{figure}

Then we consider Nine peaks of Poisson problems where the centers $(x_0^i,y_0^i)=(-0.5,0.5)+(\frac{mod(i,3)}{2},0)+(0, \frac{\lfloor i/3\rfloor}{2})$, $i=0,...,8$ (\cite{Jiao2023:GAS}) are equally distributed in the domain. The exact solution is showed in Figure \ref{fig:PS9Pexact}.

Similarly, using Monte-Carlo integration for $\int\Q(\x)\d\x$, we uniformly sample 100k points inside $\Omega$. For the total points in the domain we set 2k points of which 500 points are sampled according to each sampling method. 500 points are sampled on the boundary.

\begin{figure}[htbp]
    \centering
    \begin{subfigure}{.5\textwidth}
    \centering
    \includegraphics[scale=0.25]{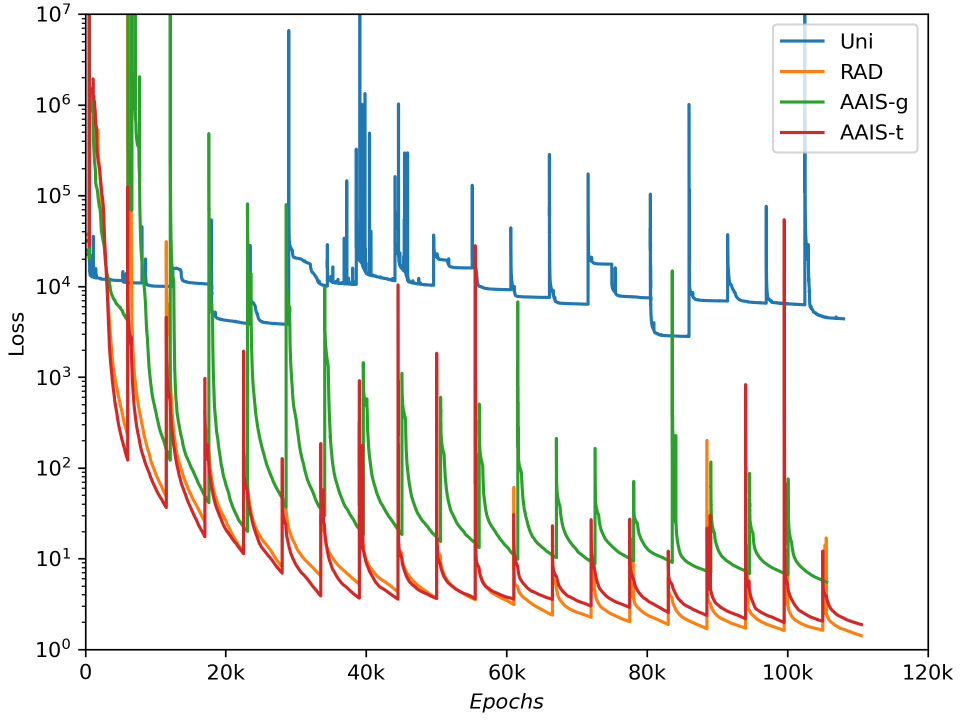}
    \end{subfigure}%
    \begin{subfigure}{.5\textwidth}
    \centering
    \includegraphics[scale=0.25]{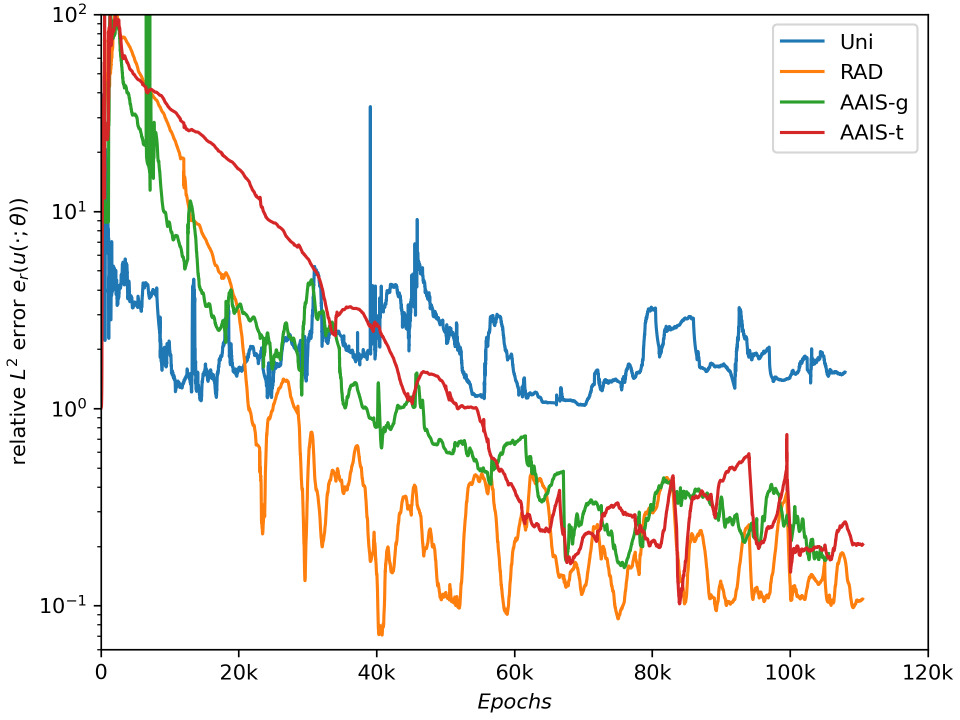}
    \end{subfigure}%
    \caption{Loss and relative errors during training for Poisson equation with nine peaks with four sampling methods. Left: the loss function. Right: the relative $L^2$ error $e_r(u(\cdot;\theta))$.
    }
    \label{fig:PS9PLossErr5000e}
\end{figure}

Firstly we train the PINNs with Adam 500 epochs and lbfgs 5000 epochs  at pre-train and each iteration for maximum iteration $M=20$. We set the numbers $N_A$ of points for AAIS be 6k, because the multi-peaks problem demands capabilities of searching for adaptive sampling methods. In Figure \ref{fig:PS9PLossErr5000e}, we list the loss and relative errors within the epochs. The adaptive sampling methods generate better results than \textit{Uni} method. According to the parameter setting in \cite{Jiao2023:GAS}, we find that we could obtain similar results with fewer points.

\begin{figure}[htbp]
    \centering
    \begin{subfigure}{.25\textwidth}
        \centering
        \includegraphics[height=0.75\textwidth,width=1.0\textwidth]{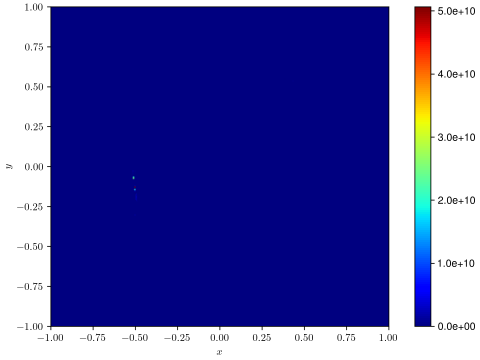}
    \end{subfigure}%
    \begin{subfigure}{.25\textwidth}
        \centering
        \includegraphics[height=0.75\textwidth,width=1.0\textwidth]{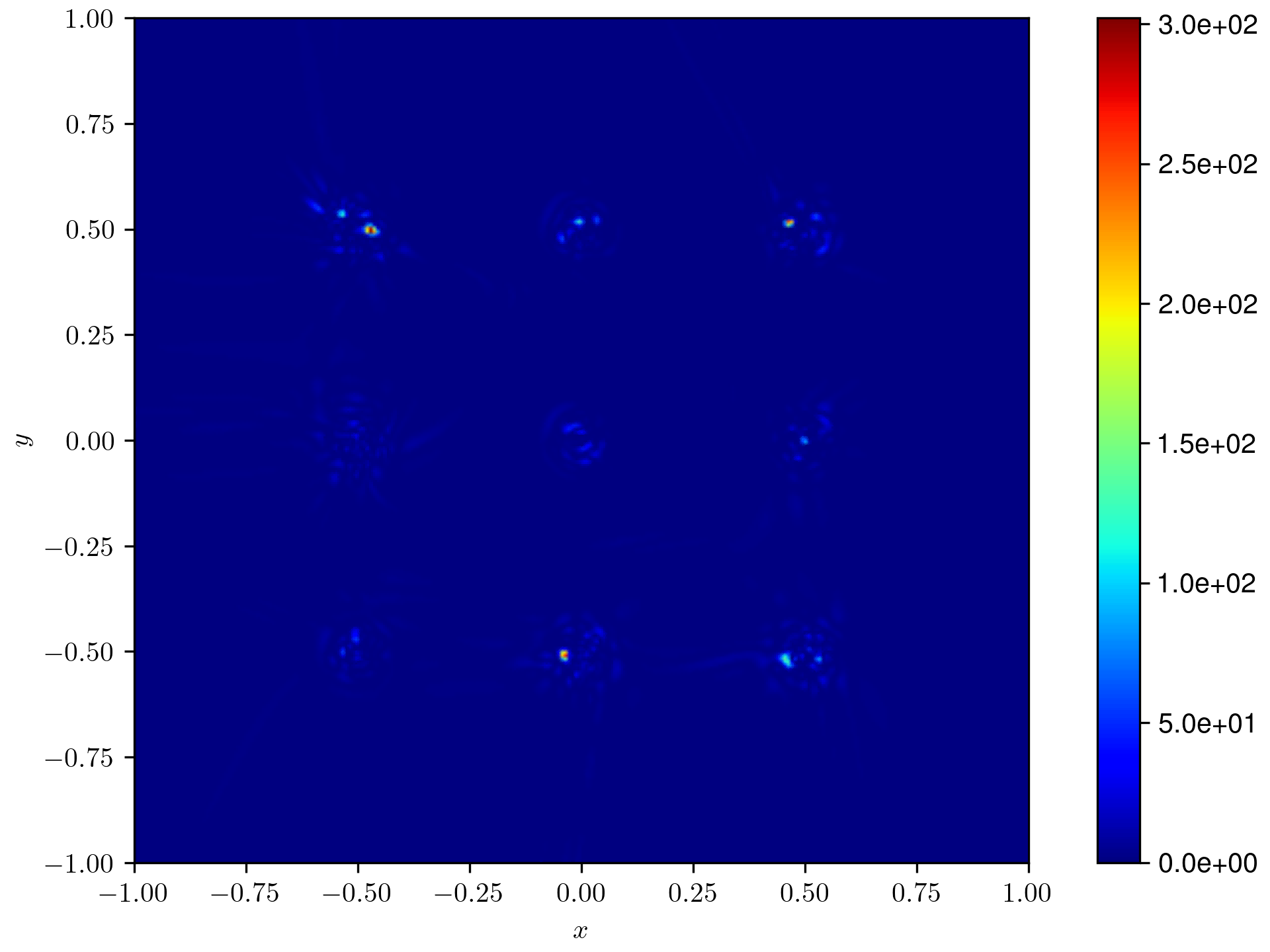}
    \end{subfigure}%
    \begin{subfigure}{.25\textwidth}
        \centering
        \includegraphics[height=0.75\textwidth,width=1.0\textwidth]{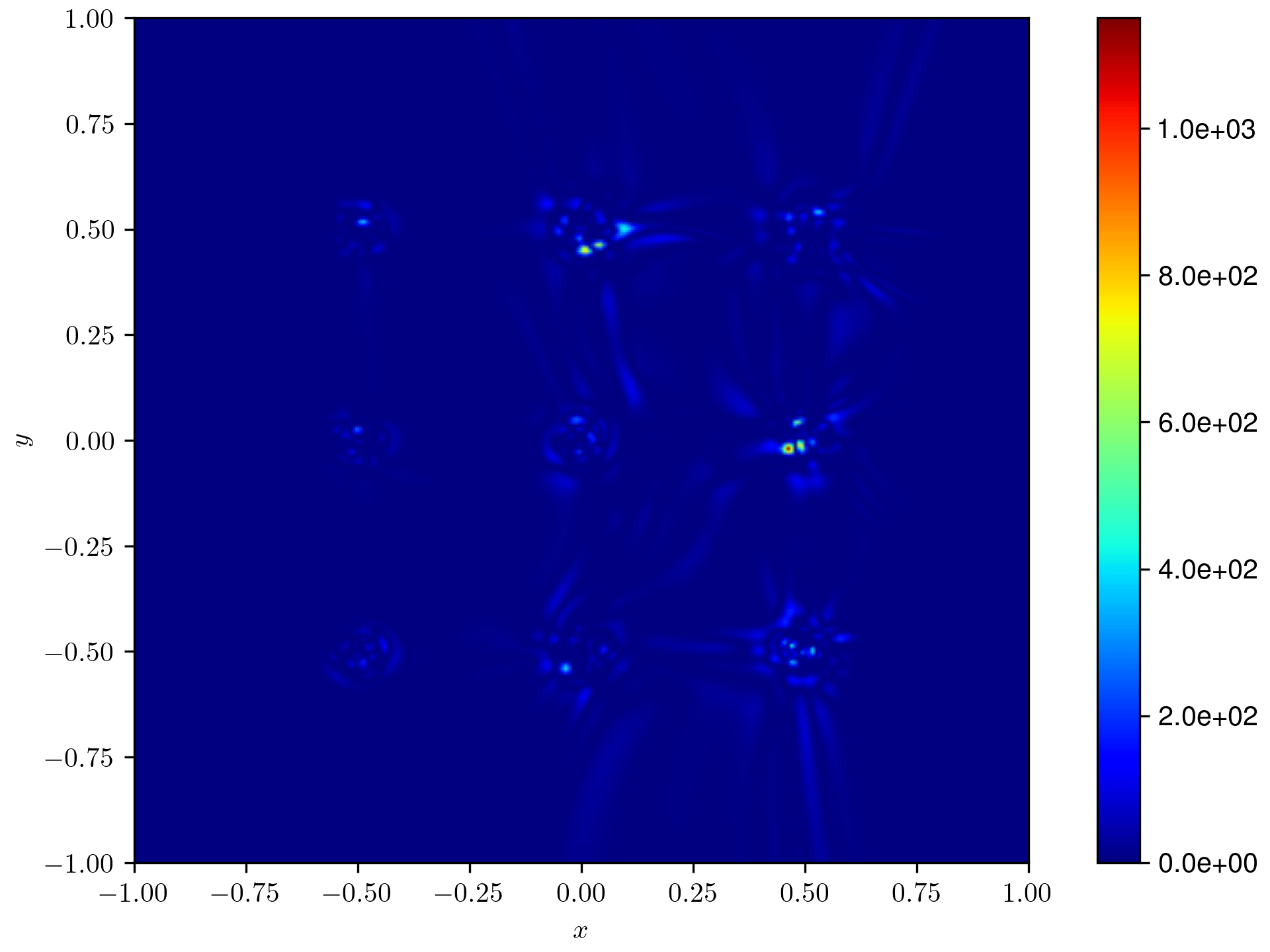}
    \end{subfigure}%
    \begin{subfigure}{.25\textwidth}
        \centering
        \includegraphics[height=0.75\textwidth,width=1.0\textwidth]{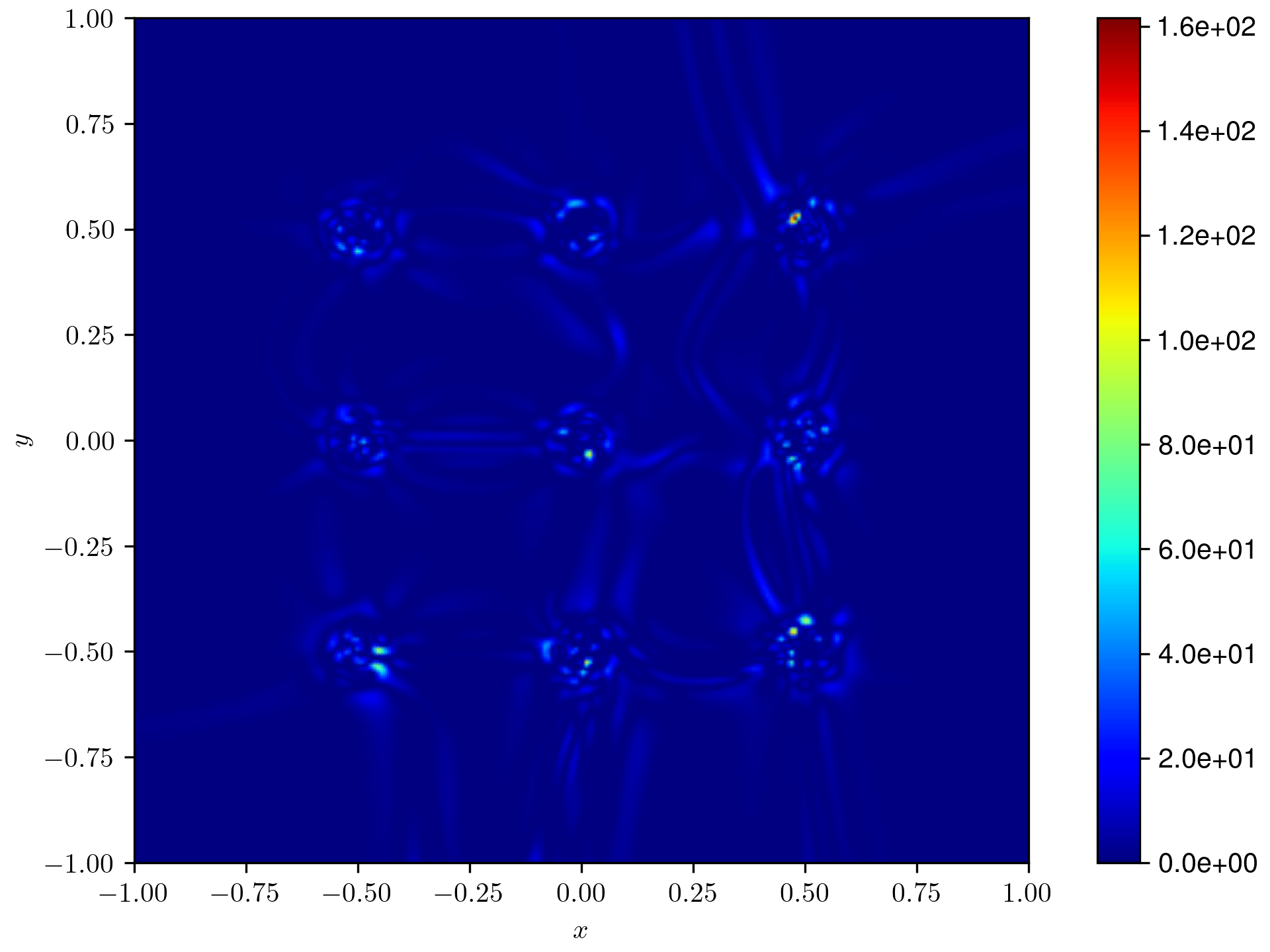}
    \end{subfigure}%
    \newline
    \raggedleft
    \begin{subfigure}{.25\textwidth}
        \centering
        \includegraphics[height=0.75\textwidth,width=1.0\textwidth]{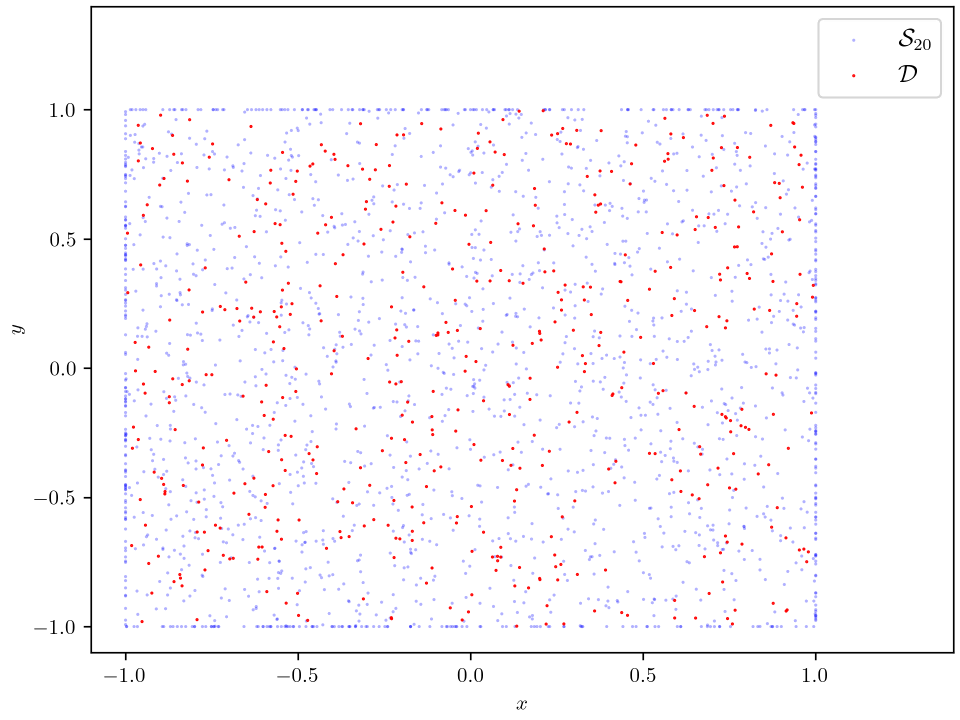}
        \caption{\textit{Uni}}
    \end{subfigure}%
    \begin{subfigure}{.25\textwidth}
        \centering
        \includegraphics[height=0.75\textwidth,width=1.0\textwidth]{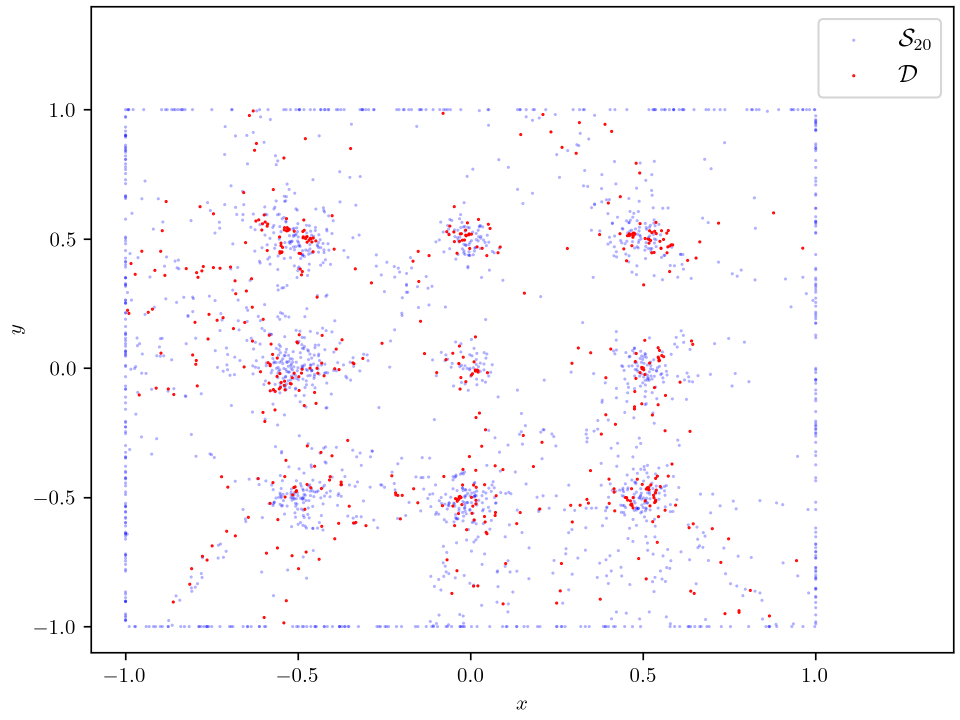}
        \caption{\textit{RAD}}
    \end{subfigure}%
    \begin{subfigure}{.25\textwidth}
        \centering
        \includegraphics[height=0.75\textwidth,width=1.0\textwidth]{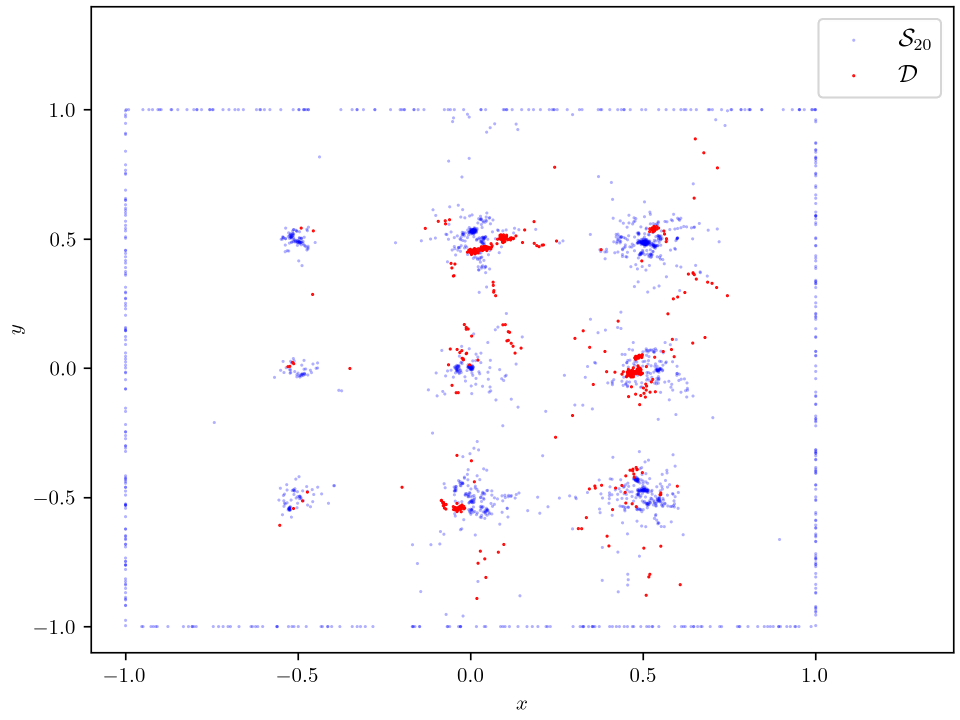}
        \caption{\textit{AAIS-g}}
    \end{subfigure}%
    \begin{subfigure}{.25\textwidth}
        \centering
        \includegraphics[height=0.75\textwidth,width=1.0\textwidth]{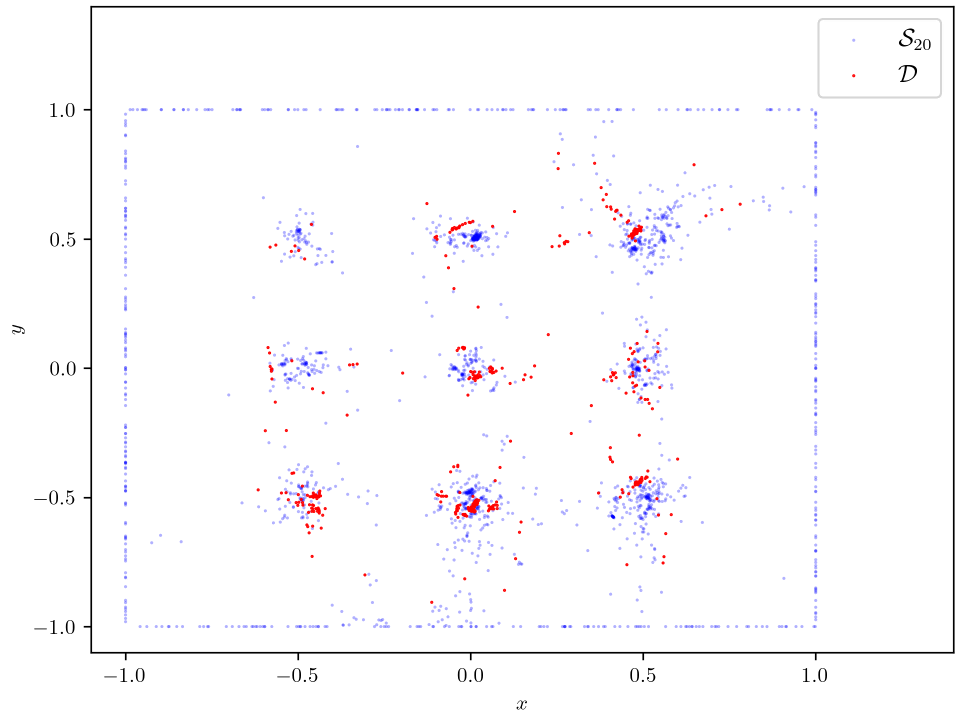}
        \caption{\textit{AAIS-t}}
    \end{subfigure}%
    \caption{Profiles of residual and nodes  for nine peaks Poisson equation at 20-th iteration. }
    \label{fig:Ps9PNode5000e}
\end{figure}

\begin{figure}[htbp]
    \centering
    \begin{subfigure}{.25\textwidth}
        \centering
        \includegraphics[height=0.75\textwidth,width=1.0\textwidth]{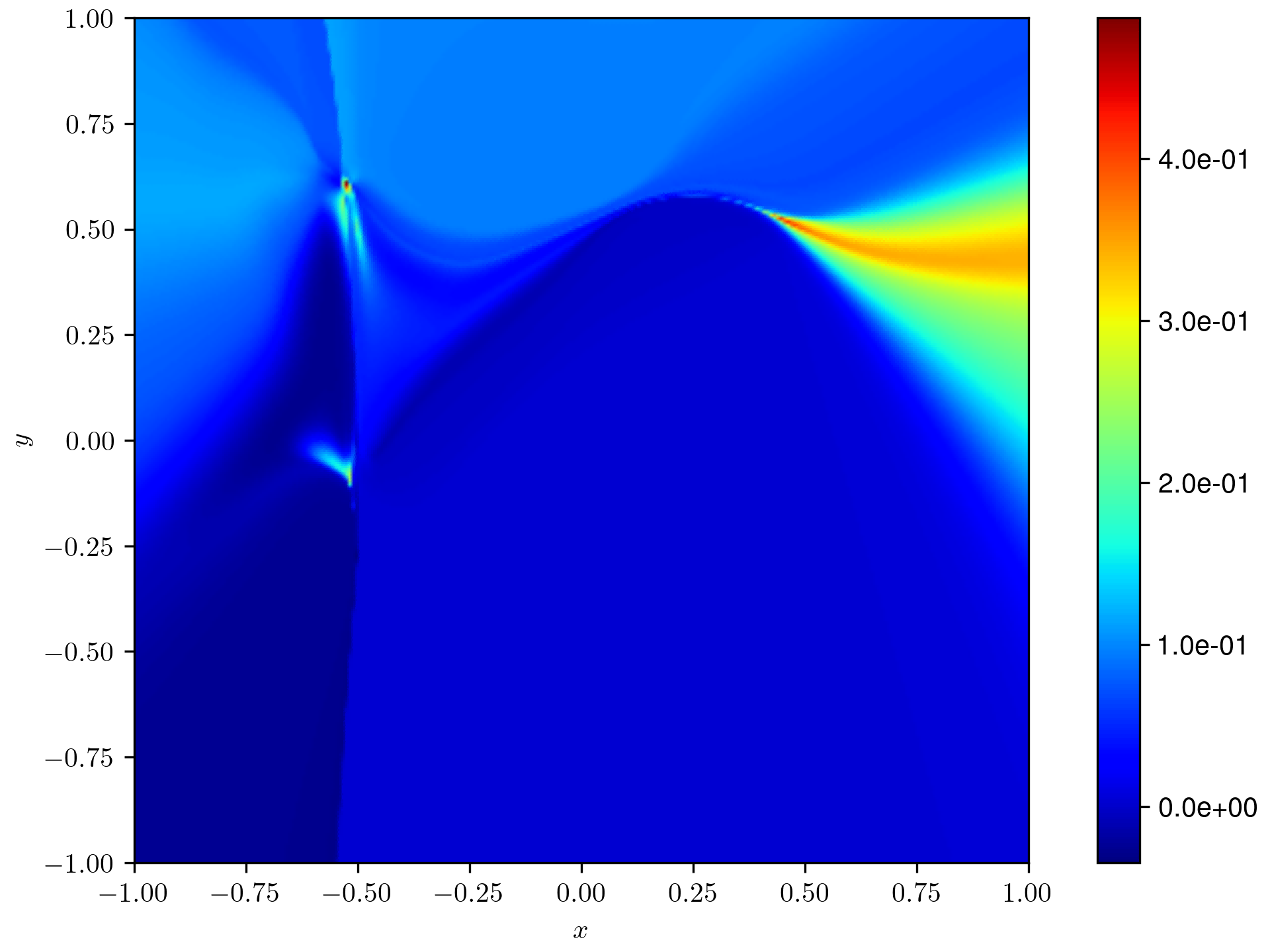}
    \end{subfigure}%
    \begin{subfigure}{.25\textwidth}
        \centering
        \includegraphics[height=0.75\textwidth,width=1.0\textwidth]{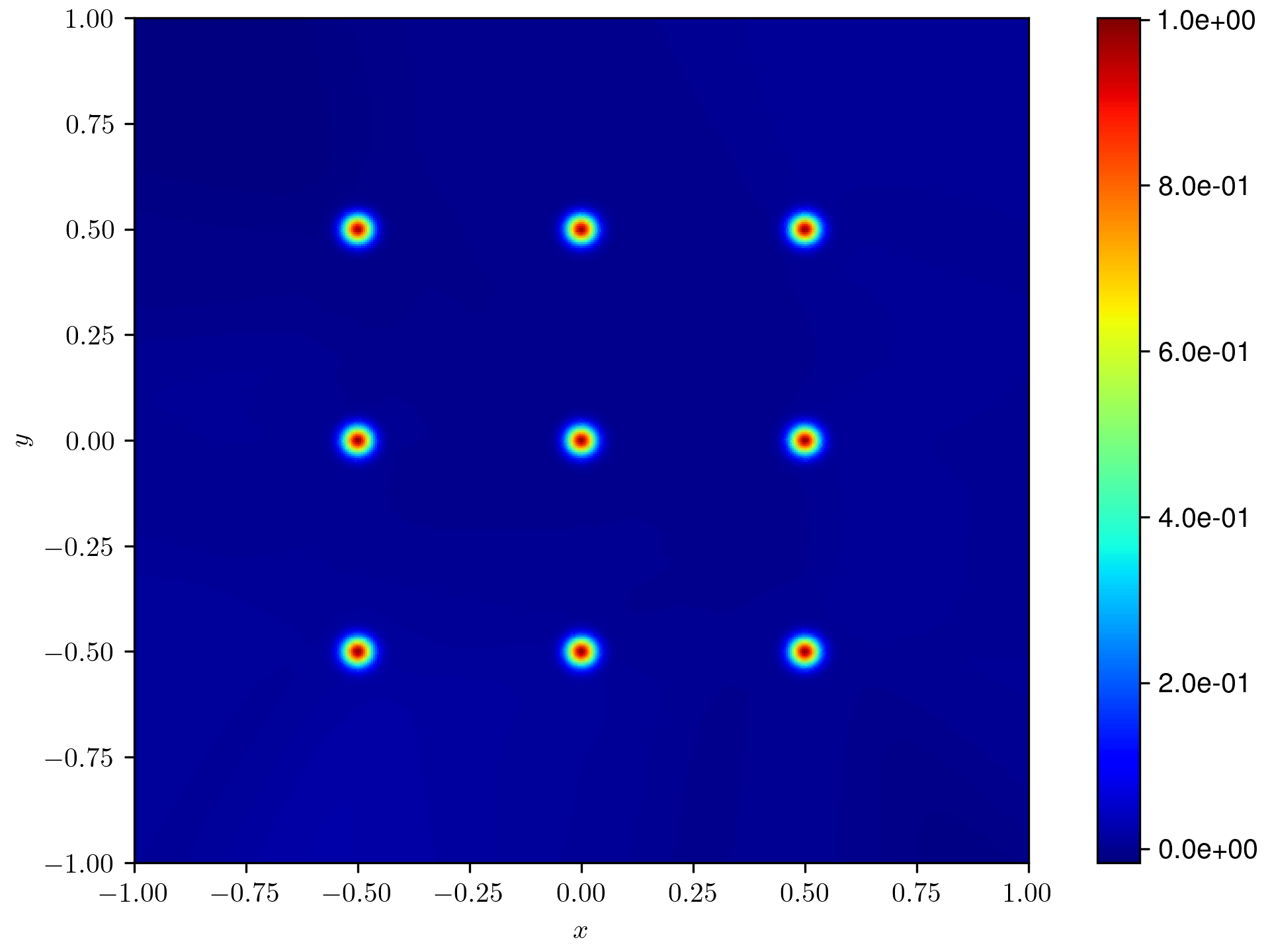}
    \end{subfigure}%
    \begin{subfigure}{.25\textwidth}
        \centering
        \includegraphics[height=0.75\textwidth,width=1.0\textwidth]{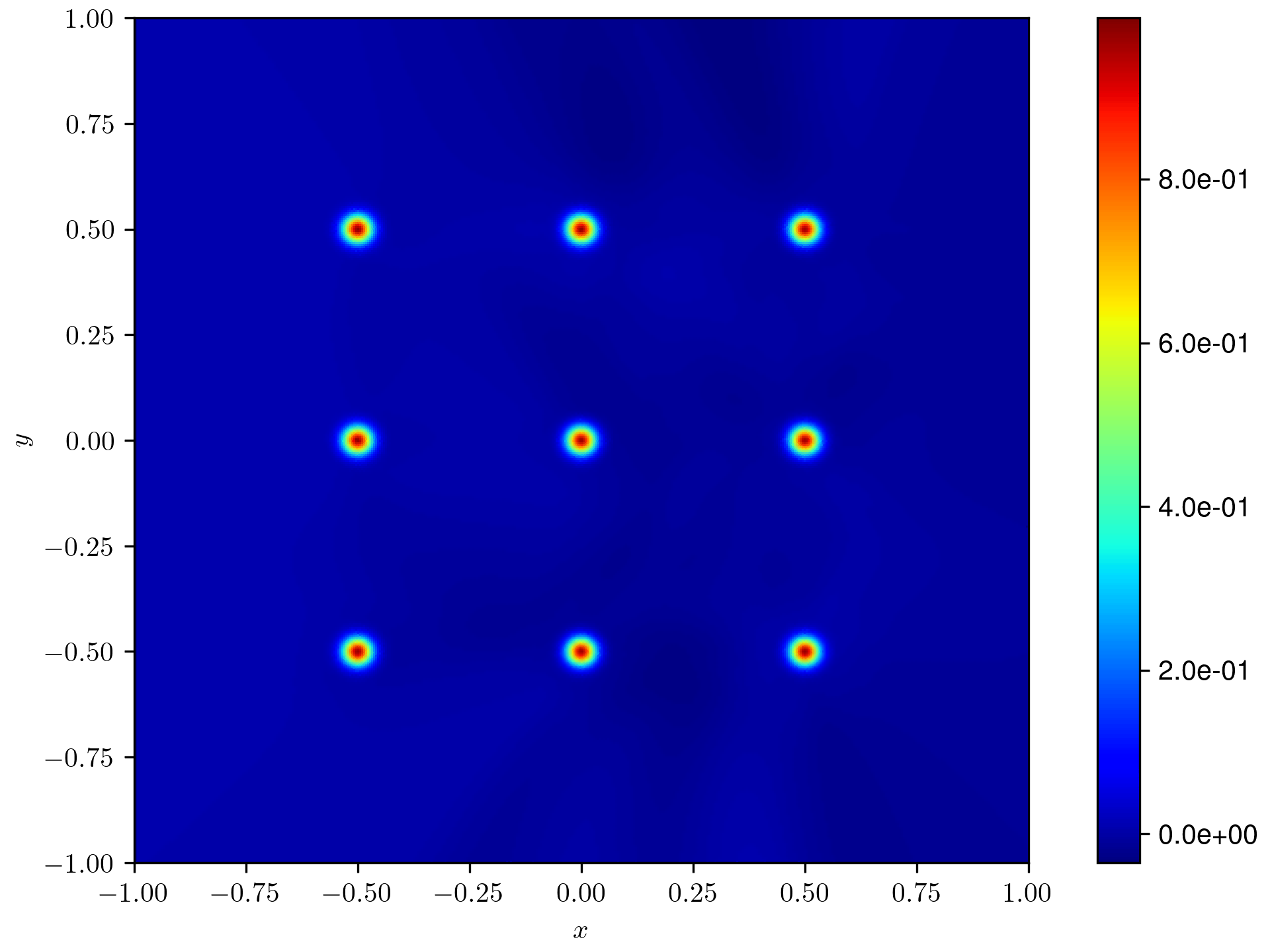}
    \end{subfigure}%
    \begin{subfigure}{.25\textwidth}
        \centering
        \includegraphics[height=0.75\textwidth,width=1.0\textwidth]{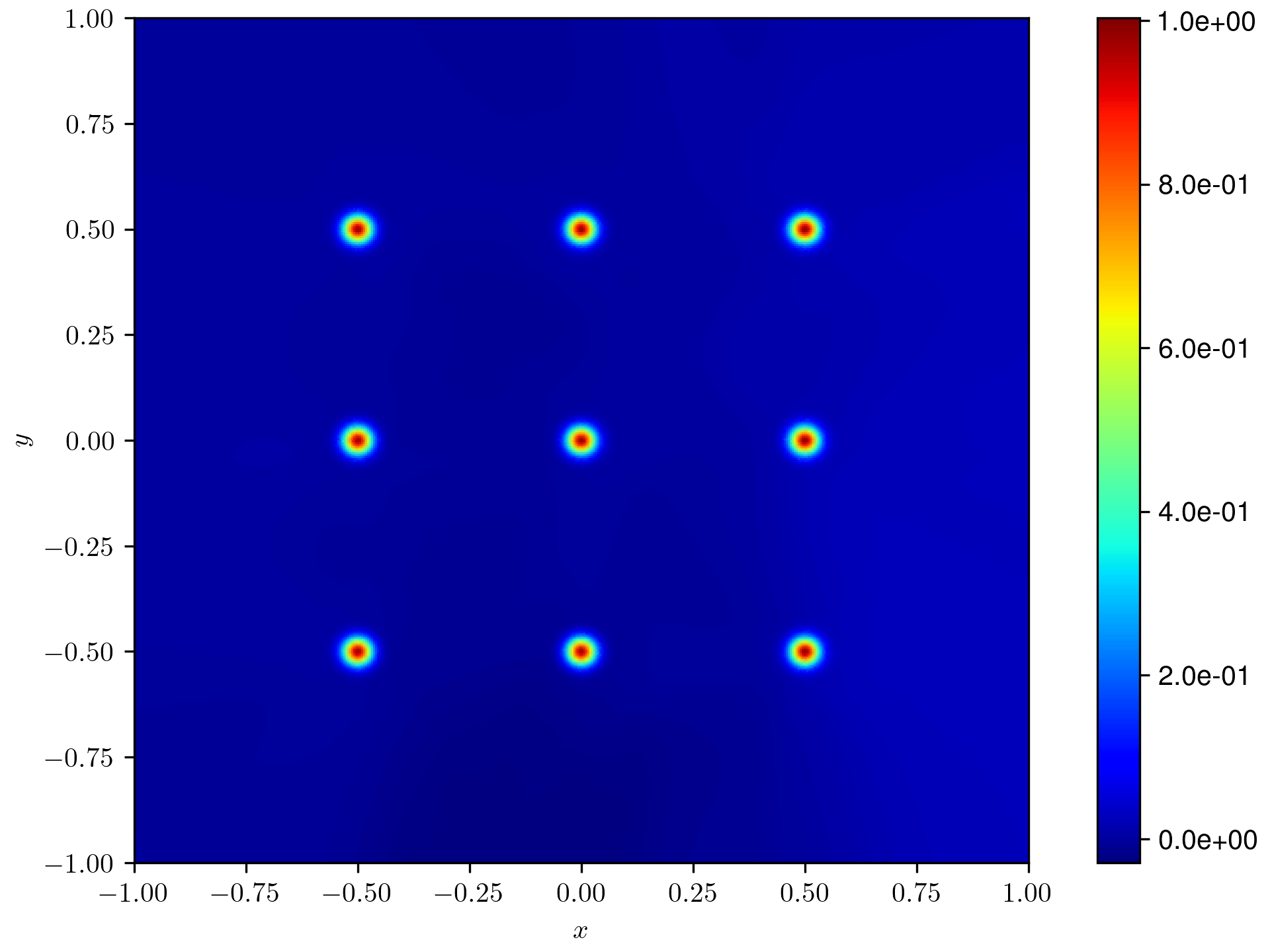}
    \end{subfigure}%
    \newline
    \raggedleft
    \begin{subfigure}{.25\textwidth}
        \centering
        \includegraphics[height=0.75\textwidth,width=1.0\textwidth]{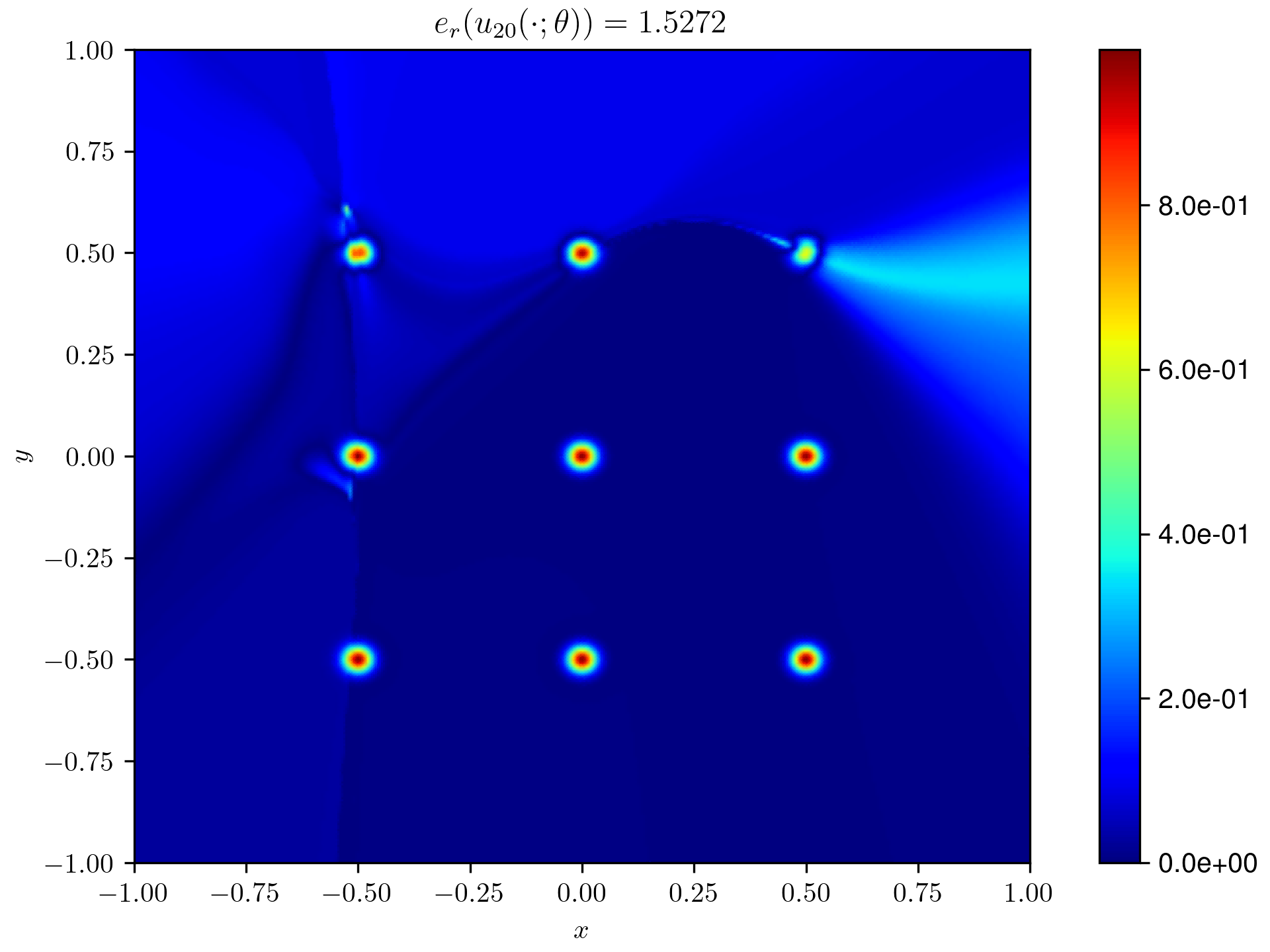}
        \caption{\textit{Uni}}
    \end{subfigure}%
    \begin{subfigure}{.25\textwidth}
        \centering
        \includegraphics[height=0.75\textwidth,width=1.0\textwidth]{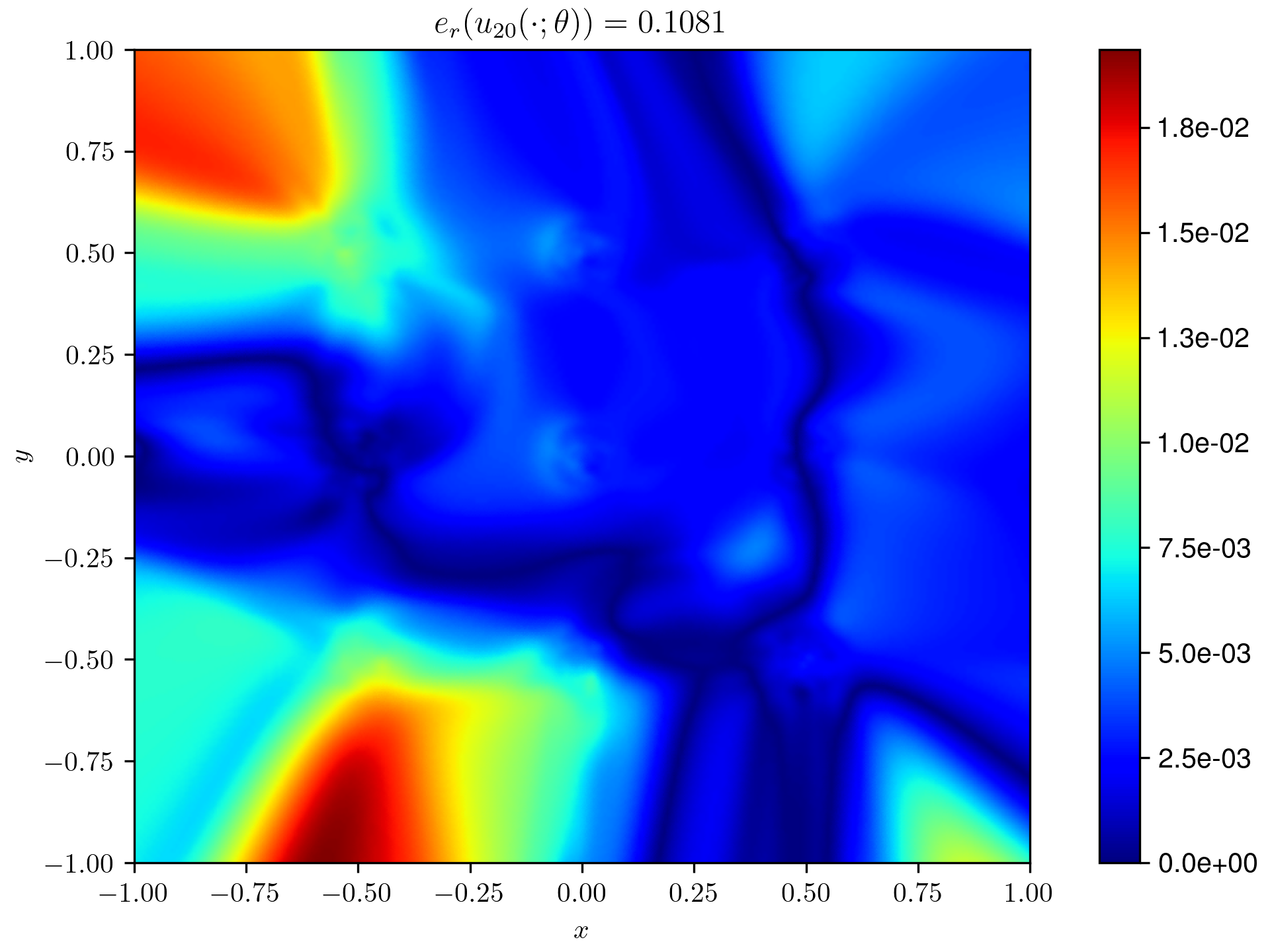}
        \caption{\textit{RAD}}
    \end{subfigure}%
    \begin{subfigure}{.25\textwidth}
        \centering
        \includegraphics[height=0.75\textwidth,width=1.0\textwidth]{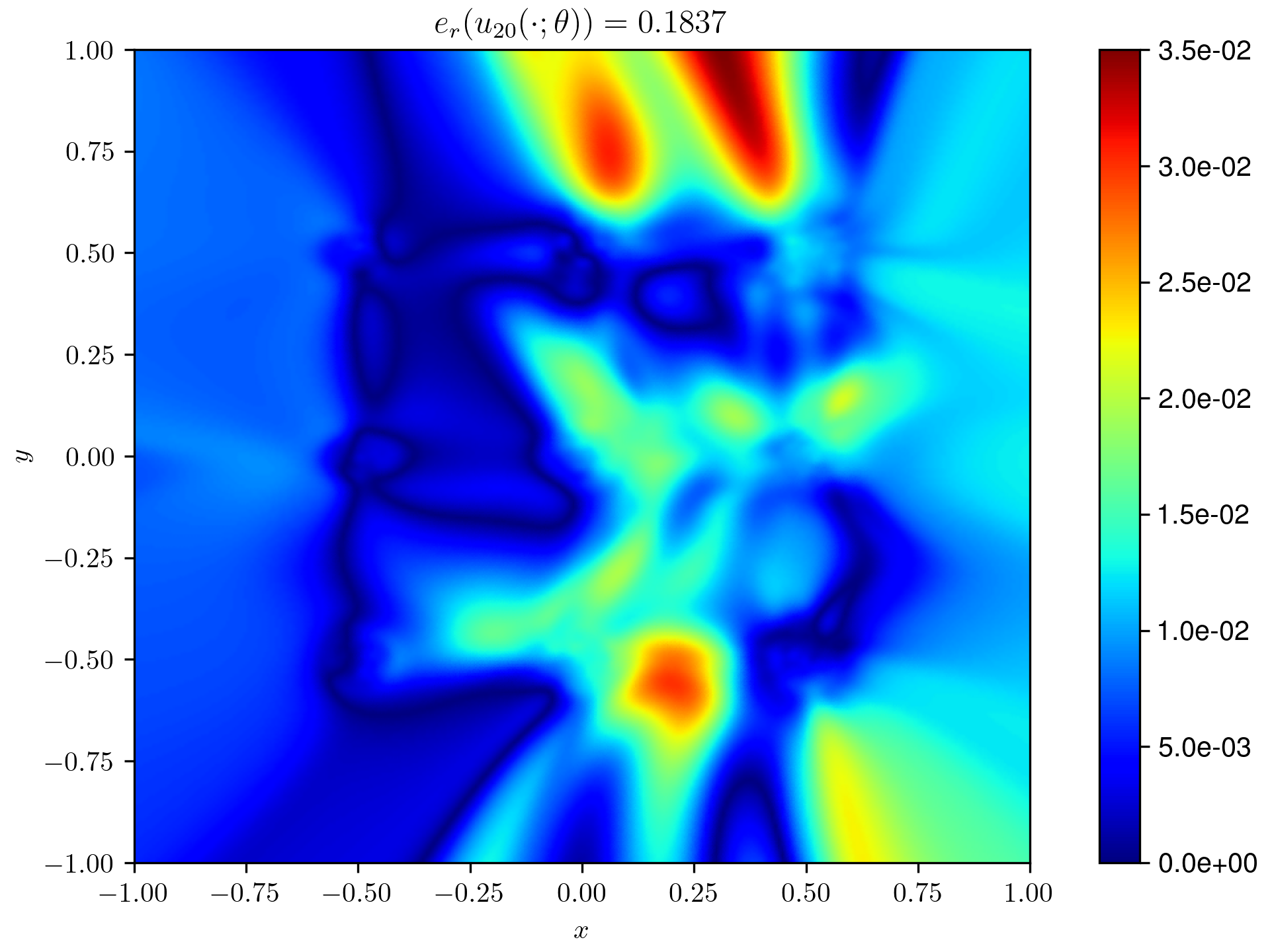}
        \caption{\textit{AAIS-g}}
    \end{subfigure}%
    \begin{subfigure}{.25\textwidth}
        \centering
        \includegraphics[height=0.75\textwidth,width=1.0\textwidth]{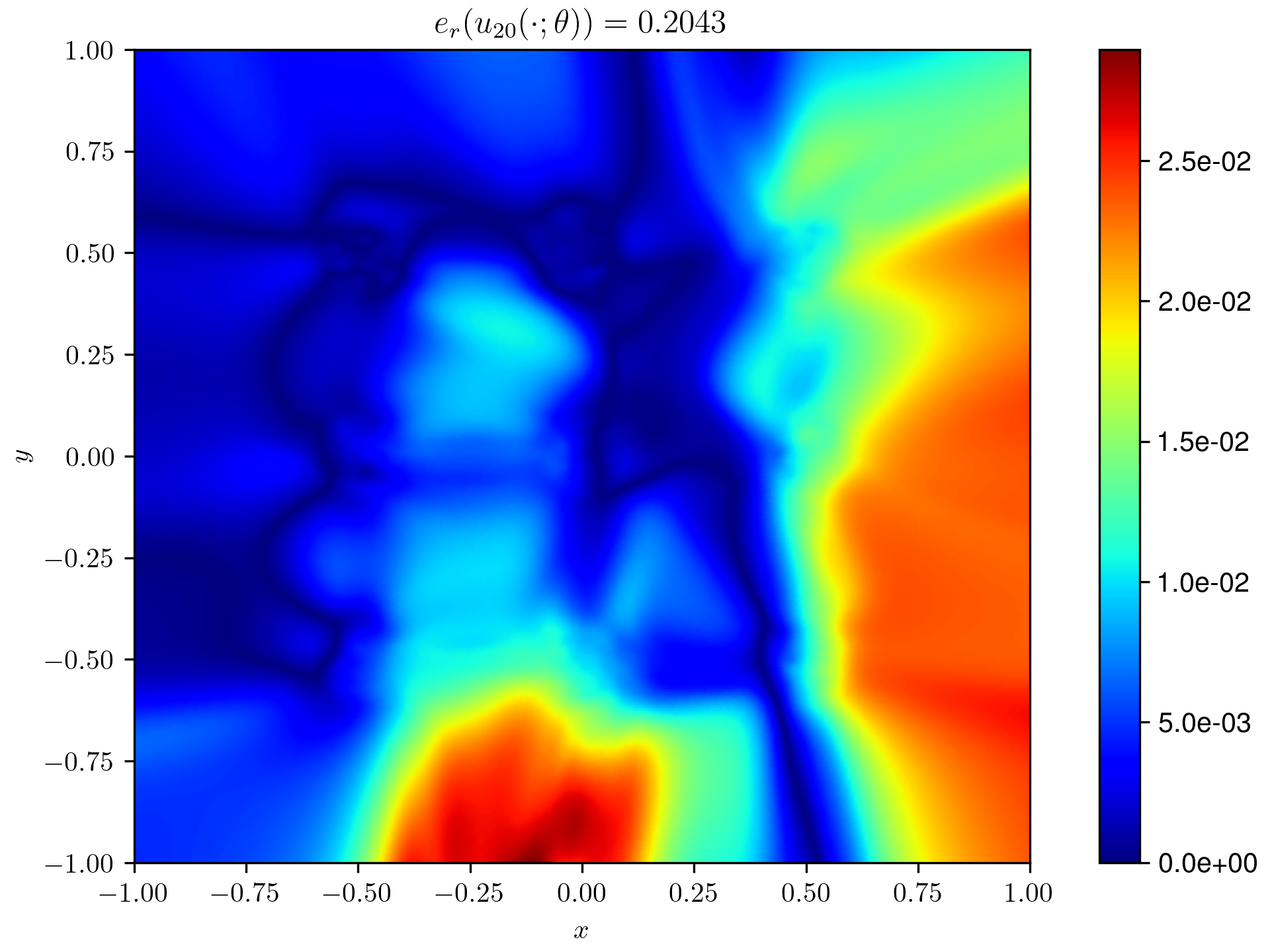}
        \caption{\textit{AAIS-t}}
    \end{subfigure}%
    \caption{Profiles of absolute error and neural network solutions for nine peak Poisson equation. First row: numeral solutions. Second row: absolute error.}
    \label{fig:Ps9PErr5000e}
\end{figure}

In Figure \ref{fig:Ps9PNode5000e} we show the residual $\Q$ and nodes respectively. Like one peak problem, the nodes of adaptive sampling methods focus around the high singularities. The residual $\Q$ would decay away from singularities with increasing frequency. The absolute error and solution profiles are listed in Figure \ref{fig:Ps9PErr5000e}, the singularities hide from the absolute error implying being well-trained there.

\subsection{High-dimensional Poisson problems}
In this part, we focus on the high dimensional multi-peaks Poisson problems. Above experiments are showed in two dimensions, \textit{RAD} method stands out of our proposed AAIS algorithm. The main reason is that the size of searching nodes for \textit{RAD} is 100k, i.e. $N_S=100000$, the simple Monte-Carlo method would obtain satisfactory results in 2D. However, in high dimensional problems, 100k searching nodes would be insufficient, and the Monte-Carlo method would fail. However, our proposed AAIS algorithm work due to the importance sampling. Therefore, in this subsection we choose \textit{RAD} and \textit{AAIS-t} to test $5, 9, 20$ dimensional Poisson problems with multi-peaks.

Given $\Omega = (-1,1)^d$ and the source term $f$ and boundary term $g$ of equations \eqref{pde:Poisson1Peak} are defined by the exact solution
\begin{equation*}
    u^*(x,y) = \sum_{i=1}^c\sum_{j=1^d}\exp\left[-K((x_j-x_j^i)^2)\right], (x_1,x_2,...,x_d)\in \Omega,
\end{equation*}
for 5D and 9D problems we let $K=100$ and 15D problems we let $K=10$.

Moreover, for high dimensional Poisson problems, it is difficult to efficiently compute numerical errors, in the following we firstly sample points uniformly in the domain then combine them with samplings from Gaussians, whose mean and covariance are determined by each part of the solution  leading by one of centers, then use these points to compute numerical errors.
\subsection{Five-dimension two-peaks problem}
\begin{figure}[htbp]
    \centering
    \includegraphics[scale=0.125]{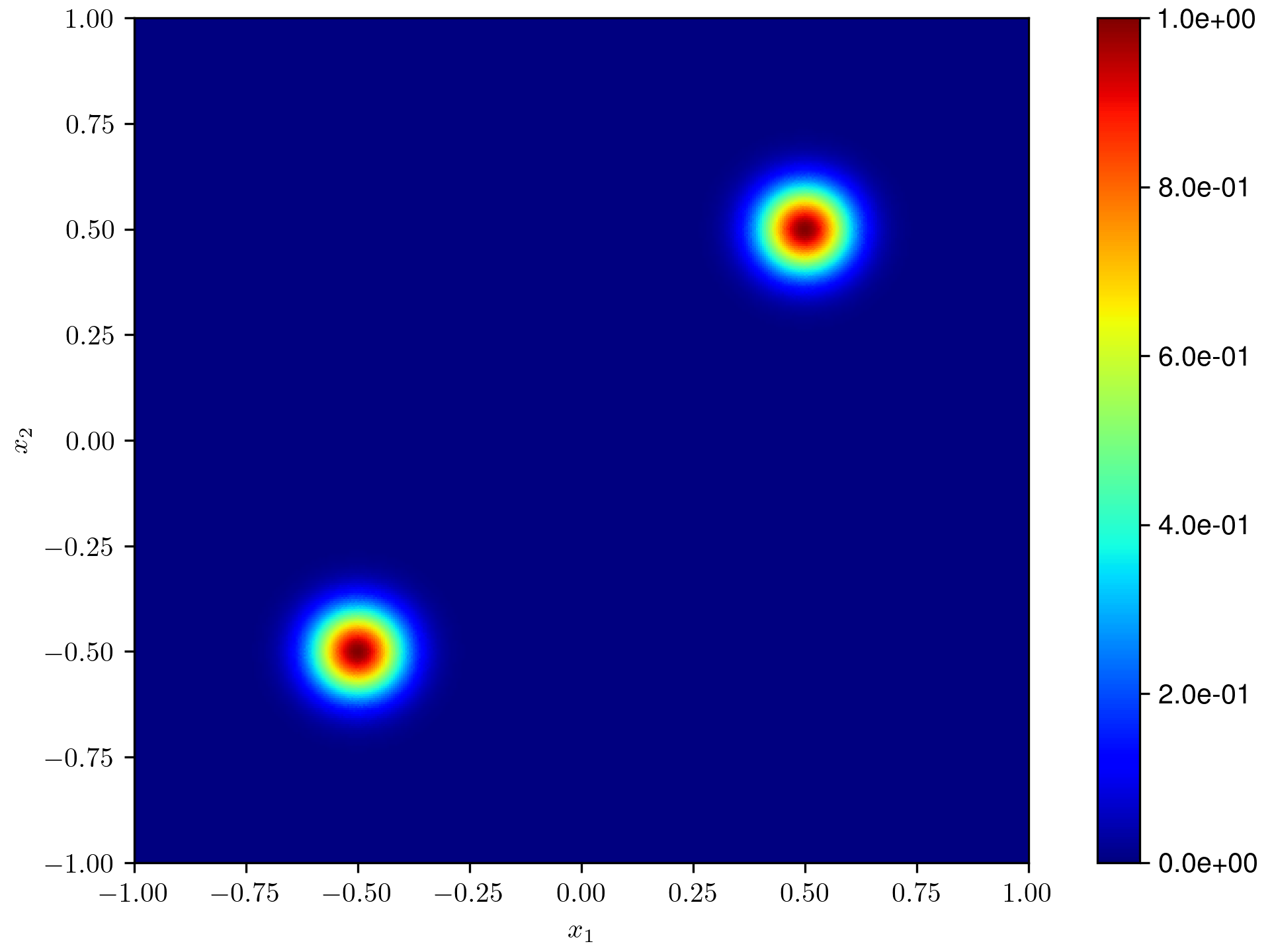}
    \caption{Projection on $x_1x_2$ plane of exact solution for Poisson 5D problems with two peaks in \eqref{pde:Poisson1Peak}.}
    \label{fig:PS5Dexact}
\end{figure}

We set the centers $(x_1^i,x_2^i,...x_5^i)=(0.5*(-1)^i, 0.5*(-1)^i,0,0,0)$ for $i=1,2$. The exact solution is showed in Figure \ref{fig:PS5Dexact}. For testing the relative errors, due to the limitation of storage, we uniformly sample 100k points in the domain from, 15k Gaussian samples for each mode and 10k points on the boundary. The neural network structure has 6 hidden layers with 64 neurons in each layer. The max iteration is $20$. The training schedule is Adam 500 epochs and  lbfgs of 2000 epochs for pre-training, Adam of 500 epochs and lbfgs of 10000 epochs for adaptive training.

In the following we plot the solution under the projection on $x_1x_2$-plane at the hypersurface $(x_1,x_2, 0,...,0)$.

\begin{figure}[htbp]
    \centering
    \begin{subfigure}{.5\textwidth}
        \centering
        \includegraphics[scale=0.25]{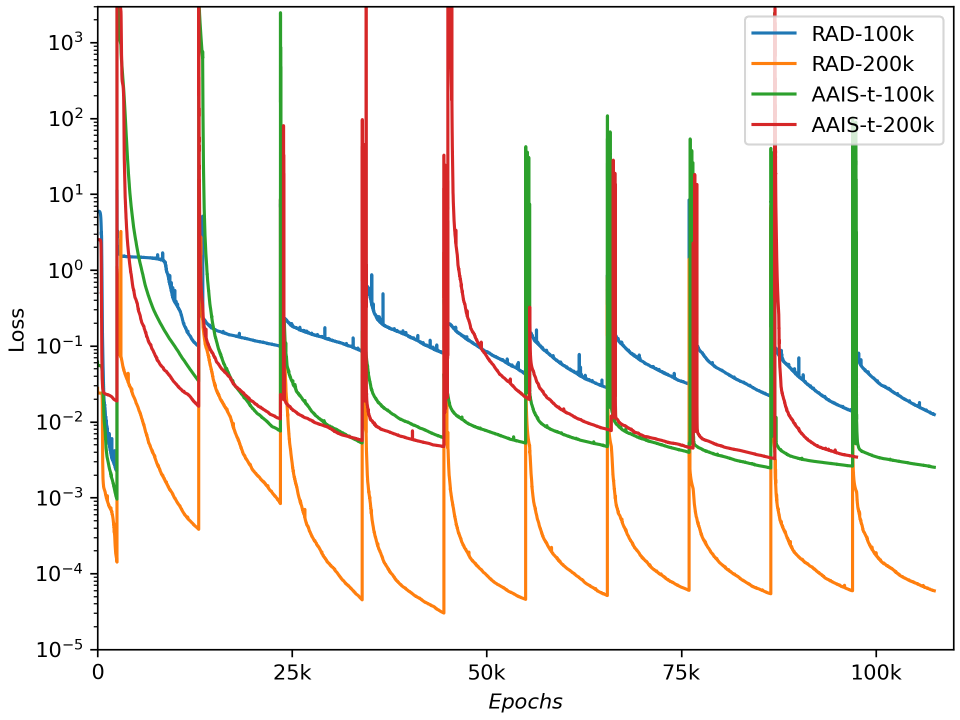}
        \end{subfigure}%
        \begin{subfigure}{.5\textwidth}
        \centering
        \includegraphics[scale=0.25]{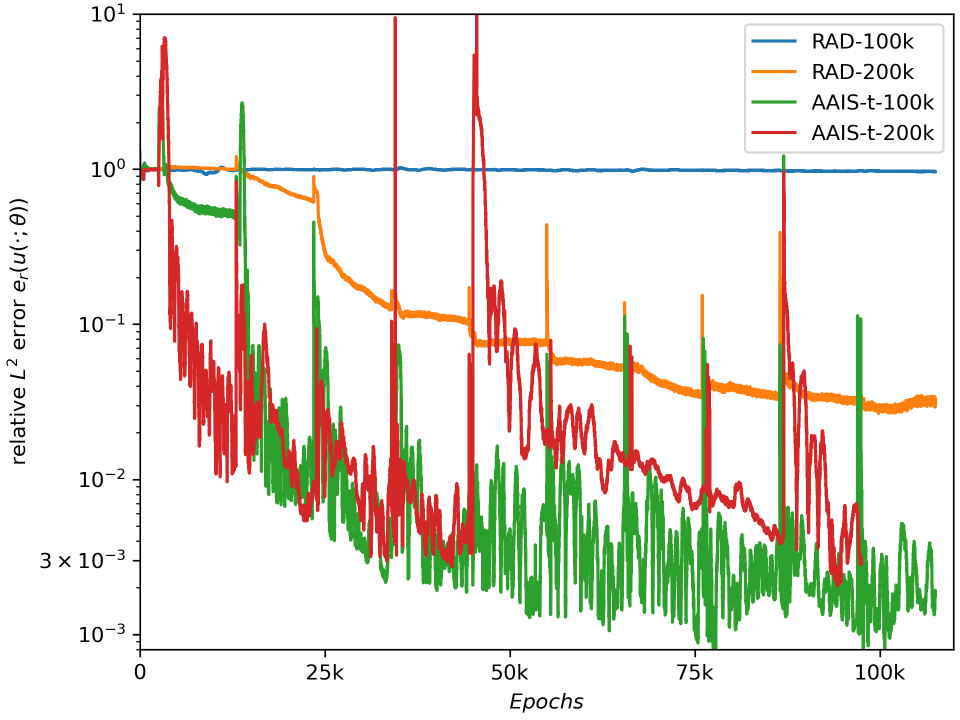}
        \end{subfigure}%
        \caption{Loss and relative errors for 5D Poisson equation. Left: the loss function. Right: the relative $L^2$ error $e_r(u(\cdot;\theta))$.
        }
    \label{fig:PS5DErr}
\end{figure}

\begin{figure}[htbp]
    \centering
    \begin{subfigure}{.25\textwidth}
        \centering
        \includegraphics[height=0.75\textwidth,width=1.0\textwidth]{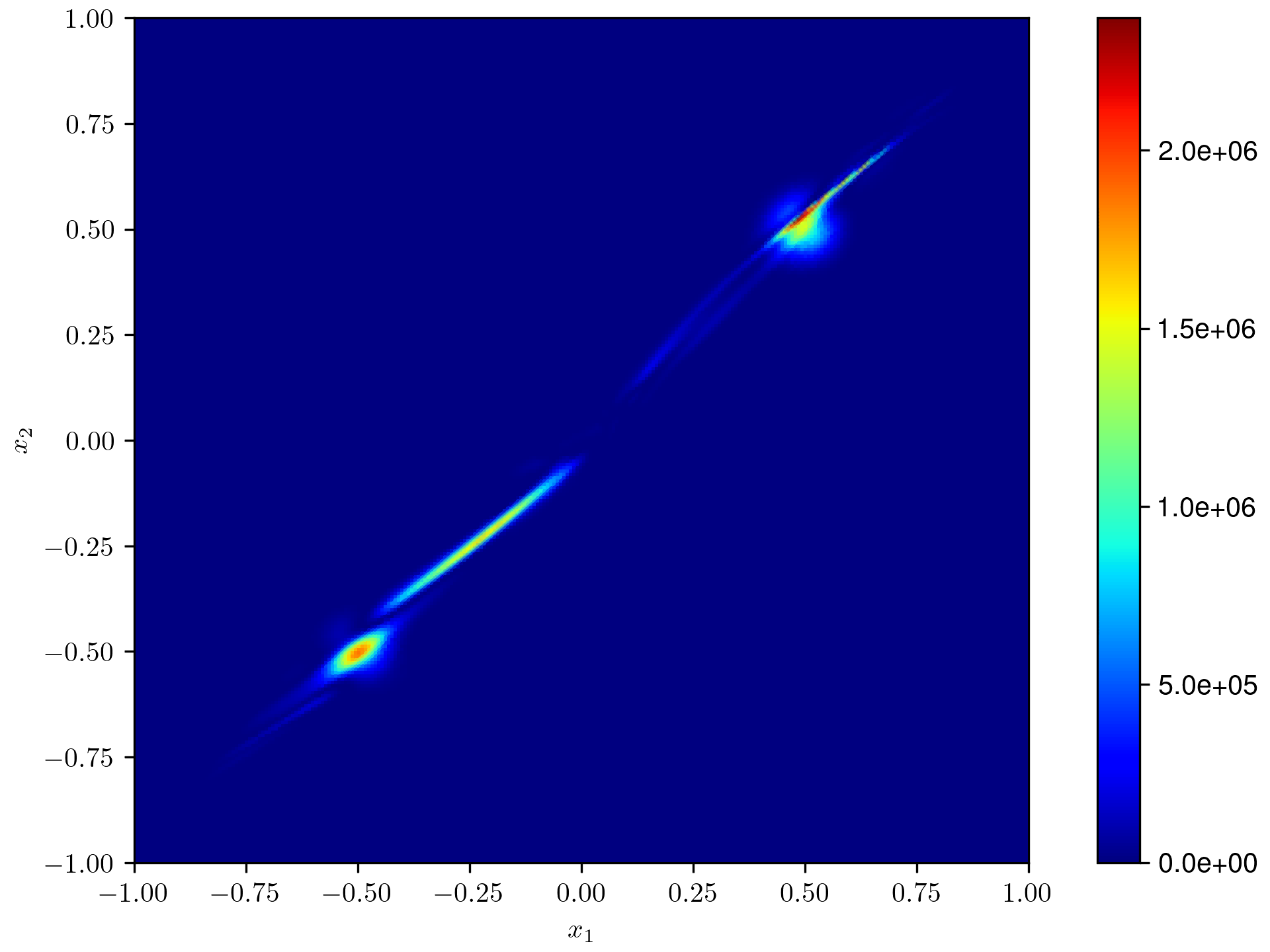}
    \end{subfigure}%
    \begin{subfigure}{.25\textwidth}
        \centering
        \includegraphics[height=0.75\textwidth,width=1.0\textwidth]{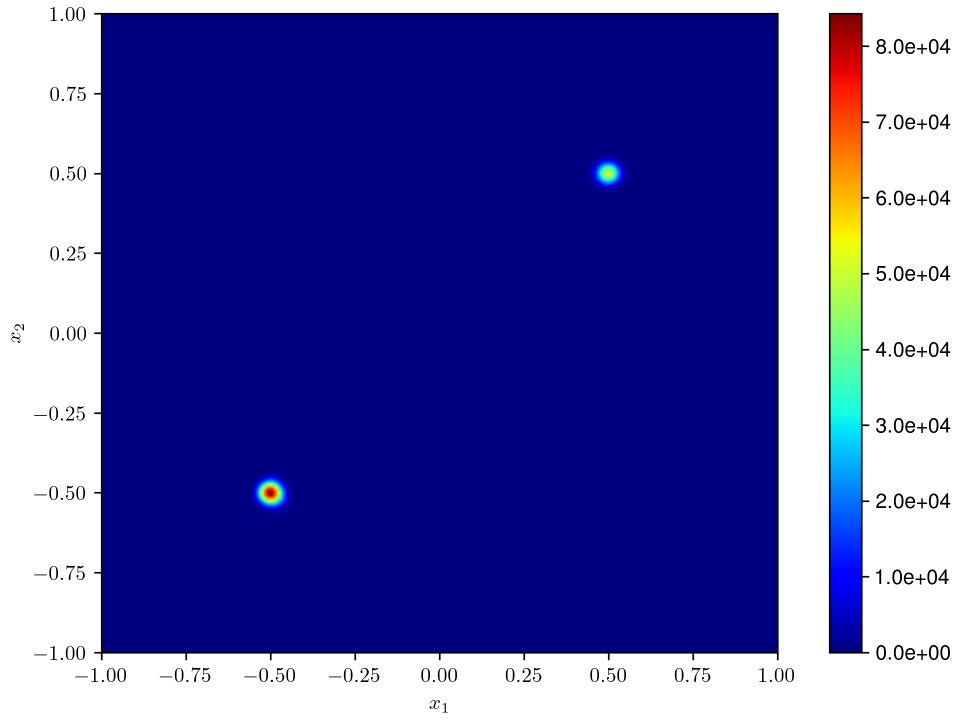}
    \end{subfigure}%
    \begin{subfigure}{.25\textwidth}
        \centering
        \includegraphics[height=0.75\textwidth,width=1.0\textwidth]{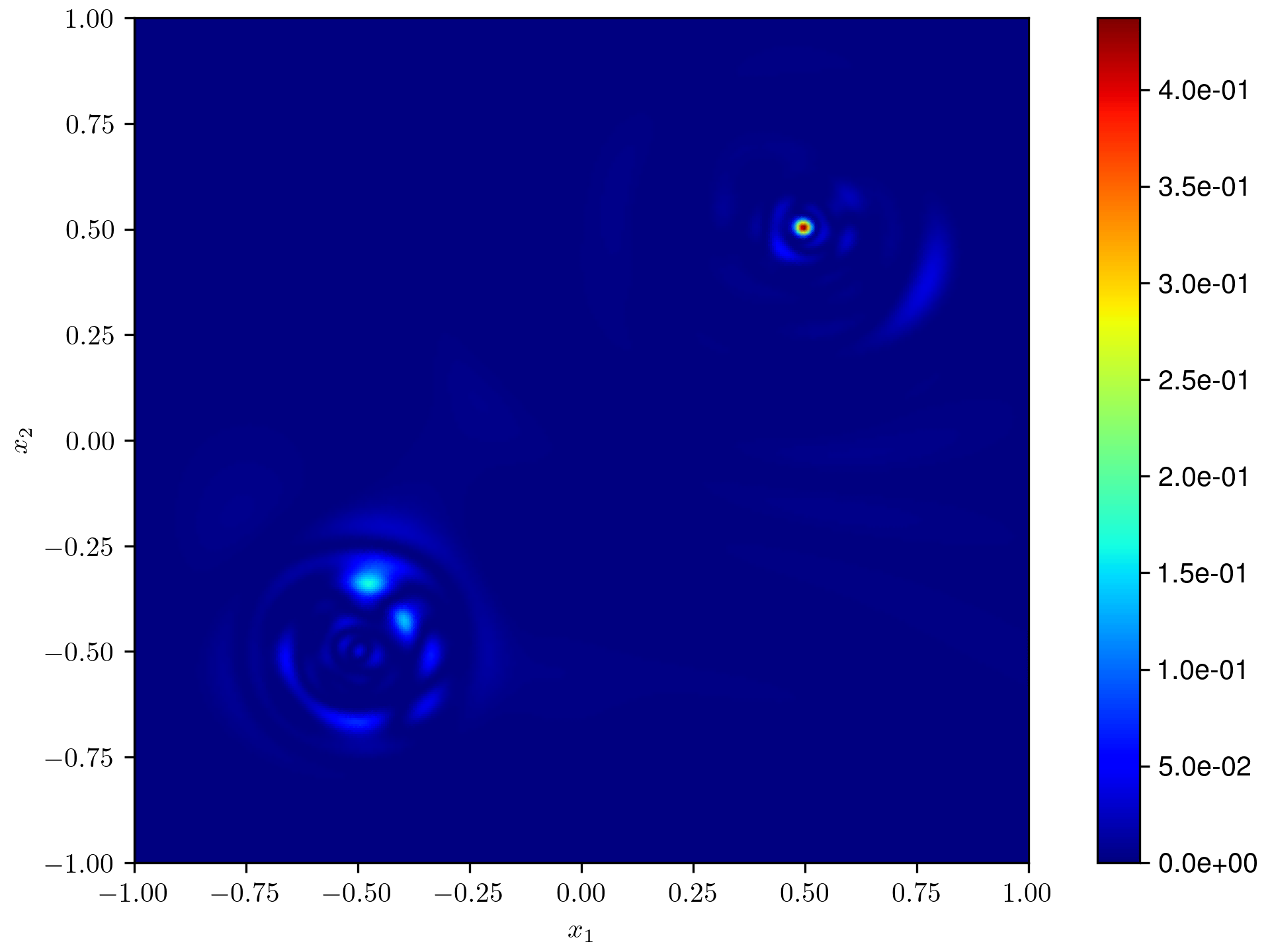}
    \end{subfigure}%
    \begin{subfigure}{.25\textwidth}
        \centering
        \includegraphics[height=0.75\textwidth,width=1.0\textwidth]{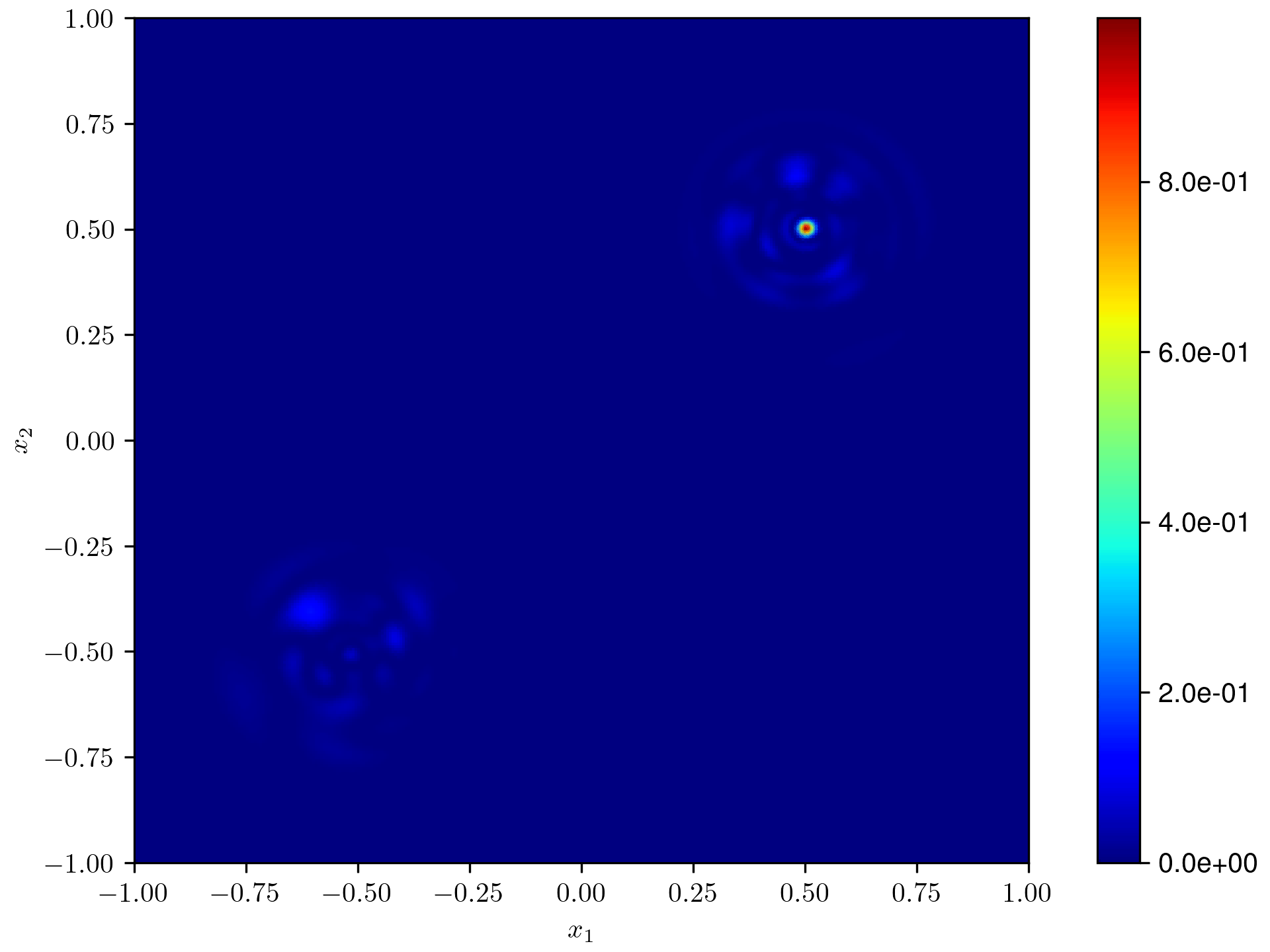}
    \end{subfigure}%
    \newline
    \raggedleft
    \begin{subfigure}{.25\textwidth}
        \centering
        \includegraphics[height=0.75\textwidth,width=1.0\textwidth]{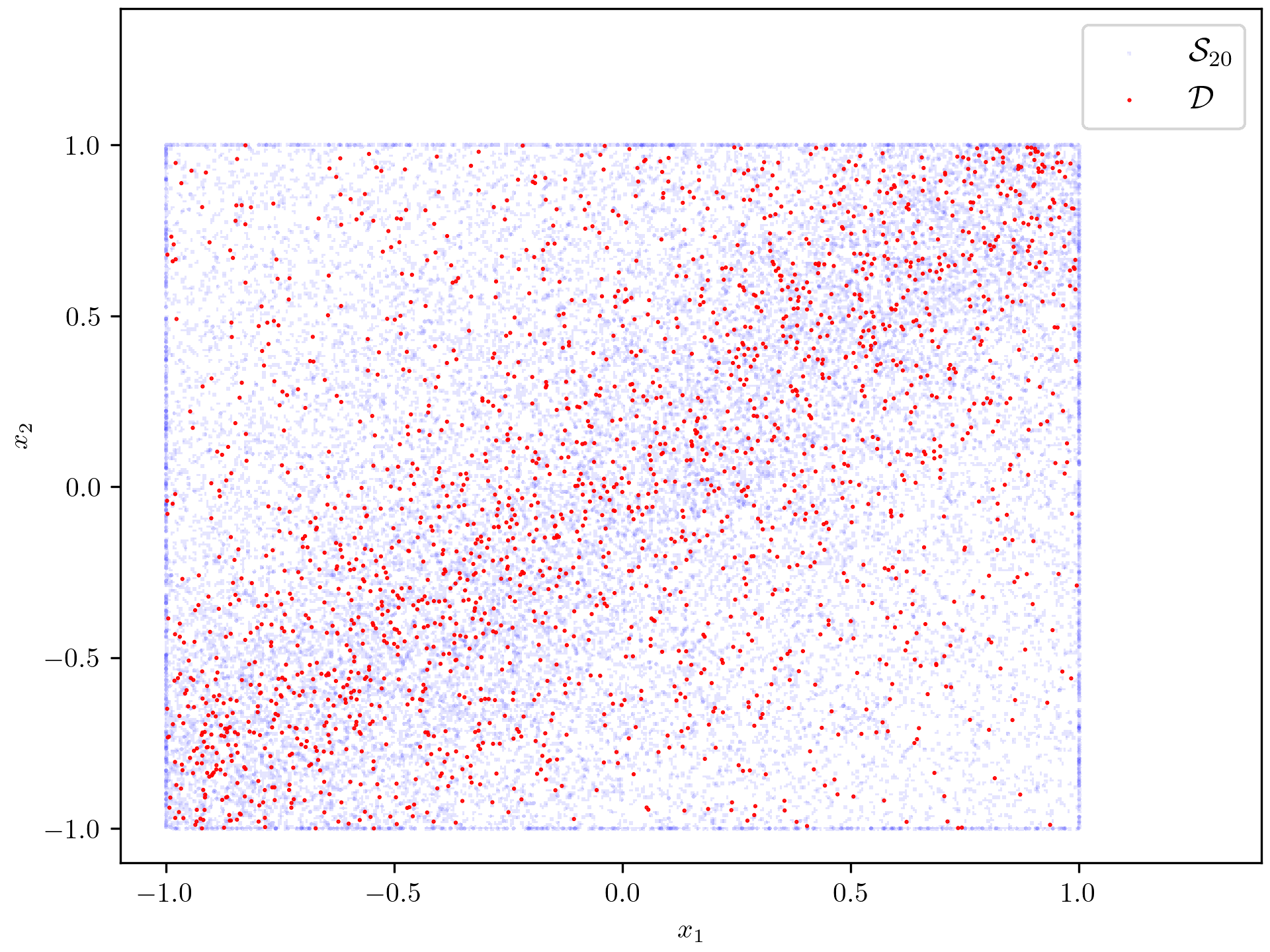}
        \caption{\textit{RAD}, $N_S=100k$.}
    \end{subfigure}%
    \begin{subfigure}{.25\textwidth}
        \centering
        \includegraphics[height=0.75\textwidth,width=1.0\textwidth]{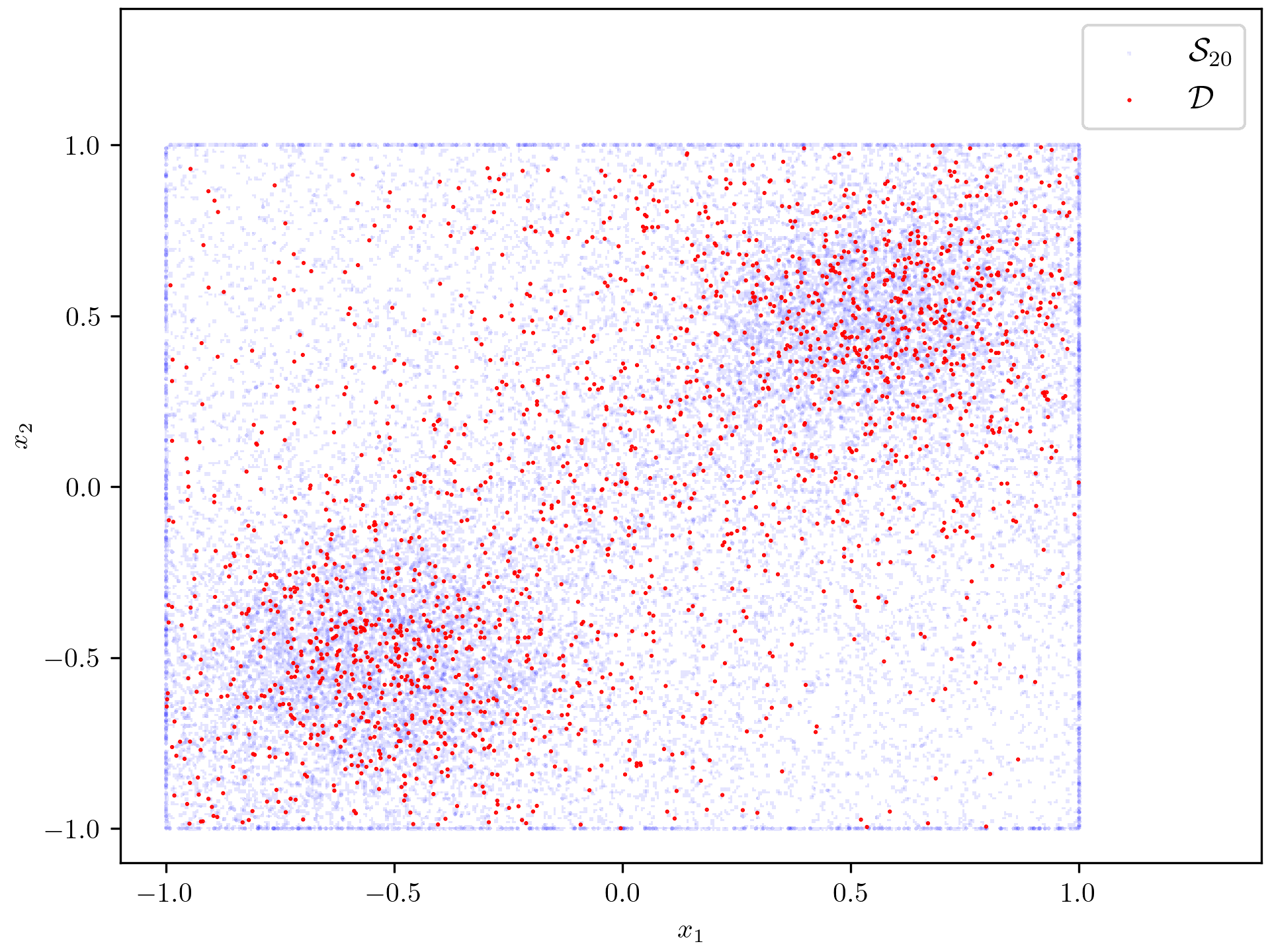}
        \caption{\textit{RAD}, $N_S=200k$.}
    \end{subfigure}%
    \begin{subfigure}{.25\textwidth}
        \centering
        \includegraphics[height=0.75\textwidth,width=1.0\textwidth]{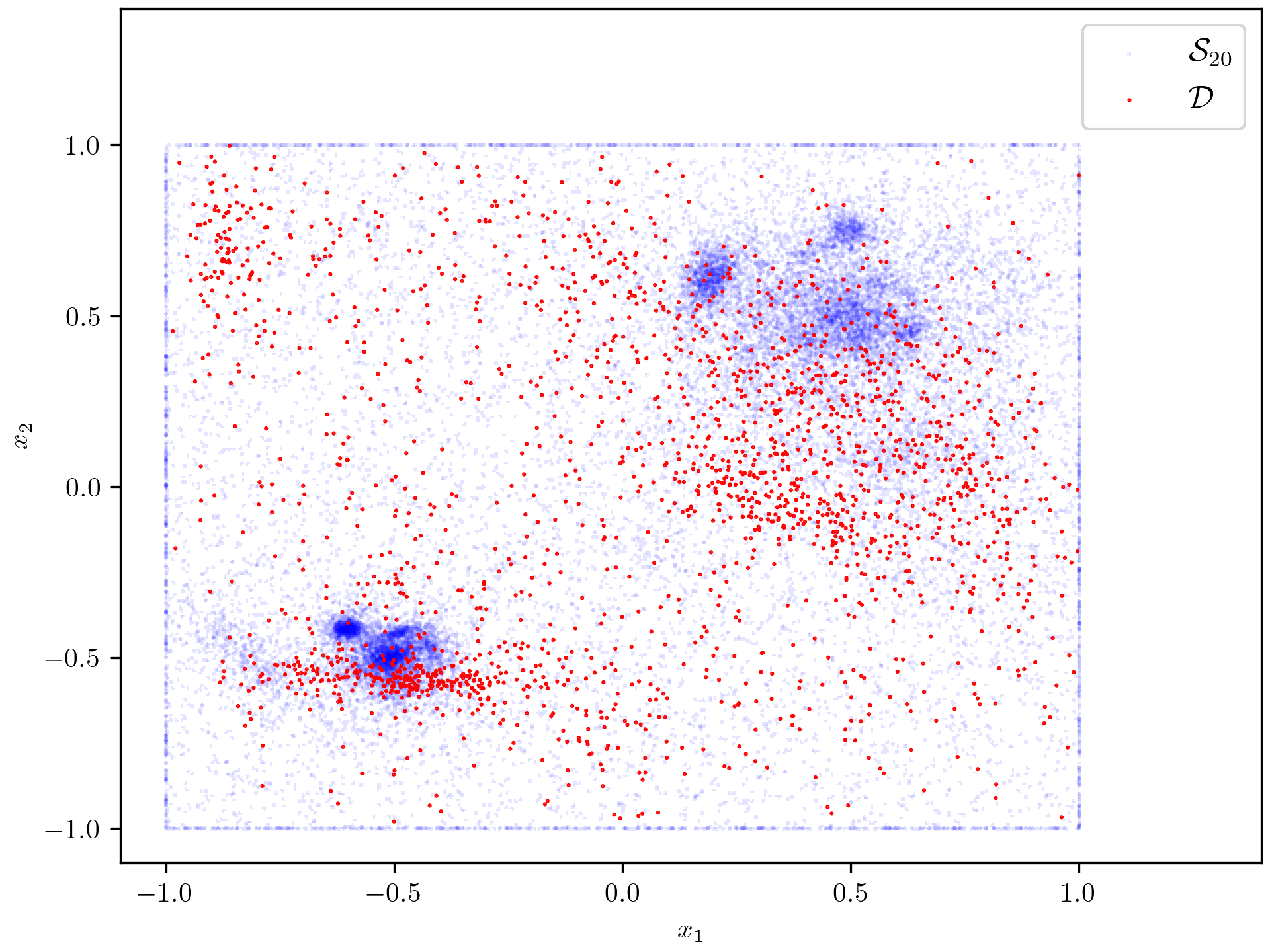}
        \caption{\textit{AAIS-t}, $N_S=100k$.}
    \end{subfigure}%
    \begin{subfigure}{.25\textwidth}
        \centering
        \includegraphics[height=0.75\textwidth,width=1.0\textwidth]{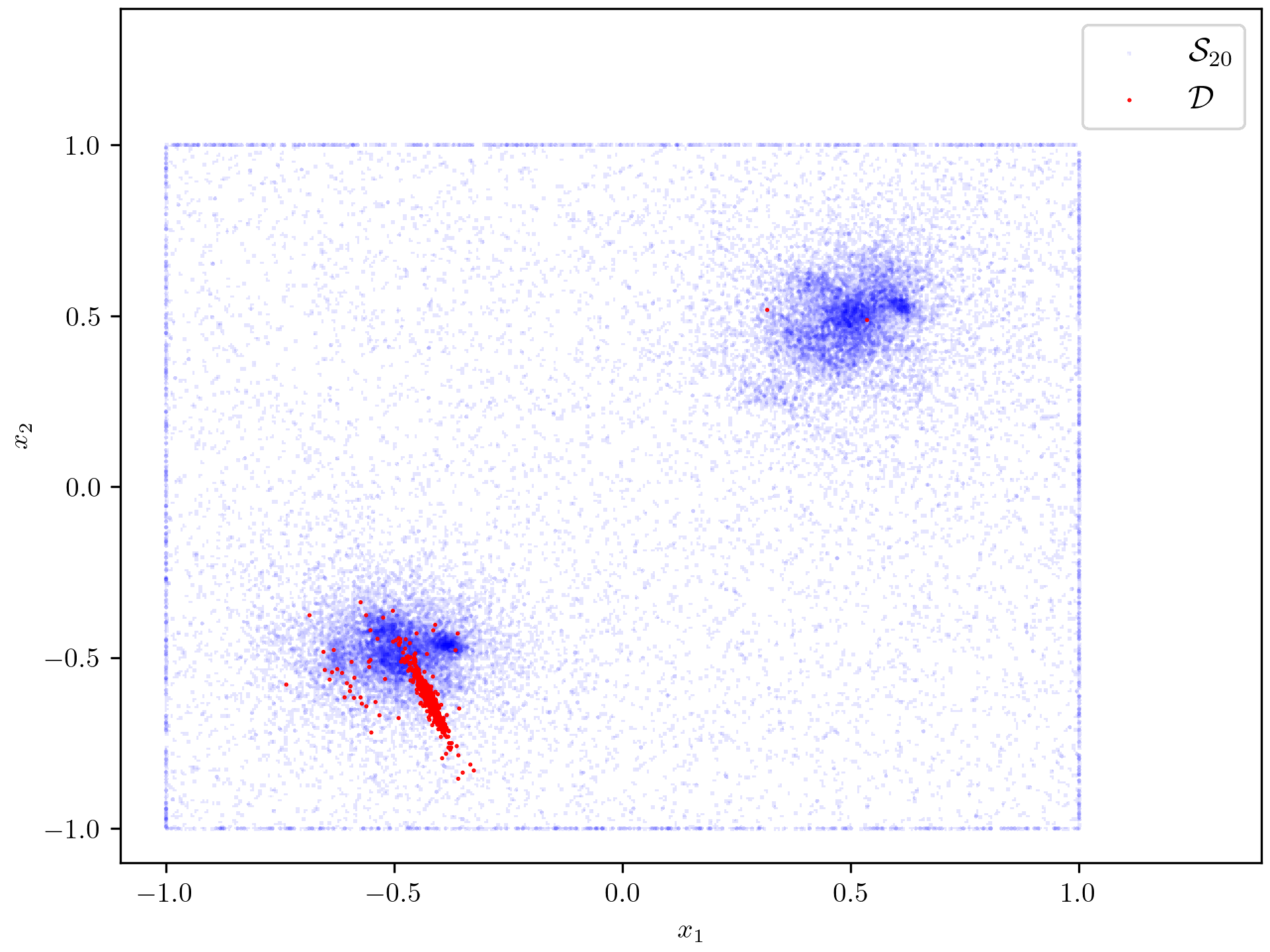}
        \caption{\textit{AAIS-t}, $N_S=200k$.}
    \end{subfigure}%
    \caption{$x_1x_2$-plane profiles of residual and nodes for 5D Poisson problem at 20-th iteration.}
    \label{fig:Ps5DNode}
\end{figure}

$N_S$ for both \textit{RAD} and \textit{AAIS-t} is 100k and 200k. The loss and errors are showed in Figure \ref{fig:PS5DErr}. The \textit{RAD} method fails to solve the 5D problem with 100k searching points due to the poor searching ability, but for 200k searching points \textit{RAD} could solve limitedly. \textit{AAIS-t} could solve the problem more accurately from importance sampling.
Figure \ref{fig:Ps5DNode} shows the loss function and node distributions at the last iteration. The \textit{RAD} method makes the nodes cluster around the singularities, and this effect becomes more pronounced as the number of searching nodes increases. This may explain why the \textit{RAD} method fails with 100k searching points but finds the solution with 200k searching points. For \textit{AAIS-t} method, firstly we could see that nodes focus on the singularities which could lead to a better solution, but since the number of searching points is still relatively small, AAIS algorithm may fail to mimic the loss function properly when the frequency of the loss become higher.

\begin{figure}[htbp]
    \centering
    \begin{subfigure}{.25\textwidth}
        \centering
        \includegraphics[height=0.75\textwidth,width=1.0\textwidth]{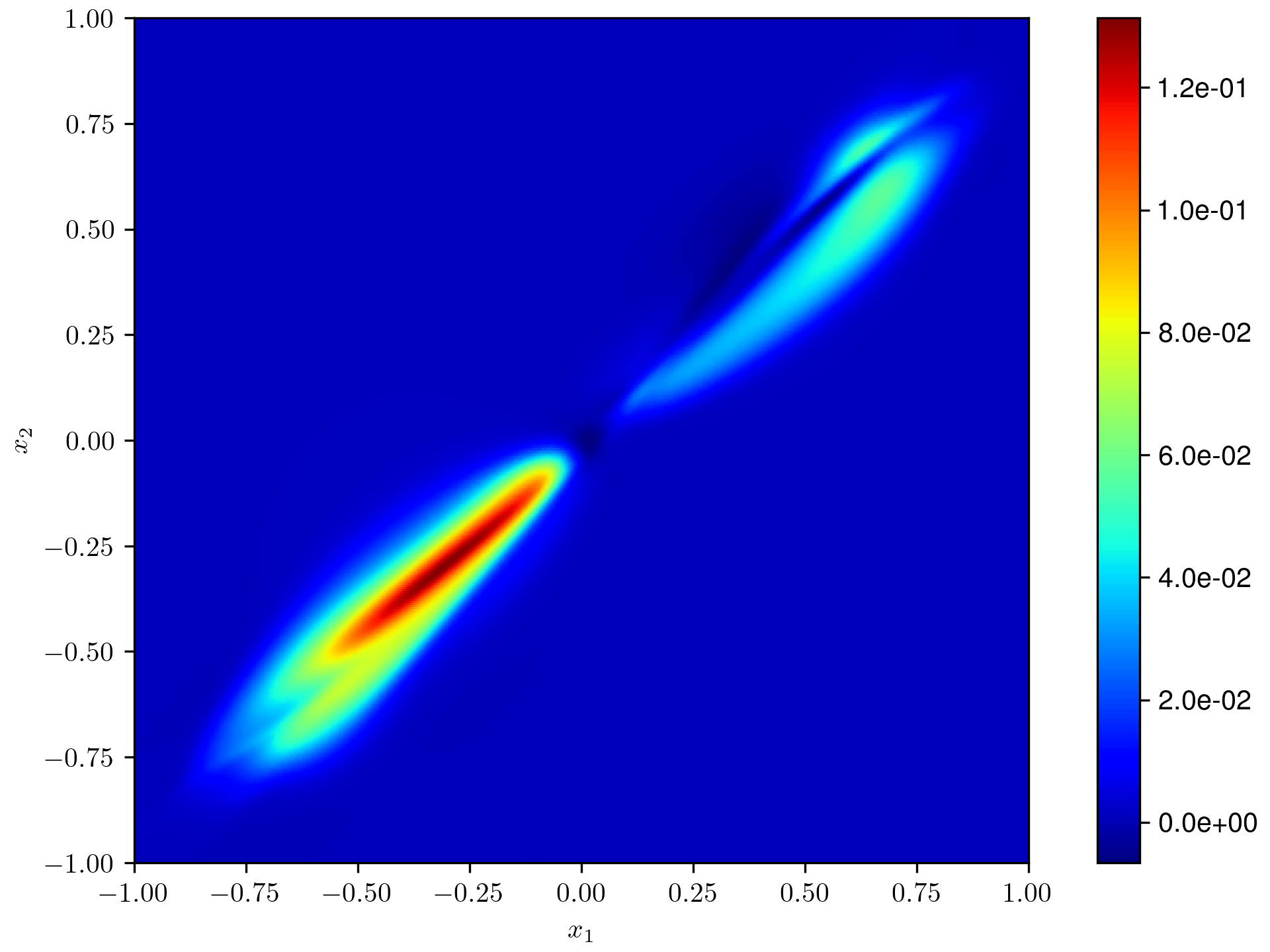}
    \end{subfigure}%
    \begin{subfigure}{.25\textwidth}
        \centering
        \includegraphics[height=0.75\textwidth,width=1.0\textwidth]{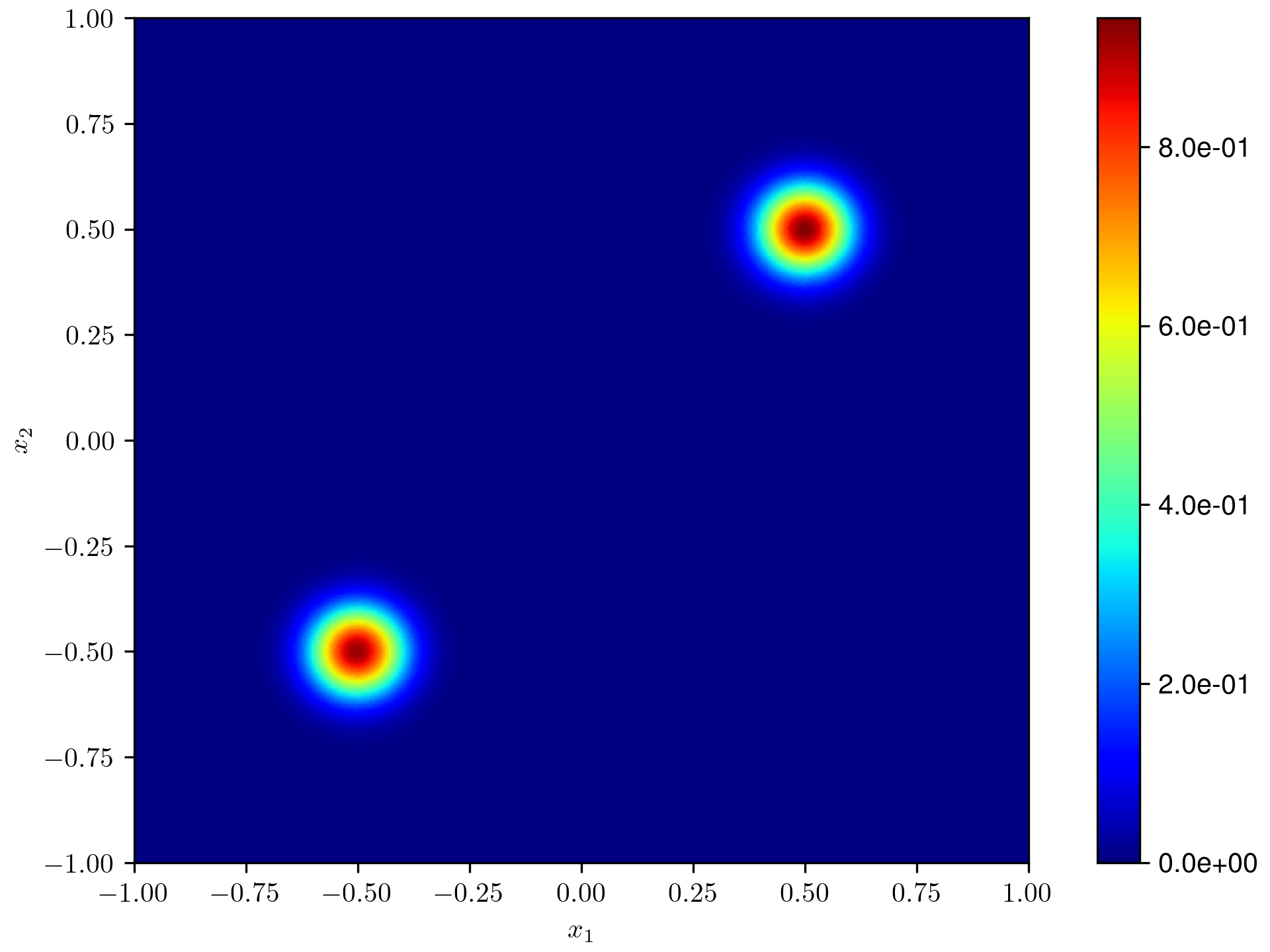}
    \end{subfigure}%
    \begin{subfigure}{.25\textwidth}
        \centering
        \includegraphics[height=0.75\textwidth,width=1.0\textwidth]{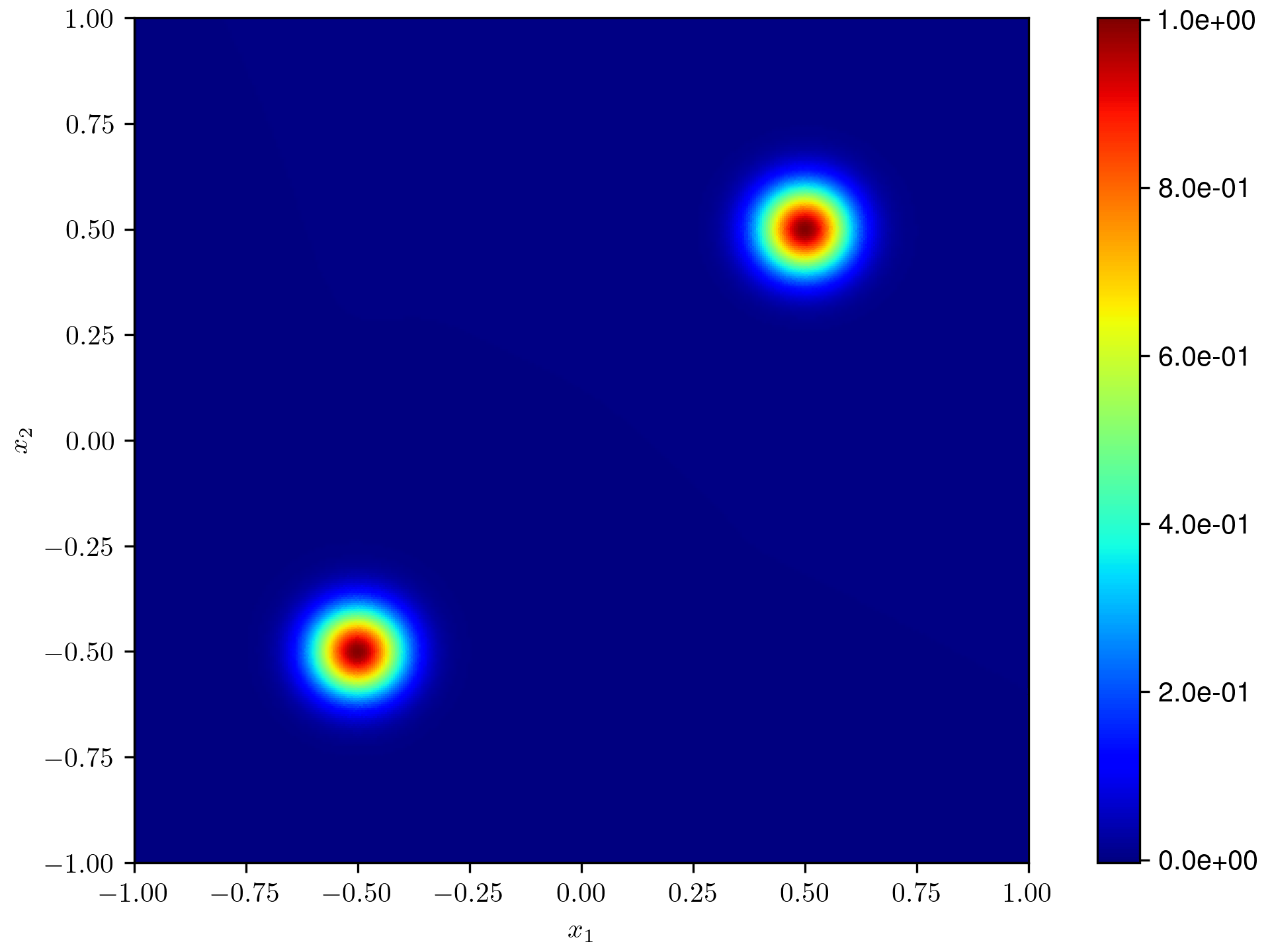}
    \end{subfigure}%
    \begin{subfigure}{.25\textwidth}
        \centering
        \includegraphics[height=0.75\textwidth,width=1.0\textwidth]{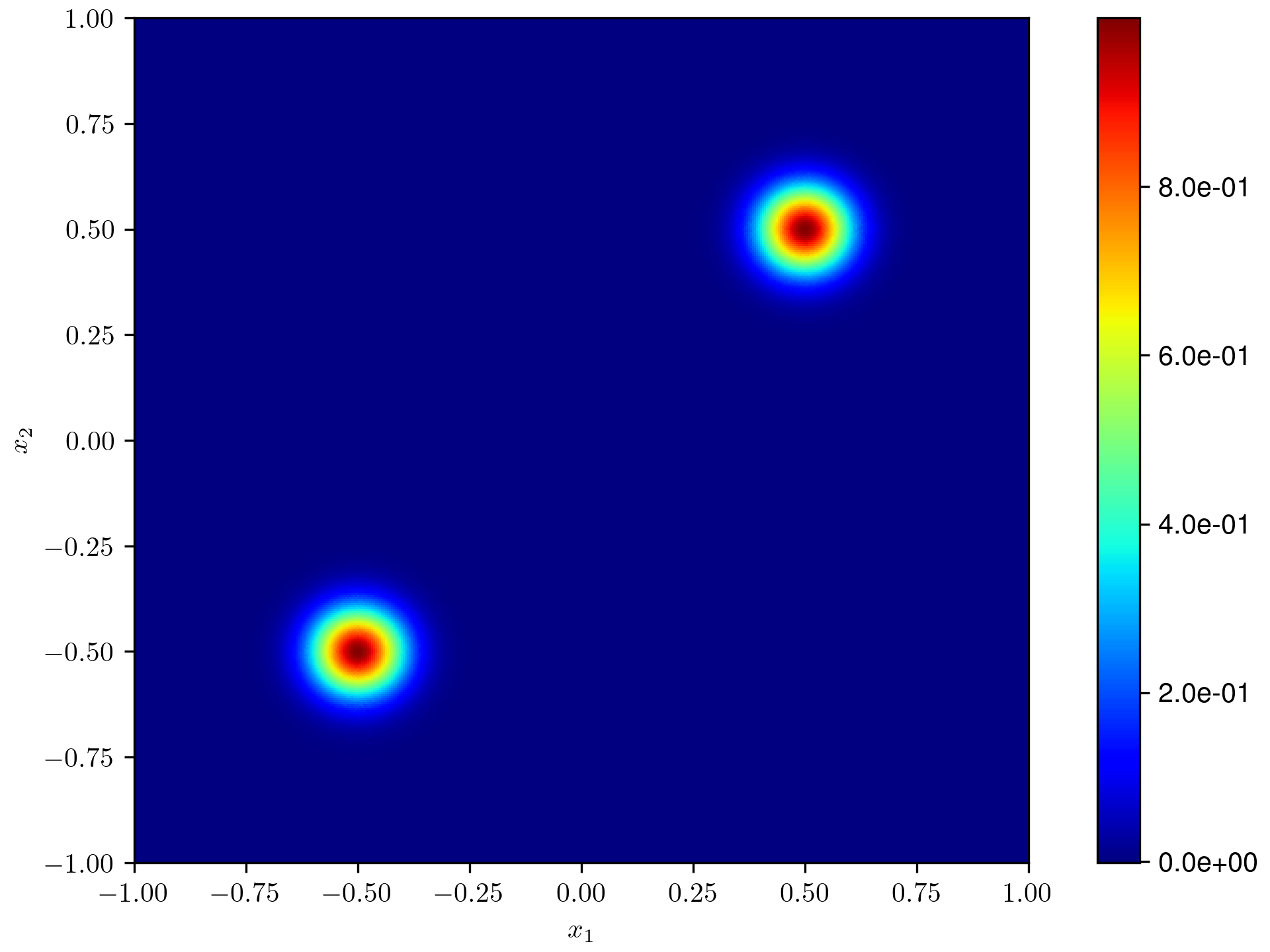}
    \end{subfigure}%
    \newline
    \raggedleft
    \begin{subfigure}{.25\textwidth}
        \centering
        \includegraphics[height=0.75\textwidth,width=1.0\textwidth]{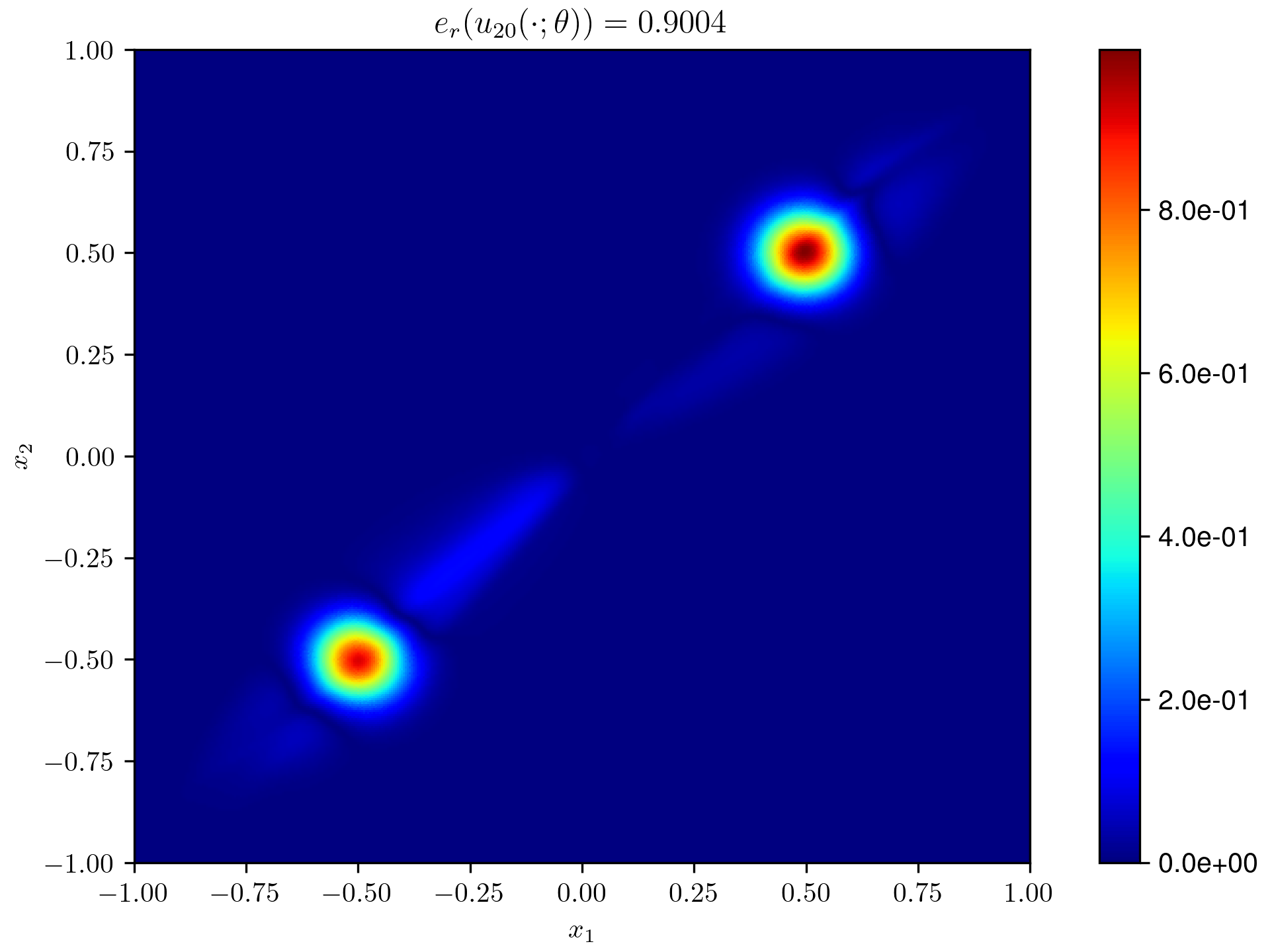}
        \caption{\textit{RAD}, $N_S=100k$.}
    \end{subfigure}%
    \begin{subfigure}{.25\textwidth}
        \centering
        \includegraphics[height=0.75\textwidth,width=1.0\textwidth]{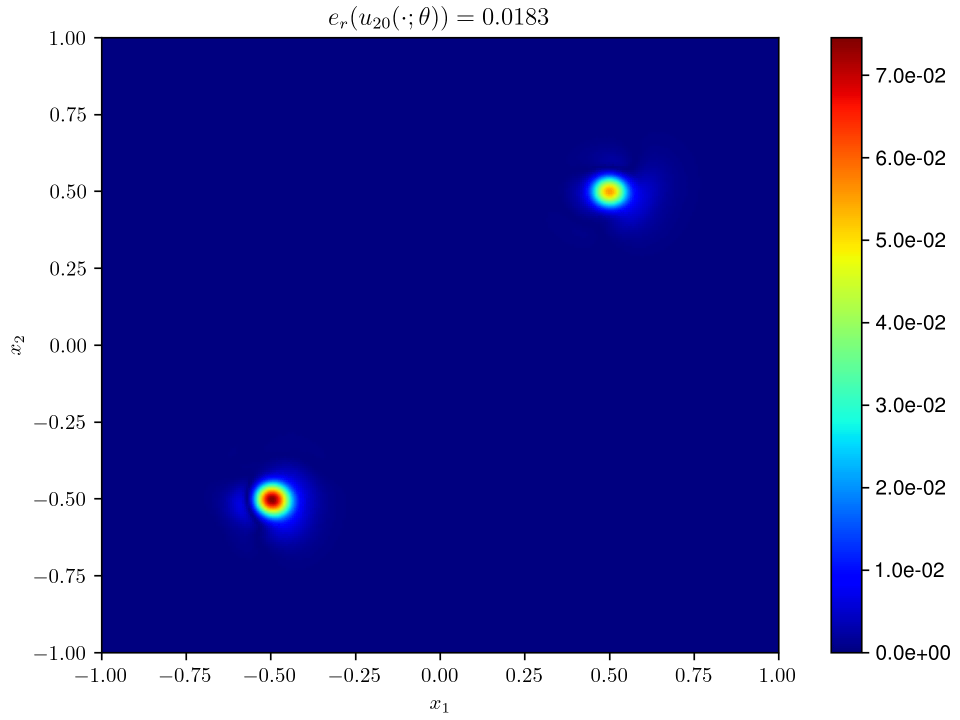}
        \caption{\textit{RAD}, $N_S=200k$.}
    \end{subfigure}%
    \begin{subfigure}{.25\textwidth}
        \centering
        \includegraphics[height=0.75\textwidth,width=1.0\textwidth]{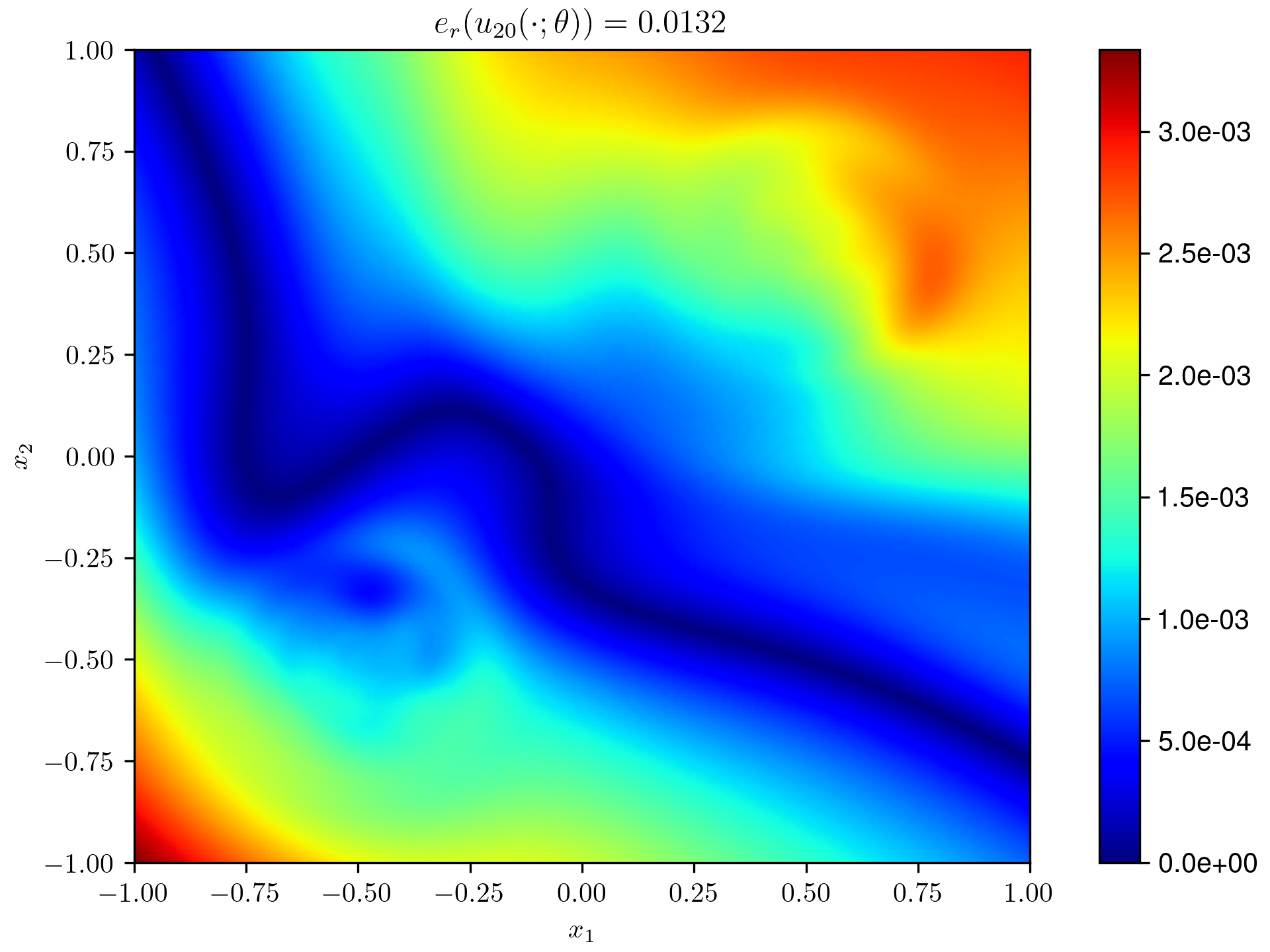}
        \caption{\textit{AAIS-t}, $N_S=100k$.}
    \end{subfigure}%
    \begin{subfigure}{.25\textwidth}
        \centering
        \includegraphics[height=0.75\textwidth,width=1.0\textwidth]{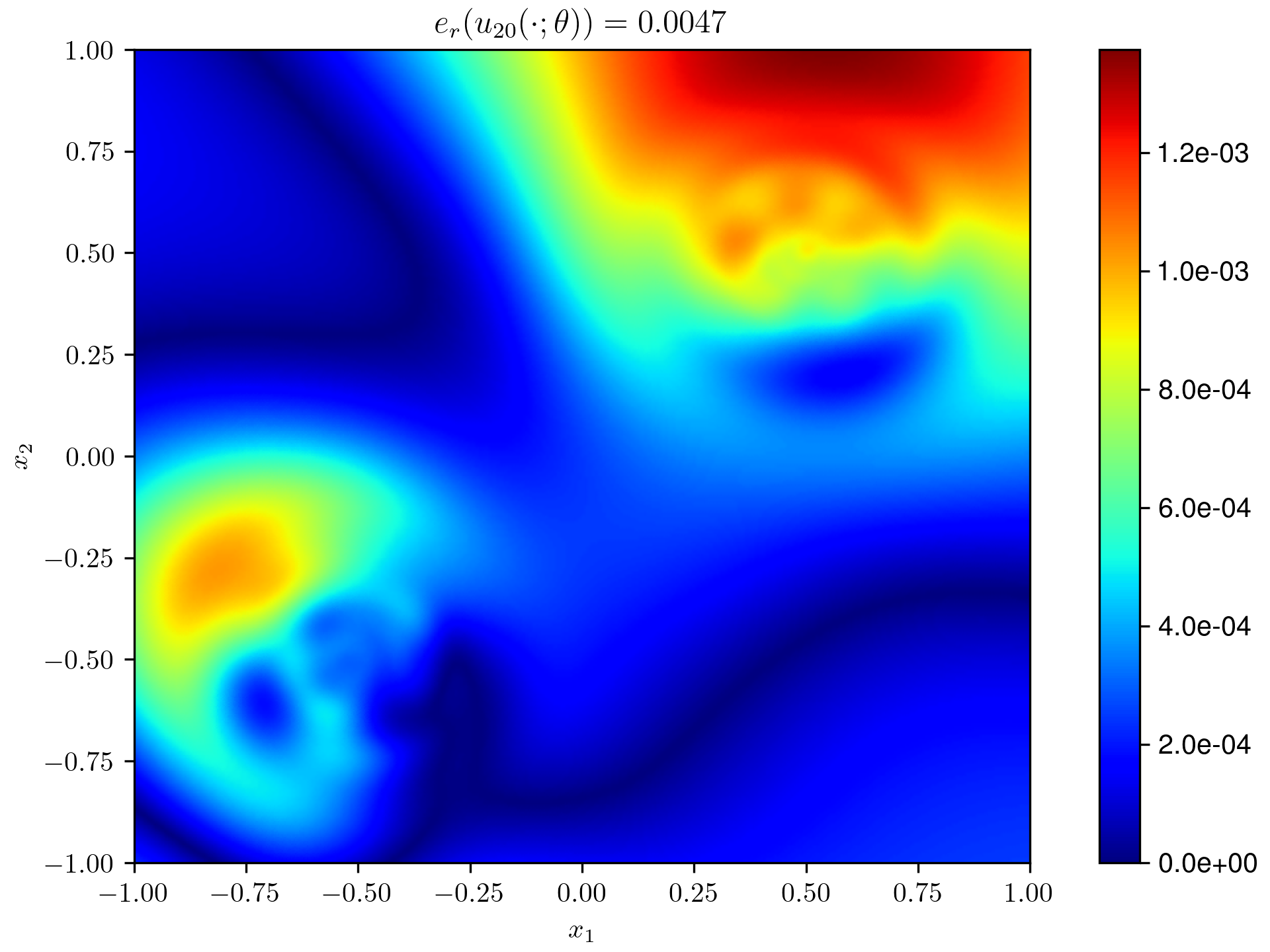}
        \caption{\textit{AAIS-t}, $N_S=200k$.}
    \end{subfigure}%
    \caption{Profiles of absolute error and neural network solutions for 5D Poisson equation. First row: numeral solutions. Second row: absolute error.}
    \label{fig:Ps5Dsol}
\end{figure}

Moreover, in Figure \ref{fig:Ps5Dsol}, we could see that when $N_S=100000$, \textit{RAD} could not solve the problem but \textit{AAIS-t} could. For $N_S=200000$, it could be seen that  both two adaptive methods succeed to solve the problem because of more focus on the singularities. However, \textit{RAD} method could not solve very well, reflected in the low frequency of absolute error, which implies that the solution does not fit very well at the singularities. For \textit{AAIS-t} method, similarly as the solution behaviors in 2D problem, the frequency of the absolute error is relatively high, and the singularities hide from the error, implying the better solvability of our proposed AAIS algorithms.

\subsubsection{Nine-dimensional two peaks problem}
\begin{figure}[htbp]
    \centering
    \includegraphics[scale=0.125]{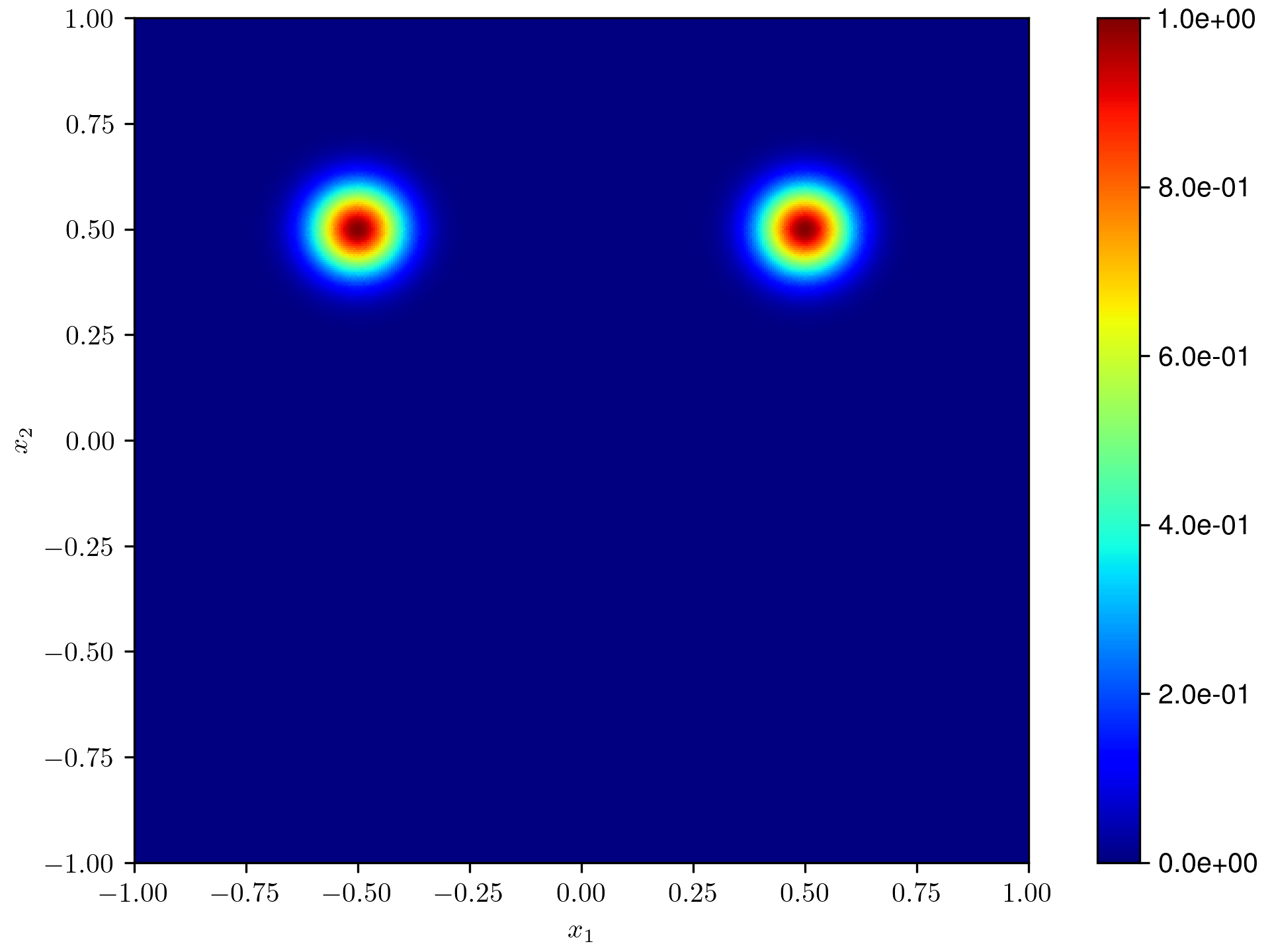}
    \caption{Projection on $x_1x_2$ plane of exact solution for Poisson 9D problems with two peaks in \eqref{pde:Poisson1Peak}.}
    \label{fig:PS9Dexact}
\end{figure}

In this part we consider a nine dimensional Poisson problems with two peaks where the centers are $(\pm 0.5, 0.5, 0,...0)$. The exact solution is showed in Figure \ref{fig:PS9Dexact}. We let $N_{in}=20000$, $\D=2000$, $N_b=5000$ in Algorithm \ref{alg:PINNSamplingResample}. The neural network structure has 64 neurons with 6 hidden layers. We train the PINNs with 500 epochs Adam, 2000 epochs lbfgs in the pre-training and 500 epochs Adam, 10000 lbfgs in the adaptive training.

\begin{figure}[htbp]
    \centering
    \begin{subfigure}{.5\textwidth}
        \centering
        \includegraphics[scale=0.25]{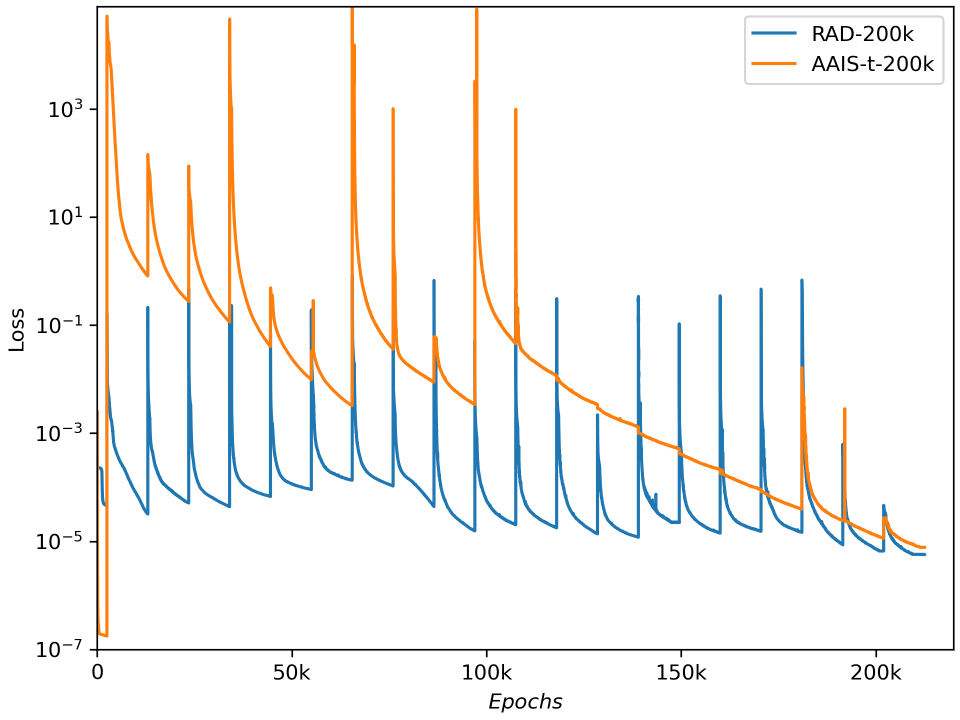}
        \end{subfigure}%
        \begin{subfigure}{.5\textwidth}
        \centering
        \includegraphics[scale=0.25]{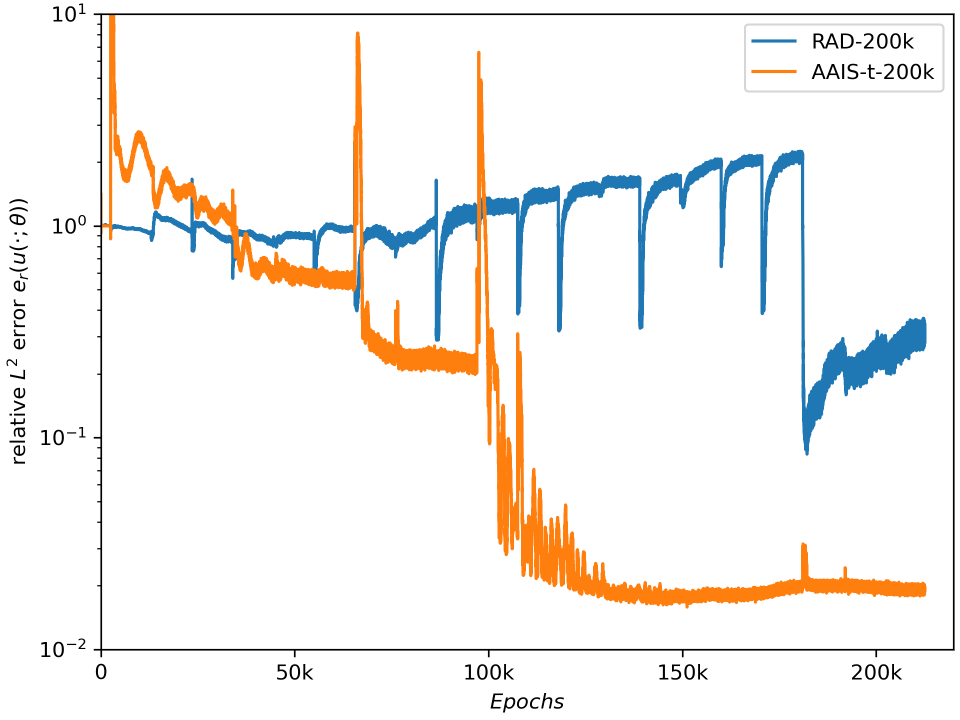}
        \end{subfigure}%
        \caption{Loss and relative errors for 9D Poisson equation. Left: the loss function. Right: the relative $L^2$ error $e_r(u(\cdot;\theta))$.
        }
    \label{fig:PS9DErr}
\end{figure}

Here we only choose $N_S=200000$, the loss and relative errors are showed in Figure \ref{fig:PS9DErr}. 
The \textit{RAD} method fails to solve the 9D problem but \textit{AAIS-t} would solve the problem accurately.

\begin{figure}[htbp]
    \centering
    \begin{subfigure}{.25\textwidth}
        \centering
        \includegraphics[height=0.75\textwidth,width=1.0\textwidth]{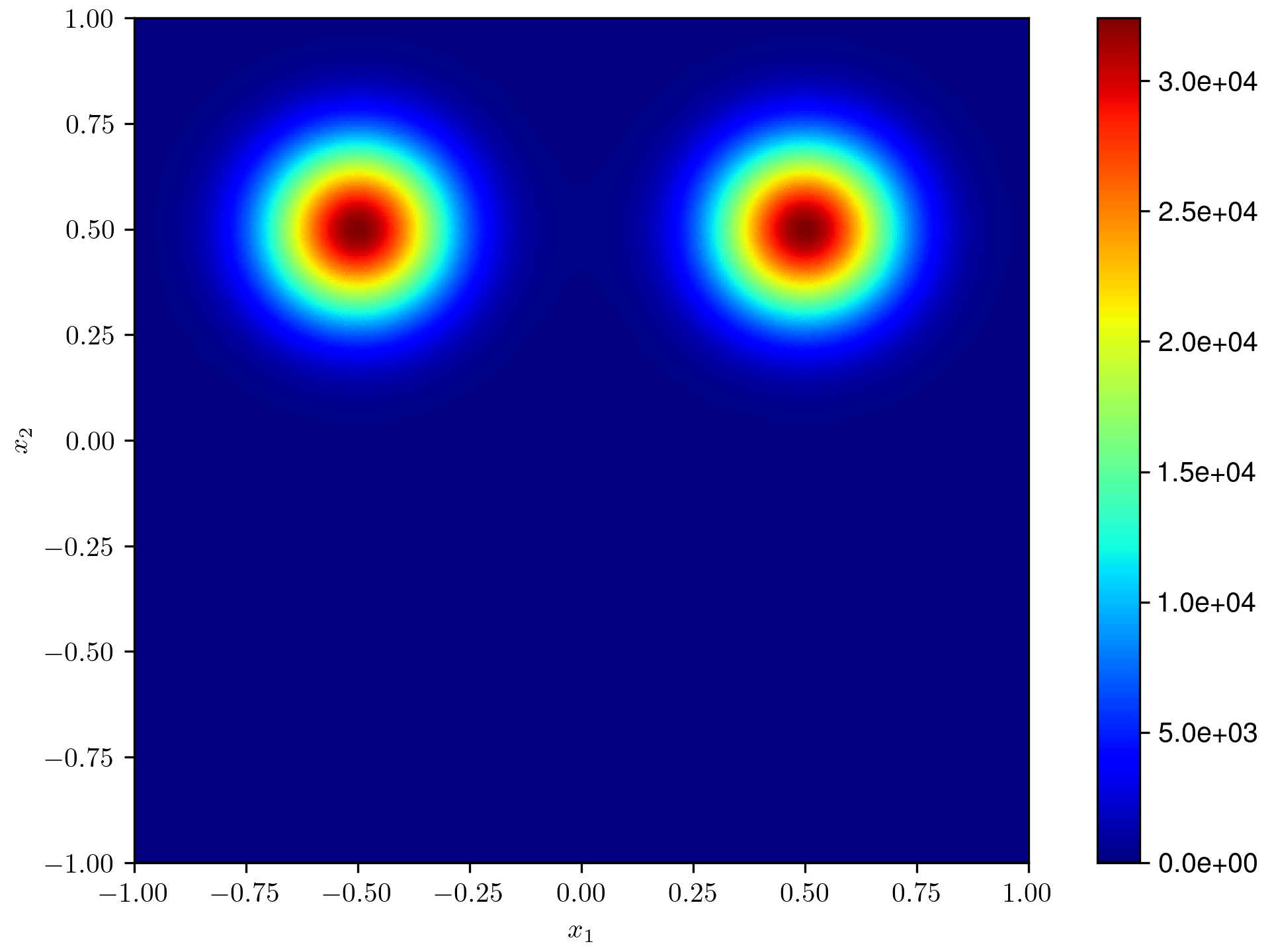}
      \end{subfigure}%
    \begin{subfigure}{.25\textwidth}
        \centering
        \includegraphics[height=0.75\textwidth,width=1.0\textwidth]{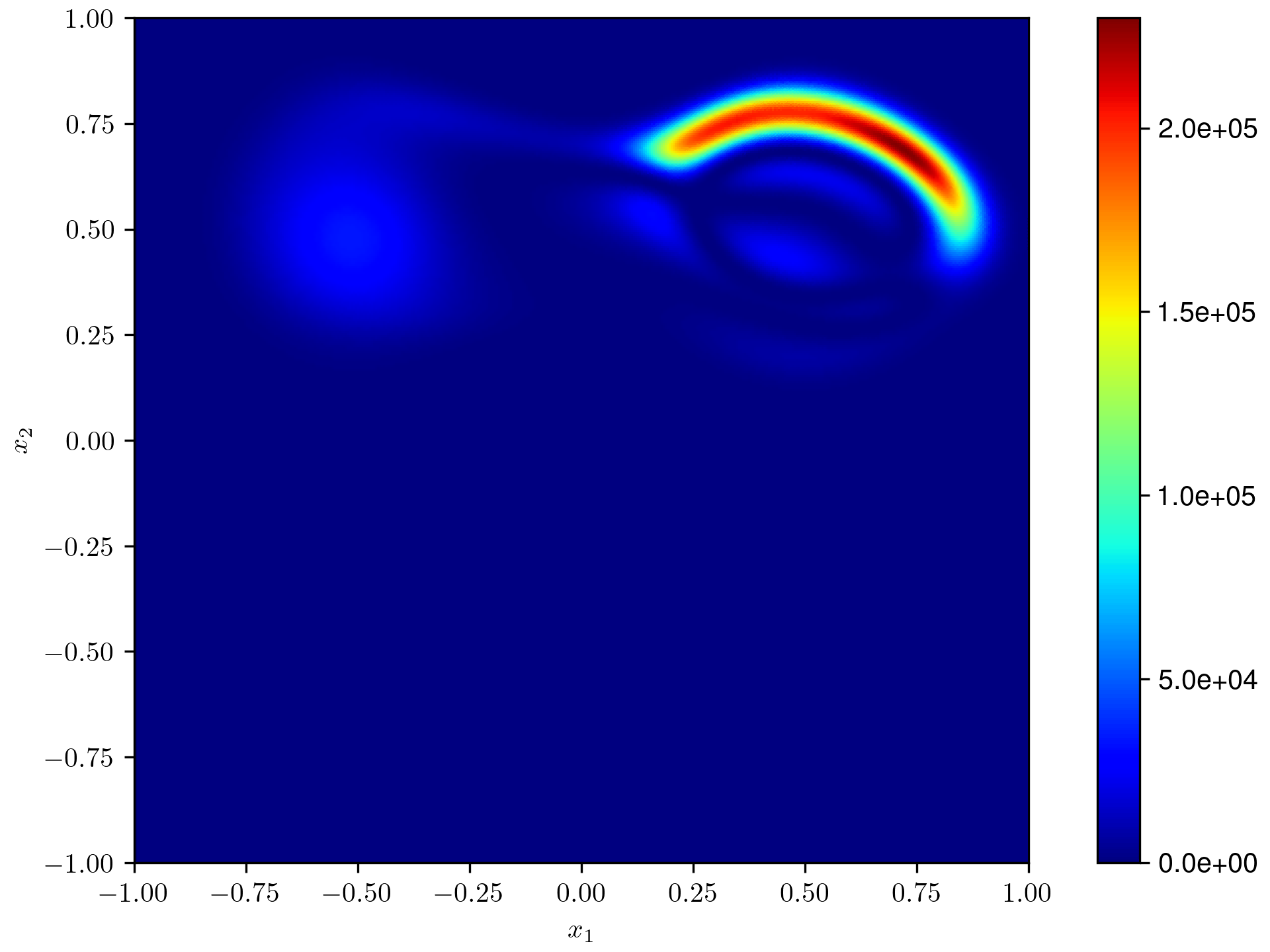}
    \end{subfigure}%
    \begin{subfigure}{.25\textwidth}
        \centering
        \includegraphics[height=0.75\textwidth,width=1.0\textwidth]{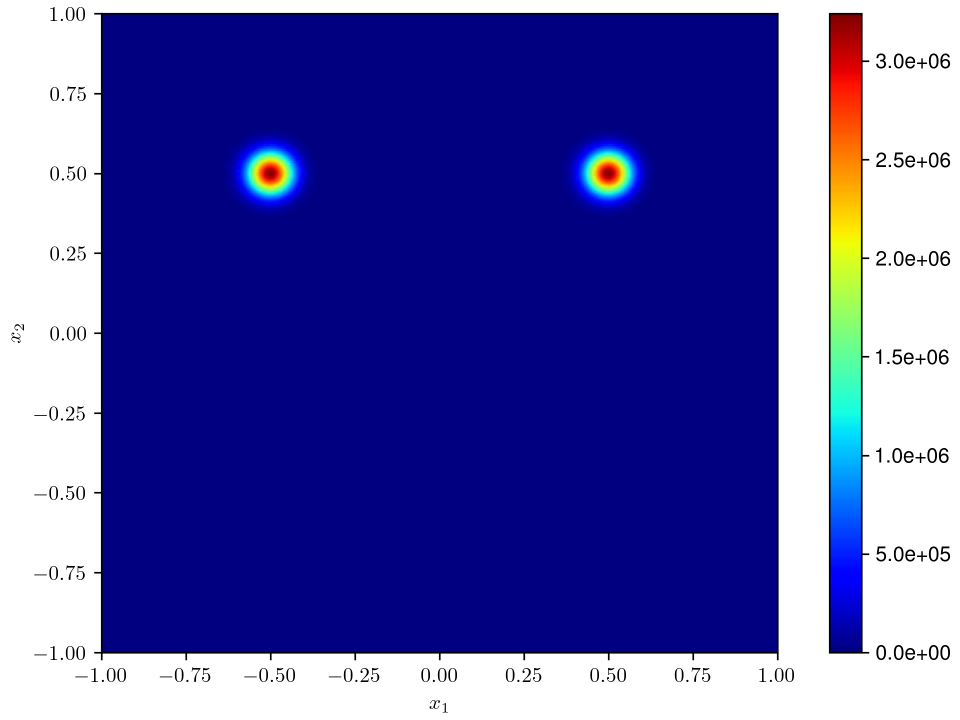}
    \end{subfigure}%
    \begin{subfigure}{.25\textwidth}
        \centering
        \includegraphics[height=0.75\textwidth,width=1.0\textwidth]{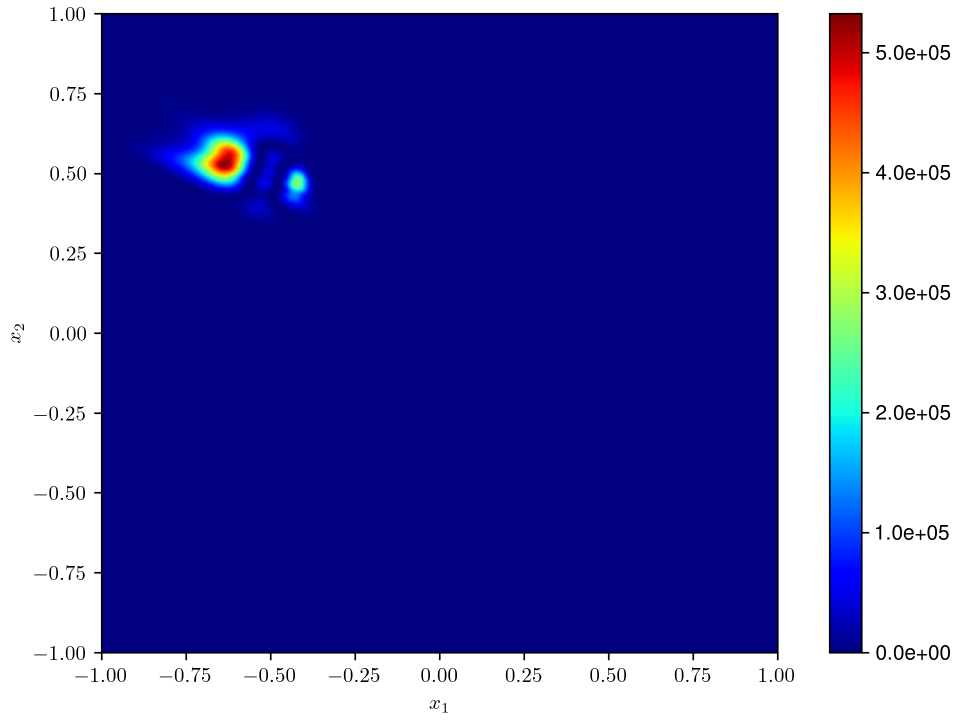}
    \end{subfigure}%
    \newline
    \raggedleft
    \begin{subfigure}{.25\textwidth}
        \centering
        \includegraphics[height=0.75\textwidth,width=1.0\textwidth]{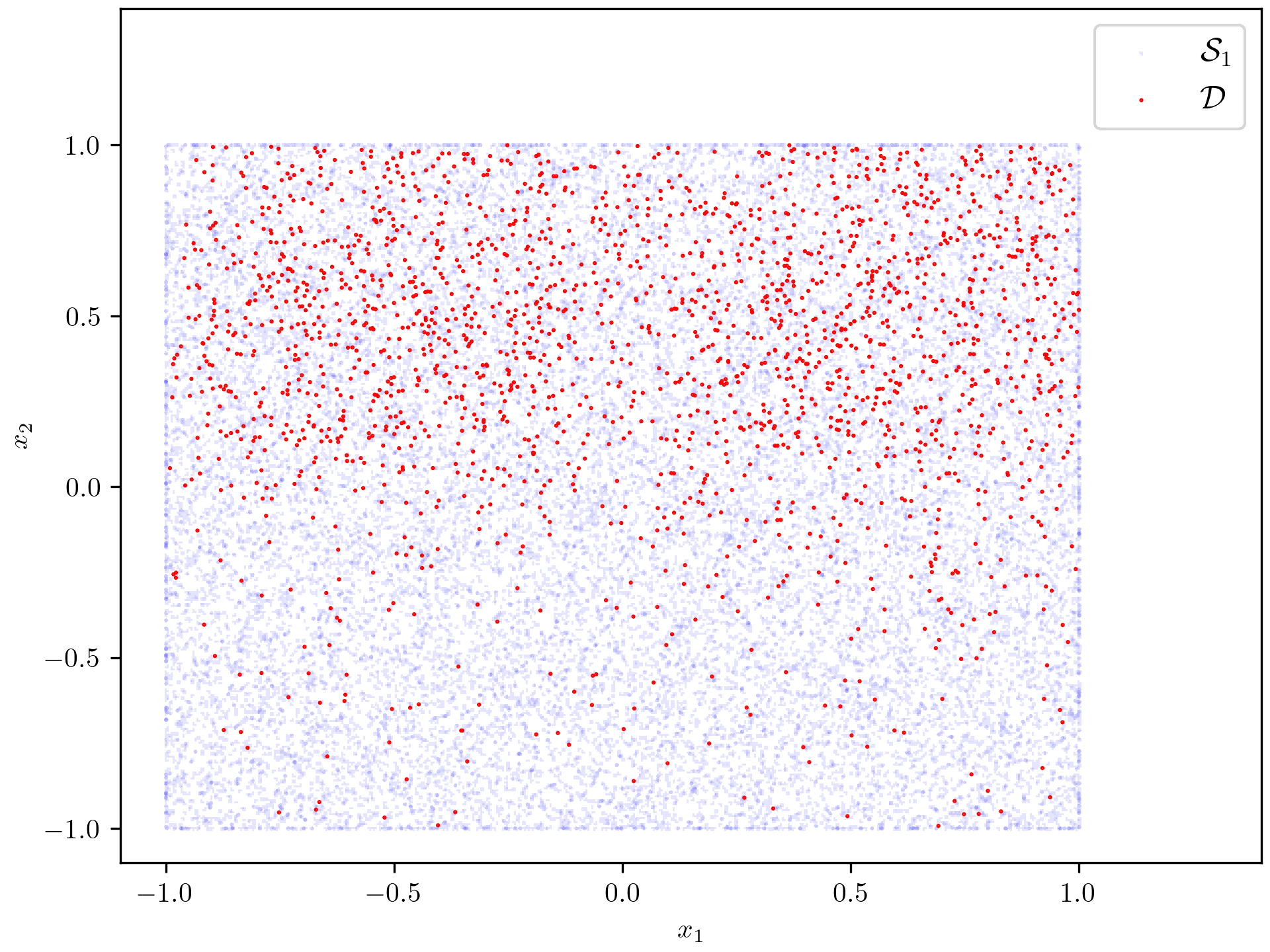}
        \caption{\textit{RAD}, 1st iteration.}
    \end{subfigure}%
    \begin{subfigure}{.25\textwidth}
        \centering
        \includegraphics[height=0.75\textwidth,width=1.0\textwidth]{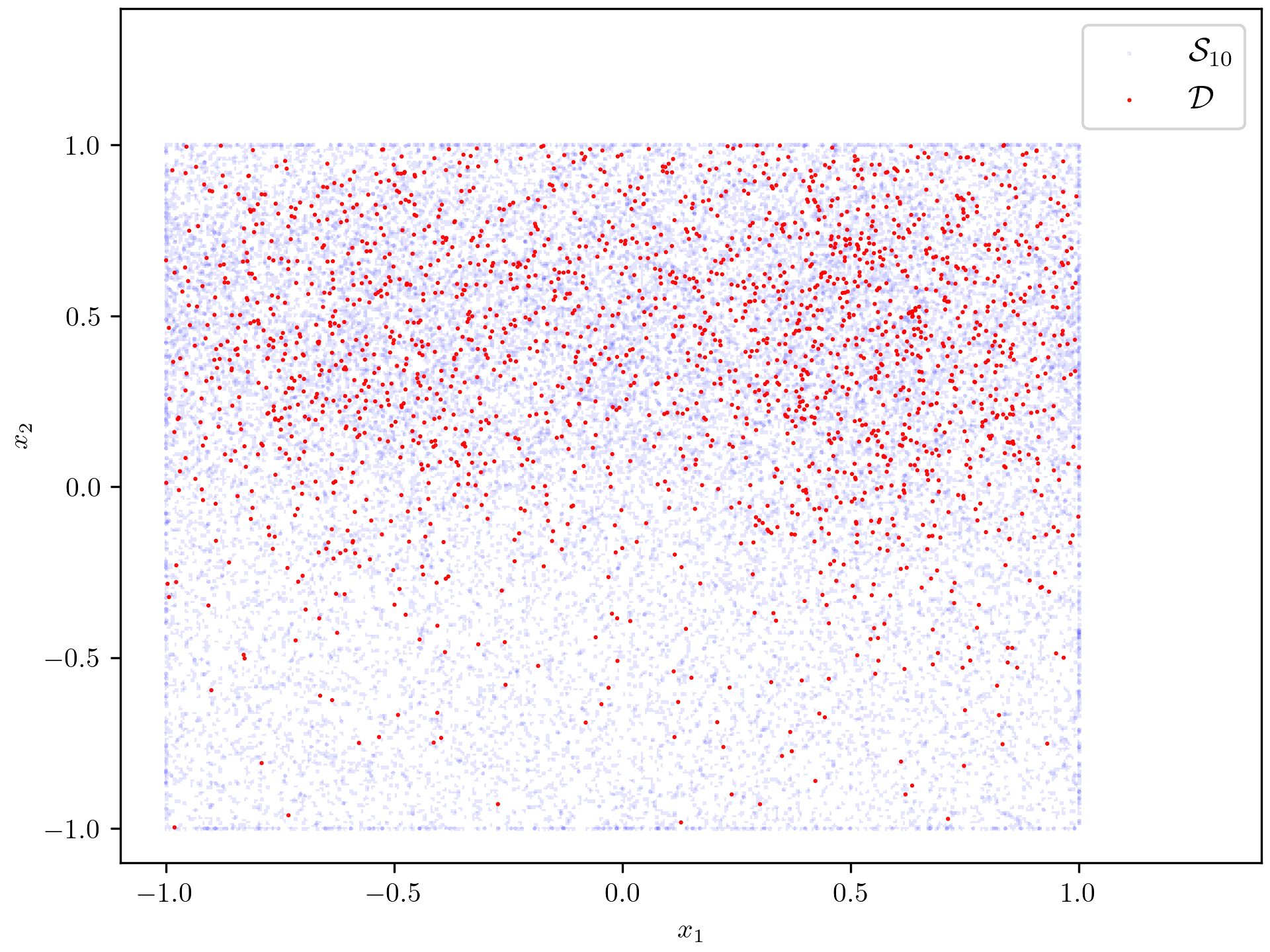}
        \caption{\textit{RAD}, 10th iteration.}
    \end{subfigure}%
    \begin{subfigure}{.25\textwidth}
        \centering
        \includegraphics[height=0.75\textwidth,width=1.0\textwidth]{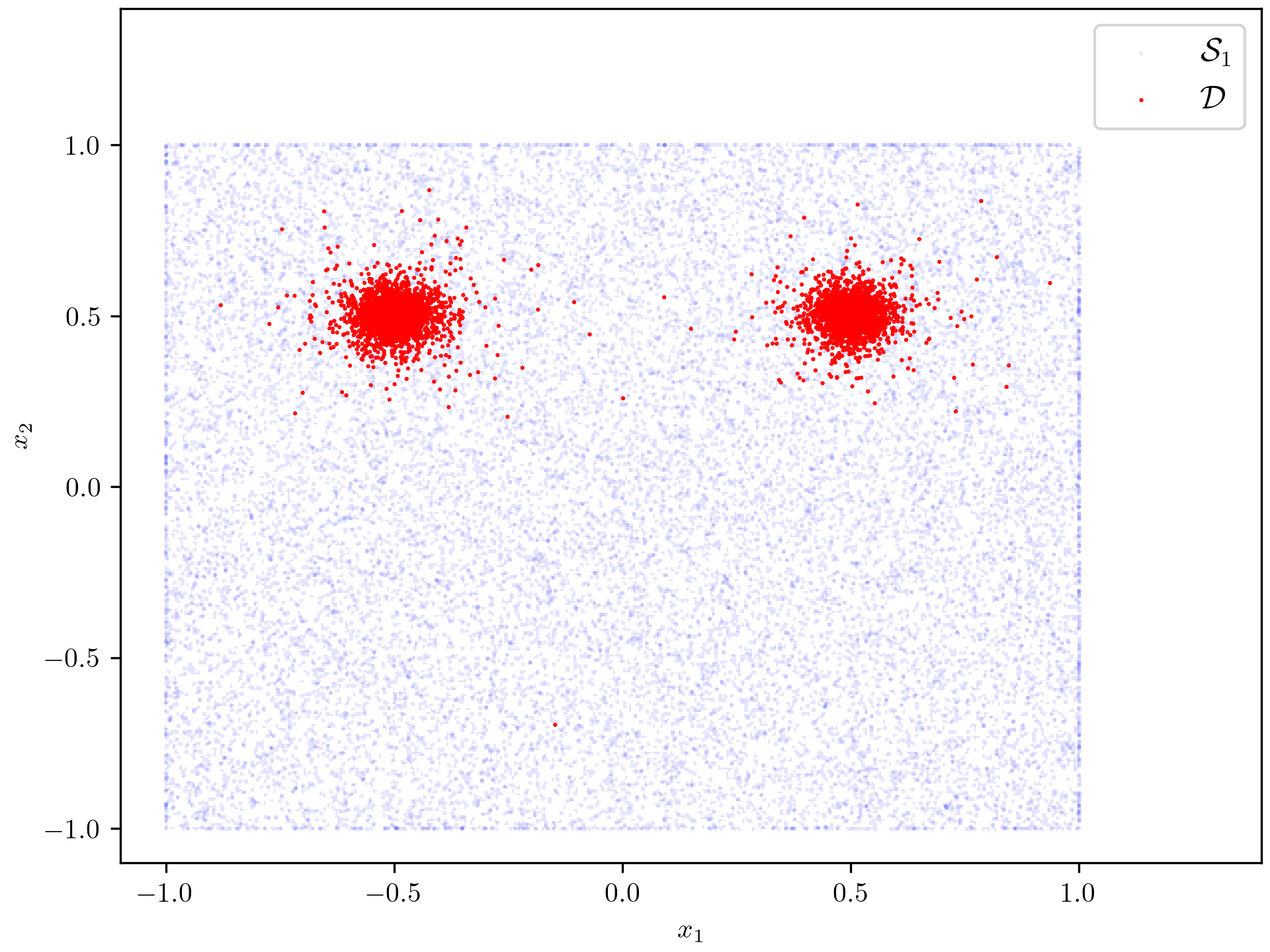}
        \caption{\textit{AAIS-t}, 1st iteration.}
    \end{subfigure}%
    \begin{subfigure}{.25\textwidth}
        \centering
        \includegraphics[height=0.75\textwidth,width=1.0\textwidth]{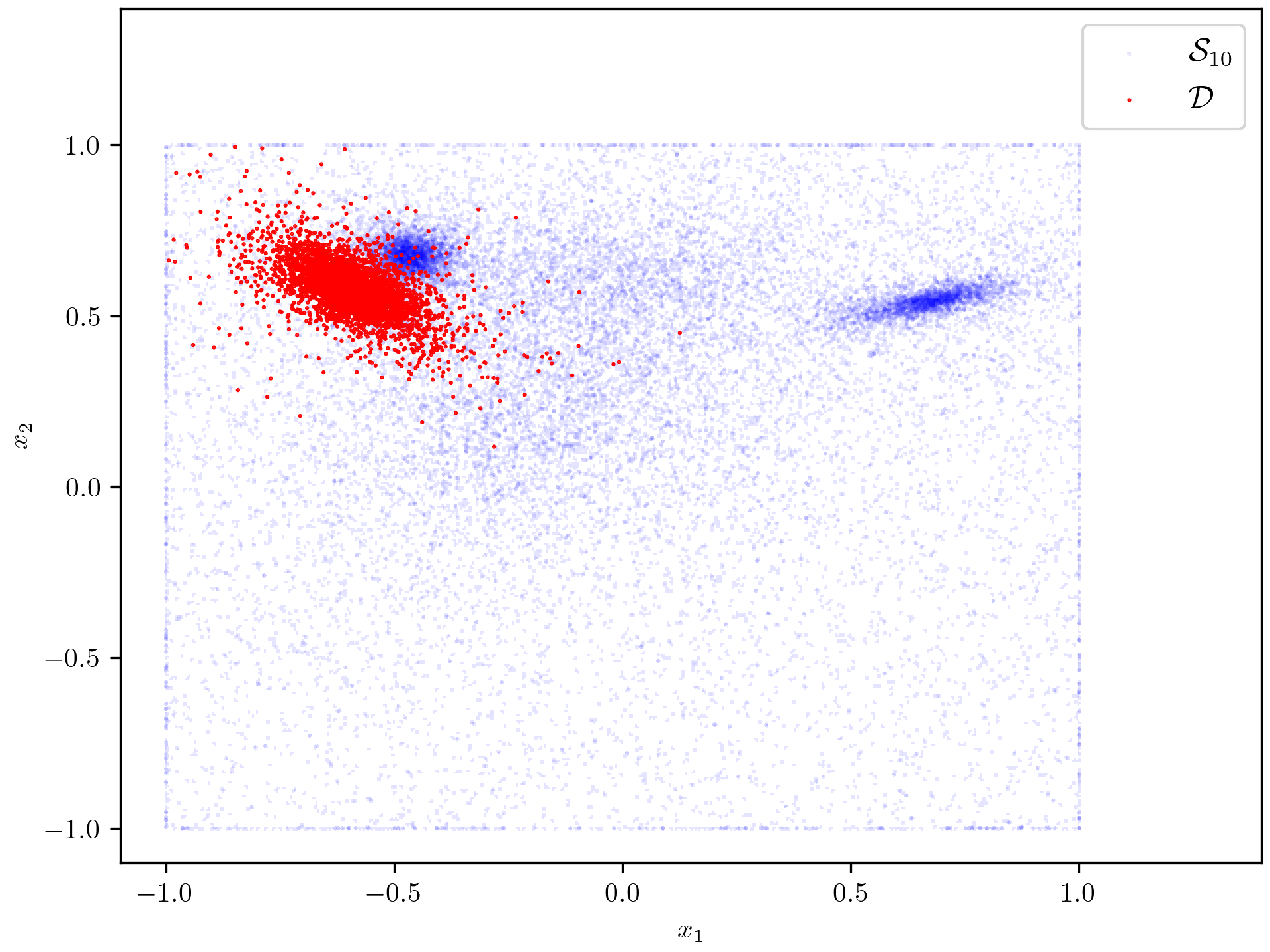}
        \caption{\textit{AAIS-t}, 10th iteration.}
    \end{subfigure}%
    \caption{$x_1x_2$-plane profiles of residual and nodes for 9D Poisson problem at 1st and 10th iteration.}
    \label{fig:Ps9DNode}
\end{figure}

Take a step forward, the residual and the node at the 1st and 10th iteration are plotted in Figure \ref{fig:Ps9DNode}. \textit{RAD} method just realize  the area of bigger loss, but fail to locate precisely. However, our proposed AAIS algorithm make the nodes clustering around the singularities after pre-training, implying the efficient training of PINNs. But at the last iteration, we could see that the \textit{AAIS-t} algorithm also fails to mimic the residual due to the sparsity of 200k searching points in the $9D$ domain.

\begin{figure}[htbp]
    \centering
    \begin{subfigure}{.25\textwidth}
        \centering
        \includegraphics[height=0.75\textwidth,width=1.0\textwidth]{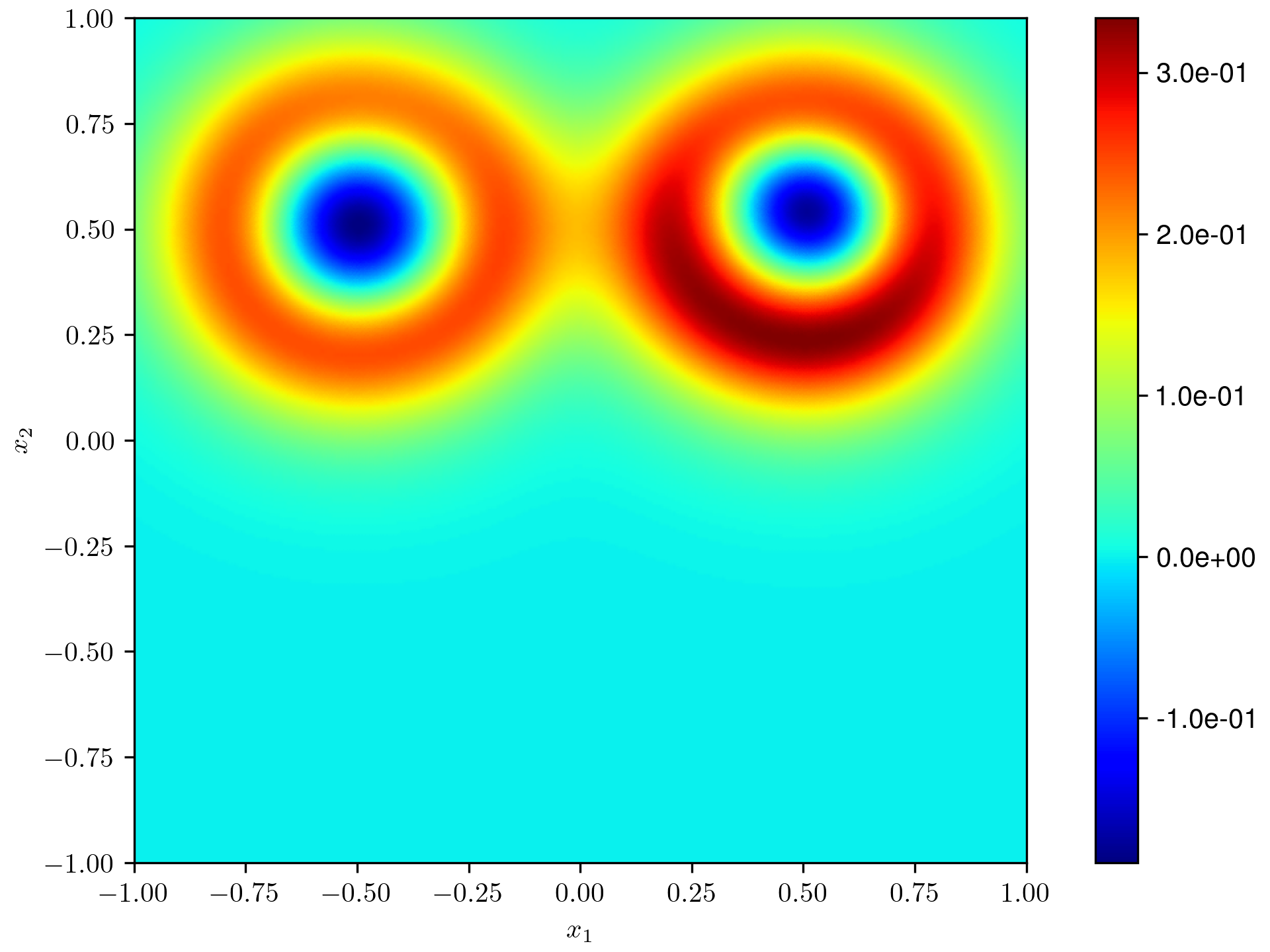}
        \caption{\textit{RAD}, solution.}
    \end{subfigure}%
    \begin{subfigure}{.25\textwidth}
        \centering
        \includegraphics[height=0.75\textwidth,width=1.0\textwidth]{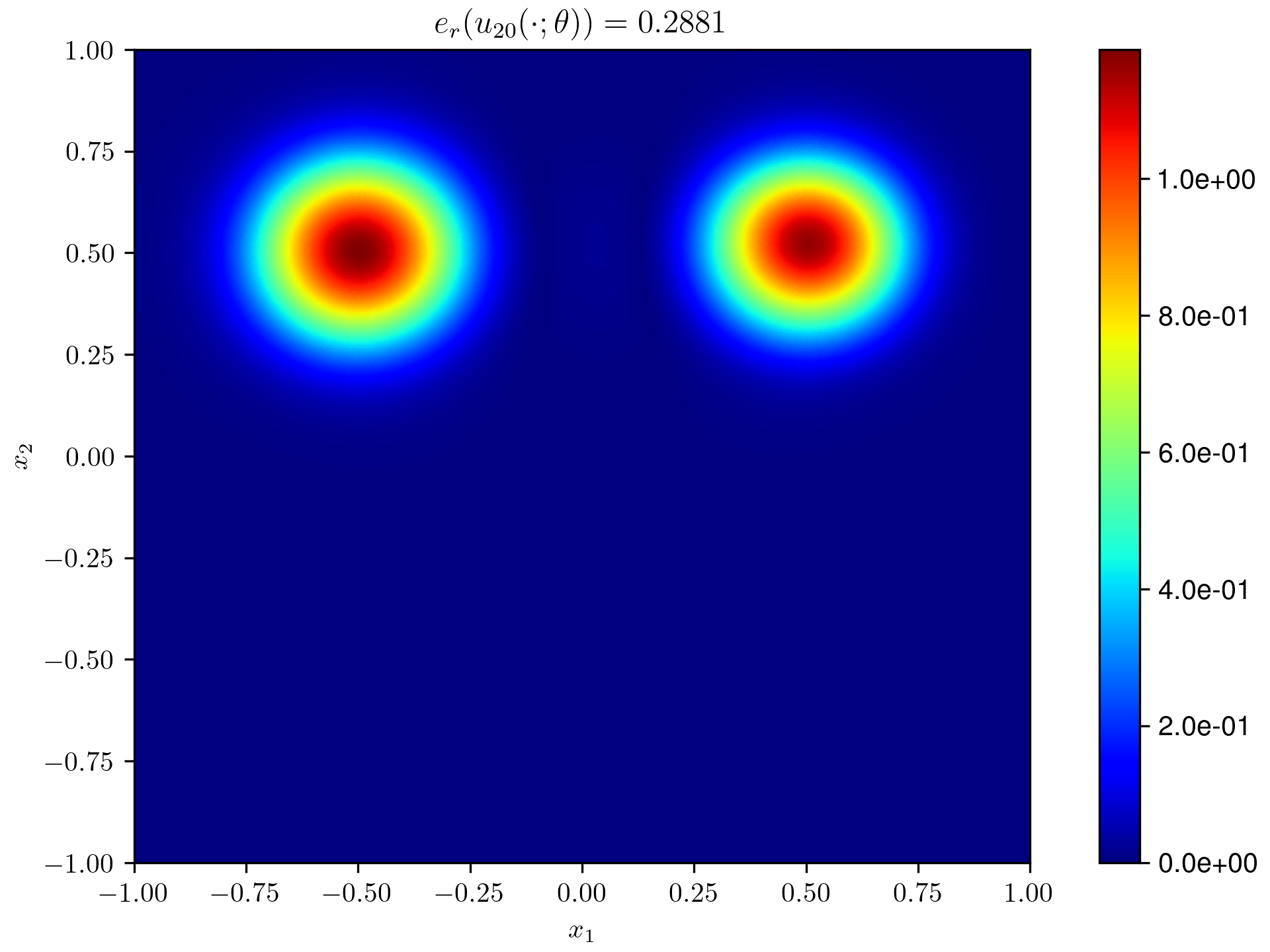}
        \caption{\textit{RAD}, absolute error.}
    \end{subfigure}%
    \begin{subfigure}{.25\textwidth}
        \centering
        \includegraphics[height=0.75\textwidth,width=1.0\textwidth]{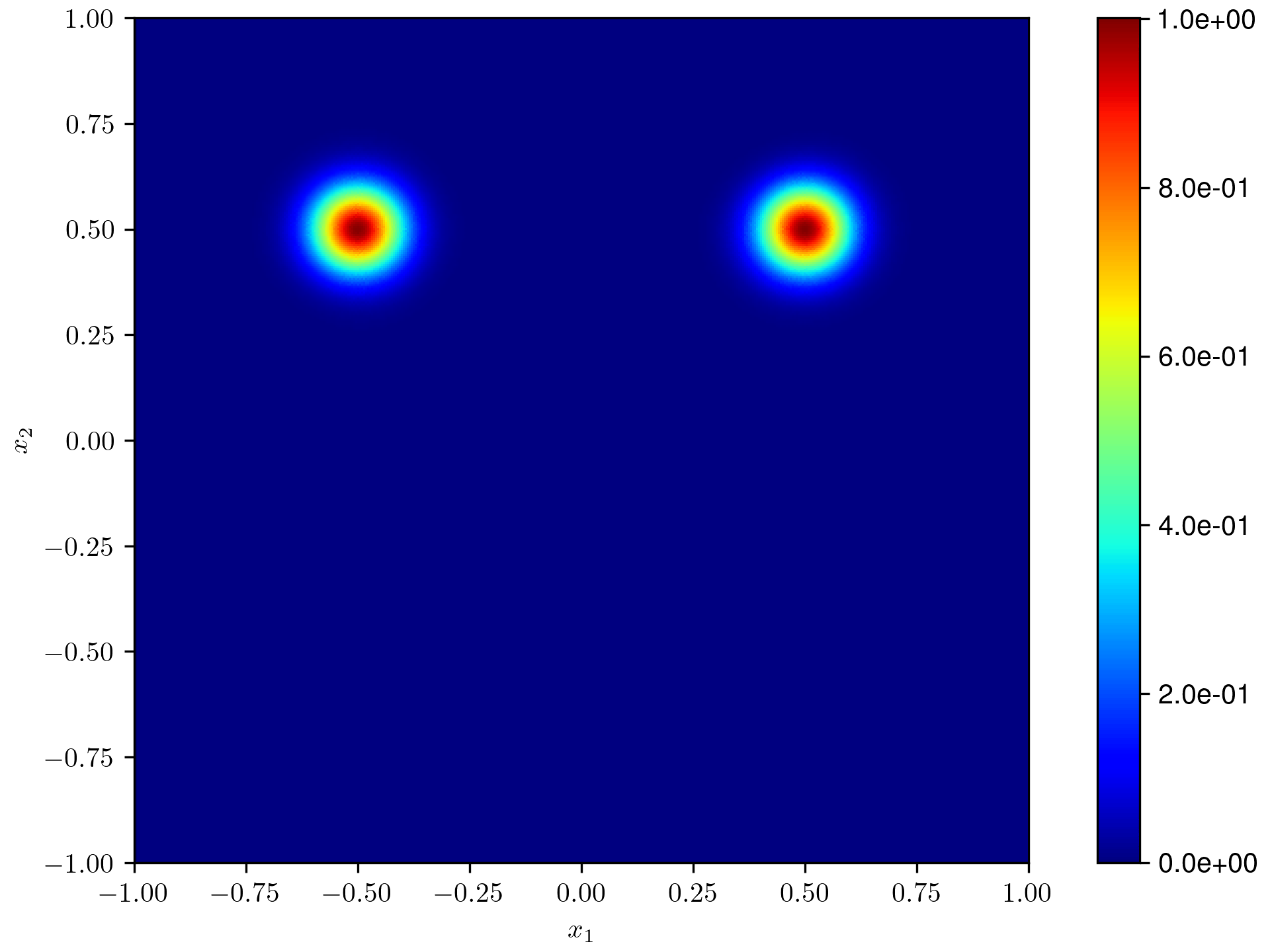}
        \caption{\textit{AAIS-t}, solution.}
    \end{subfigure}%
    \begin{subfigure}{.25\textwidth}
        \centering
        \includegraphics[height=0.75\textwidth,width=1.0\textwidth]{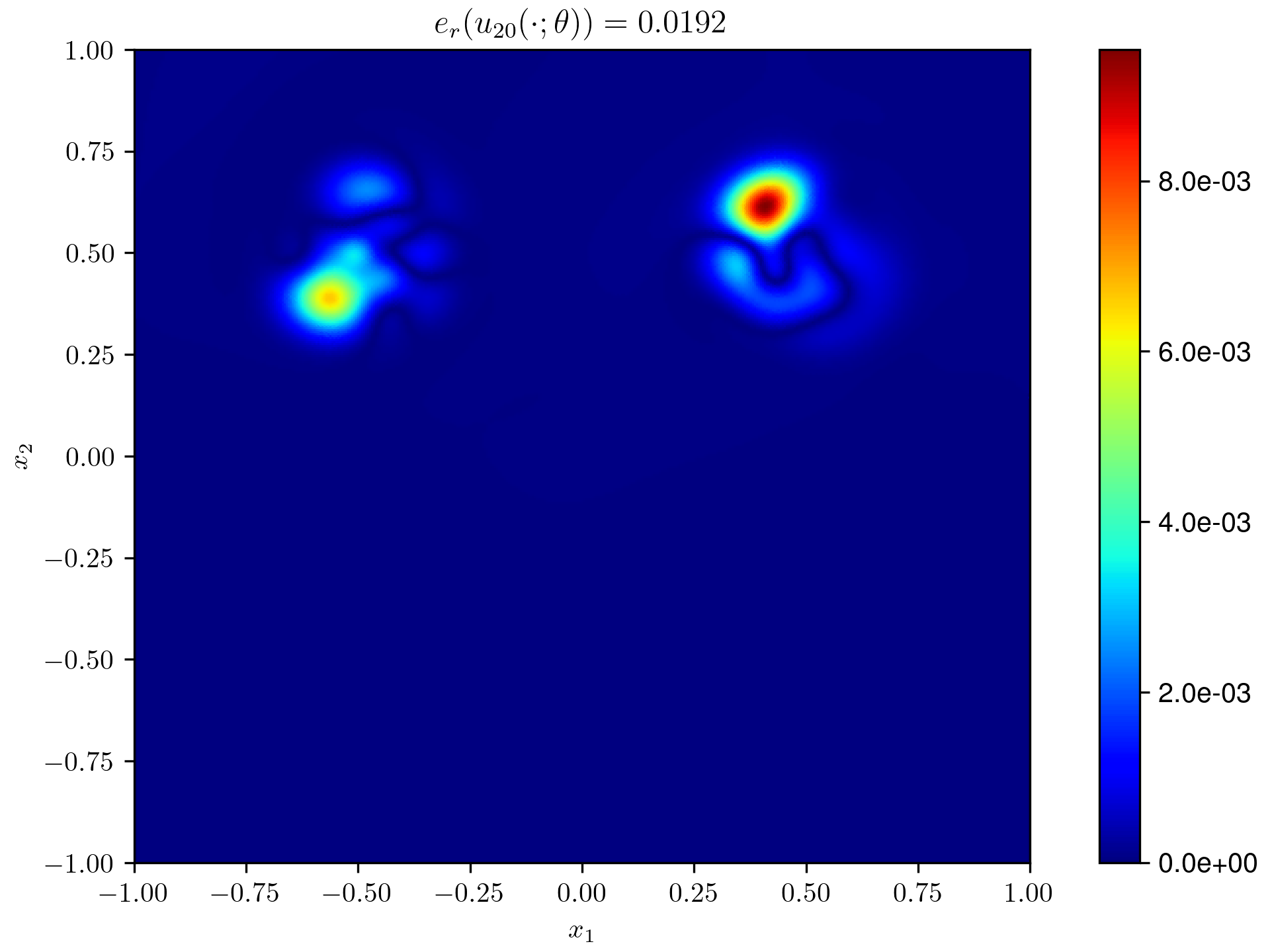}
        \caption{\textit{AAIS-t}, absolute error.}
    \end{subfigure}%
    \caption{Profiles of absolute error and neural network solutions for 9D Poisson equation.}
    \label{fig:Ps9Dsol}
\end{figure}

The solution and the absolute error are showed in Figure \ref{fig:Ps9Dsol}, we could see that the \textit{RAD} method could not solve the problem accurately but realize the locations of singularities. \textit{AAIS-t} algorithm could perfectly solve the problem with increasing frequency of absolute error. Therefore, as far as we know, it is the first time adaptive methods could solve high dimensional Poisson problems with multiple singularities. Our proposed AAIS algorithm based importance sampling could obtain satisfactory numerical results in adaptive PINNs.
\subsubsection{Fifteen dimensional one peak problem}
\begin{figure}[htbp]
    \centering
    \includegraphics[scale=0.125]{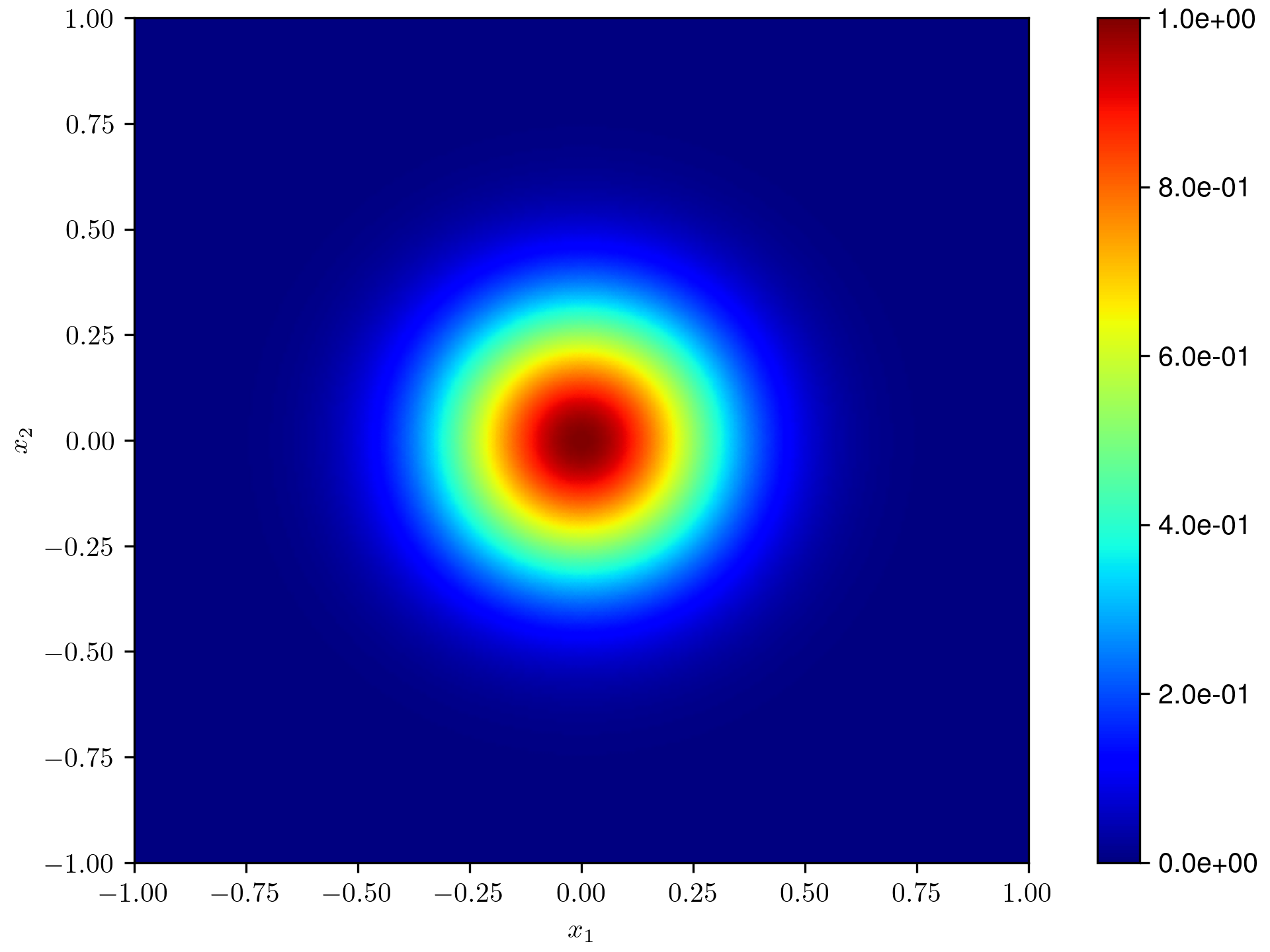}
    \caption{Projection on $x_1x_2$ plane of exact solution for Poisson 15D problems with two peaks in \eqref{pde:Poisson1Peak}.}
    \label{fig:PS15Dexact}
\end{figure}

In this part we test our proposed AAIS algorithm in a very high dimension, the center is $(0,0,...,0)$(see in Figure \ref{fig:PS15Dexact}). The grid points are the same in 9D case and the neural network structure has 20 neurons with 7 hidden layers due to the single peak solution. Here we let $N_S=100000$ because of the memory limitation of GPU.

\begin{figure}[htbp]
    \centering
    \begin{subfigure}{.5\textwidth}
        \centering
        \includegraphics[scale=0.25]{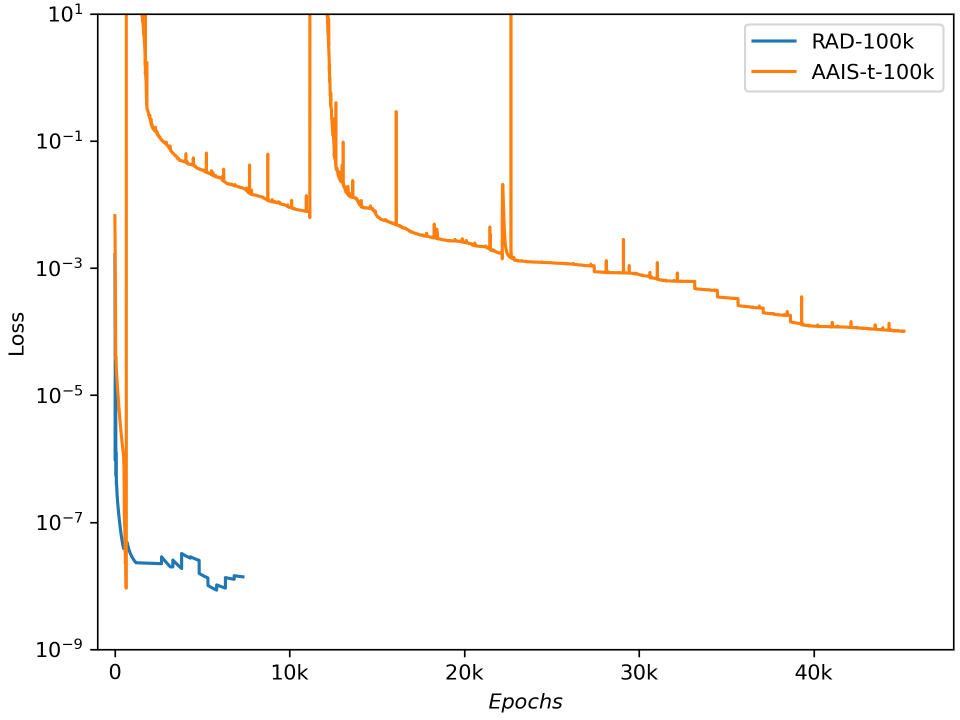}
        \end{subfigure}%
        \begin{subfigure}{.5\textwidth}
        \centering
        \includegraphics[scale=0.25]{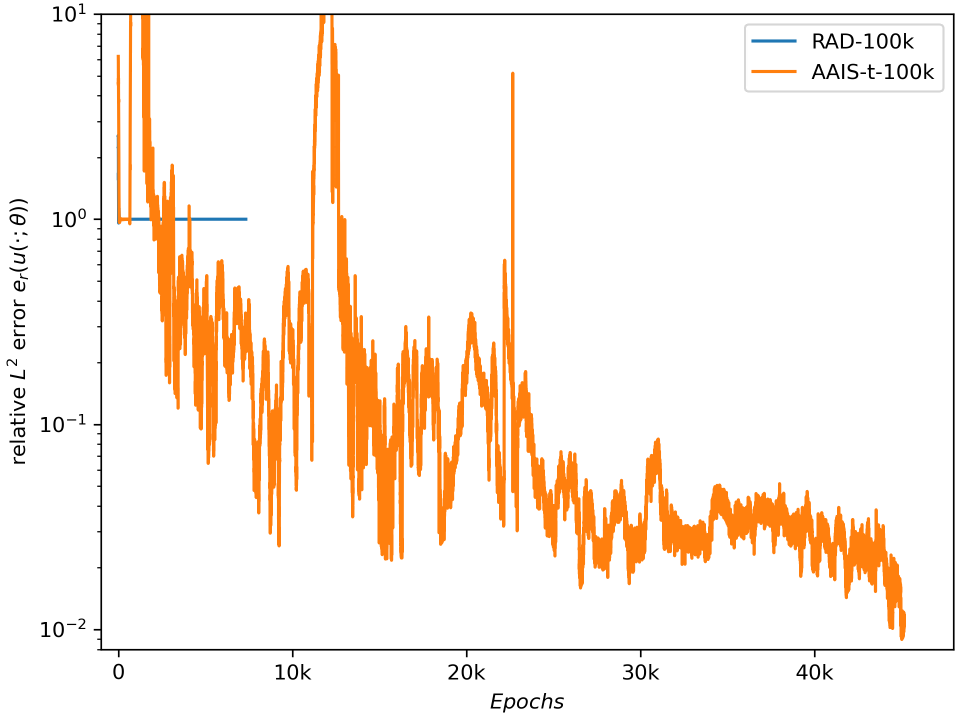}
        \end{subfigure}%
        \caption{Loss and relative errors for 15D Poisson equation. Left: the loss function. Right: the relative $L^2$ error $e_r(u(\cdot;\theta))$.
        }
    \label{fig:PS15DErr}
\end{figure}

\begin{figure}[htbp]
    \centering
    \begin{subfigure}{.25\textwidth}
        \centering
        \includegraphics[height=0.75\textwidth,width=1.0\textwidth]{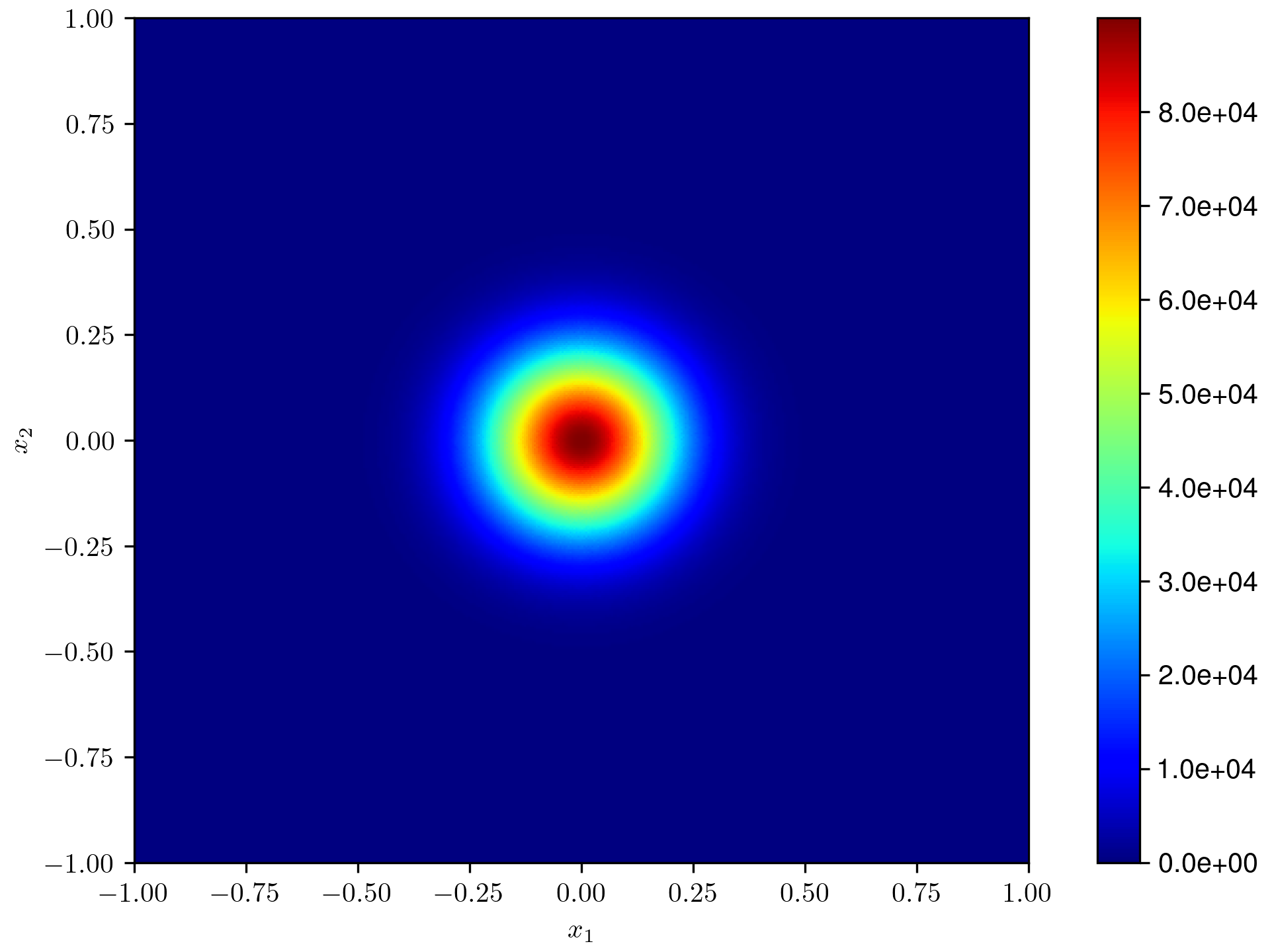}
      \end{subfigure}%
    \begin{subfigure}{.25\textwidth}
        \centering
        \includegraphics[height=0.75\textwidth,width=1.0\textwidth]{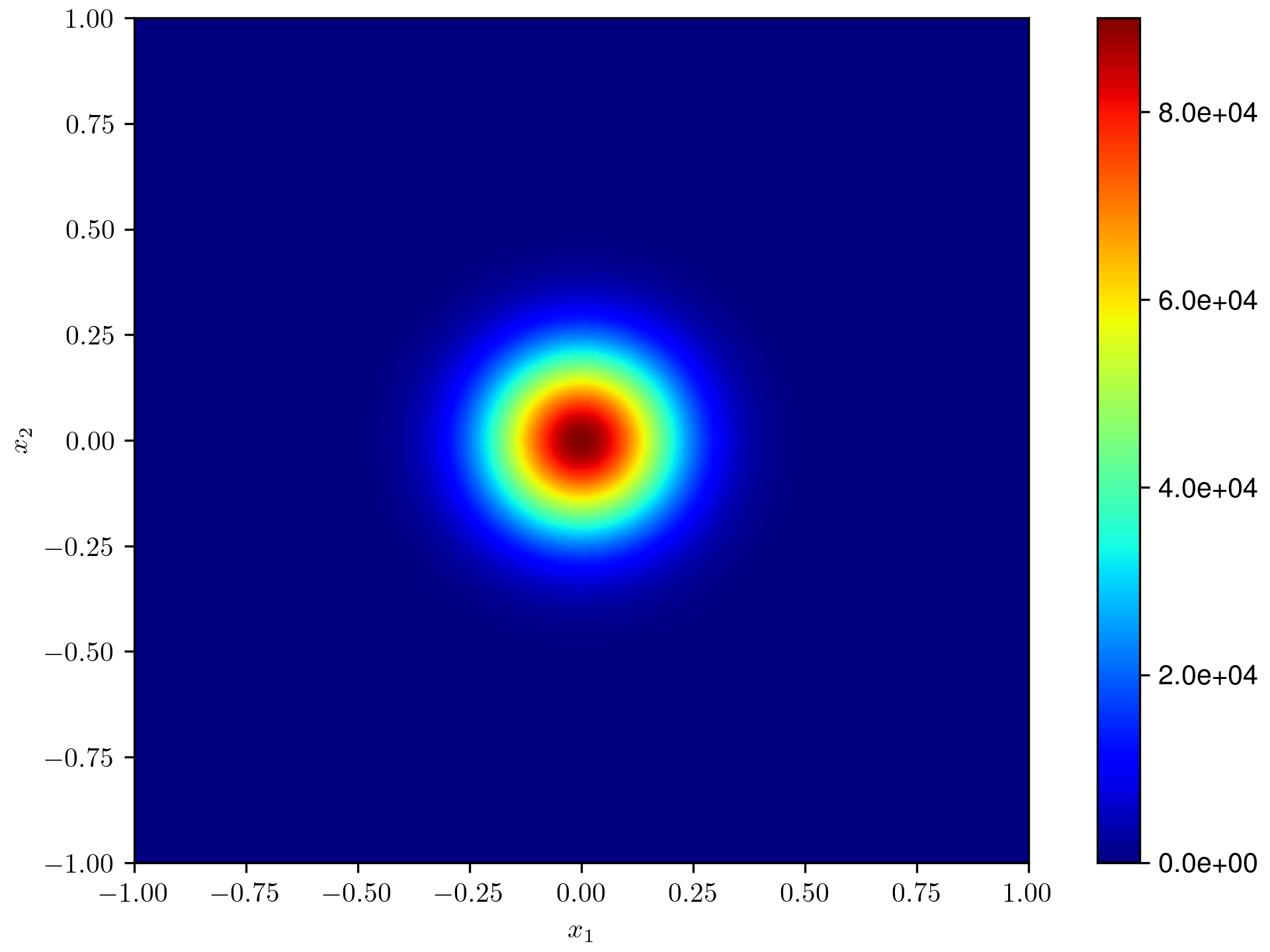}
    \end{subfigure}%
    \begin{subfigure}{.25\textwidth}
        \centering
        \includegraphics[height=0.75\textwidth,width=1.0\textwidth]{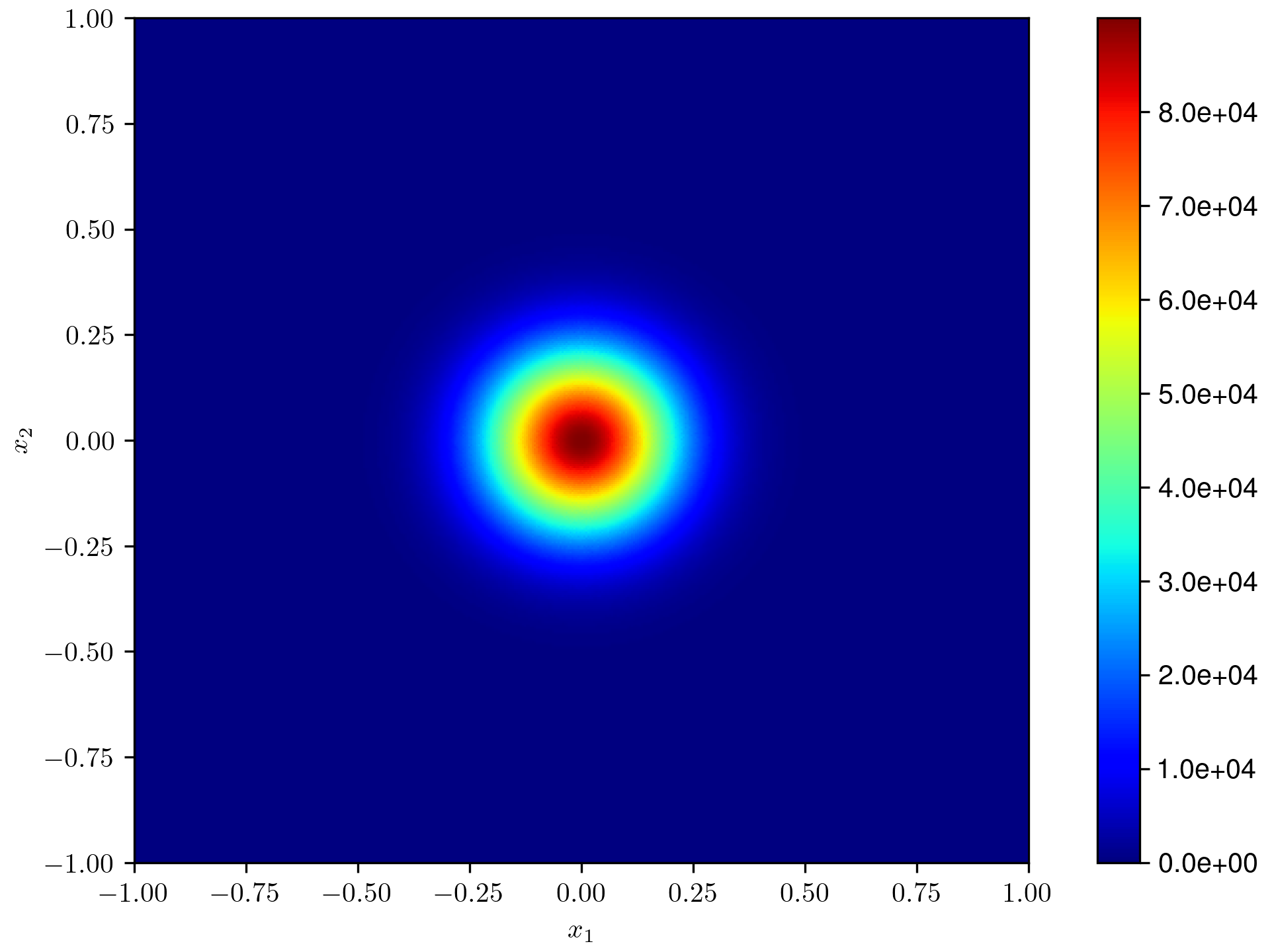}
    \end{subfigure}%
    \begin{subfigure}{.25\textwidth}
        \centering
        \includegraphics[height=0.75\textwidth,width=1.0\textwidth]{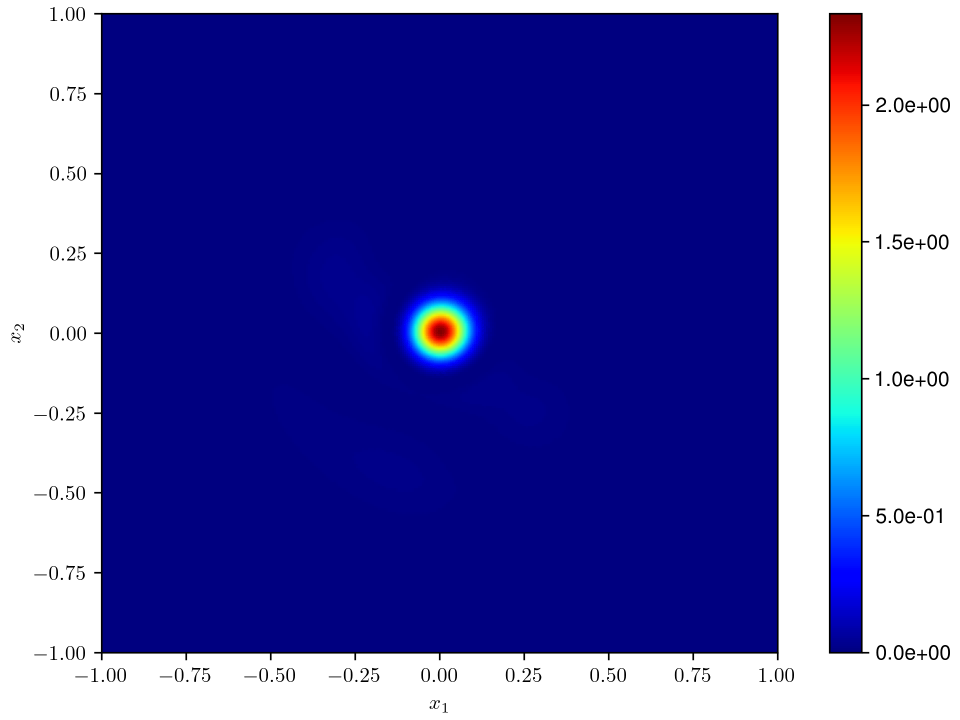}
    \end{subfigure}%
    \newline
    \raggedleft
    \begin{subfigure}{.25\textwidth}
        \centering
        \includegraphics[height=0.75\textwidth,width=1.0\textwidth]{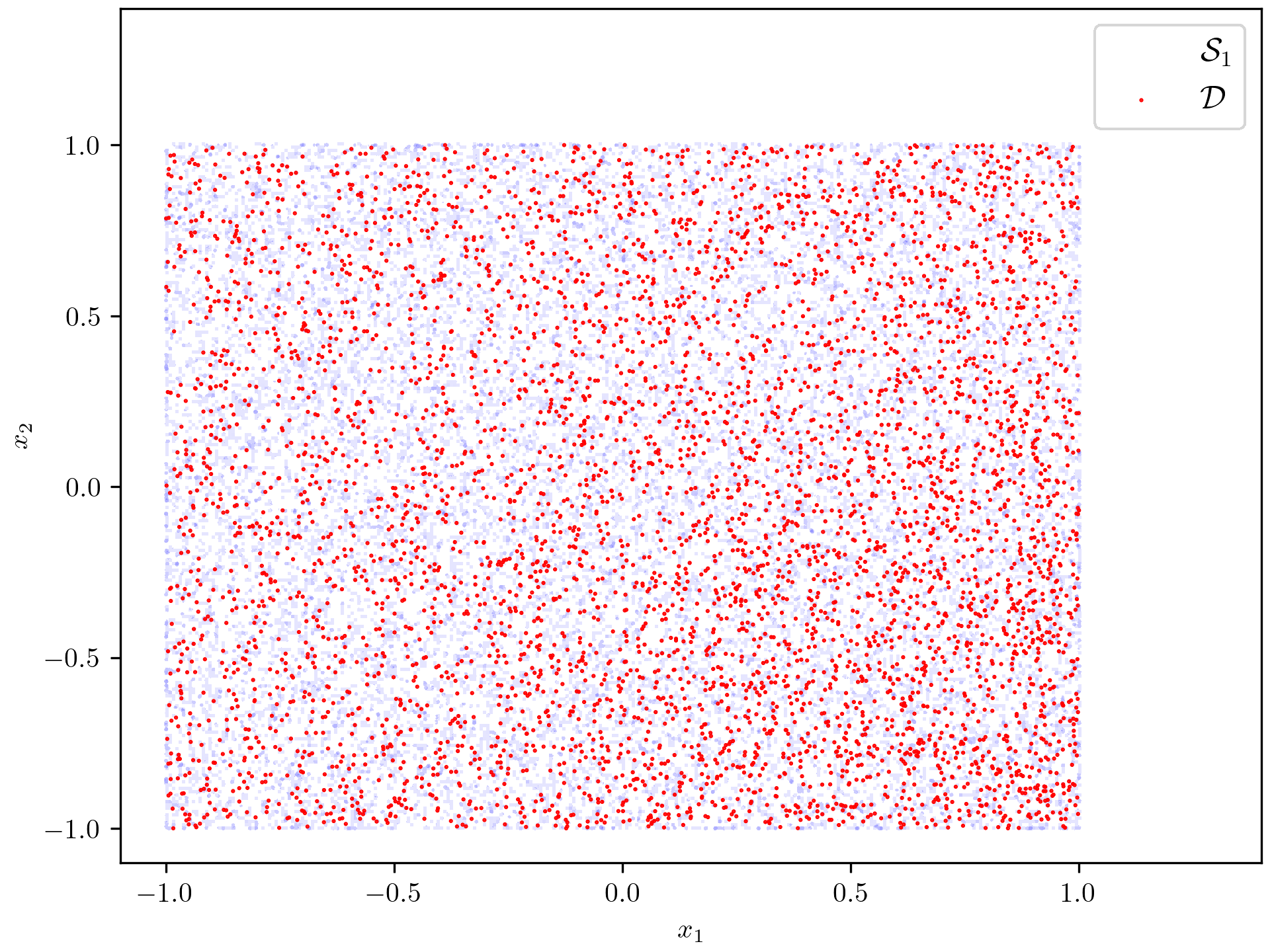}
        \caption{\textit{RAD}, 1st iteration.}
    \end{subfigure}%
    \begin{subfigure}{.25\textwidth}
        \centering
        \includegraphics[height=0.75\textwidth,width=1.0\textwidth]{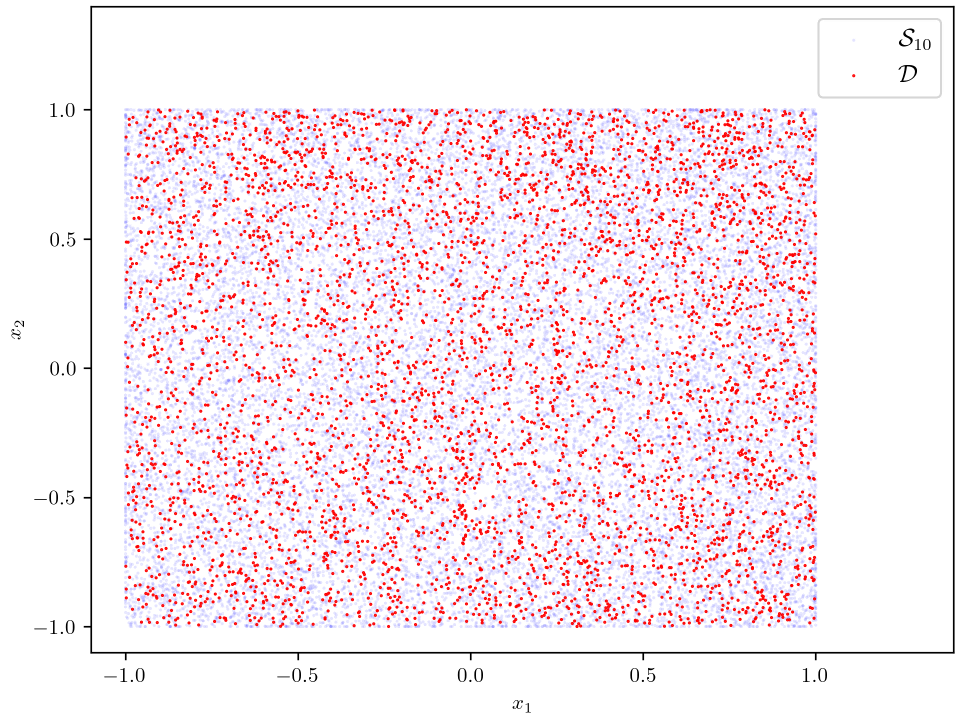}
        \caption{\textit{RAD}, 10th iteration.}
    \end{subfigure}%
    \begin{subfigure}{.25\textwidth}
        \centering
        \includegraphics[height=0.75\textwidth,width=1.0\textwidth]{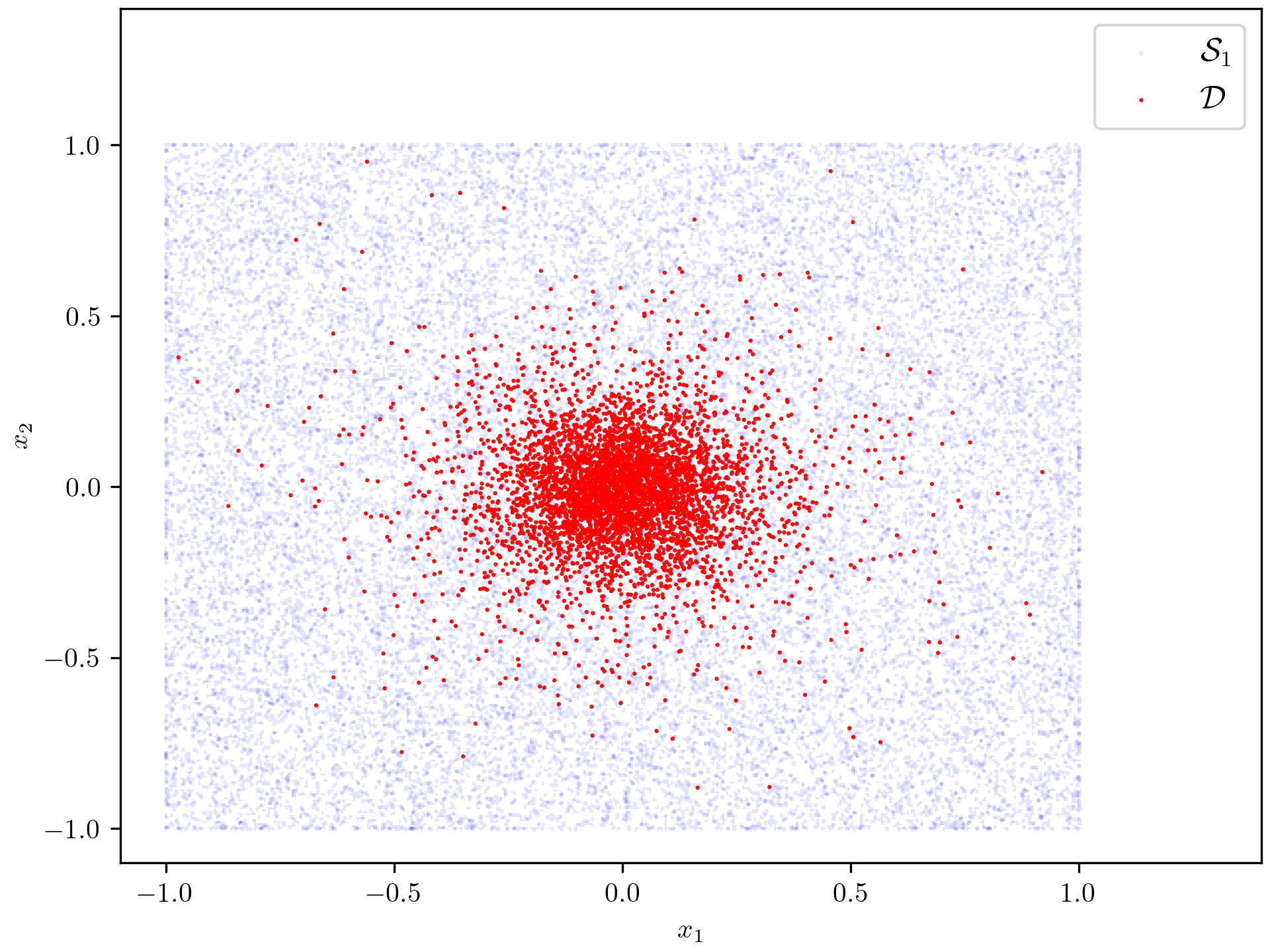}
        \caption{\textit{AAIS-t}, 1st iteration.}
    \end{subfigure}%
    \begin{subfigure}{.25\textwidth}
        \centering
        \includegraphics[height=0.75\textwidth,width=1.0\textwidth]{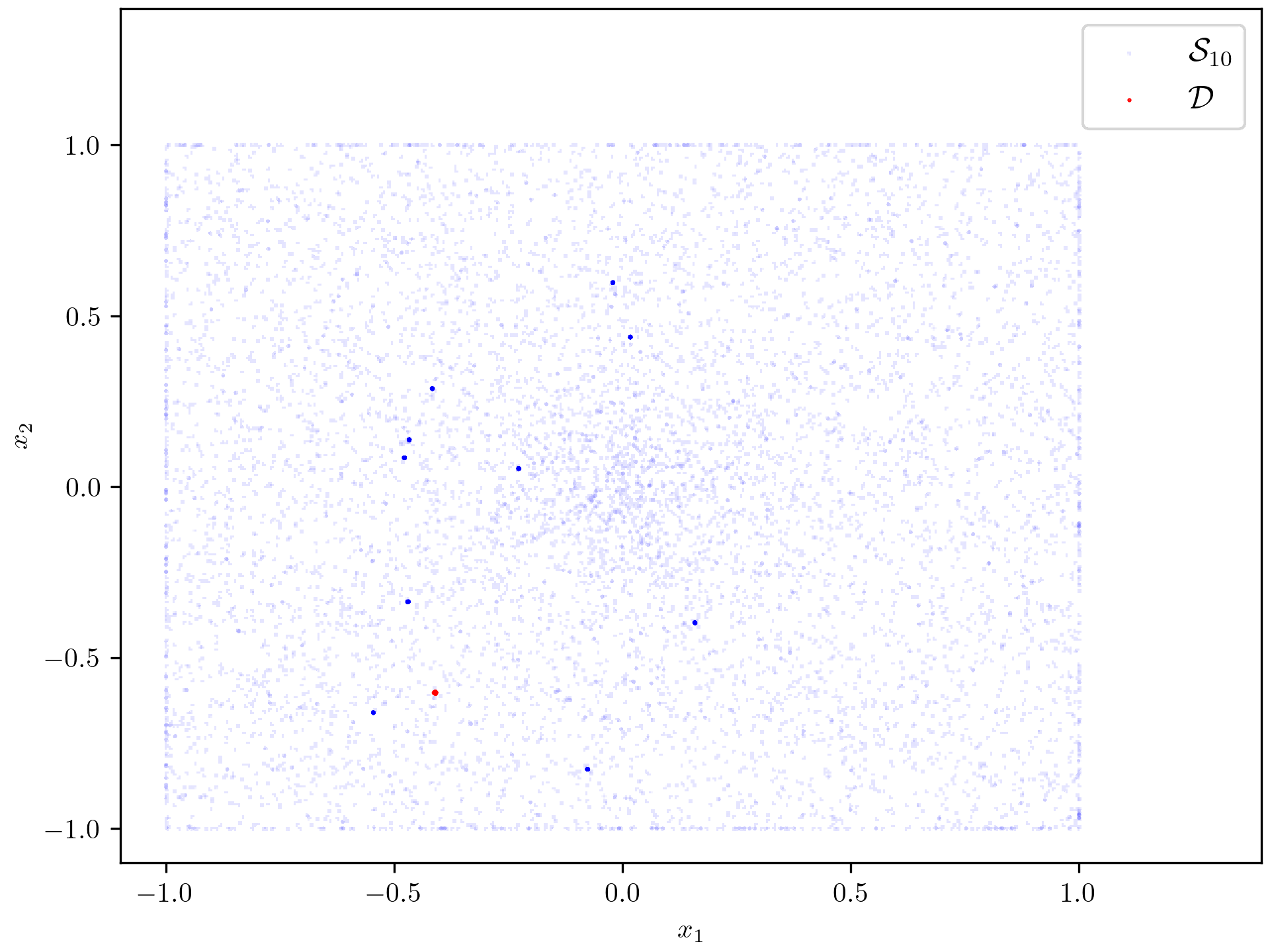}
        \caption{\textit{AAIS-t}, 10th iteration.}
    \end{subfigure}%
    \caption{$x_1x_2$-plane profiles of residual and nodes for 15D Poisson problem at 1st and 10th iteration.}
    \label{fig:Ps15DNode}
\end{figure}

\begin{figure}[htbp]
    \centering
    \begin{subfigure}{.25\textwidth}
        \centering
        \includegraphics[height=0.75\textwidth,width=1.0\textwidth]{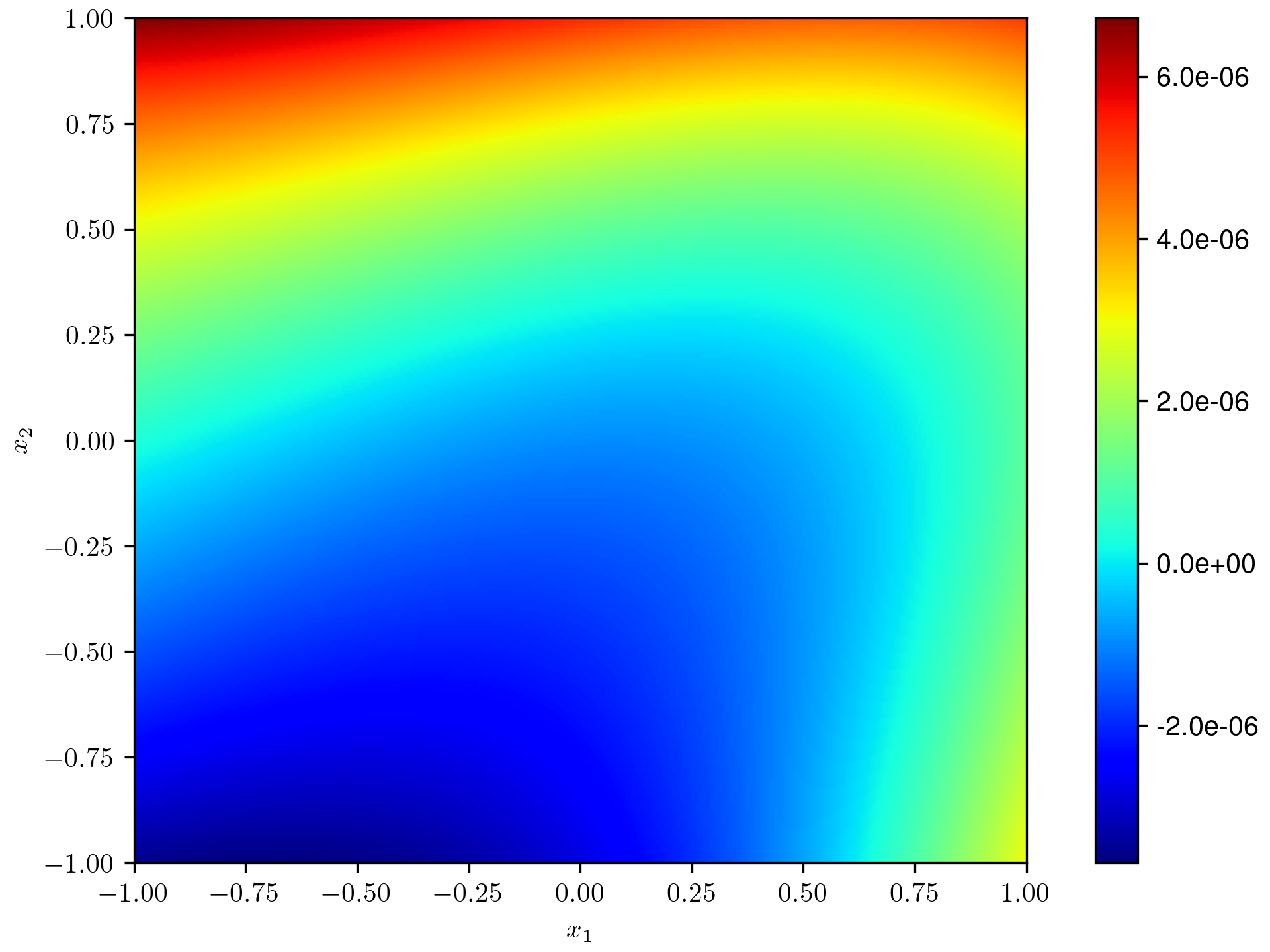}
        \caption{\textit{RAD}, solution.}
    \end{subfigure}%
    \begin{subfigure}{.25\textwidth}
        \centering
        \includegraphics[height=0.75\textwidth,width=1.0\textwidth]{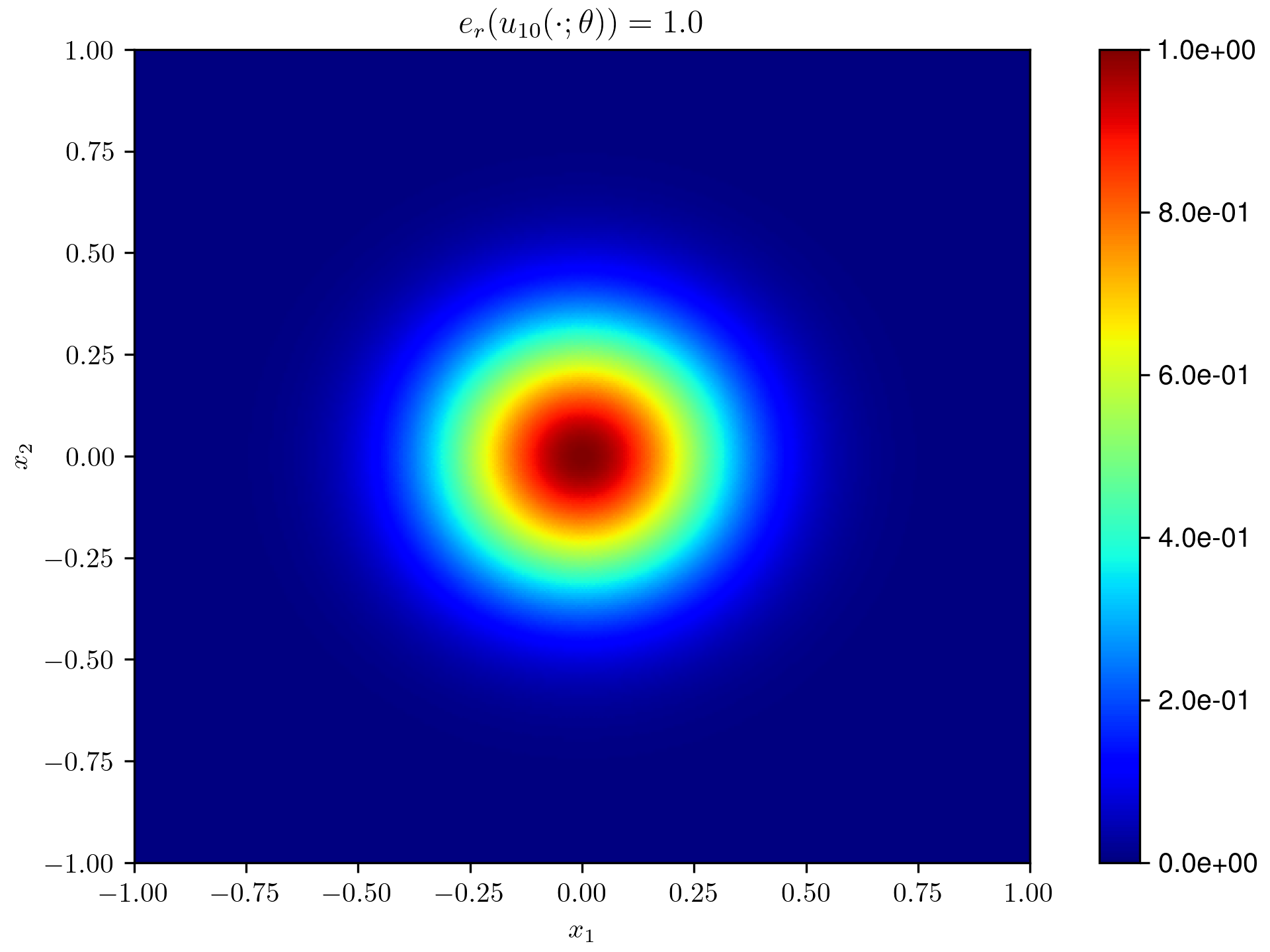}
        \caption{\textit{RAD}, absolute error.}
    \end{subfigure}%
    \begin{subfigure}{.25\textwidth}
        \centering
        \includegraphics[height=0.75\textwidth,width=1.0\textwidth]{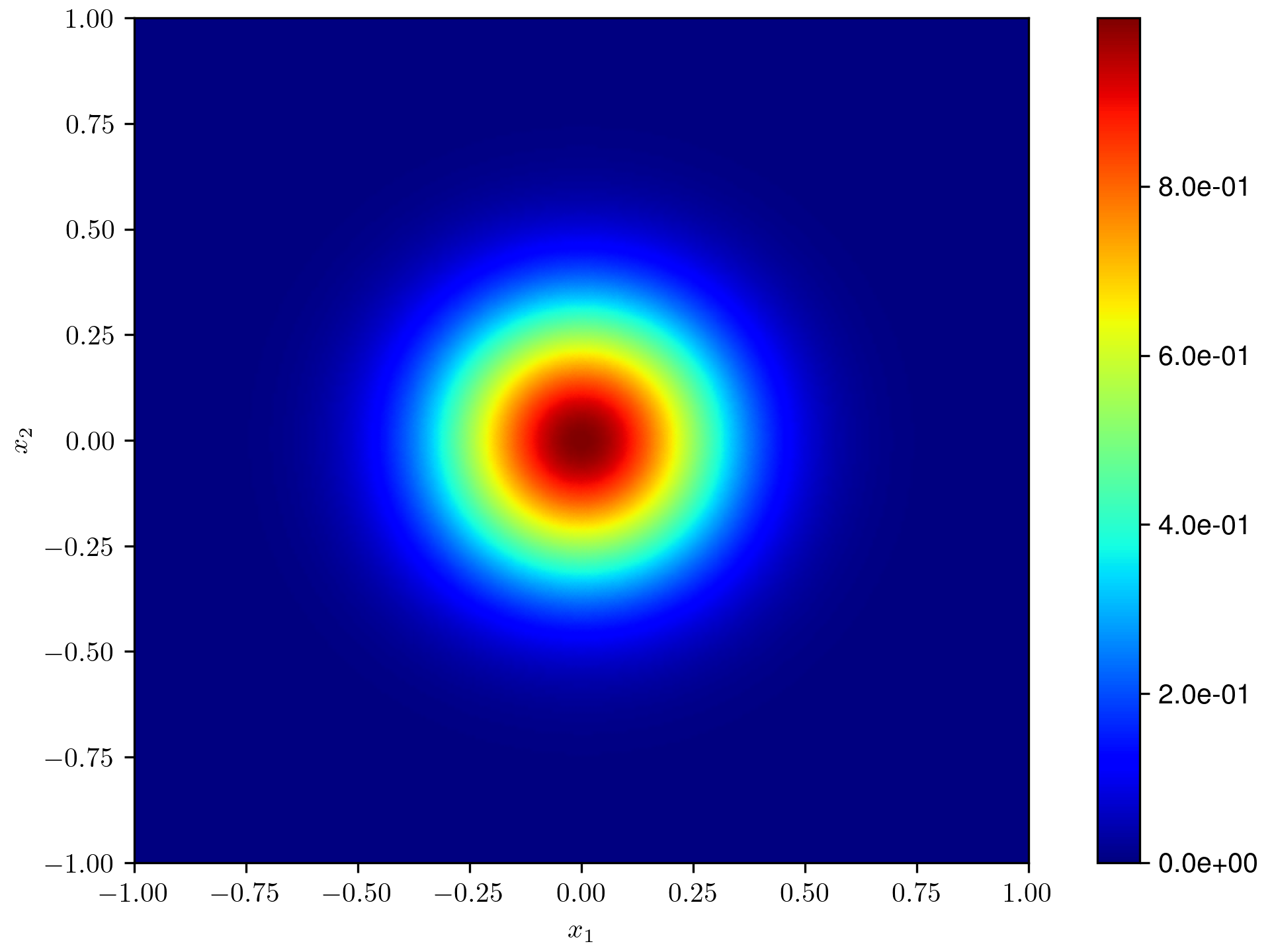}
        \caption{\textit{AAIS-t}, solution.}
    \end{subfigure}%
    \begin{subfigure}{.25\textwidth}
        \centering
        \includegraphics[height=0.75\textwidth,width=1.0\textwidth]{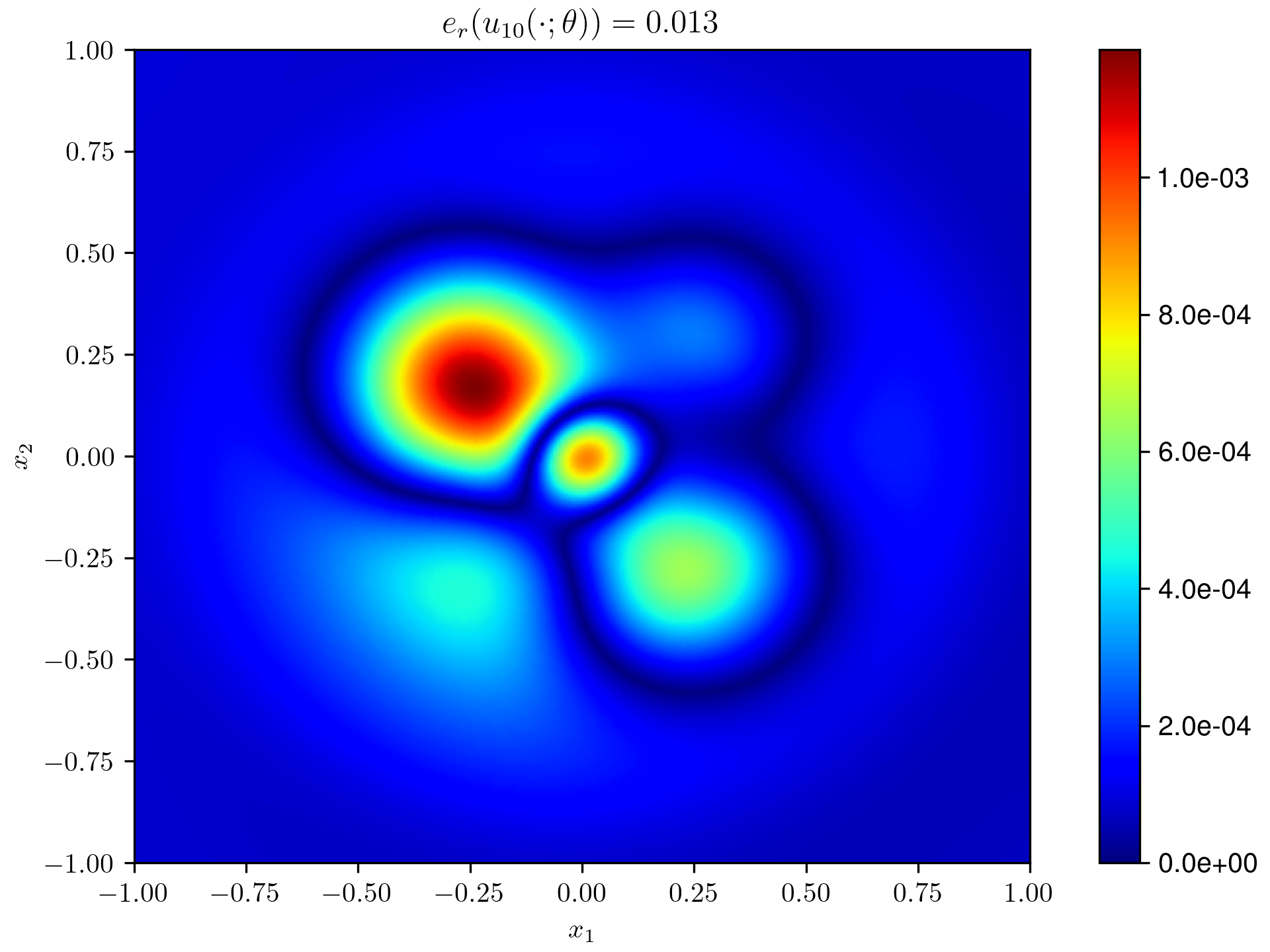}
        \caption{\textit{AAIS-t}, absolute error.}
    \end{subfigure}%
    \caption{Profiles of absolute error and neural network solutions for 15D Poisson equation.}
    \label{fig:Ps15Dsol}
\end{figure}

The relative error, loss with the training process are showed in Figure \ref{fig:PS15DErr}, the residual and the nodes are plotted in Figure \ref{fig:Ps15DNode}, and the solution behaviors followed by absolute errors are listed in Figure \ref{fig:Ps15Dsol}. We could see that \textit{RAD} fails even it could not find the area of large residuals. Our proposed AAIS algorithm could succeed to solve the problem and could sample based on the residual after pre-training. However, when the loss become more sparse, our proposed algorithm also fails to find the residual due to the sparsity of searching points.
\section{Conclusion and future works}\label{sec:conclusion}
In this work, we introduce an Annealed Adaptive Importance Sampling (AAIS) methodology for Physics-Informed Neural Networks (PINNs), which includes Gaussian mixture (\textit{AAIS-g}) and Student's t-distribution (\textit{AAIS-t}) variants. Inspired by the Expectation Maximization (EM) algorithm for finite mixtures, AAIS algorithms are adept at replicating the target function under specific parameter configurations. We have also integrated the Residual-based Adaptive Distribution (RAD) method \cite{LuLu:2023:RAD} and a uniform sampling approach into the PINNs resampling framework, offering a total of four distinct sampling strategies.

Furthermore, we investigate the performance of the four sampling methods within the resampling framework using various forward two-dimensional partial differential equations (PDEs). By observing the increased frequency of residual and absolute error, we directly witness the effectiveness of adaptive sampling compared to the uniformly sampling method (\textit{Uni}). These observations align with the assertion of the empirical Neural Tangent Kernel (NTK) theory \cite{WangYuParis:2022:NTK,Lau2024:pinnacle}, which suggests that PINNs tend to learn the solution of low-frequency parts firstly. Additionally, our proposed AAIS algorithms effectively capture the singularity and sharpness of the residual, yielding results comparable to those obtained with the RAD method.

Moreover, our proposed AAIS algorithms hold significant potential for implementations and applications, a direction we intend to explore in our future studies. They can be readily extended to high-dimensional PDEs and integrated with other training methods outlined in \cite{WangParis:2023:arxiv:PINN}. Additionally, The performance of the AAIS algorithms in inverse problems needs to be validated. The reason for the limitation of adaptive sampling methods in some PDEs (such as the KdV equation \eqref{pde:KdV}) also requires further investigation. Moving forward, we aim to exploit the advantages of the AAIS algorithm, including parameter tuning to strike a balance between high-quality mimicking behaviors and low computational cost, exploring parallel computing architectures, and experimenting with more complex mixtures.

\appendix
\section{Time-related PDEs}
{In this appendix, we extend our examination to time-dependent PDEs, necessitating a temporal-spatial interpretation of the domain  $\Omega$.  The adaptive efficiency of our proposed methods is particularly evident in the context of the Burgers' equation and the Allen-Cahn equation, with the latter benefiting from a weighted adjustment. However, for the Korteweg-de Vries (KdV) equation, the advantages of adaptive sampling were not as pronounced which may need further research.}
\label{sec:Appendix}
\subsection{Burgers' equation}
\begin{figure}[htbp]
    \centering
    \includegraphics[scale=0.125]{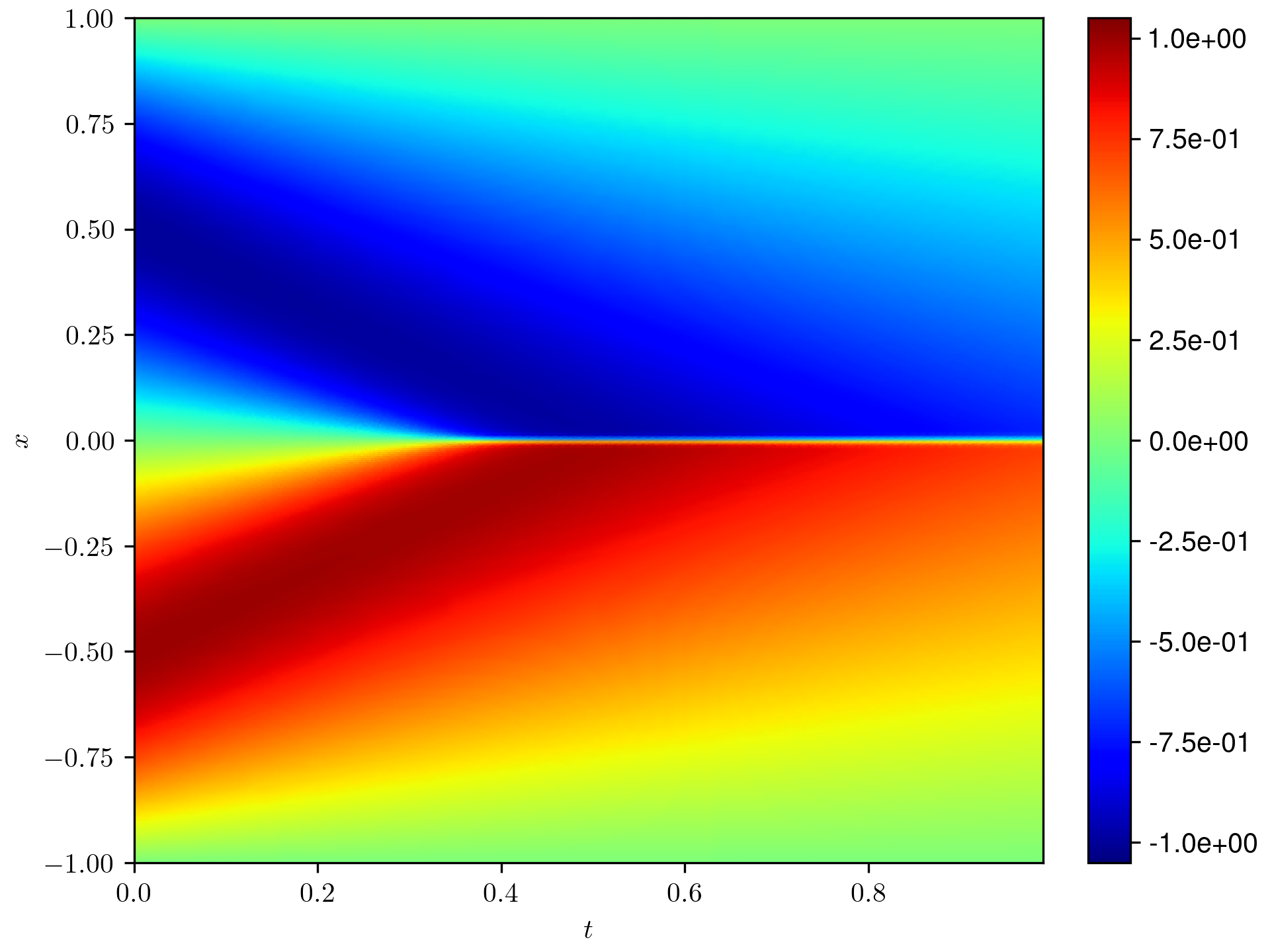}
    \caption{Reference solution for Burgers' equations in \eqref{pde:Burgers}.}
    \label{fig:Burgersref}
\end{figure}

We consider the following Burgers' equation(reference solution is given in Figure \ref{fig:Burgersref})
\begin{equation}
    \label{pde:Burgers}
    \begin{aligned}
        &\partial_t u+u\partial_x u-\frac{0.01}{\pi}\partial_{xx}u=0,~~(t,x)\in(0,1)\times(-1,1),\\
        &u(t,-1)=u(t,1)=0,~~t\in(0,1),\\
        &u(0,x) = -\sin(\pi x), ~~x\in(-1,1).
    \end{aligned}
\end{equation}
It is well-known that there is a huge gap around $x=0$ when $t\to 1^-$, the sharpness implies the efficiency of adaptive sampling methods.

\begin{figure}[htbp]
    \centering
    \begin{subfigure}{.5\textwidth}
    \centering
    \includegraphics[scale=0.25]{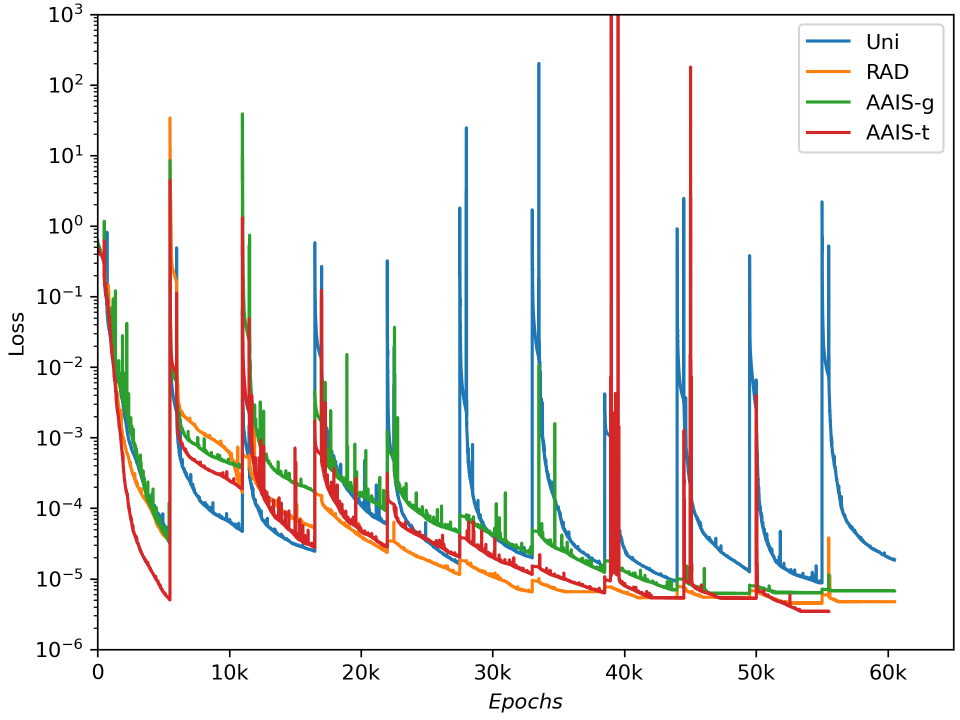}
    \end{subfigure}%
    \begin{subfigure}{.5\textwidth}
    \centering
    \includegraphics[scale=0.25]{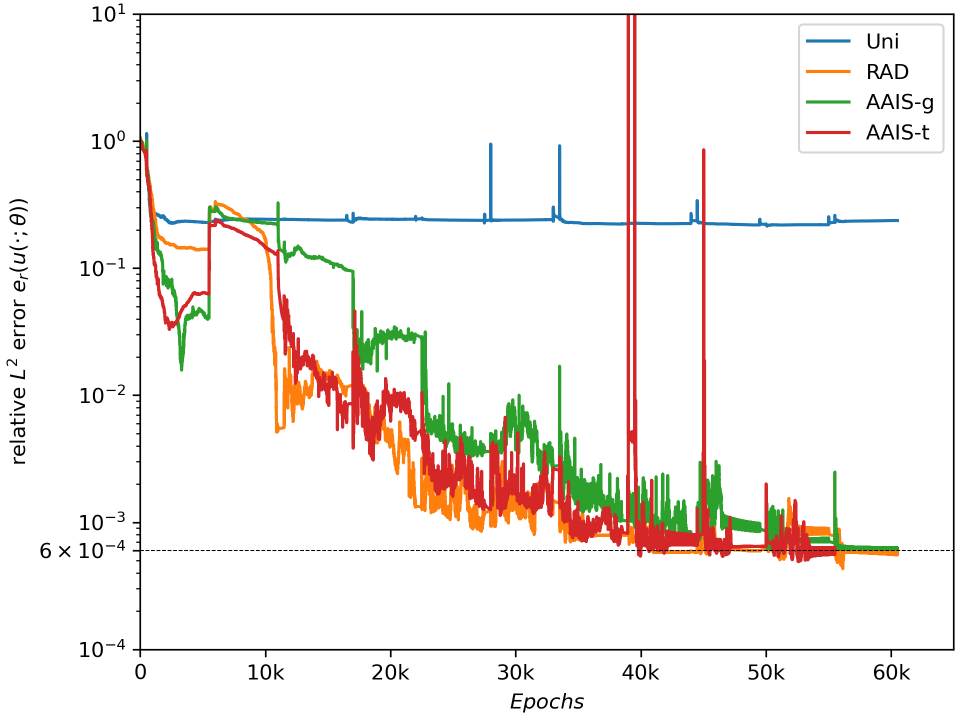}
    \end{subfigure}%
    \caption{Loss and relative errors during training for Burgers' equation with four sampling methods. Left: the loss function. Right: the relative $L^2$ error $e_r(u(\cdot;\theta))$.
    }
    \label{fig:BurgersLossErr5000e}
\end{figure}

\begin{figure}[htbp]
    \centering
    \begin{subfigure}{.25\textwidth}
        \centering
        \includegraphics[height=0.75\textwidth,width=1.0\textwidth]{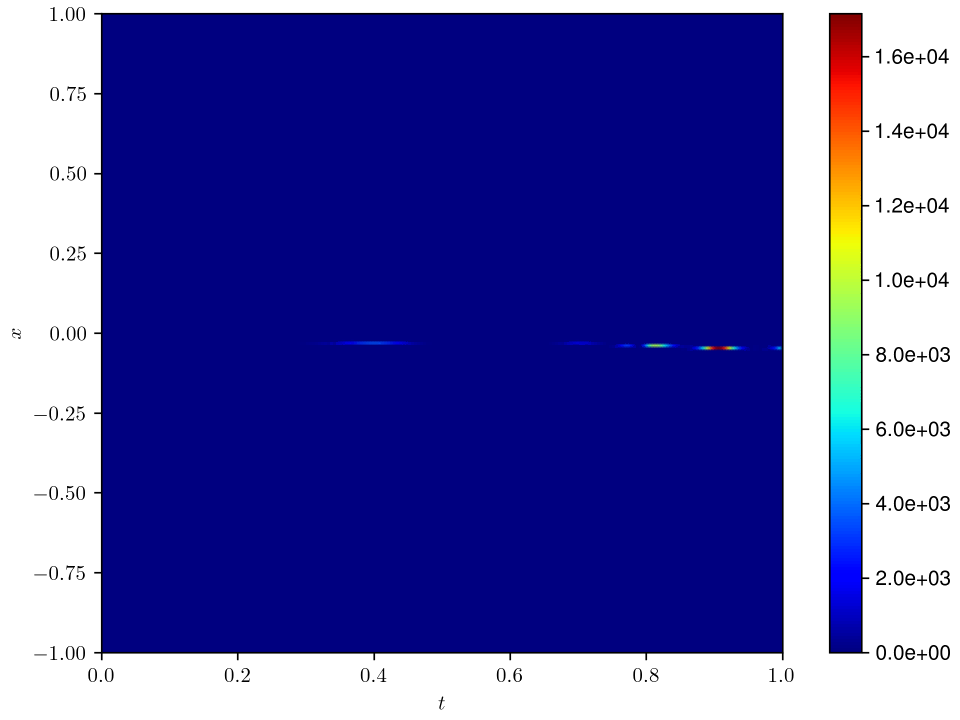}
    \end{subfigure}%
    \begin{subfigure}{.25\textwidth}
        \centering
        \includegraphics[height=0.75\textwidth,width=1.0\textwidth]{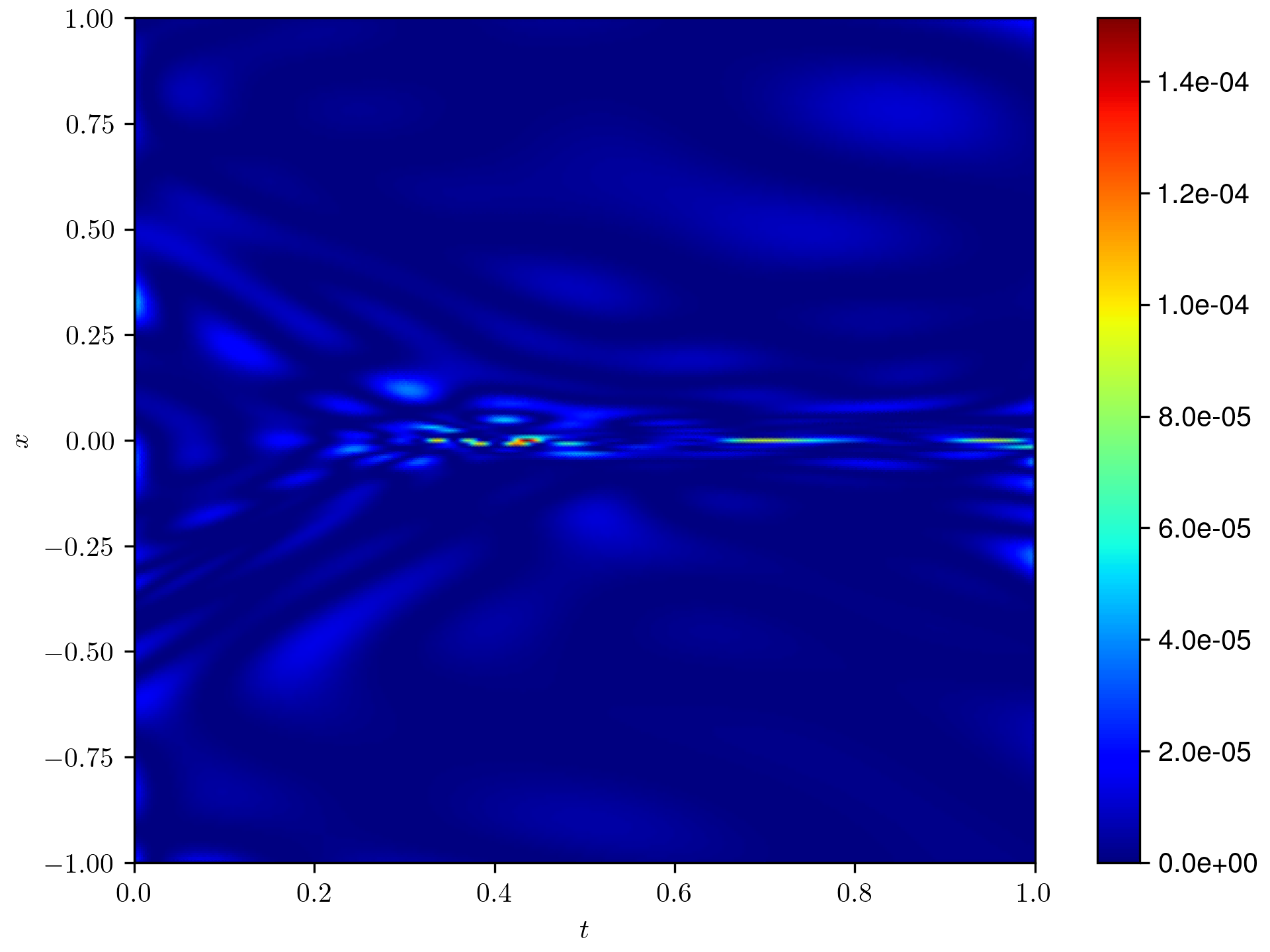}
    \end{subfigure}%
    \begin{subfigure}{.25\textwidth}
        \centering
        \includegraphics[height=0.75\textwidth,width=1.0\textwidth]{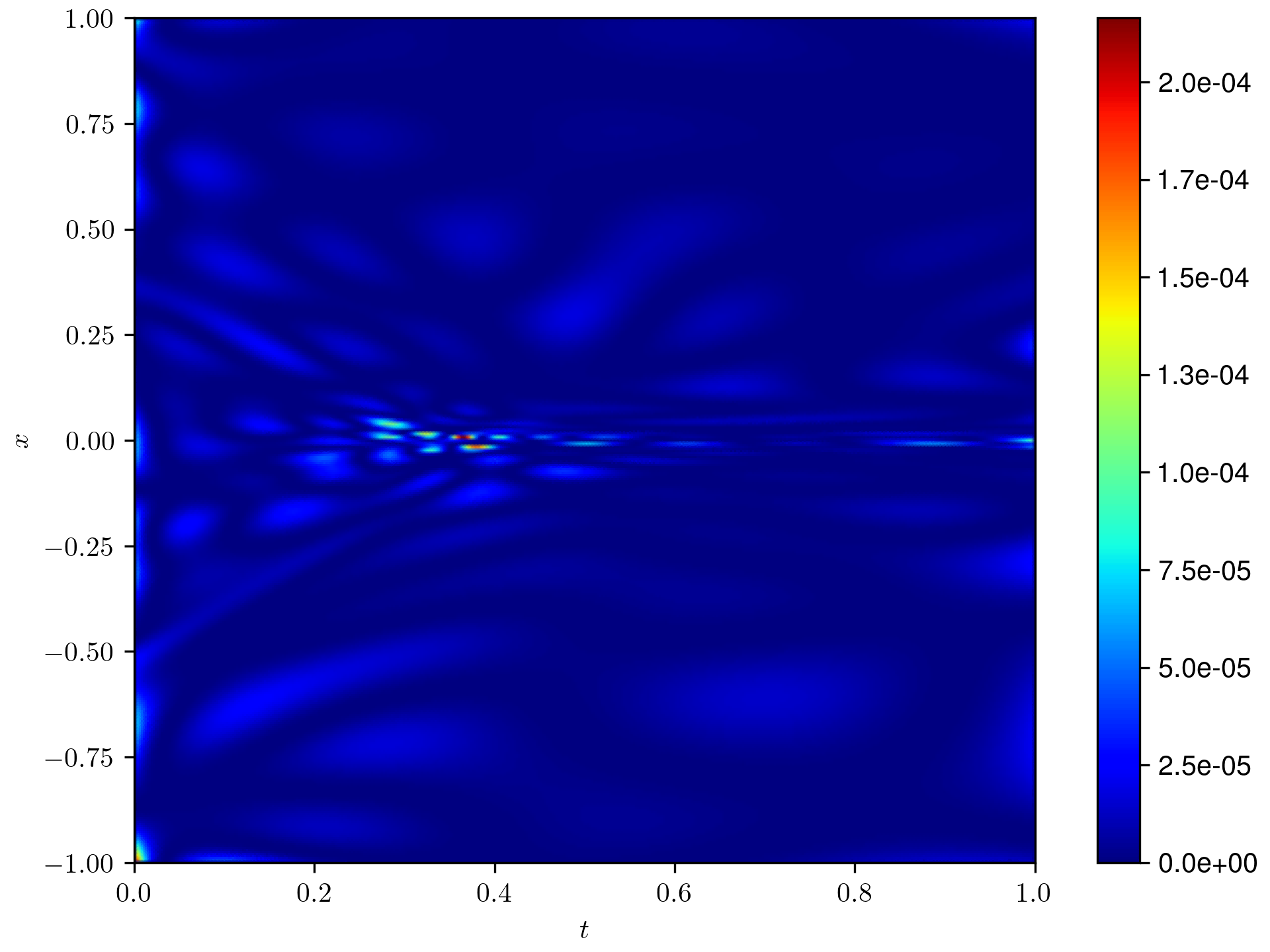}
    \end{subfigure}%
    \begin{subfigure}{.25\textwidth}
        \centering
        \includegraphics[height=0.75\textwidth,width=1.0\textwidth]{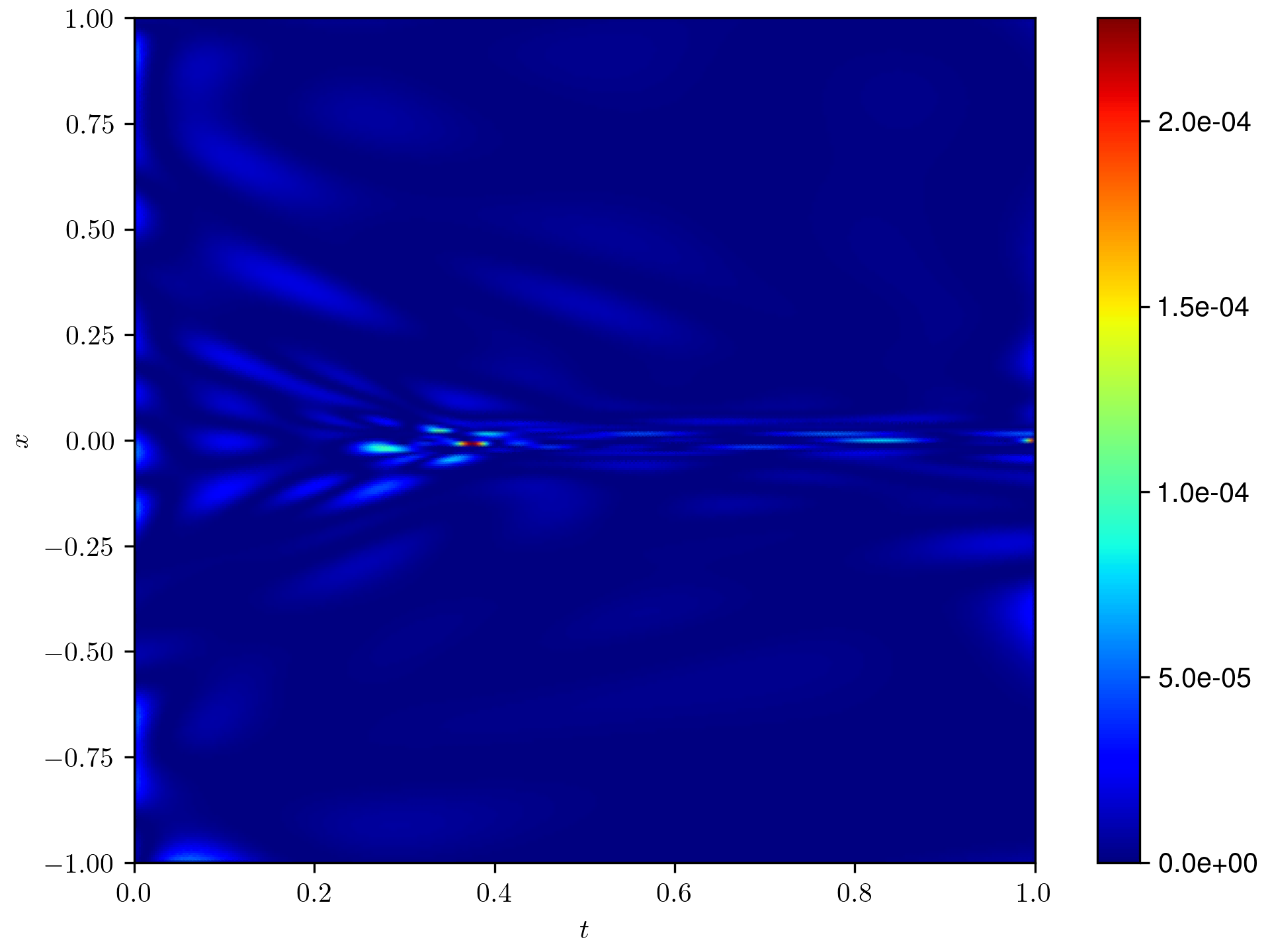}
    \end{subfigure}%
    \newline
    \raggedleft
    \begin{subfigure}{.25\textwidth}
        \centering
        \includegraphics[height=0.75\textwidth,width=1.0\textwidth]{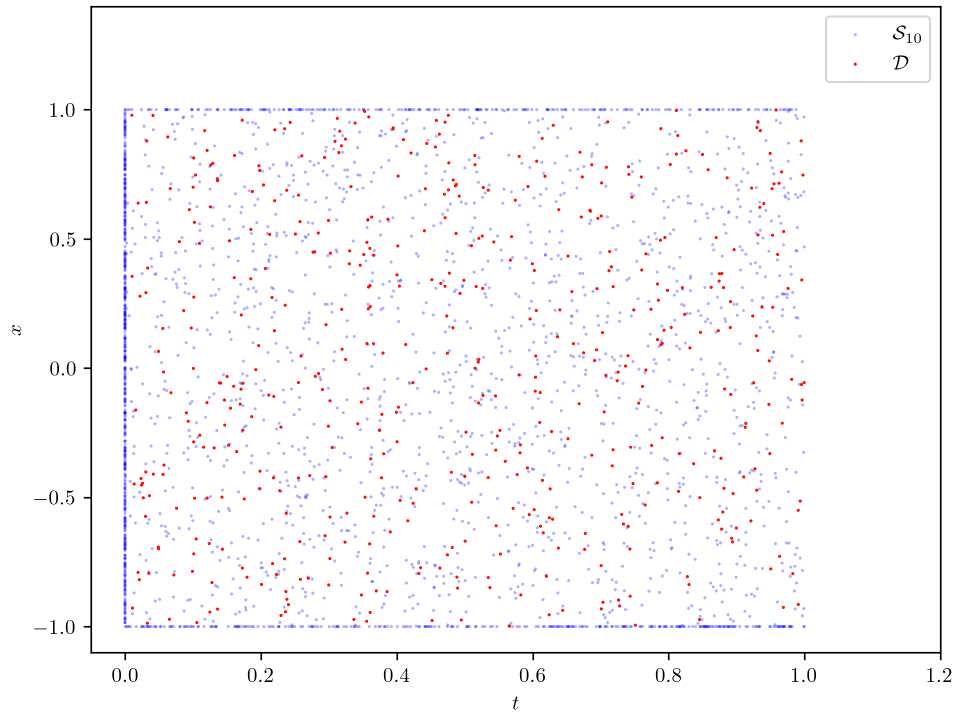}
        \caption{\textit{Uni}}
    \end{subfigure}%
    \begin{subfigure}{.25\textwidth}
        \centering
        \includegraphics[height=0.75\textwidth,width=1.0\textwidth]{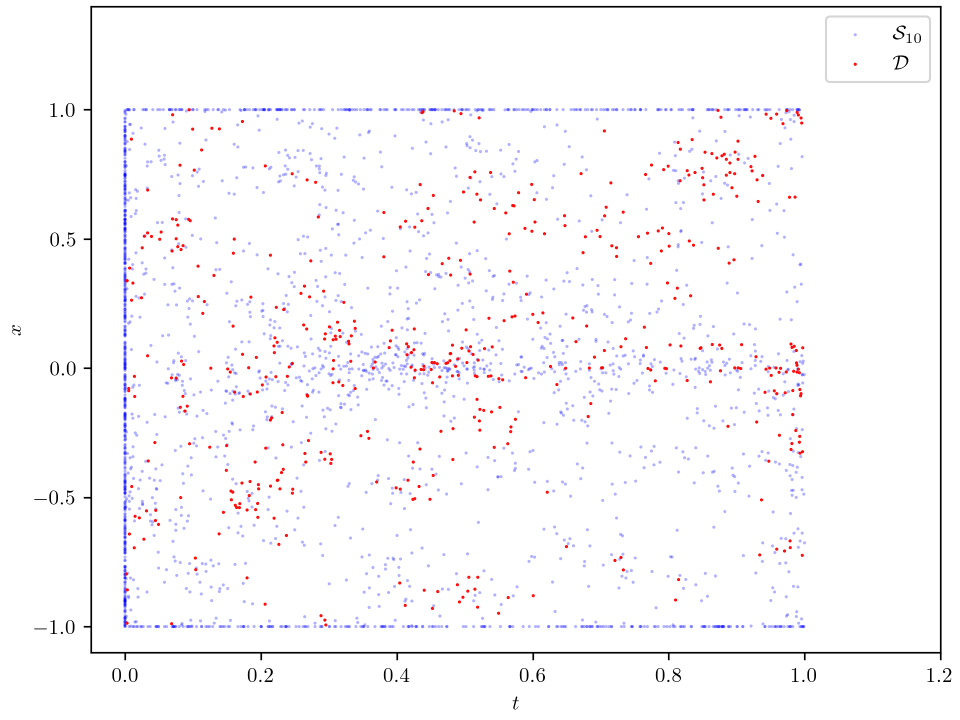}
        \caption{\textit{RAD}}
    \end{subfigure}%
    \begin{subfigure}{.25\textwidth}
        \centering
        \includegraphics[height=0.75\textwidth,width=1.0\textwidth]{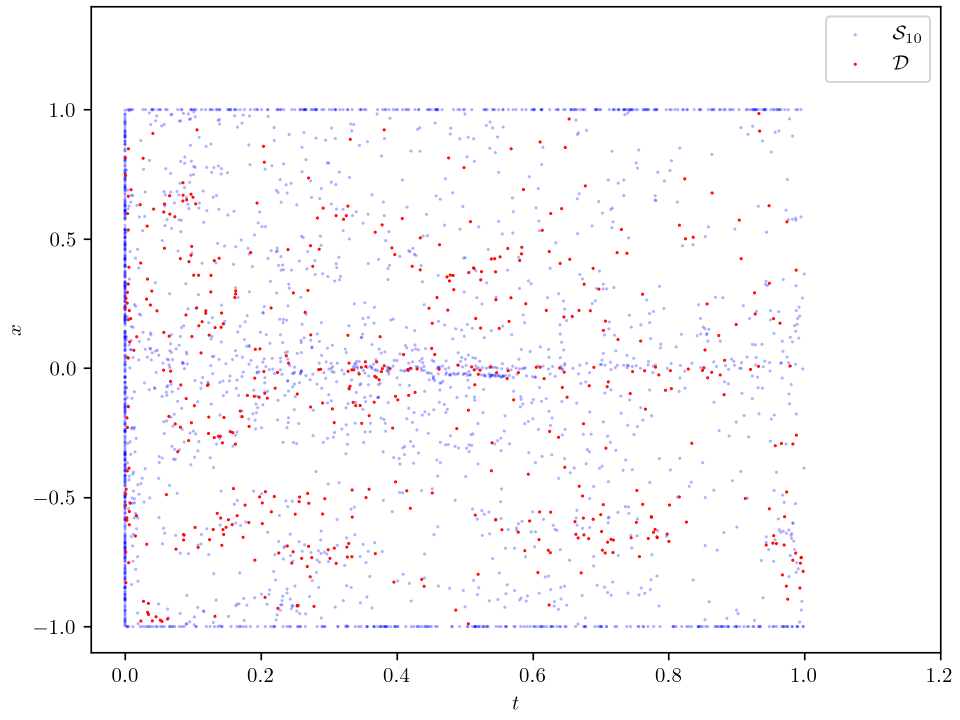}
        \caption{\textit{AAIS-g}}
    \end{subfigure}%
    \begin{subfigure}{.25\textwidth}
        \centering
        \includegraphics[height=0.75\textwidth,width=1.0\textwidth]{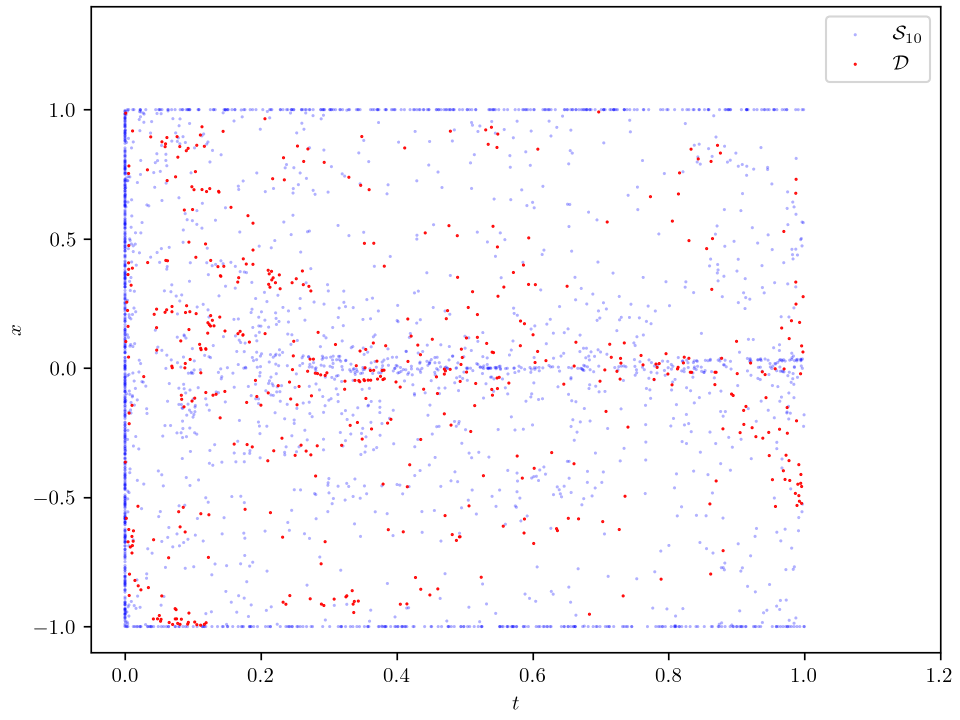}
        \caption{\textit{AAIS-t}}
    \end{subfigure}%
    \caption{Profiles of residual and nodes  for Burgers' equation after 10th training.}
    \label{fig:BurgersNode5000e}
\end{figure}

\begin{figure}[htbp]
    \centering
    \begin{subfigure}{.25\textwidth}
        \centering
        \includegraphics[height=0.75\textwidth,width=1.0\textwidth]{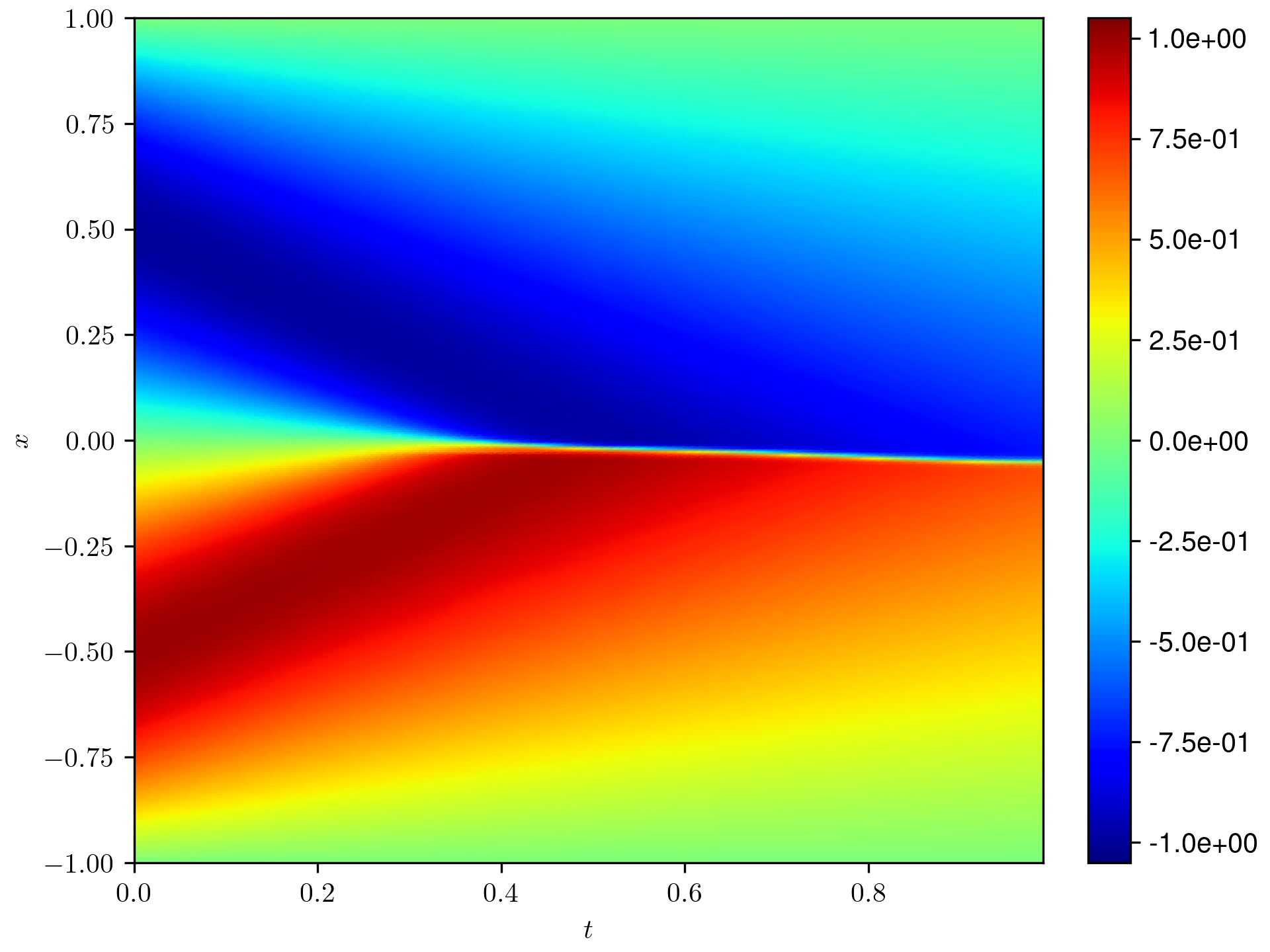}
    \end{subfigure}%
    \begin{subfigure}{.25\textwidth}
        \centering
        \includegraphics[height=0.75\textwidth,width=1.0\textwidth]{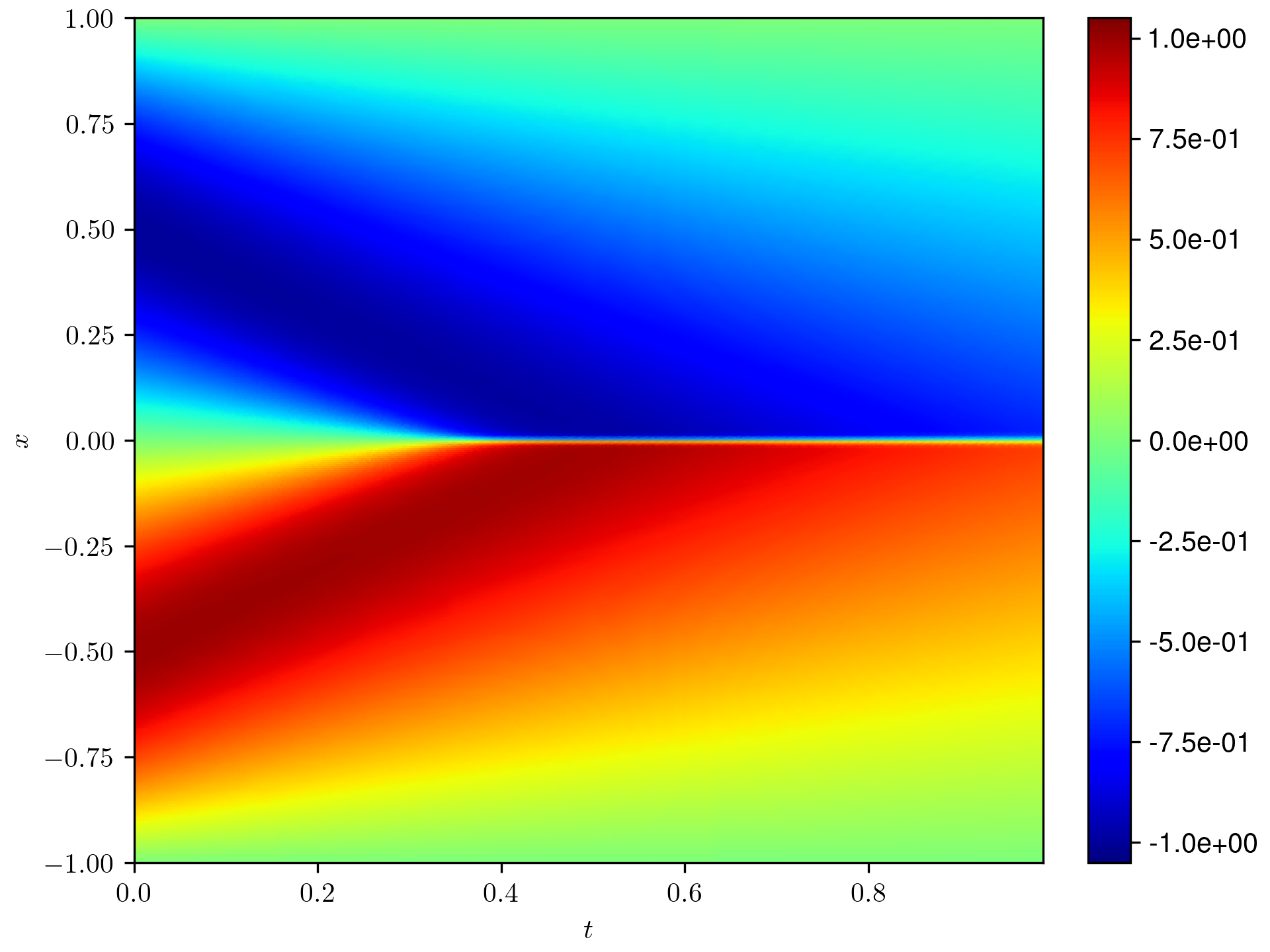}
    \end{subfigure}%
    \begin{subfigure}{.25\textwidth}
        \centering
        \includegraphics[height=0.75\textwidth,width=1.0\textwidth]{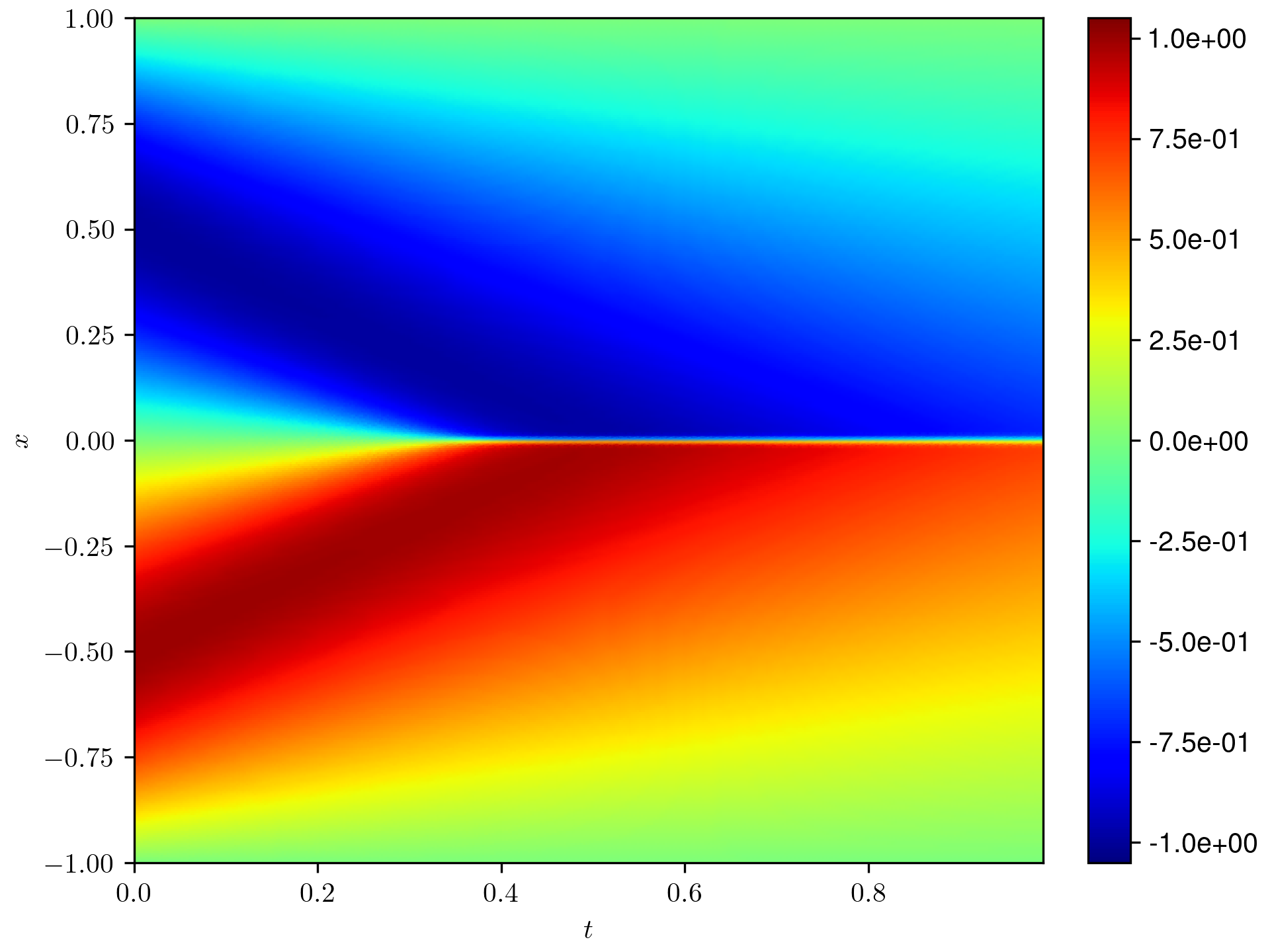}
    \end{subfigure}%
    \begin{subfigure}{.25\textwidth}
        \centering
        \includegraphics[height=0.75\textwidth,width=1.0\textwidth]{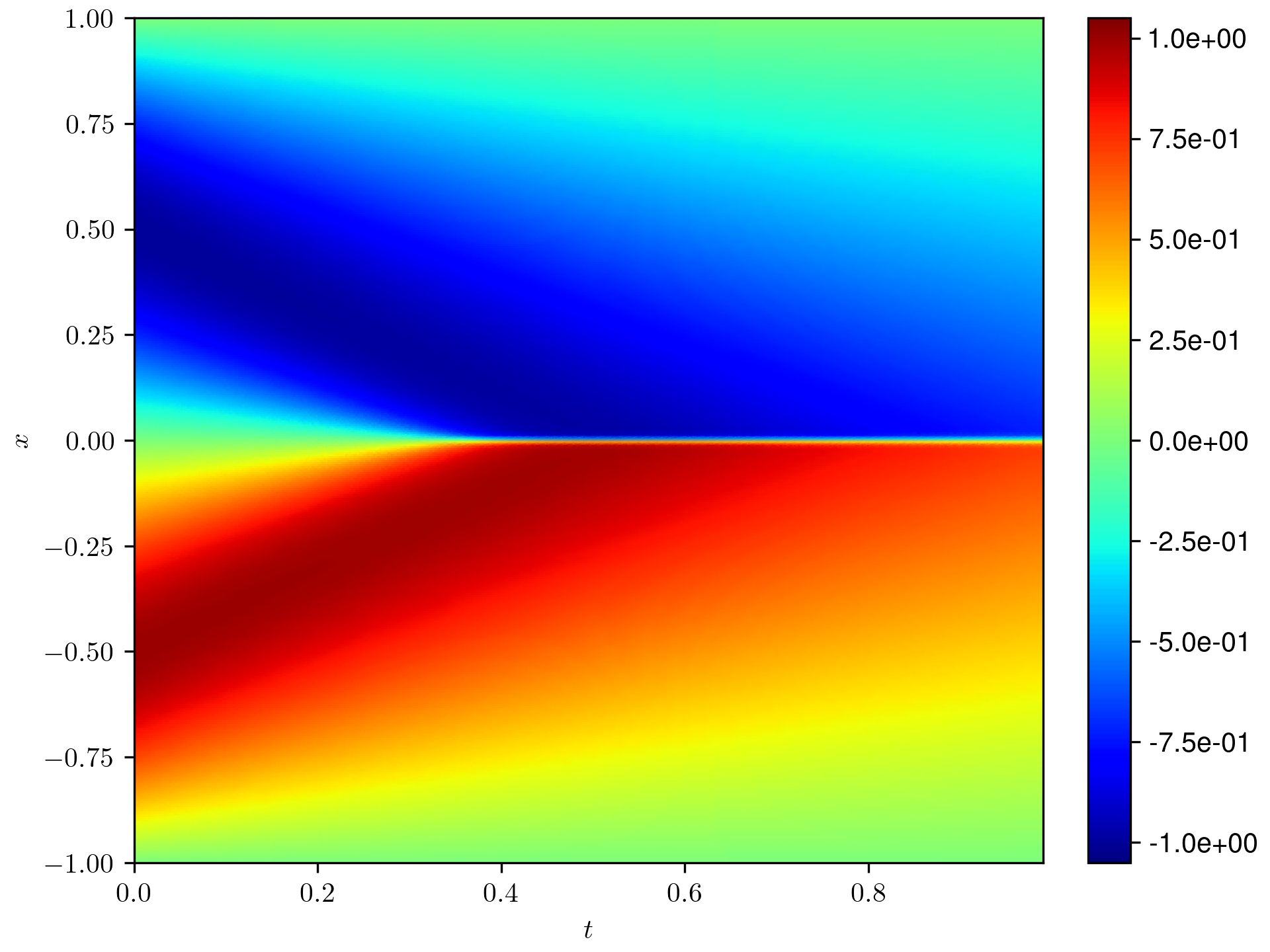}
    \end{subfigure}%
    \newline
    \raggedleft
    \begin{subfigure}{.25\textwidth}
        \centering
        \includegraphics[height=0.75\textwidth,width=1.0\textwidth]{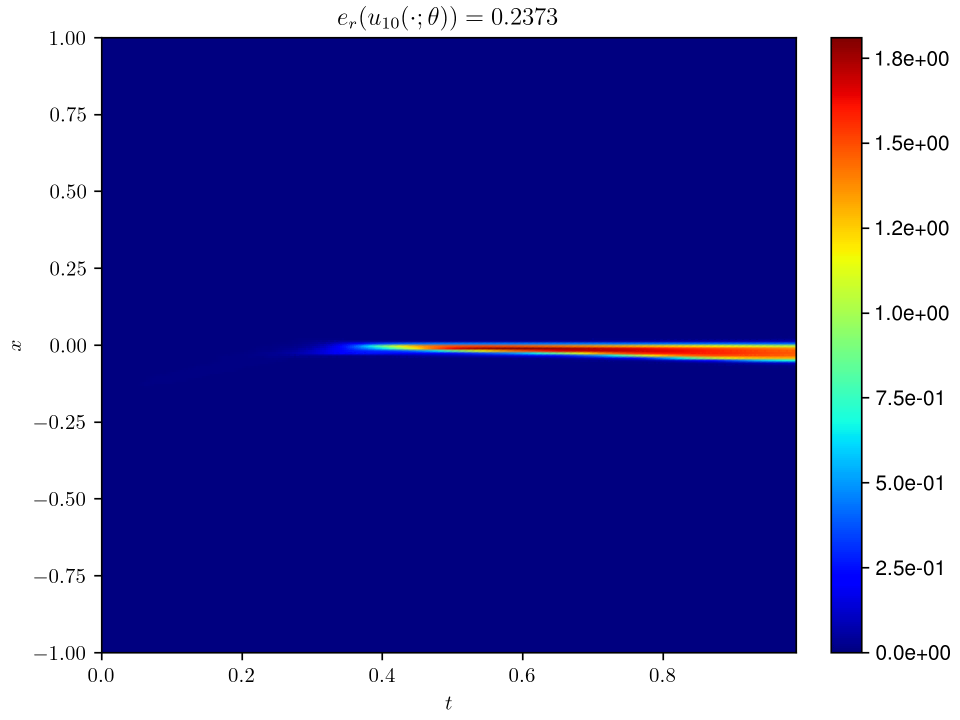}
        \caption{\textit{Uni}}
    \end{subfigure}%
    \begin{subfigure}{.25\textwidth}
        \centering
        \includegraphics[height=0.75\textwidth,width=1.0\textwidth]{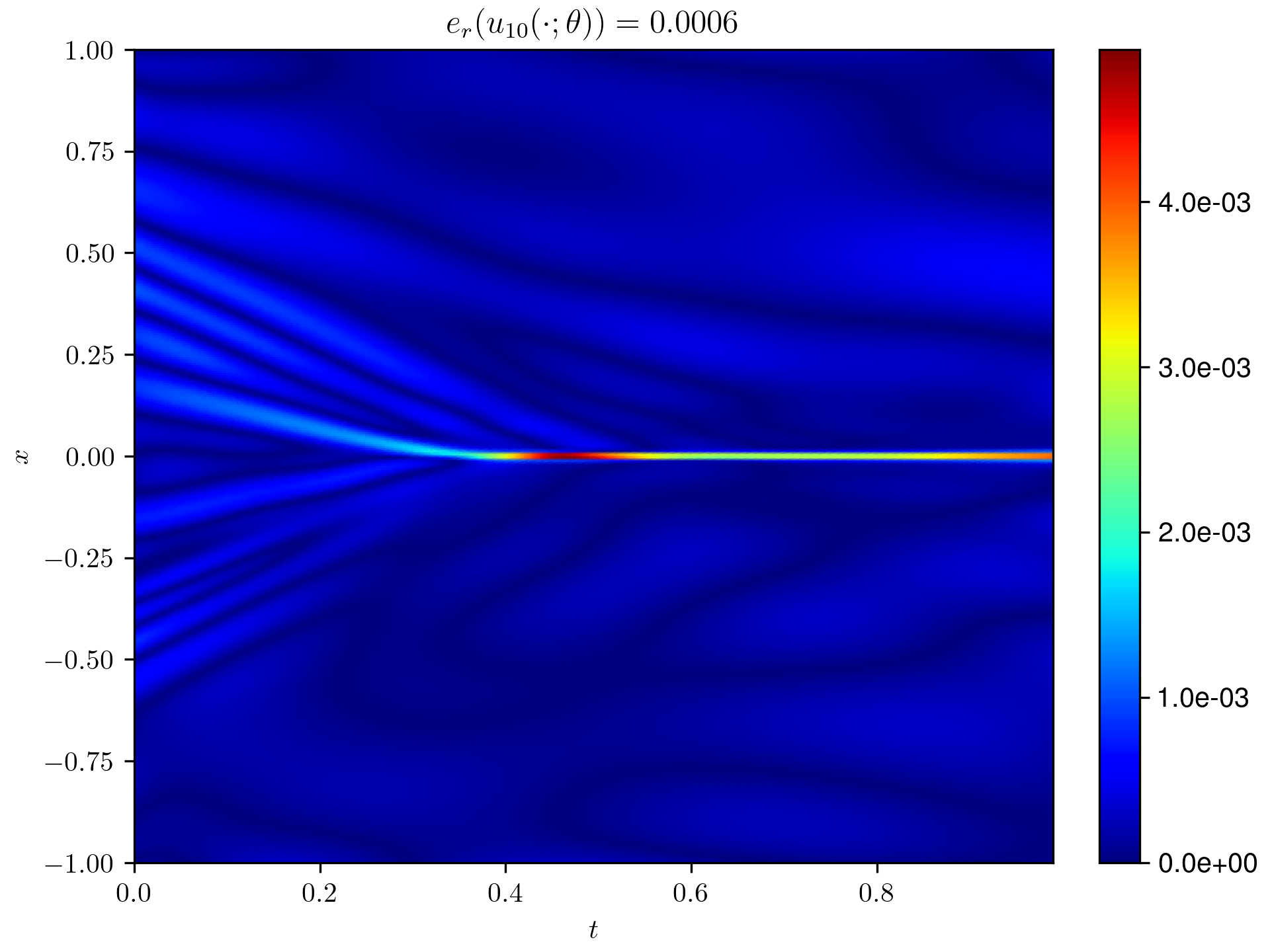}
        \caption{\textit{RAD}}
    \end{subfigure}%
    \begin{subfigure}{.25\textwidth}
        \centering
        \includegraphics[height=0.75\textwidth,width=1.0\textwidth]{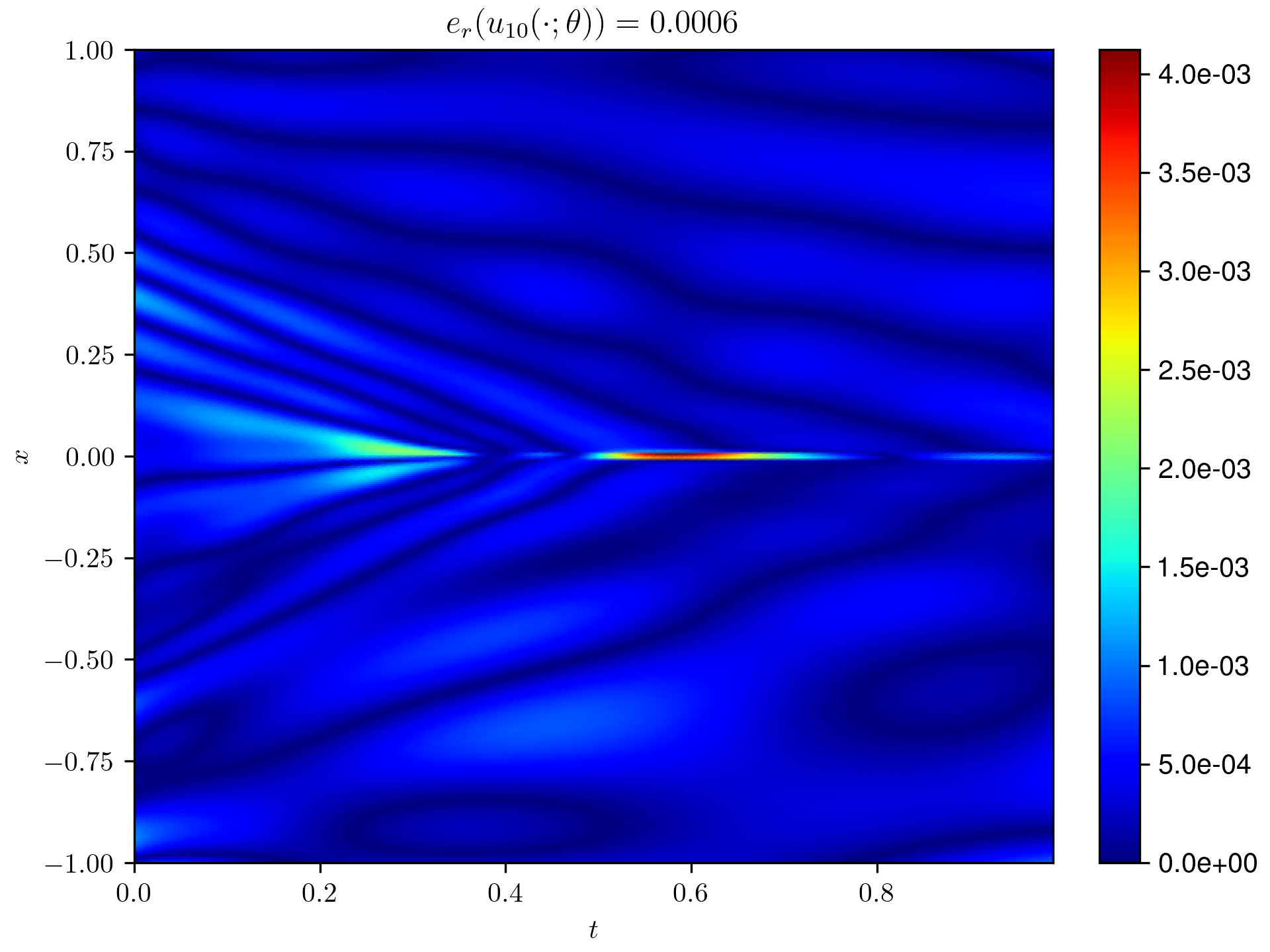}
        \caption{\textit{AAIS-g}}
    \end{subfigure}%
    \begin{subfigure}{.25\textwidth}
        \centering
        \includegraphics[height=0.75\textwidth,width=1.0\textwidth]{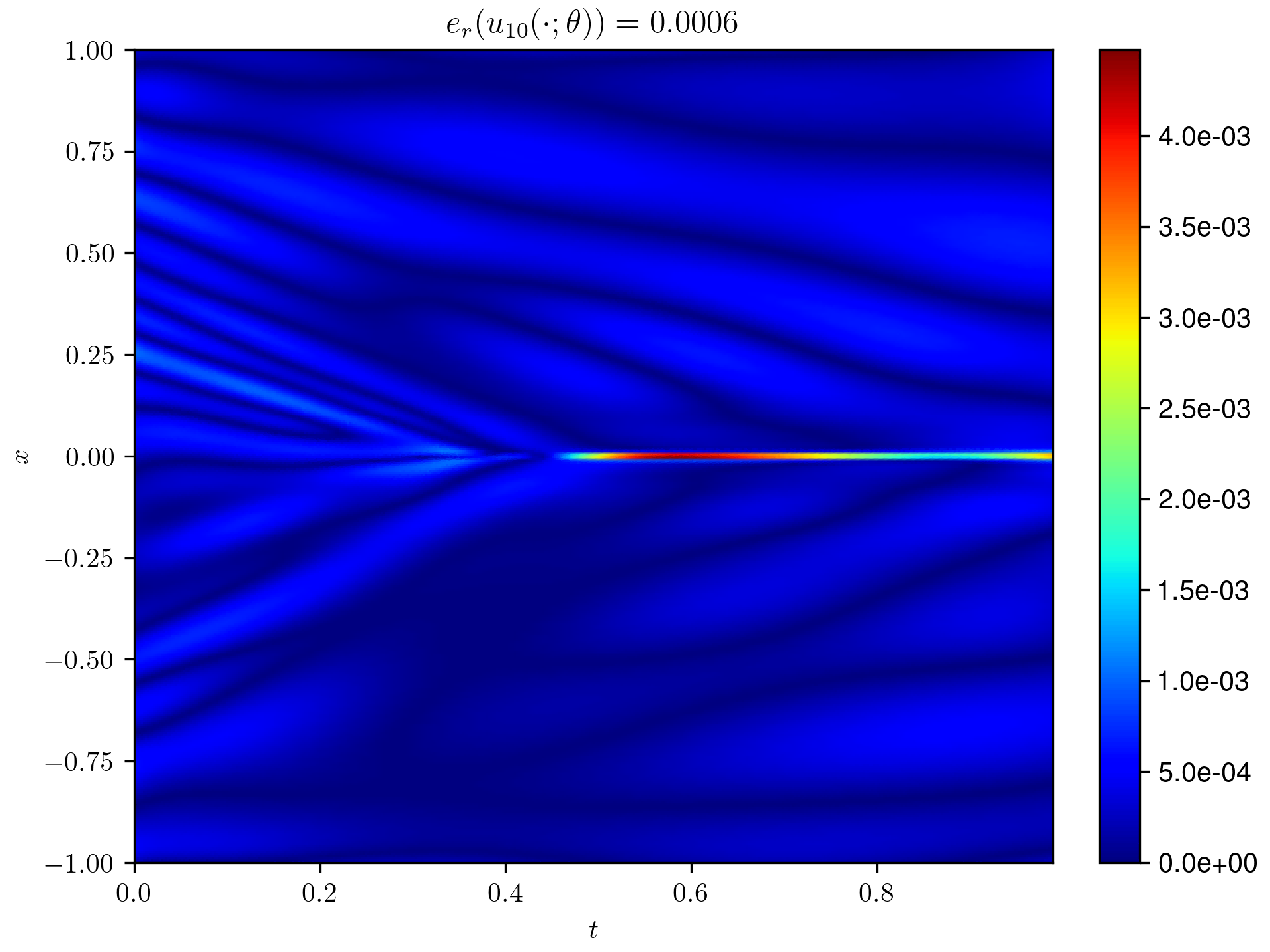}
        \caption{\textit{AAIS-t}}
    \end{subfigure}%
    \caption{Profiles of absolute error and neural network solutions for Burgers' equation after 10th training. }
    \label{fig:BurgersErr5000e}
\end{figure}

Firstly we set maximum iteration $M=10$, 500 epochs for Adam and 5000 epochs for lbfgs with learning rate 1.0 during pre-training and each iteration. $N_A$ for AAIS is 4000. The loss and relative errors during training are showed in Figure \ref{fig:BurgersLossErr5000e}, proposing that adaptive sampling methods could generate better solutions than the traditional PINNs. Figure \ref{fig:BurgersNode5000e} list the residual $\Q$ and nodes after each iteration. It is showed that adaptive sampling methods would focus more on the place $x=0$, leading to the increasing frequency of residual. The solution profiles and absolute errors in Figure \ref{fig:BurgersErr5000e} also support the statement that increasing frequency of residual implies the better solution behaviors and the increasing frequency of absolute error. The absolute error of \textit{Uni} method mainly concentrate on the sharpness due to the less training there. But the adaptive sampling methods would decrease the error there because of clustering nodes.

\subsection{Allen-Cahn equation}
\begin{figure}[htbp]
    \centering
    \includegraphics[scale=0.125]{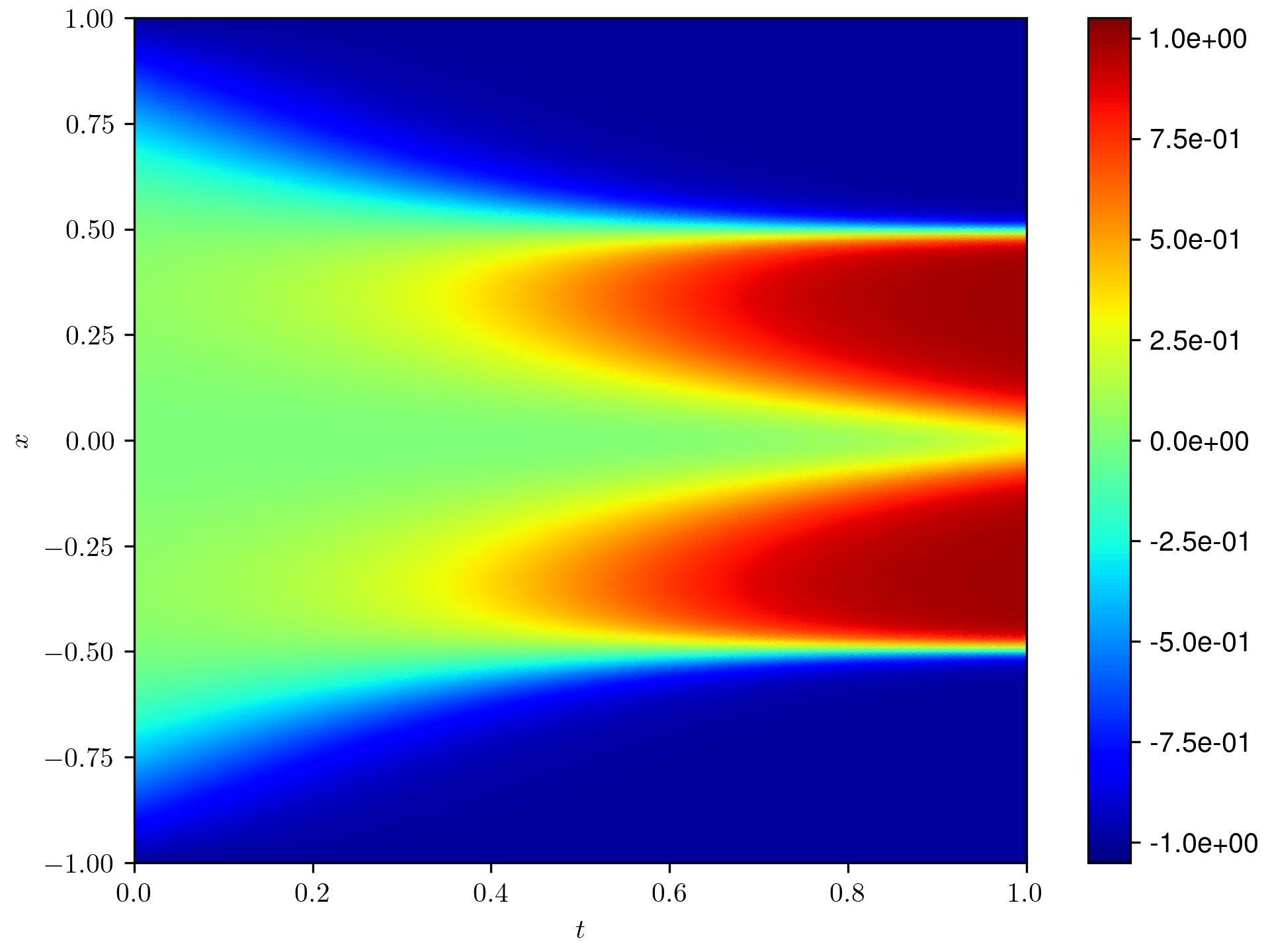}
    \caption{Reference solution for Allen-Cahn equations in \eqref{pde:AllenCahn}.}
    \label{fig:ACexact}
\end{figure}

In this subsection PINNs are used to solve two-dimensional Allen-Cahn equation(reference solution is given in Figure \ref{fig:ACexact}):
\begin{equation}
    \label{pde:AllenCahn}
    \begin{aligned}
        &\partial_t u-0.001\partial_{xx}u-5(u-u^3)=0,~~(t,x)\in(0,1)\times(-1,1),\\
        &u(t,-1)=u(t,1)=-1,~~t\in(0,1),\\
        &u(0,x) = x^2\cos(\pi x), ~~x\in(-1,1).
    \end{aligned}
\end{equation}
Firstly we set the maximum iteration $M=100$, the epochs for training in pre-train and each iteration is 500 for Adam, 2000 for lbfgs. $N_A$ for sampling of AAIS is set to be 6000. The size of training dataset is set to be 2000, where 500 nodes from which are generated from sampling methods.

\begin{figure}[htbp]
    \centering
    \begin{subfigure}{.5\textwidth}
        \centering
        \includegraphics[scale=0.25]{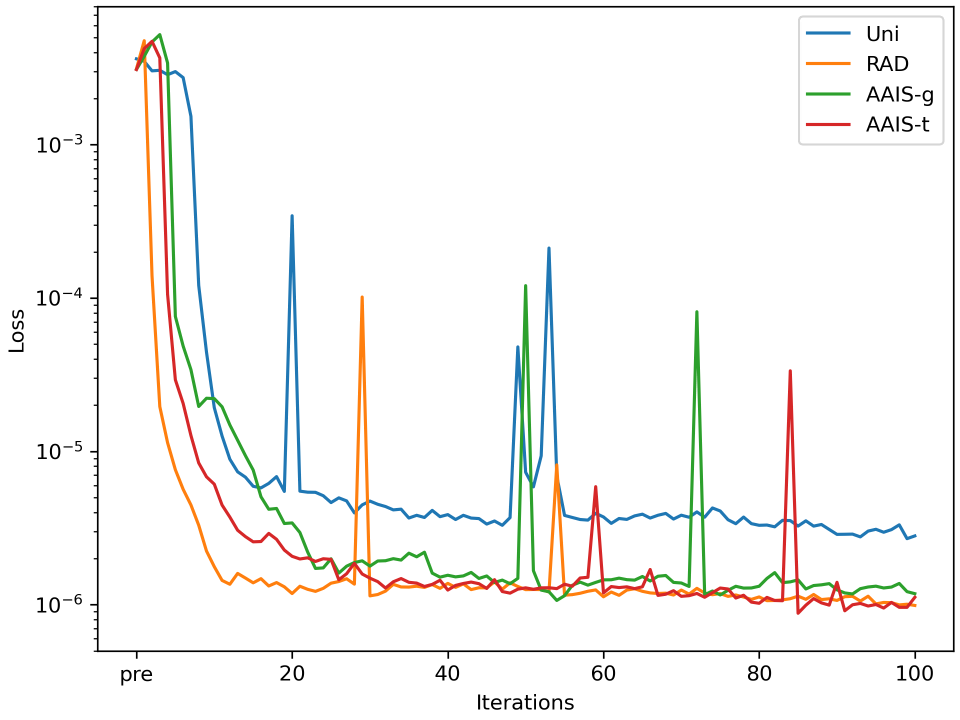}
        \end{subfigure}%
        \begin{subfigure}{.5\textwidth}
        \centering
        \includegraphics[scale=0.25]{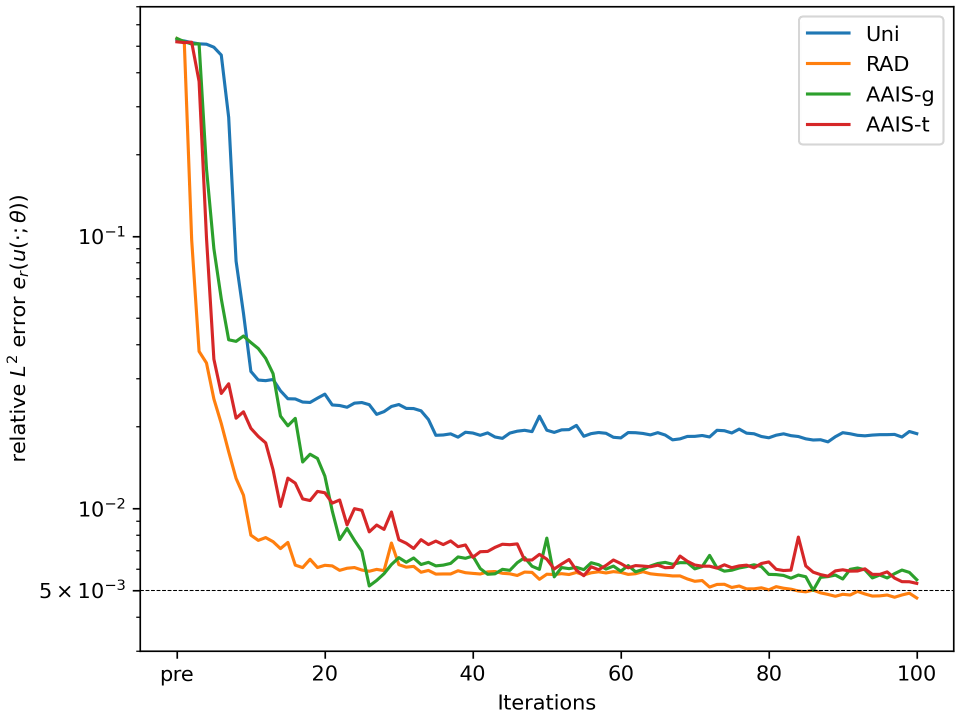}
        \end{subfigure}%
        \caption{Loss and relative errors at each iteration after training for AllenCahn equation with training schedule of lbfgs 2000 epochs each. Left: the loss function. Right: the relative $L^2$ error $e_r(u(\cdot;\theta))$.
        }
    \label{fig:ACErr2000e}
\end{figure}

Differently, all methods would fail to solve the Allen-Cahn equation otherwise we set the weight for the initial loss to be 10, i.e. $\omega_i=10$ for the nodes sampled from initial layer($t=0$). It implies that the information of initial condition should be conveyed firstly to the solution before satisfying the PDE equation. The loss and error after every iteration are plotted in Figure \ref{fig:ACErr2000e}. Unlike causality used in \cite{GaoYanZhou2023:FIPINN:I} and hard constraints used in \cite{LuLu:2023:RAD}, we could see the solution would decay very fast at a small interval then stay almost fixed. And the adaptive sampling methods would obtain better results than \textit{Uni} method.

\begin{figure}[htbp]
    \centering
    \begin{subfigure}{.33\textwidth}
        \centering
        \includegraphics[height=0.75\textwidth,width=1.0\textwidth]{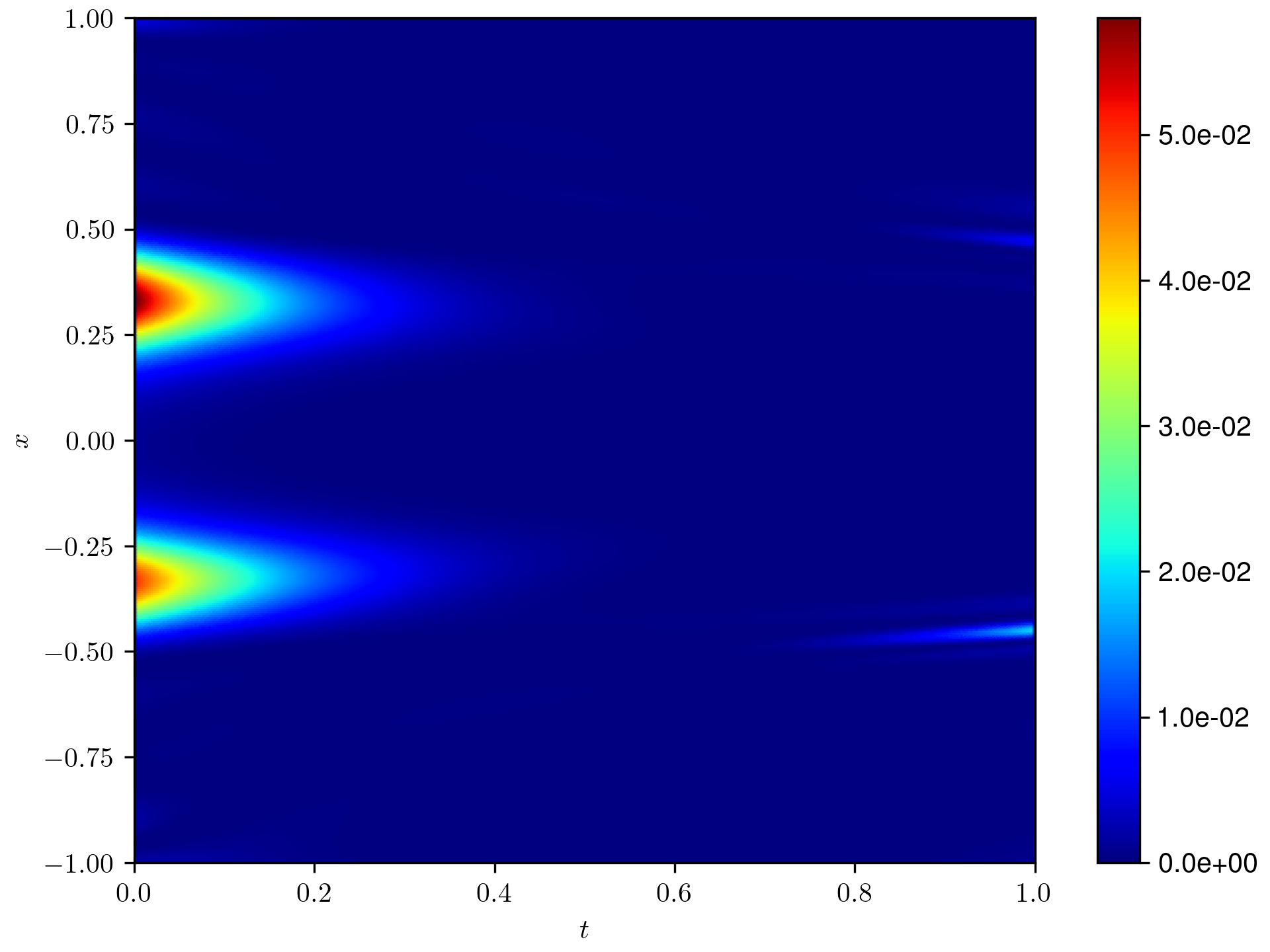}
        \caption{$\Q$ after pre-training of \textit{Uni}.}
    \end{subfigure}%
    \begin{subfigure}{.33\textwidth}
        \centering
        \includegraphics[height=0.75\textwidth,width=1.0\textwidth]{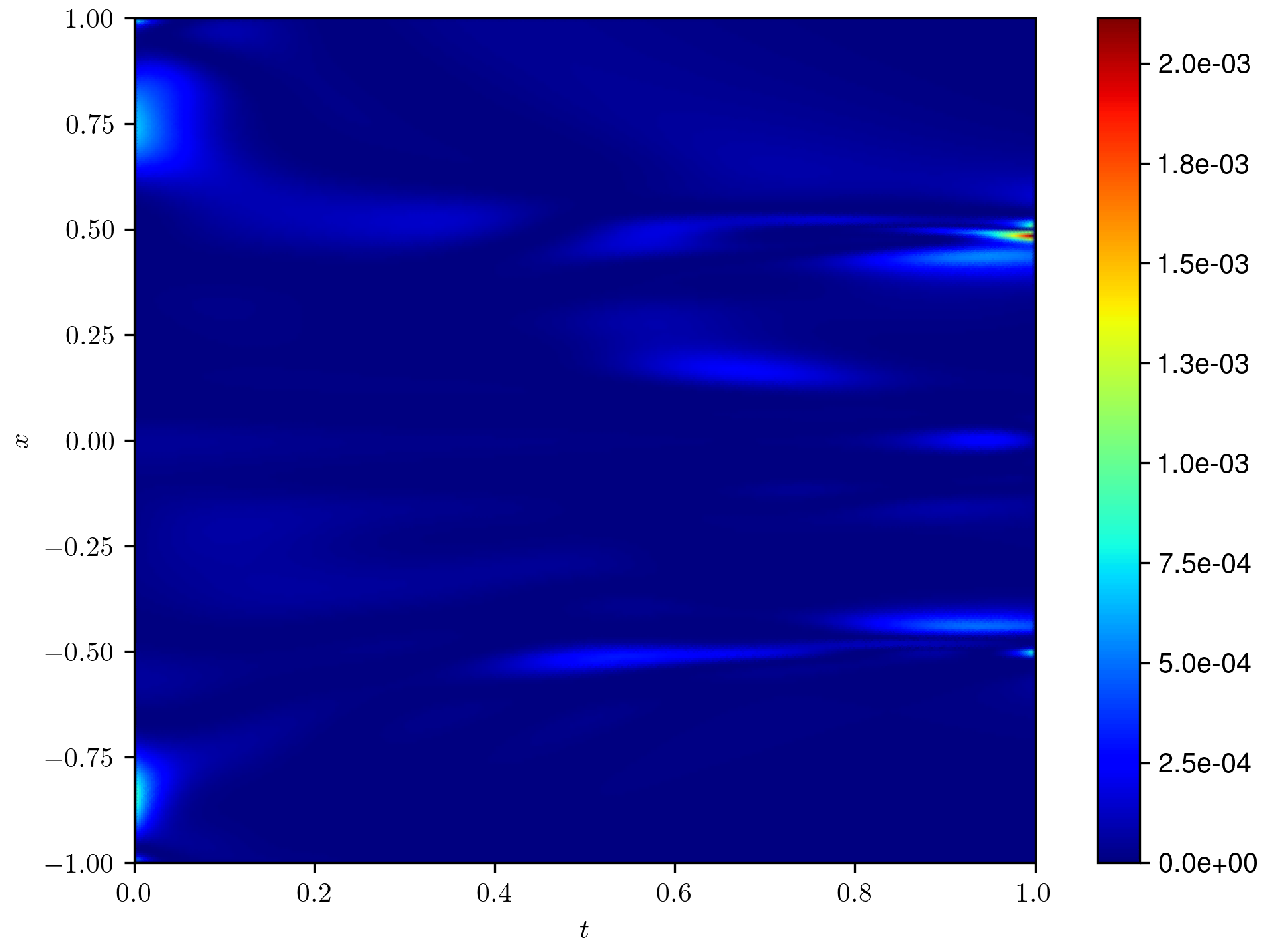}
        \caption{$\Q$ after 49th training of \textit{Uni}.}
    \end{subfigure}%
    \begin{subfigure}{.33\textwidth}
        \centering
        \includegraphics[height=0.75\textwidth,width=1.0\textwidth]{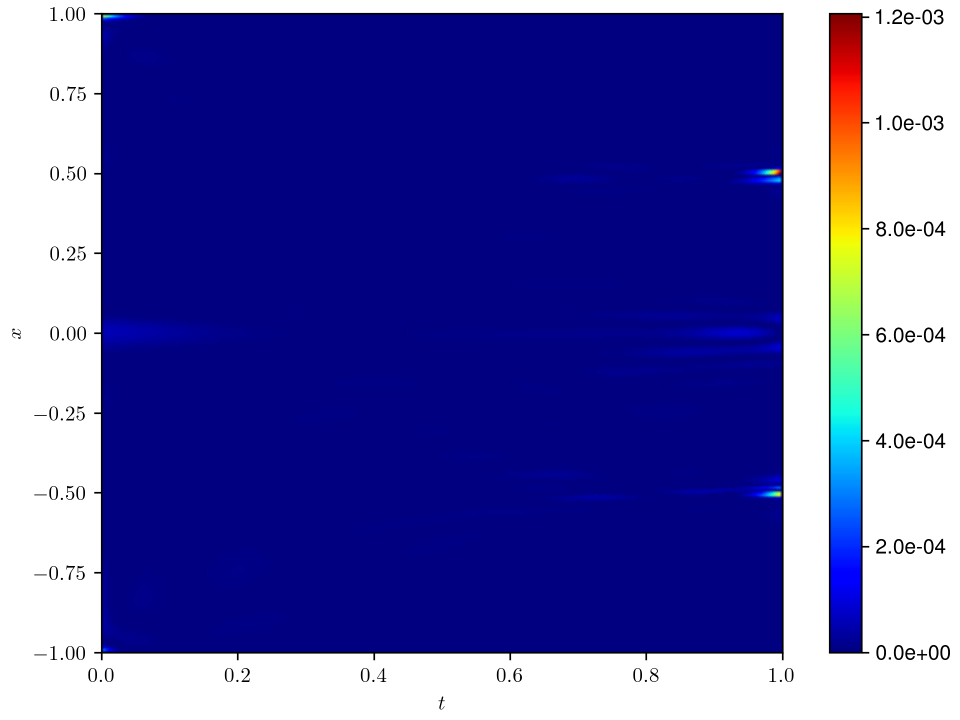}
        \caption{$\Q$ after 99th training of \textit{Uni}.}
    \end{subfigure}%
    \newline
    \raggedleft
    \begin{subfigure}{.33\textwidth}
        \centering
        \includegraphics[height=0.75\textwidth,width=1.0\textwidth]{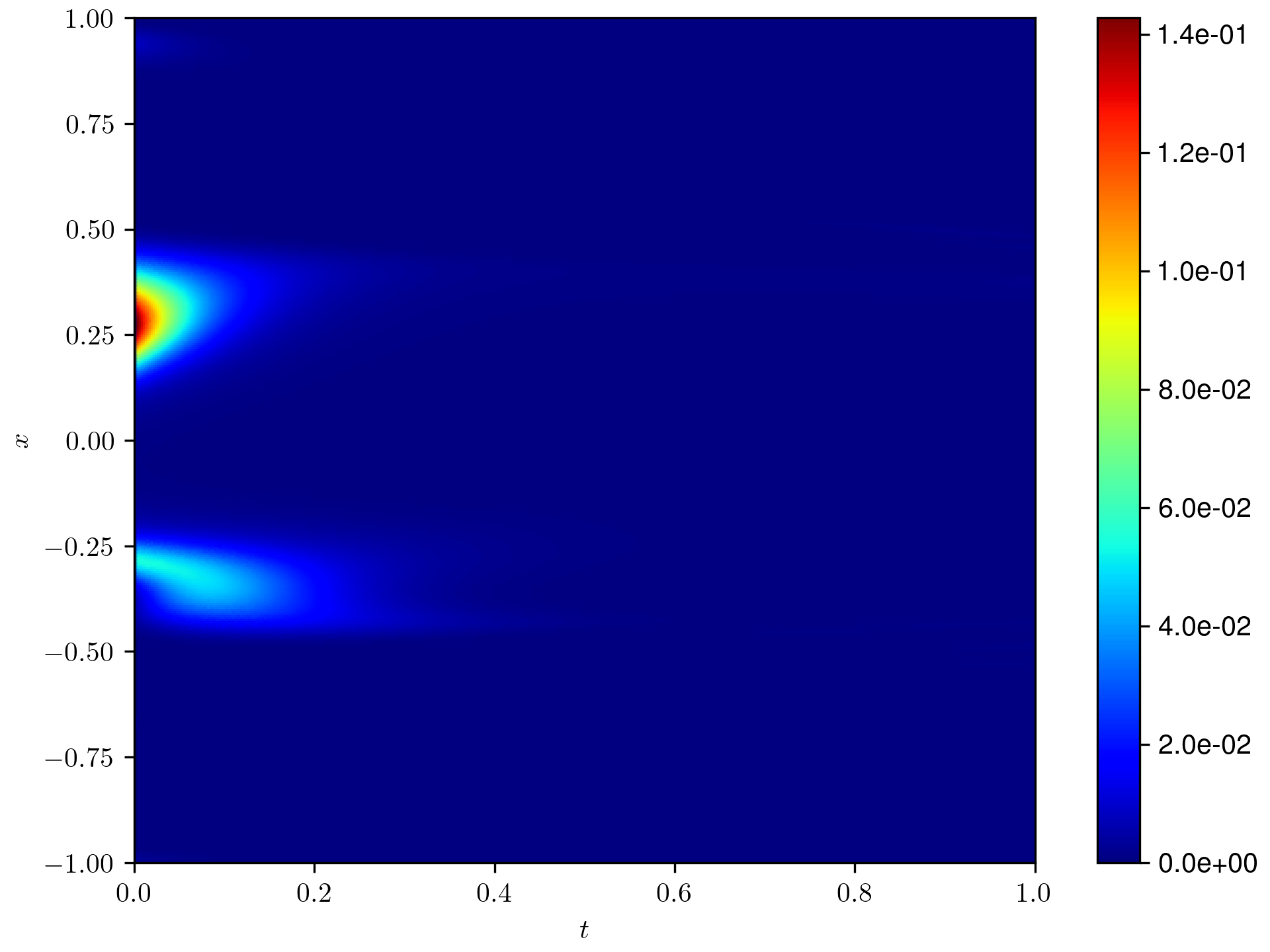}
        \caption{$\Q$ after pre-training of \textit{AAIS-t}.}
    \end{subfigure}%
    \begin{subfigure}{.33\textwidth}
        \centering
        \includegraphics[height=0.75\textwidth,width=1.0\textwidth]{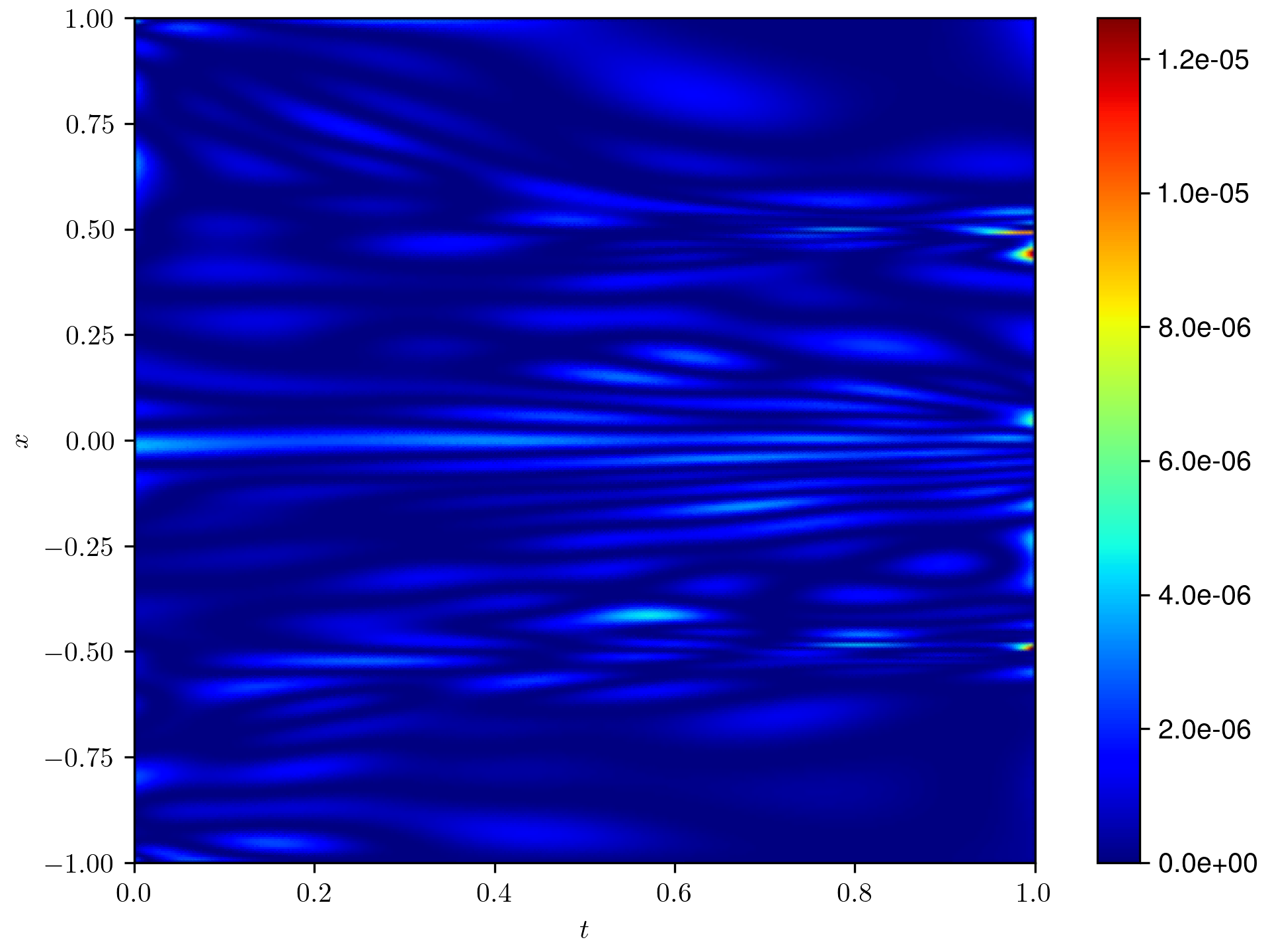}
        \caption{$\Q$ after 49th training of \textit{AAIS-t}.}
    \end{subfigure}%
    \begin{subfigure}{.33\textwidth}
        \centering
        \includegraphics[height=0.75\textwidth,width=1.0\textwidth]{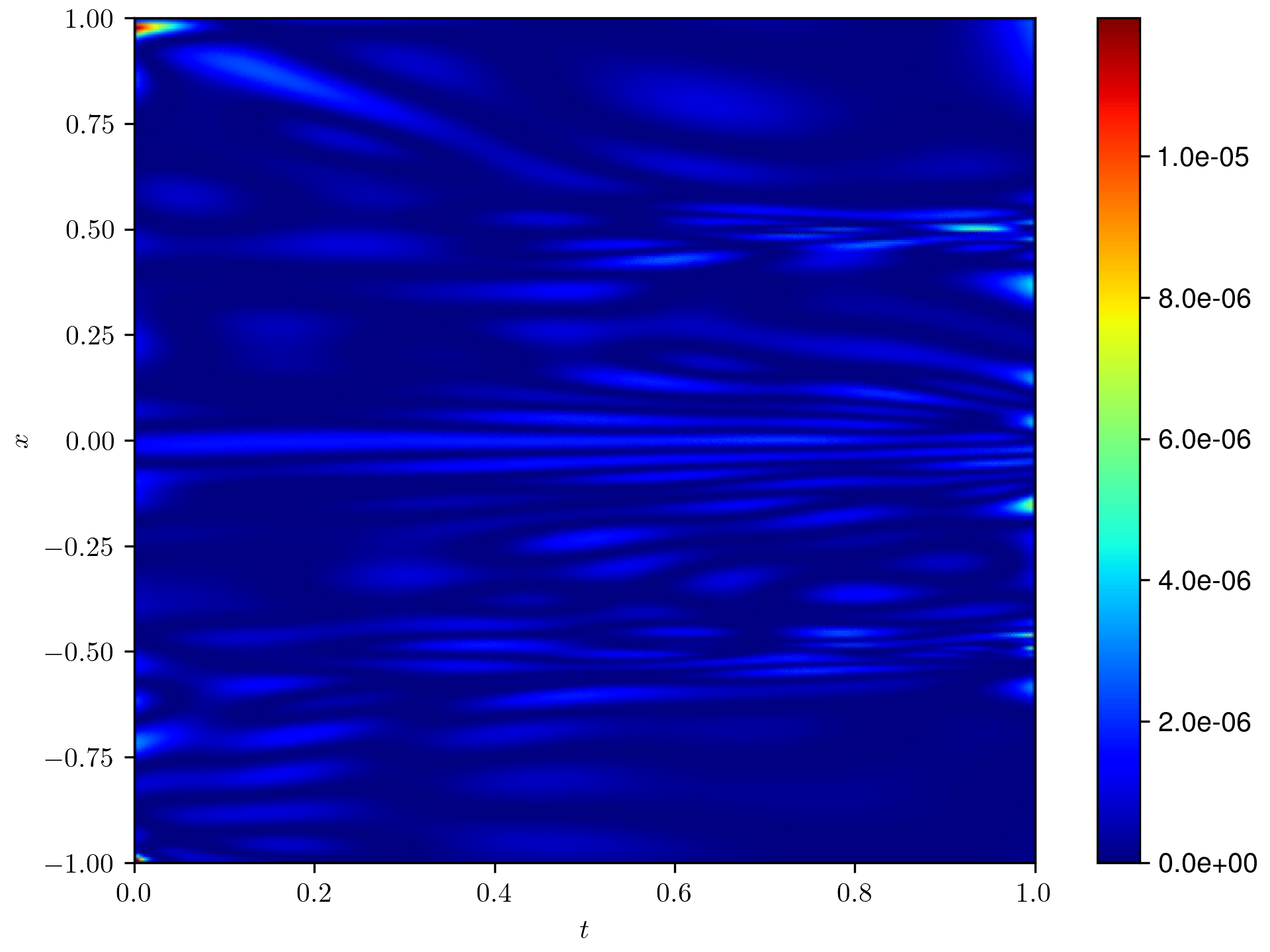}
        \caption{$\Q$ after 99th training of \textit{AAIS-t}.}
    \end{subfigure}%
    \newline
    \raggedleft
    \begin{subfigure}{.33\textwidth}
        \centering
        \includegraphics[height=0.75\textwidth,width=1.0\textwidth]{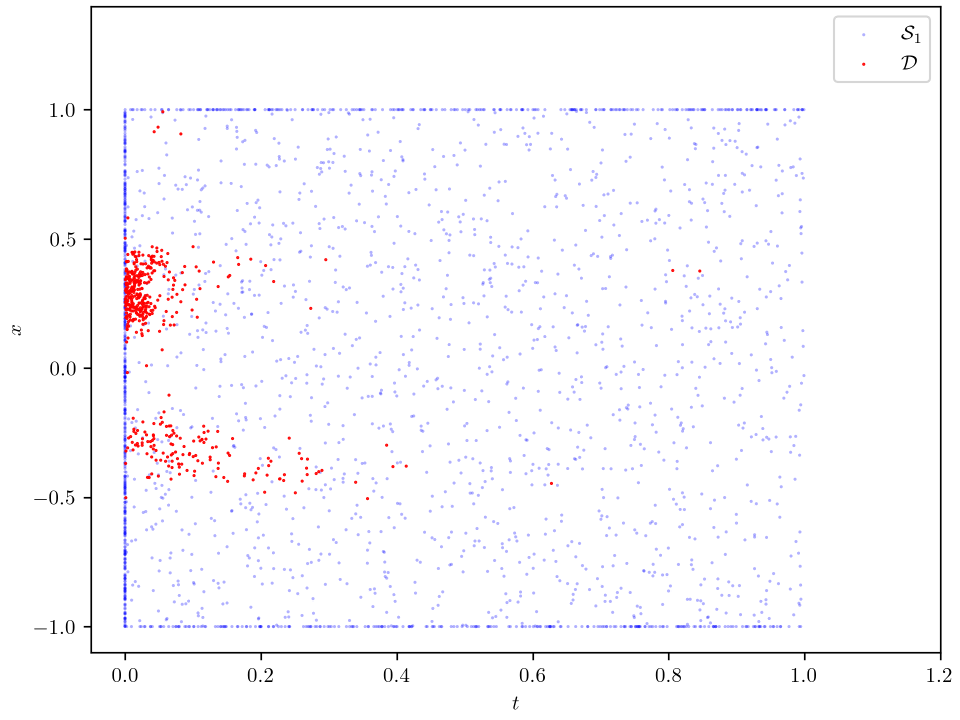}
        \caption{Nodes after pre-training of \textit{AAIS-t}.}
    \end{subfigure}%
    \begin{subfigure}{.33\textwidth}
        \centering
        \includegraphics[height=0.75\textwidth,width=1.0\textwidth]{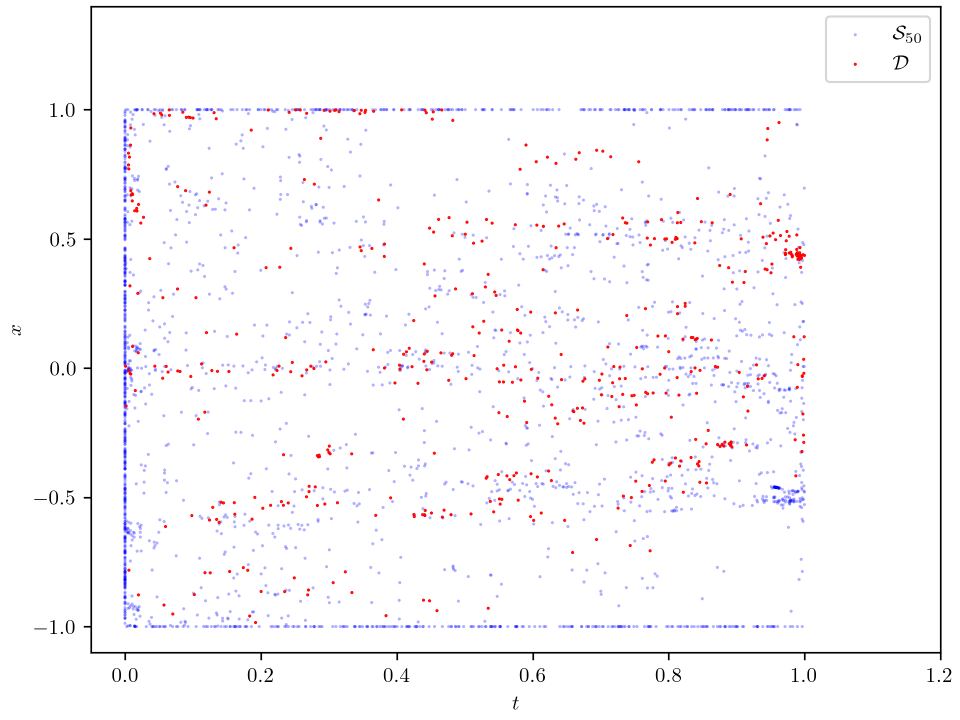}
        \caption{Nodes after 49th training of \textit{AAIS-t}.}
    \end{subfigure}%
    \begin{subfigure}{.33\textwidth}
        \centering
        \includegraphics[height=0.75\textwidth,width=1.0\textwidth]{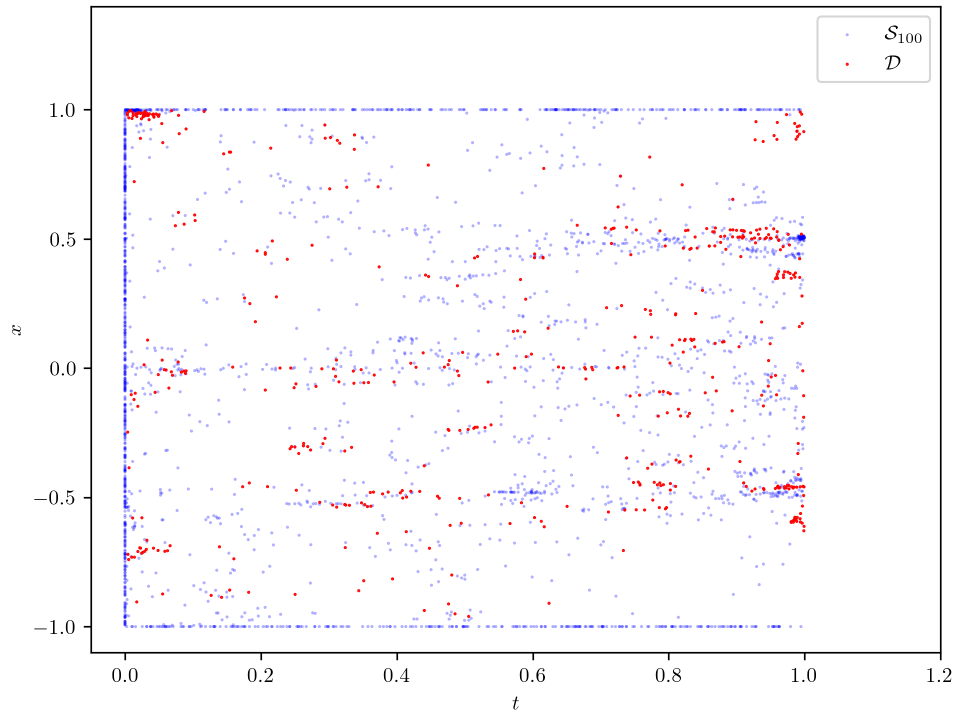}
        \caption{Nodes after 99th training of \textit{AAIS-t}.}
    \end{subfigure}%
    \caption{Profiles of residual and nodes for \textit{Uni} method and \textit{AAIS-t} method after pre-training, 49th training and 99th training. First row: residual for the \text{Uni} method. Second and Third row: residual and nodes for the \textit{AAIS-t} method. }
    \label{fig:ACloss2000e}
\end{figure}

\begin{figure}[htbp]
    \centering
    \begin{subfigure}{.25\textwidth}
        \centering
        \includegraphics[height=0.75\textwidth,width=1.0\textwidth]{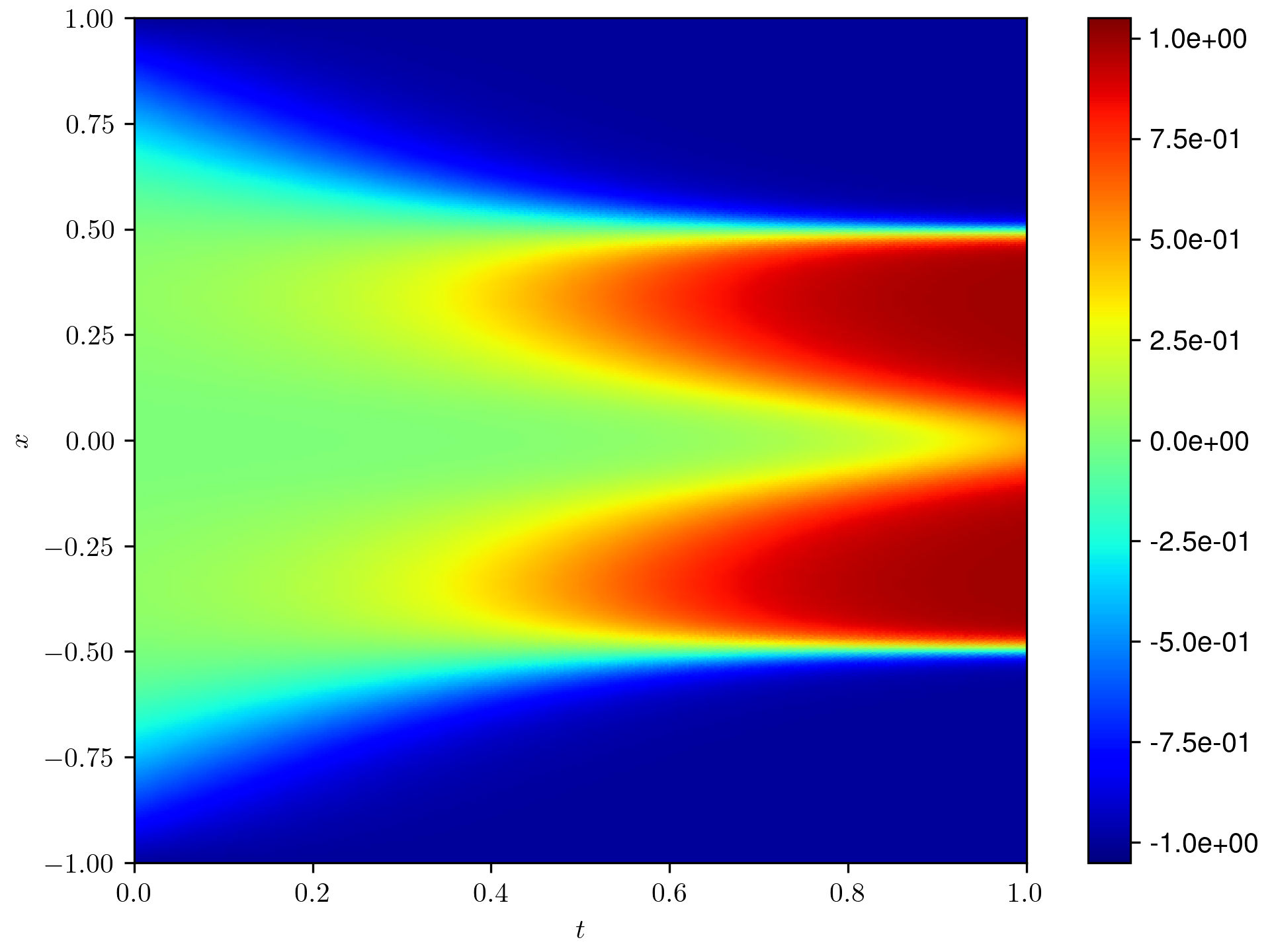}
    \end{subfigure}%
    \begin{subfigure}{.25\textwidth}
        \centering
        \includegraphics[height=0.75\textwidth,width=1.0\textwidth]{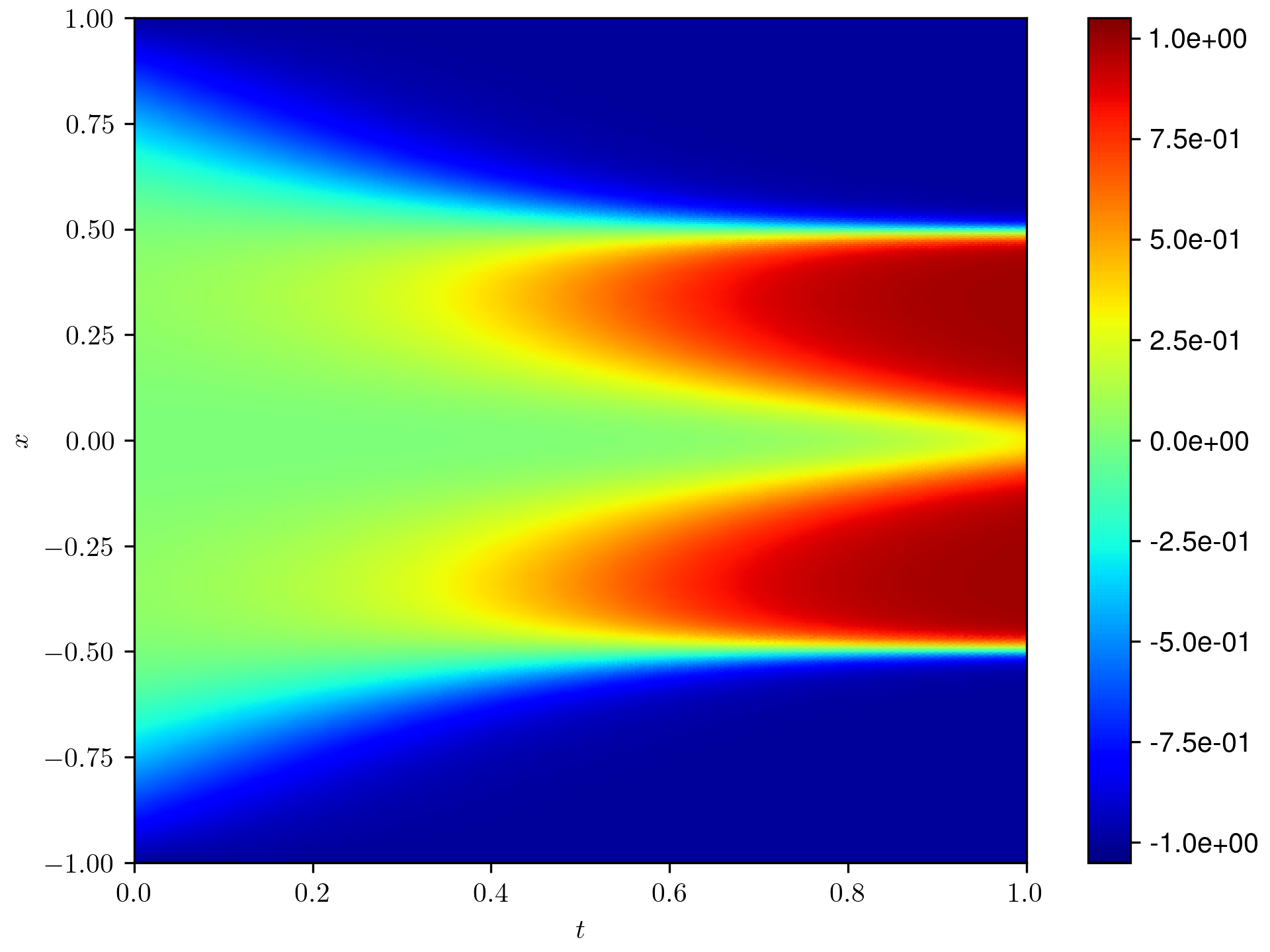}
    \end{subfigure}%
    \begin{subfigure}{.25\textwidth}
        \centering
        \includegraphics[height=0.75\textwidth,width=1.0\textwidth]{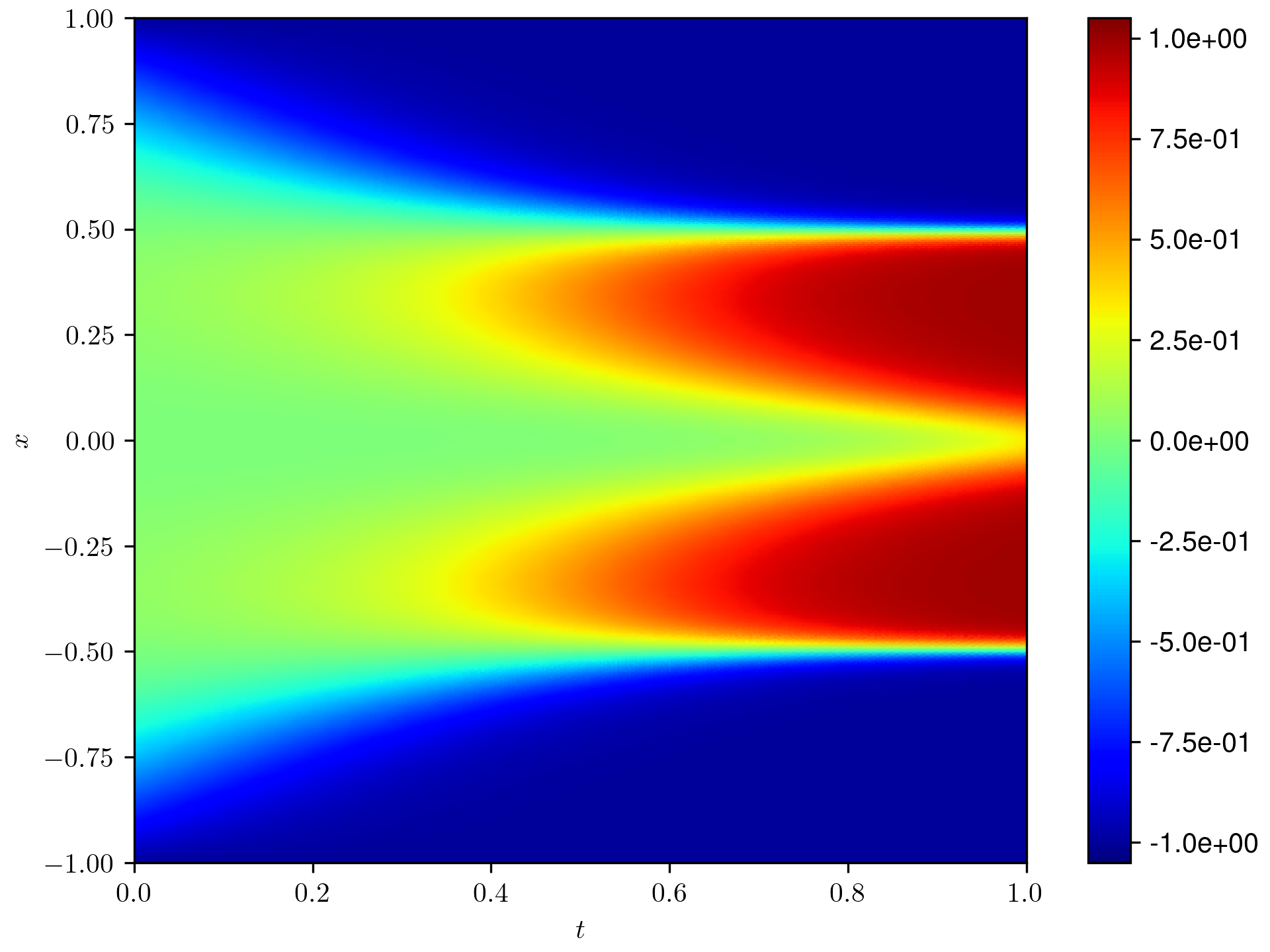}
    \end{subfigure}%
    \begin{subfigure}{.25\textwidth}
        \centering
        \includegraphics[height=0.75\textwidth,width=1.0\textwidth]{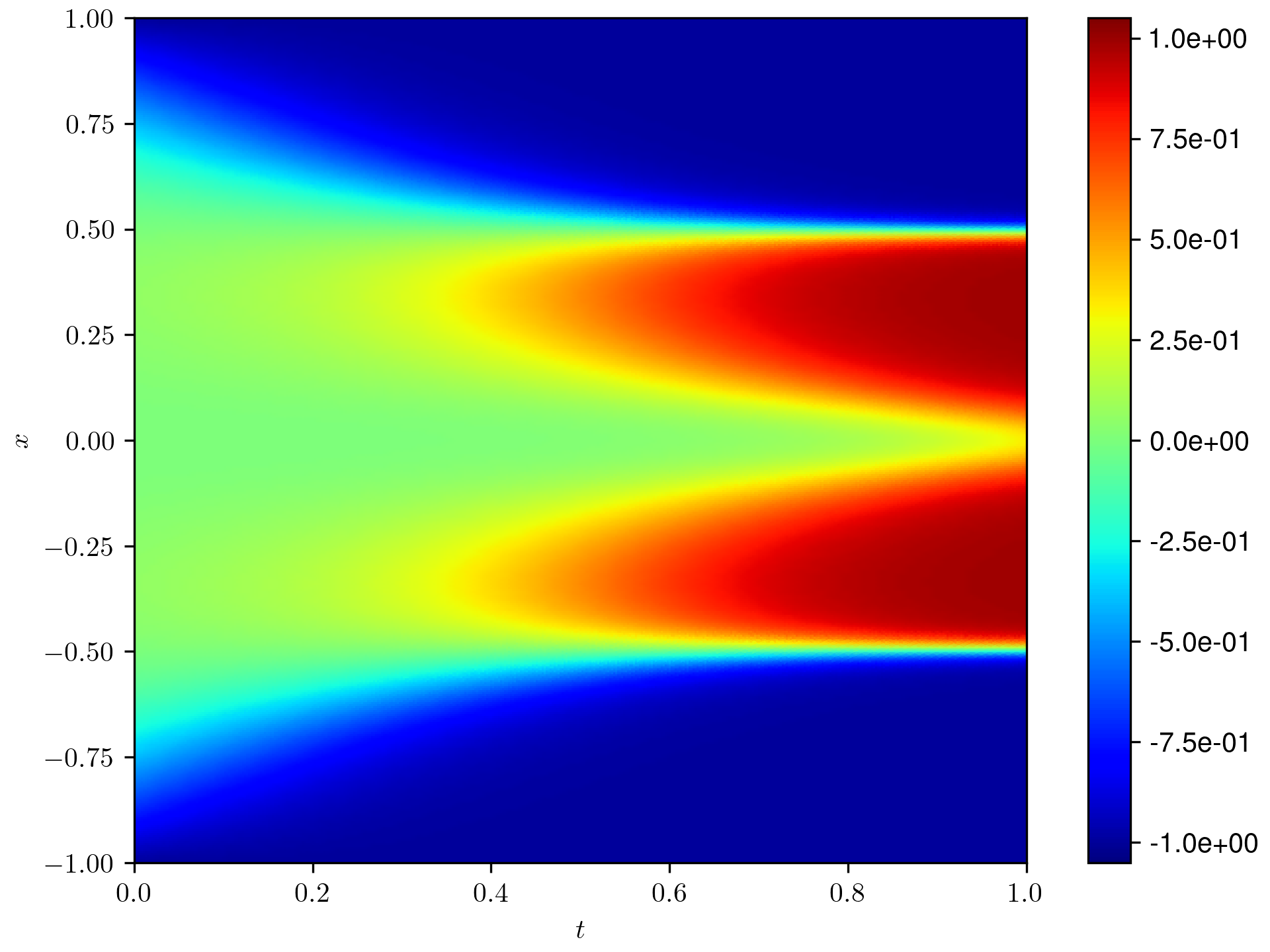}
    \end{subfigure}%
    \newline
    \raggedleft
    \begin{subfigure}{.25\textwidth}
        \centering
        \includegraphics[height=0.75\textwidth,width=1.0\textwidth]{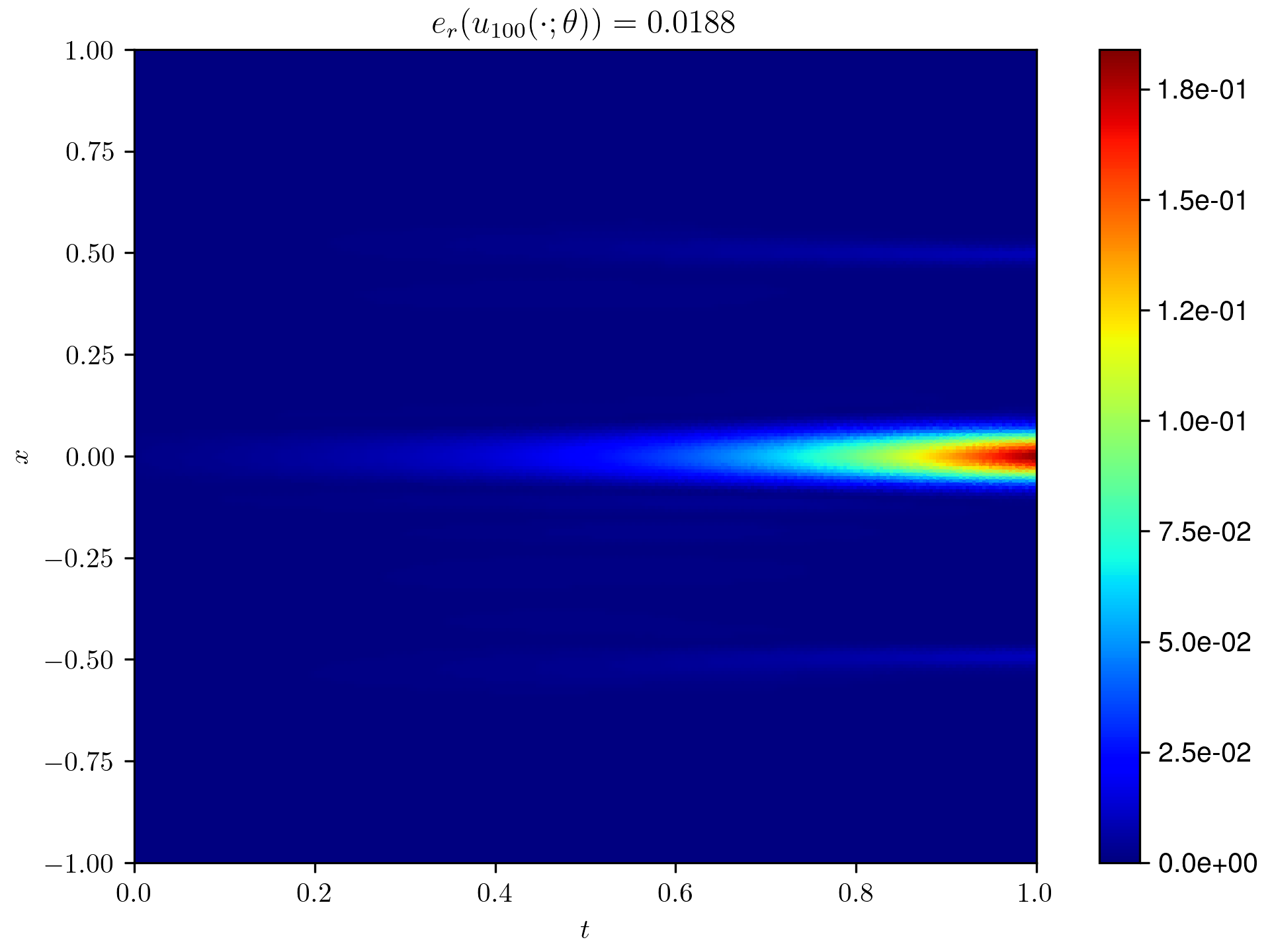}
        \caption{\textit{Uni}}
    \end{subfigure}%
    \begin{subfigure}{.25\textwidth}
        \centering
        \includegraphics[height=0.75\textwidth,width=1.0\textwidth]{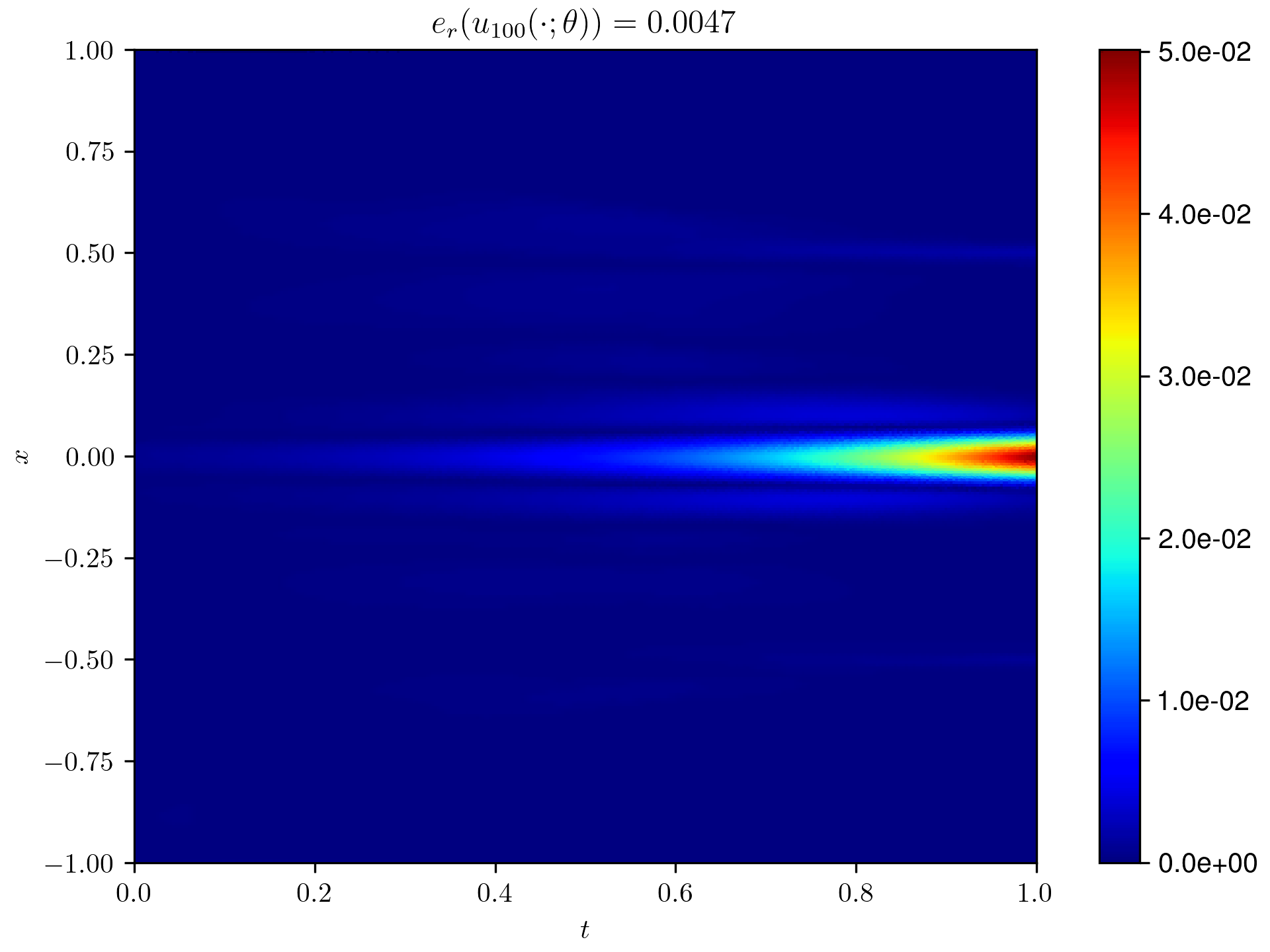}
        \caption{\textit{RAD}}
    \end{subfigure}%
    \begin{subfigure}{.25\textwidth}
        \centering
        \includegraphics[height=0.75\textwidth,width=1.0\textwidth]{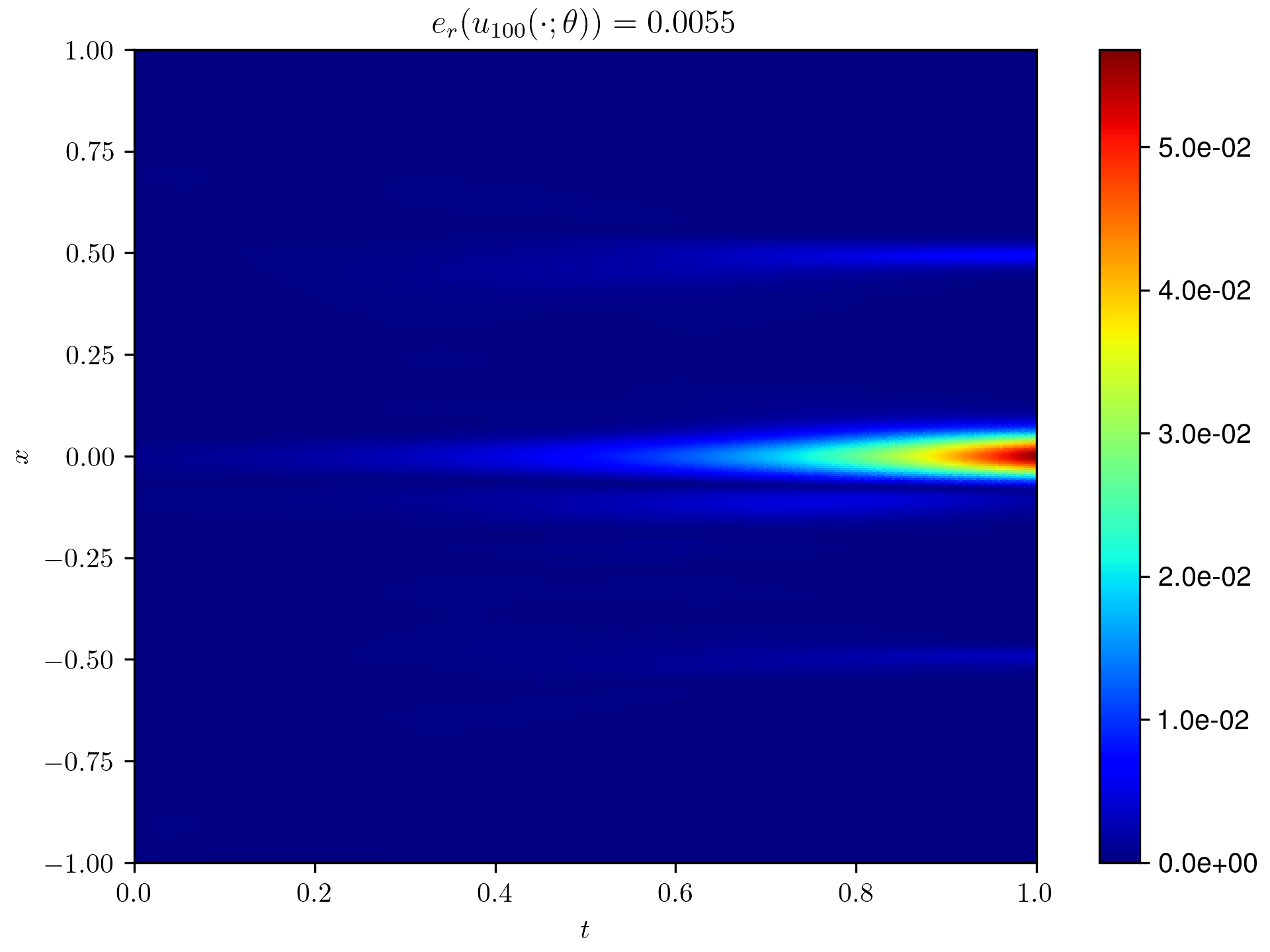}
        \caption{\textit{AAIS-g}}
    \end{subfigure}%
    \begin{subfigure}{.25\textwidth}
        \centering
        \includegraphics[height=0.75\textwidth,width=1.0\textwidth]{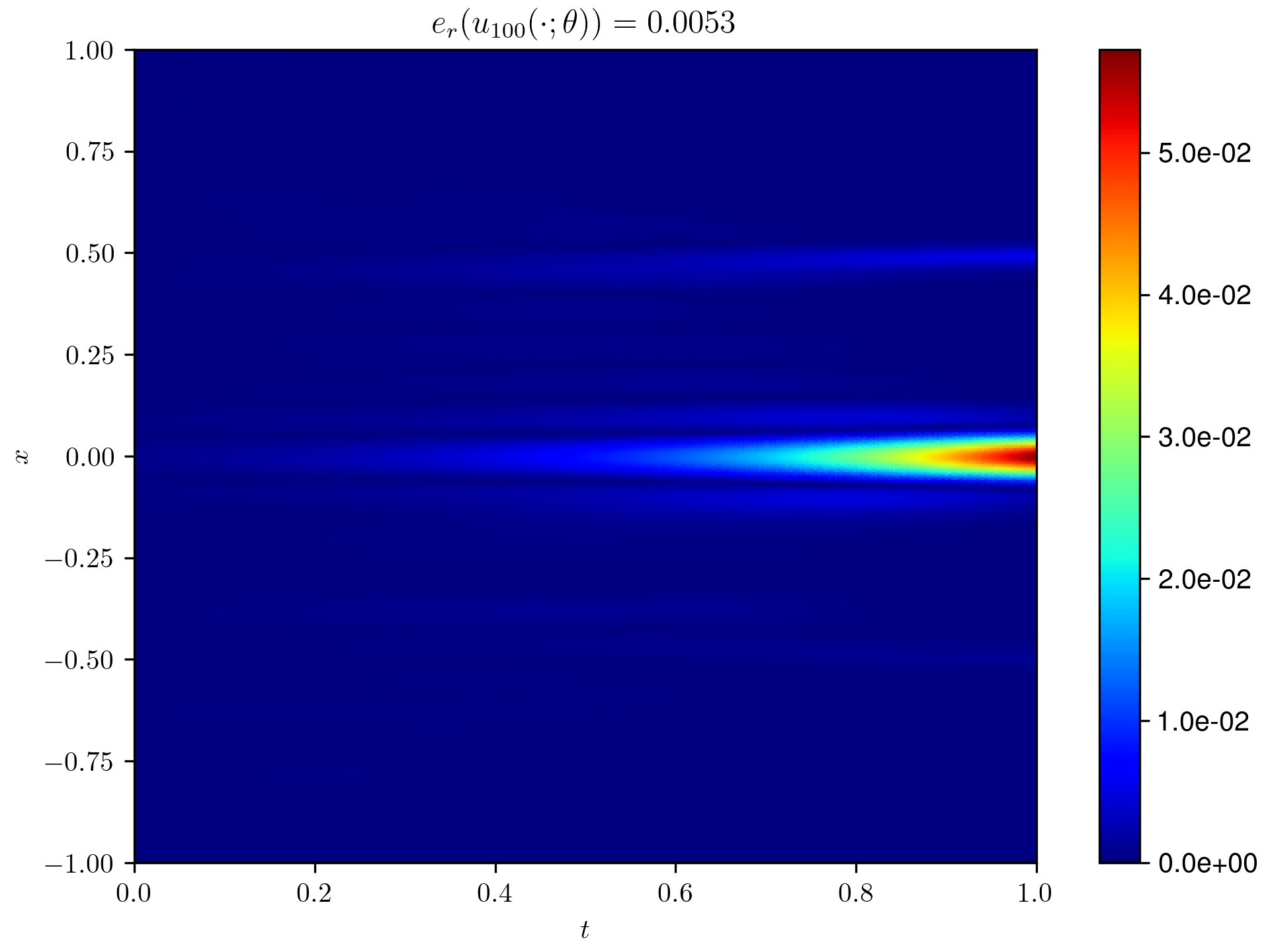}
        \caption{\textit{AAIS-t}}
    \end{subfigure}%
    \caption{Profiles of absolute error and neural network solutions for Allen-Cahn equation after 100th training. }
    \label{fig:ACAbs2000e}
\end{figure}

The profiles of residual and nodes are showed in Figure \ref{fig:ACloss2000e}. Here we use \textit{AAIS-t} as an example to illustrate the efficiency of adaptive sampling methods. After pre-training, the residual would focus on the location near $t=0$. And adaptive sampling methods would sample around there. After 49th training, the residual would focus on the place near the finial time $t=1$ and place $x=\pm 0.5$. However, the \textit{Uni} method would not train a lot there. In the end we could see \textit{Uni} method would still have a bigger residual around $x=\pm 0.5$ and $t=1$ where adaptive sampling methods do not have. And the frequency of residual increase during training for adaptive sampling methods may show a better solution.
Figure \ref{fig:ACAbs2000e} list the solution and absolute error of each method. Notice that at the end of training, the absolute error would mainly concentrate near $t=1.0$, $x=0$, but the location of high error would not reflect in the residual significantly, so all four methods could not concentrate there. It may be caused by the nonlinearity of the PDE equation \eqref{pde:AllenCahn}. 

\subsection{Korteweg-de Vries equation}
\begin{figure}[htbp]
    \centering
    \includegraphics[scale=0.125]{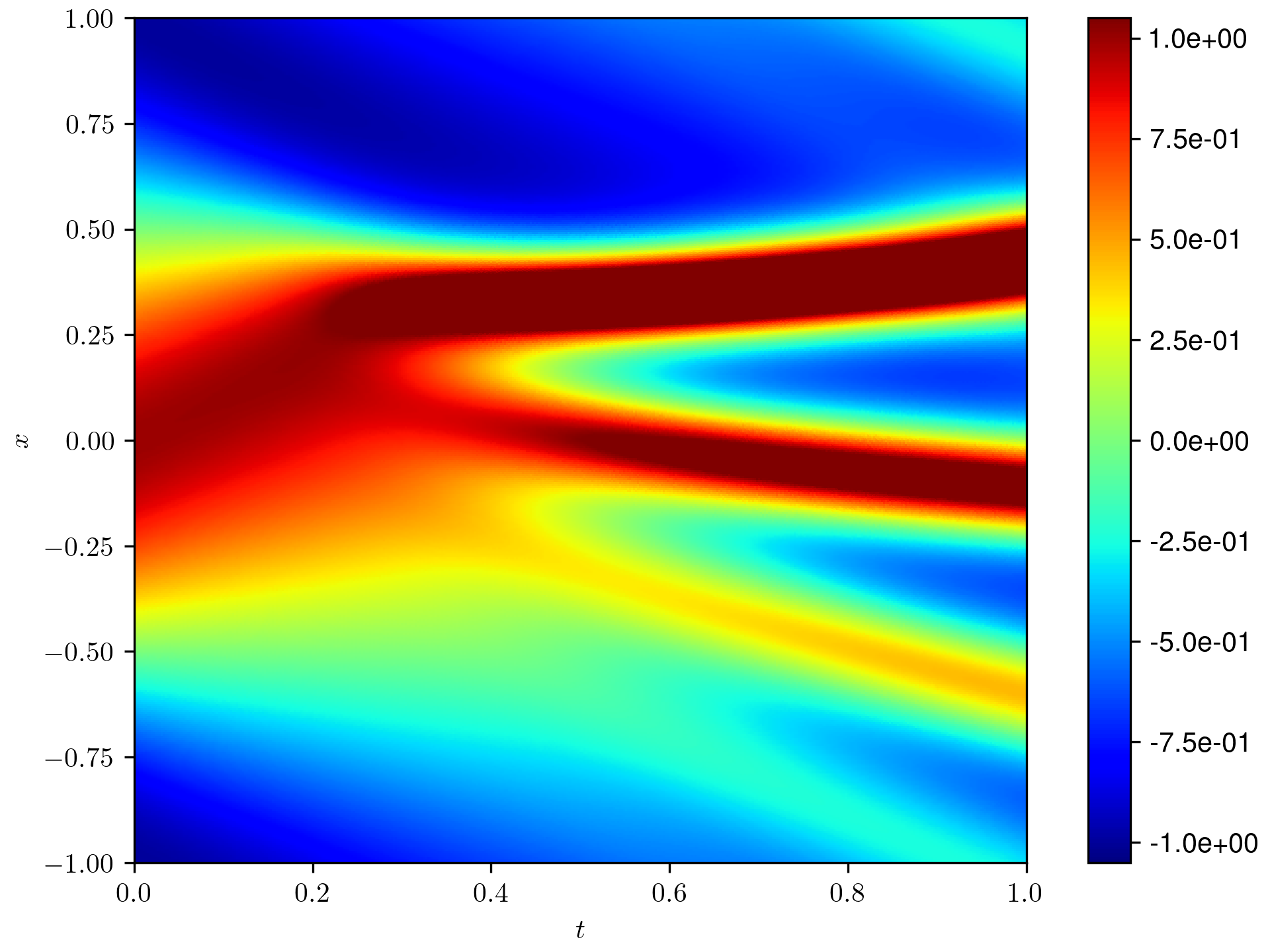}
    \caption{Reference solution for KdV equations in \eqref{pde:KdV}.}
    \label{fig:KdVexact}
\end{figure}

\begin{figure}[htbp]
    \centering
    \begin{subfigure}{.5\textwidth}
    \centering
    \includegraphics[scale=0.25]{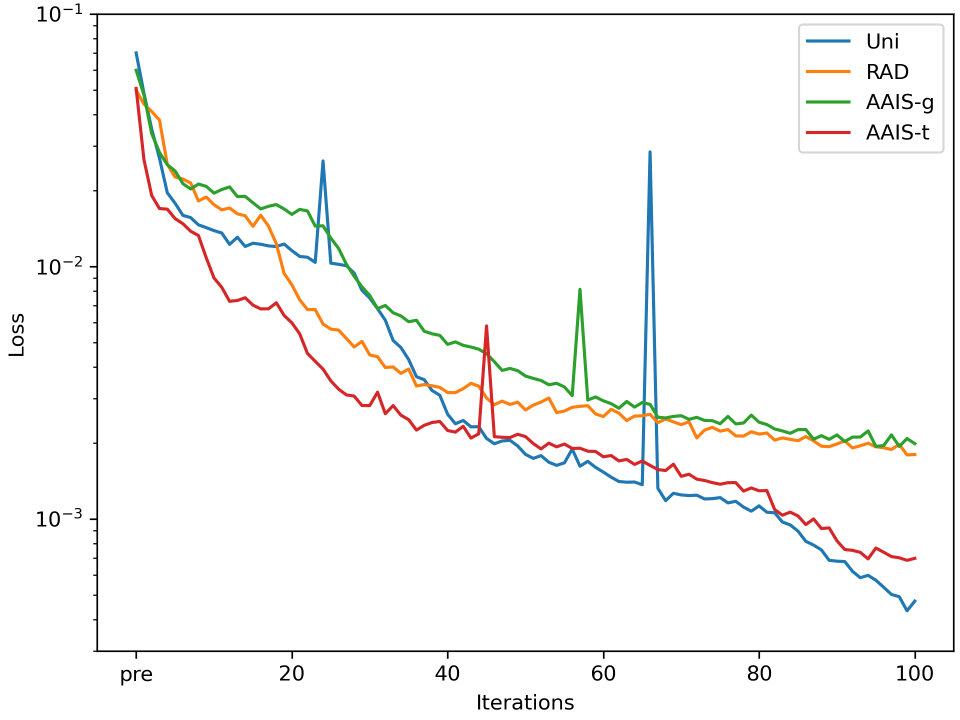}
    \end{subfigure}%
    \begin{subfigure}{.5\textwidth}
    \centering
    \includegraphics[scale=0.25]{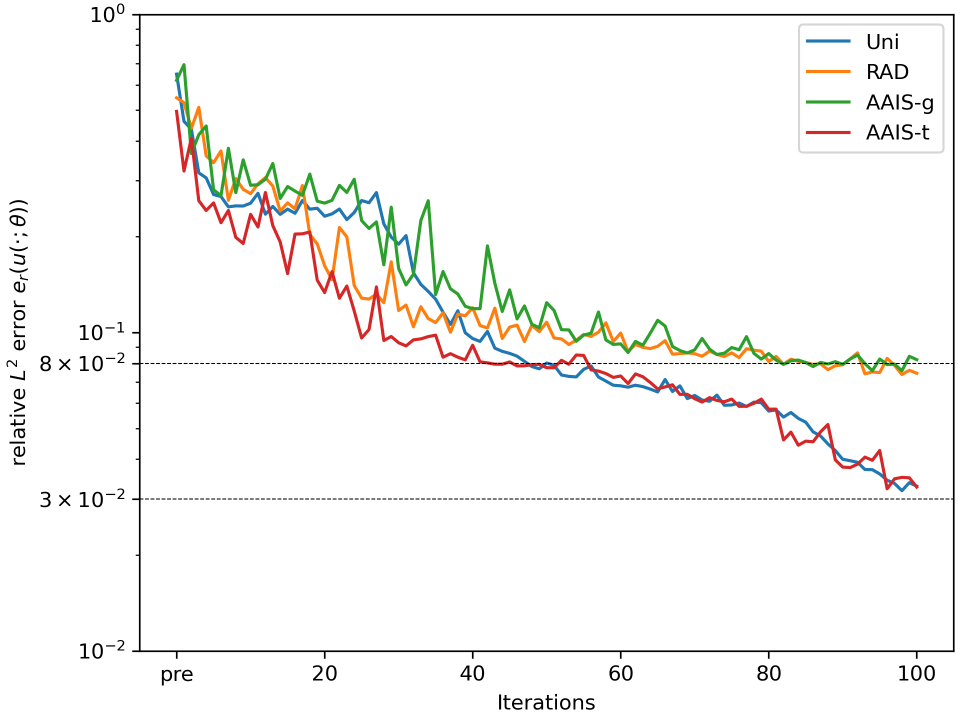}
    \end{subfigure}%
    \caption{Loss and relative errors during training for KdV equation with four sampling methods. Left: the loss function. Right: the relative $L^2$ error $e_r(u(\cdot;\theta))$.
    }
    \label{fig:KdVLossErr500e}
\end{figure}

\begin{figure}[htbp]
    \centering
    \begin{subfigure}{.25\textwidth}
        \centering
        \includegraphics[height=0.75\textwidth,width=1.0\textwidth]{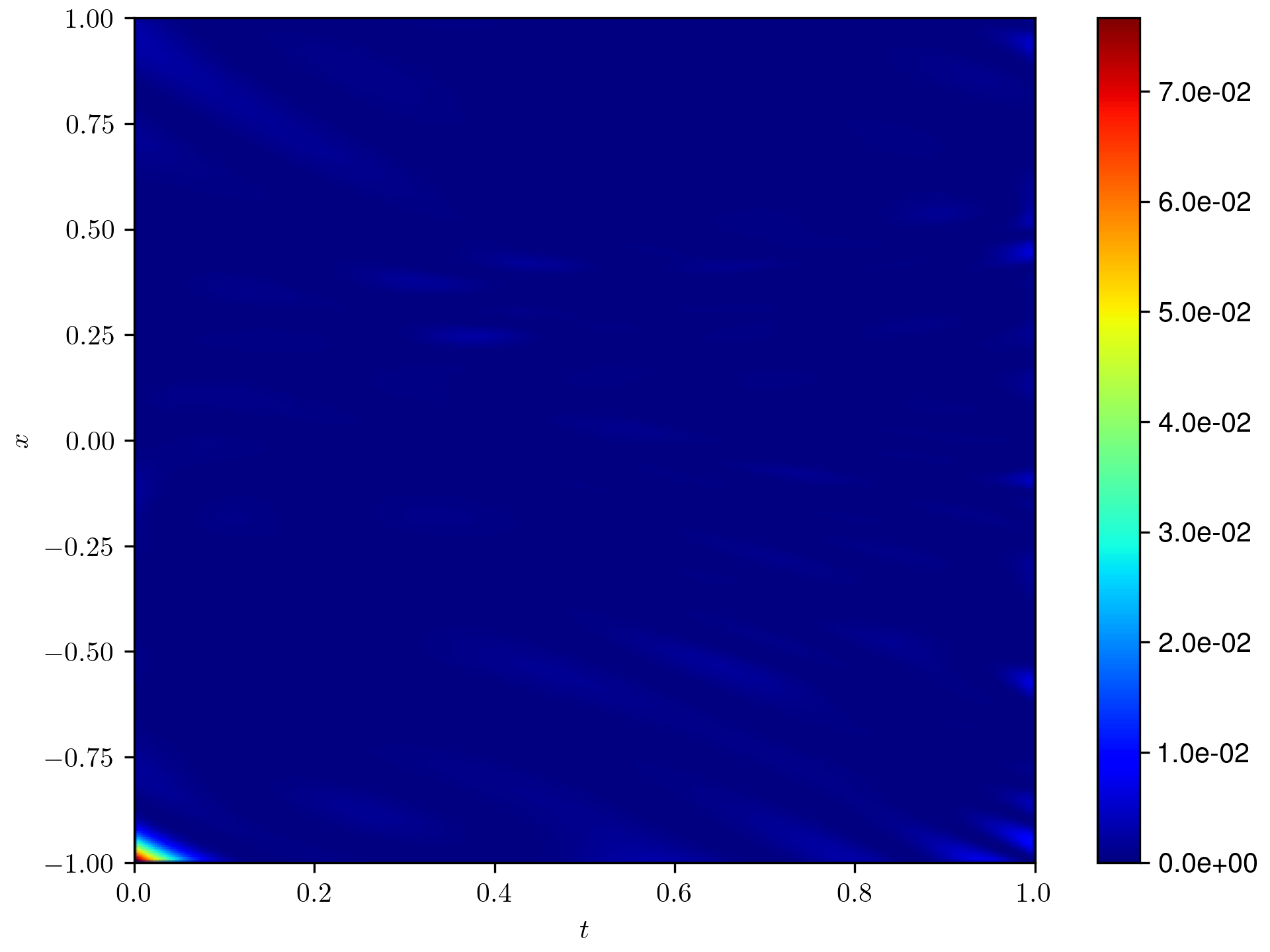}
    \end{subfigure}%
    \begin{subfigure}{.25\textwidth}
        \centering
        \includegraphics[height=0.75\textwidth,width=1.0\textwidth]{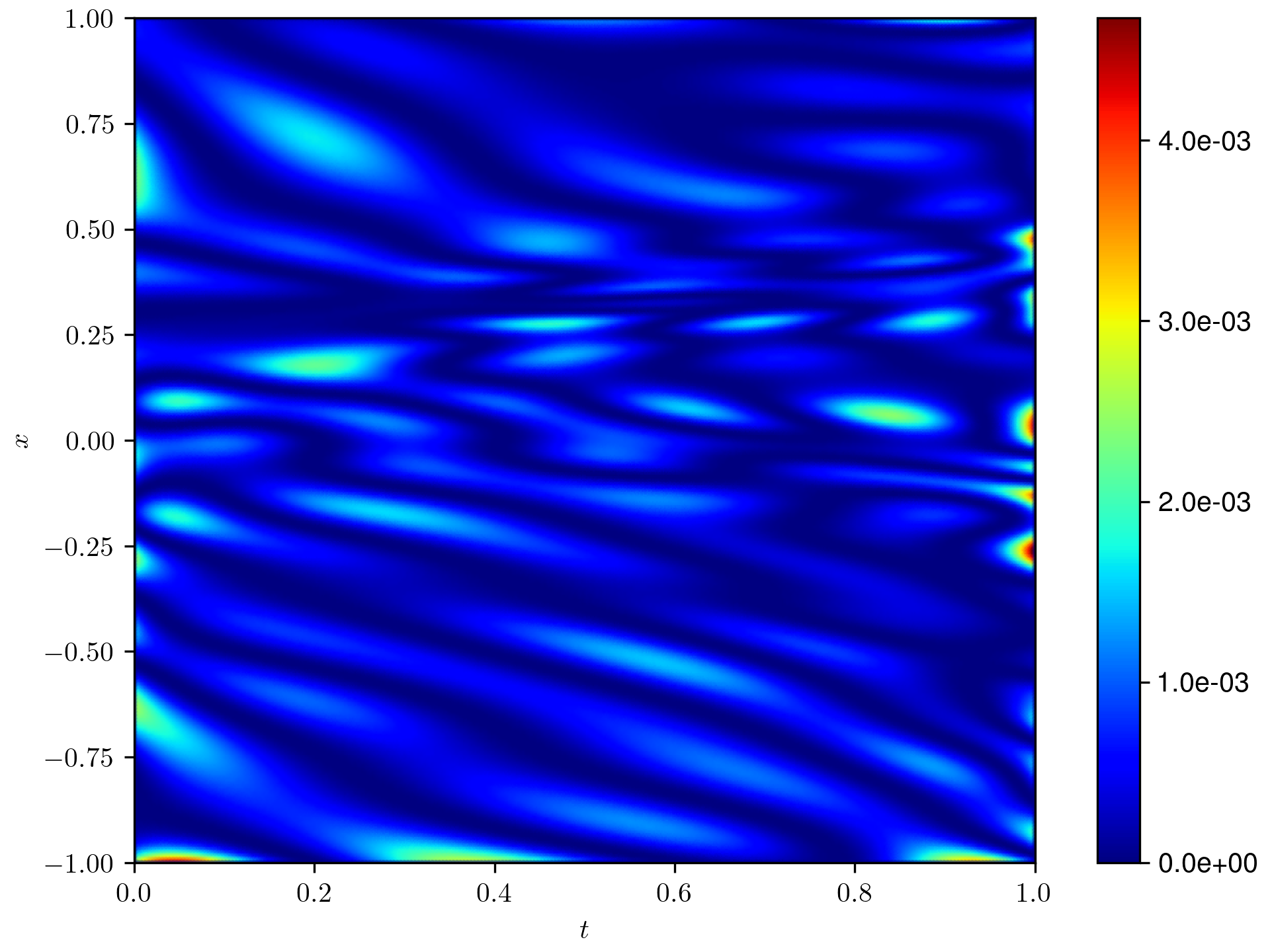}
    \end{subfigure}%
    \begin{subfigure}{.25\textwidth}
        \centering
        \includegraphics[height=0.75\textwidth,width=1.0\textwidth]{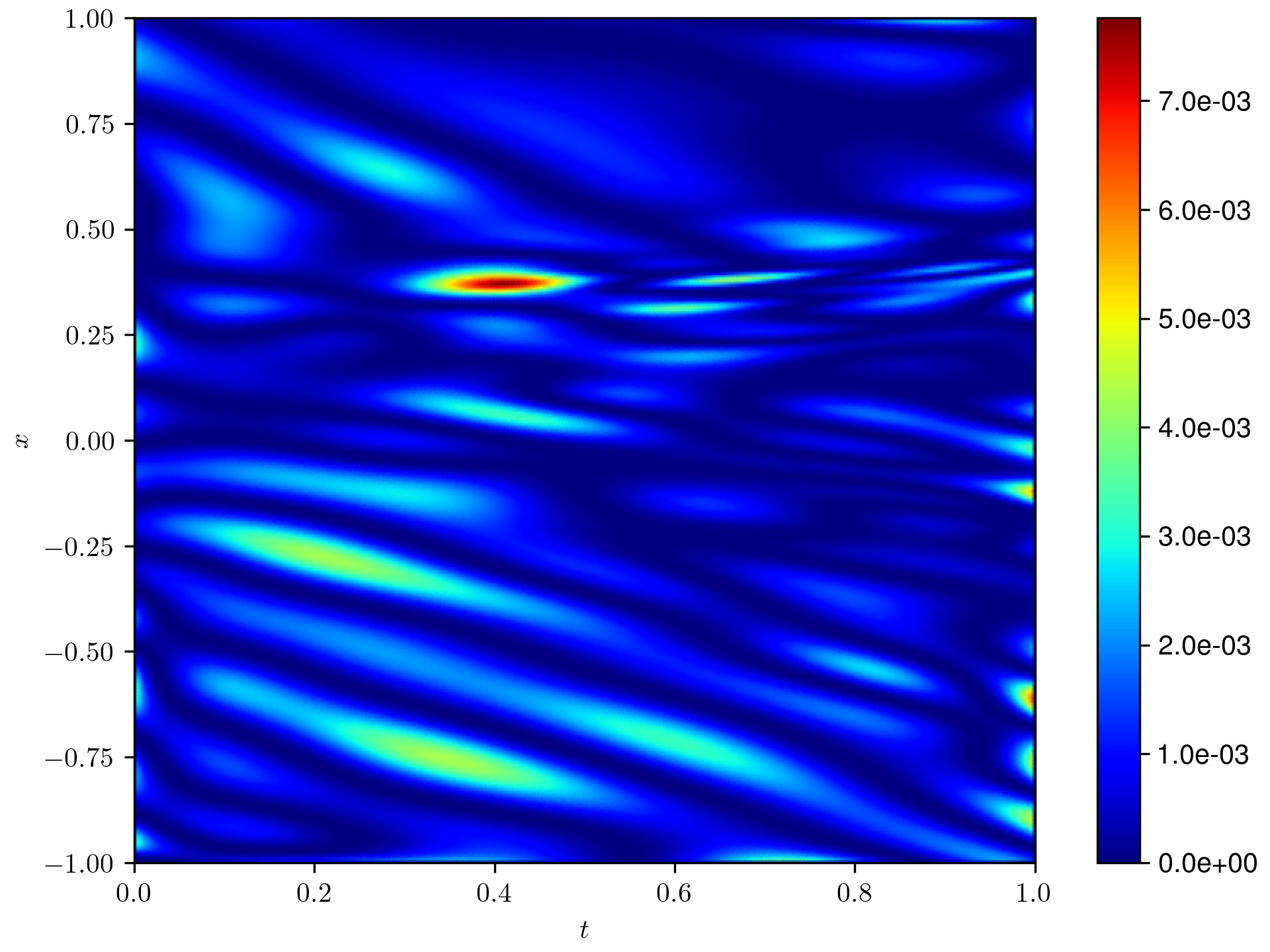}
    \end{subfigure}%
    \begin{subfigure}{.25\textwidth}
        \centering
        \includegraphics[height=0.75\textwidth,width=1.0\textwidth]{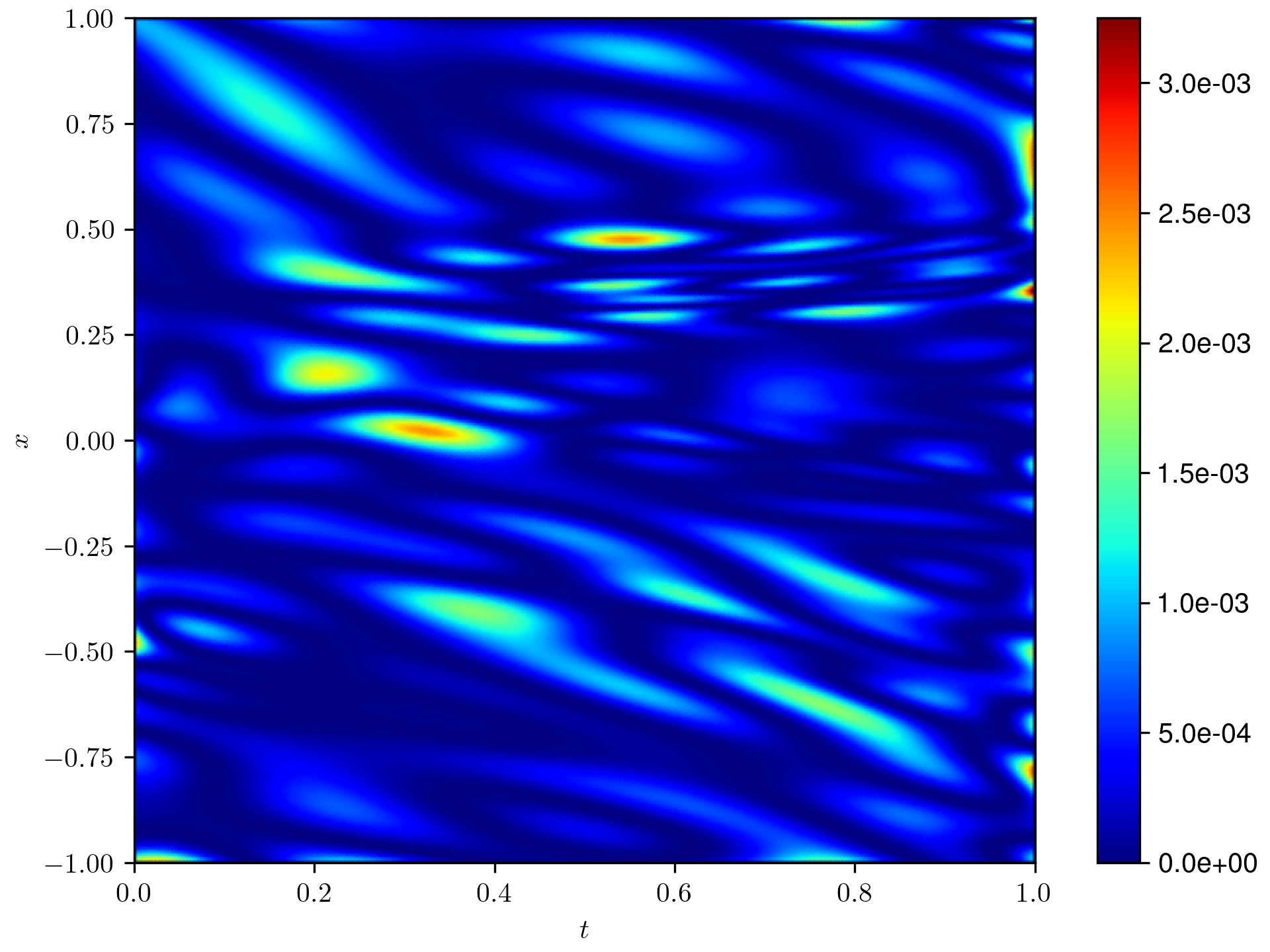}
    \end{subfigure}%
    \newline
    \raggedleft
    \begin{subfigure}{.25\textwidth}
        \centering
        \includegraphics[height=0.75\textwidth,width=1.0\textwidth]{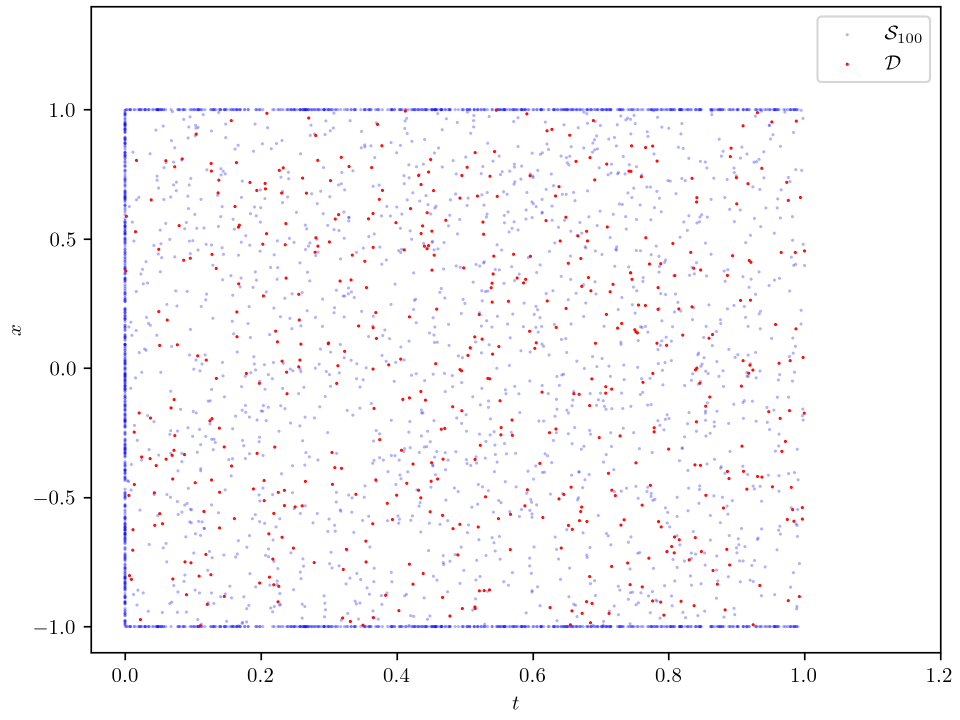}
        \caption{\textit{Uni}}
    \end{subfigure}%
    \begin{subfigure}{.25\textwidth}
        \centering
        \includegraphics[height=0.75\textwidth,width=1.0\textwidth]{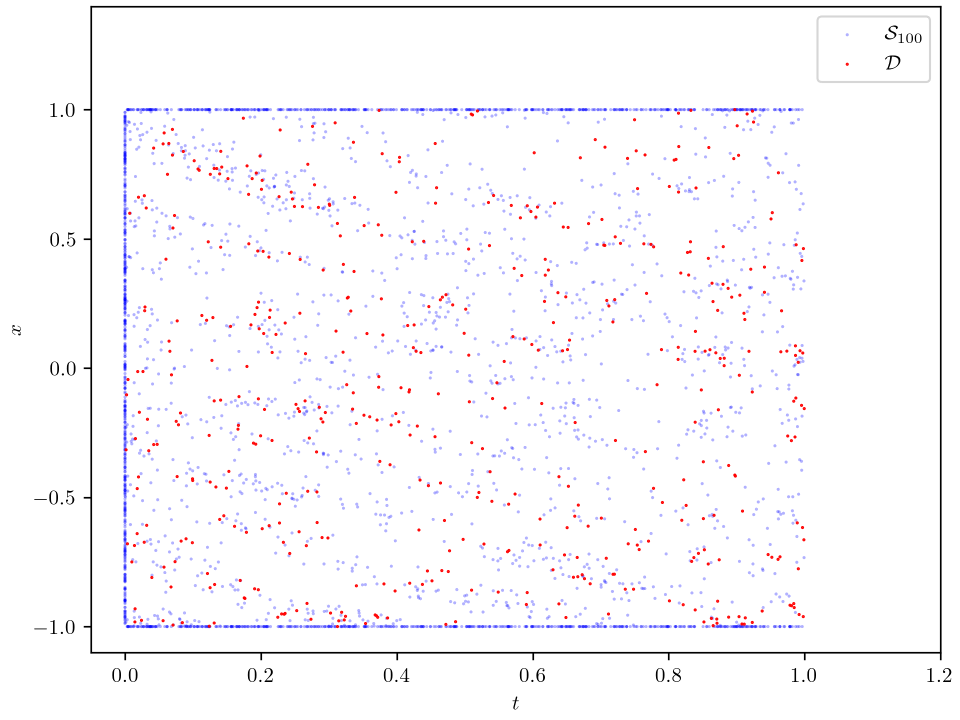}
        \caption{\textit{RAD}}
    \end{subfigure}%
    \begin{subfigure}{.25\textwidth}
        \centering
        \includegraphics[height=0.75\textwidth,width=1.0\textwidth]{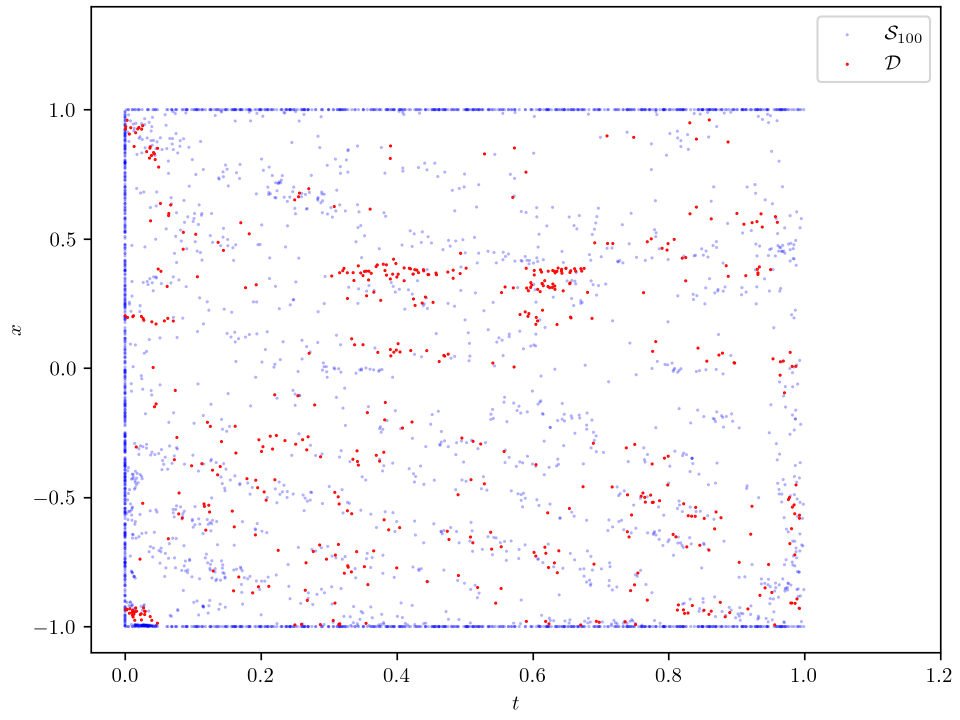}
        \caption{\textit{AAIS-g}}
    \end{subfigure}%
    \begin{subfigure}{.25\textwidth}
        \centering
        \includegraphics[height=0.75\textwidth,width=1.0\textwidth]{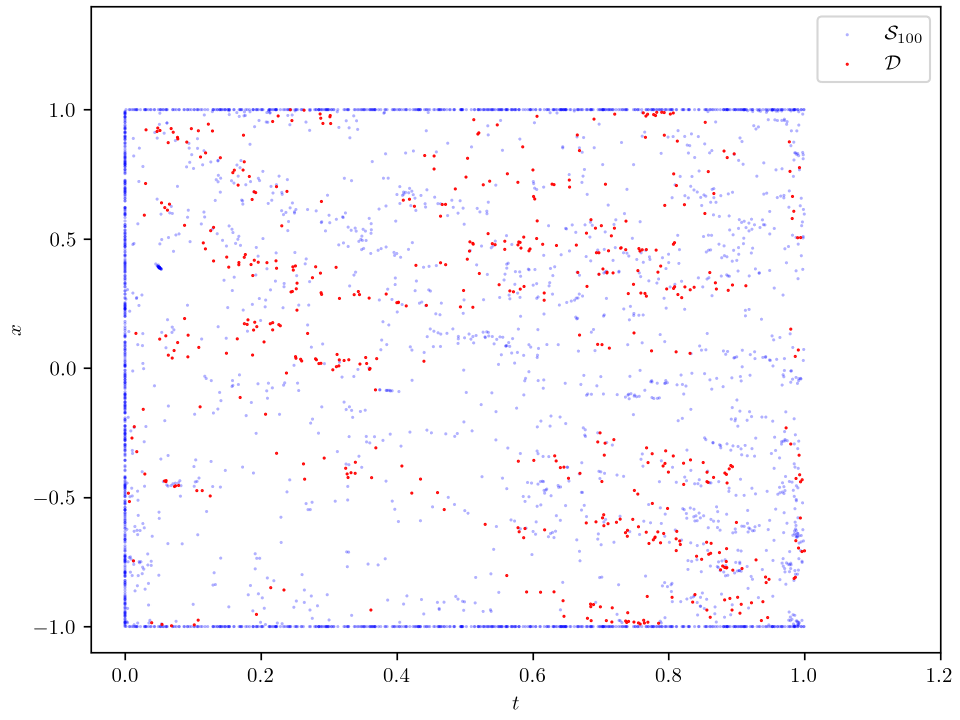}
        \caption{\textit{AAIS-t}}
    \end{subfigure}%
    \caption{Profiles of residual and nodes  for KdV equation after 100th training.}
    \label{fig:KdVNode5000e}
\end{figure}

\begin{figure}[htbp]
    \centering
    \begin{subfigure}{.25\textwidth}
        \centering
        \includegraphics[height=0.75\textwidth,width=1.0\textwidth]{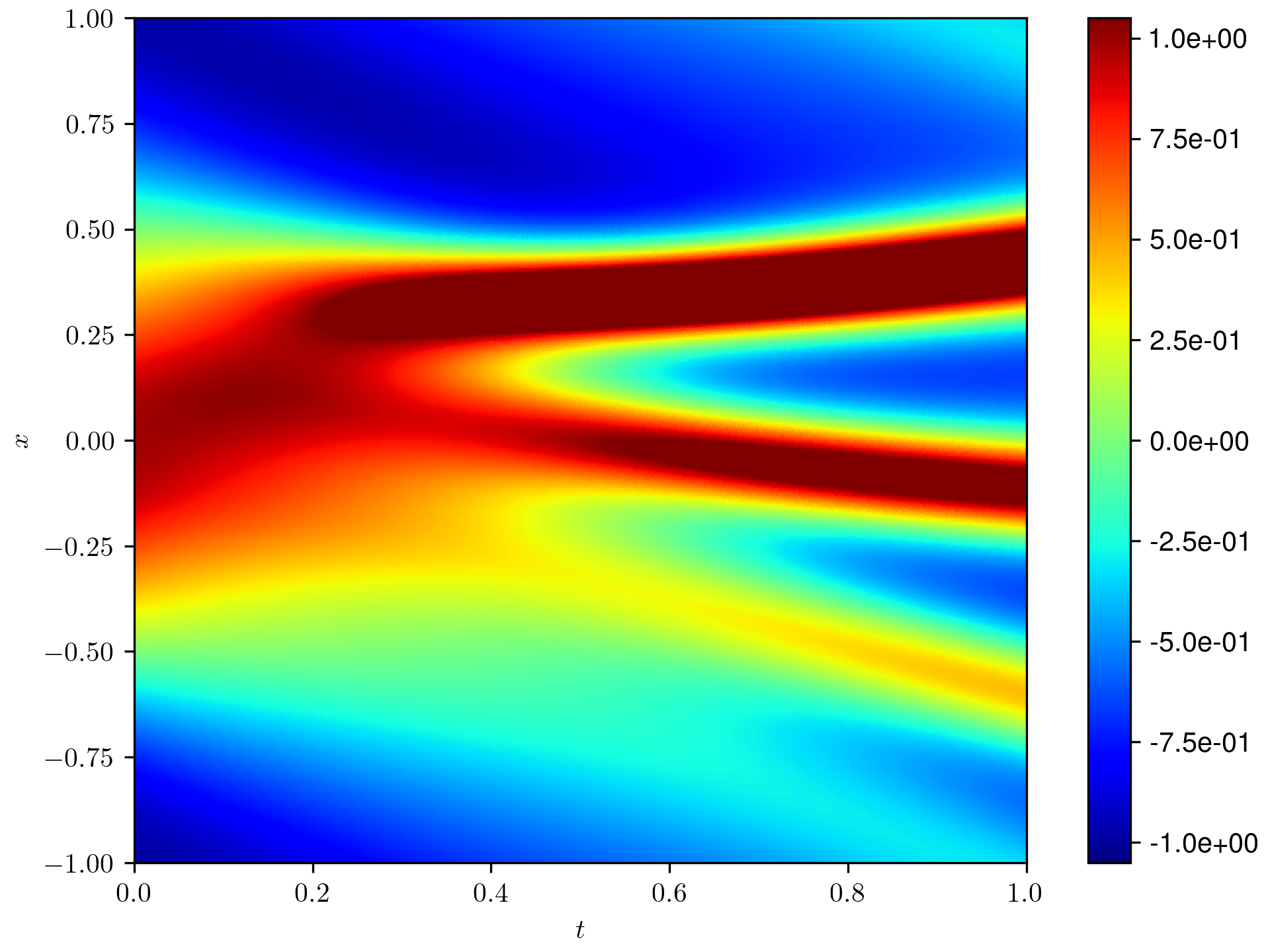}
    \end{subfigure}%
    \begin{subfigure}{.25\textwidth}
        \centering
        \includegraphics[height=0.75\textwidth,width=1.0\textwidth]{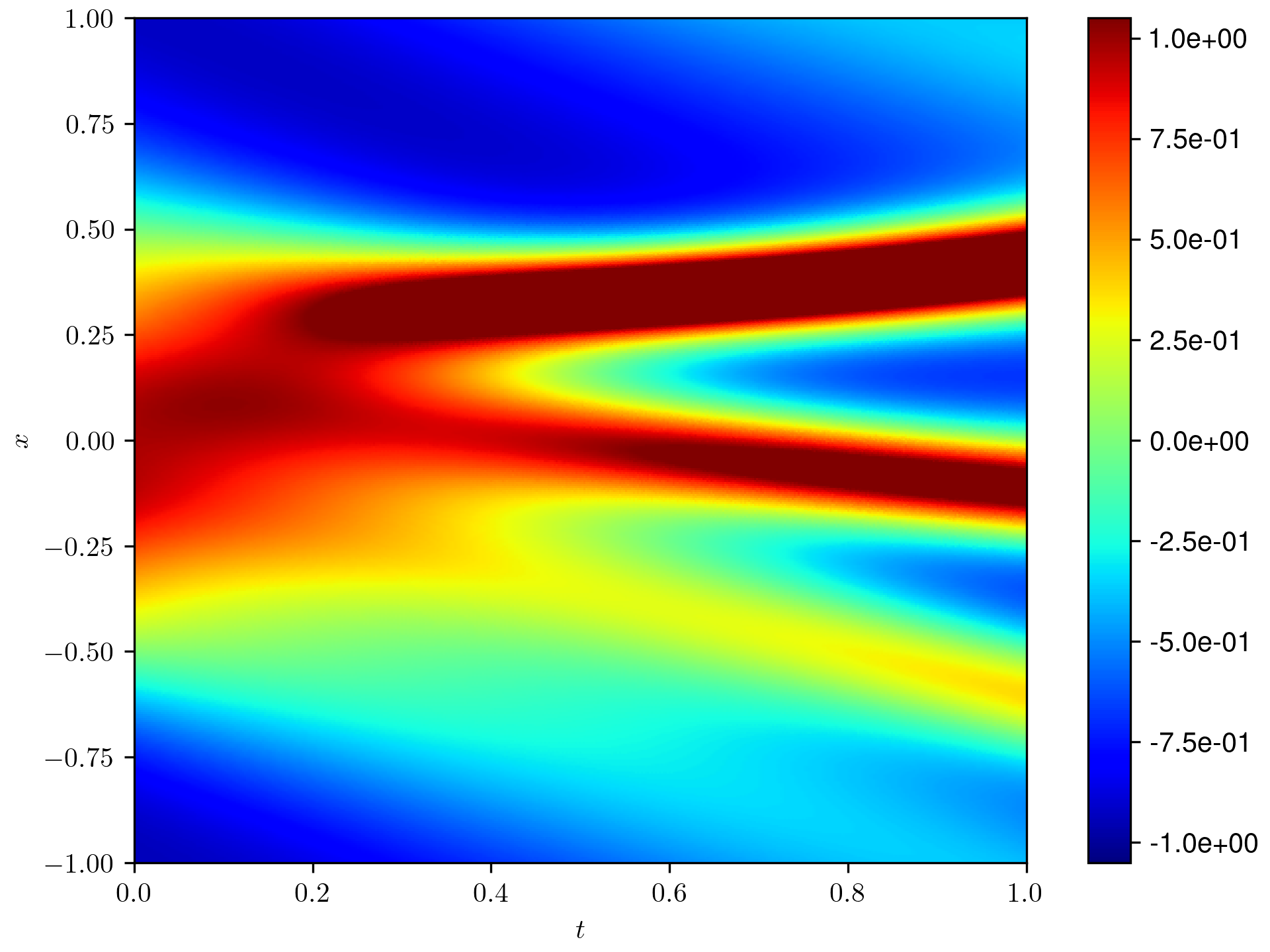}
    \end{subfigure}%
    \begin{subfigure}{.25\textwidth}
        \centering
        \includegraphics[height=0.75\textwidth,width=1.0\textwidth]{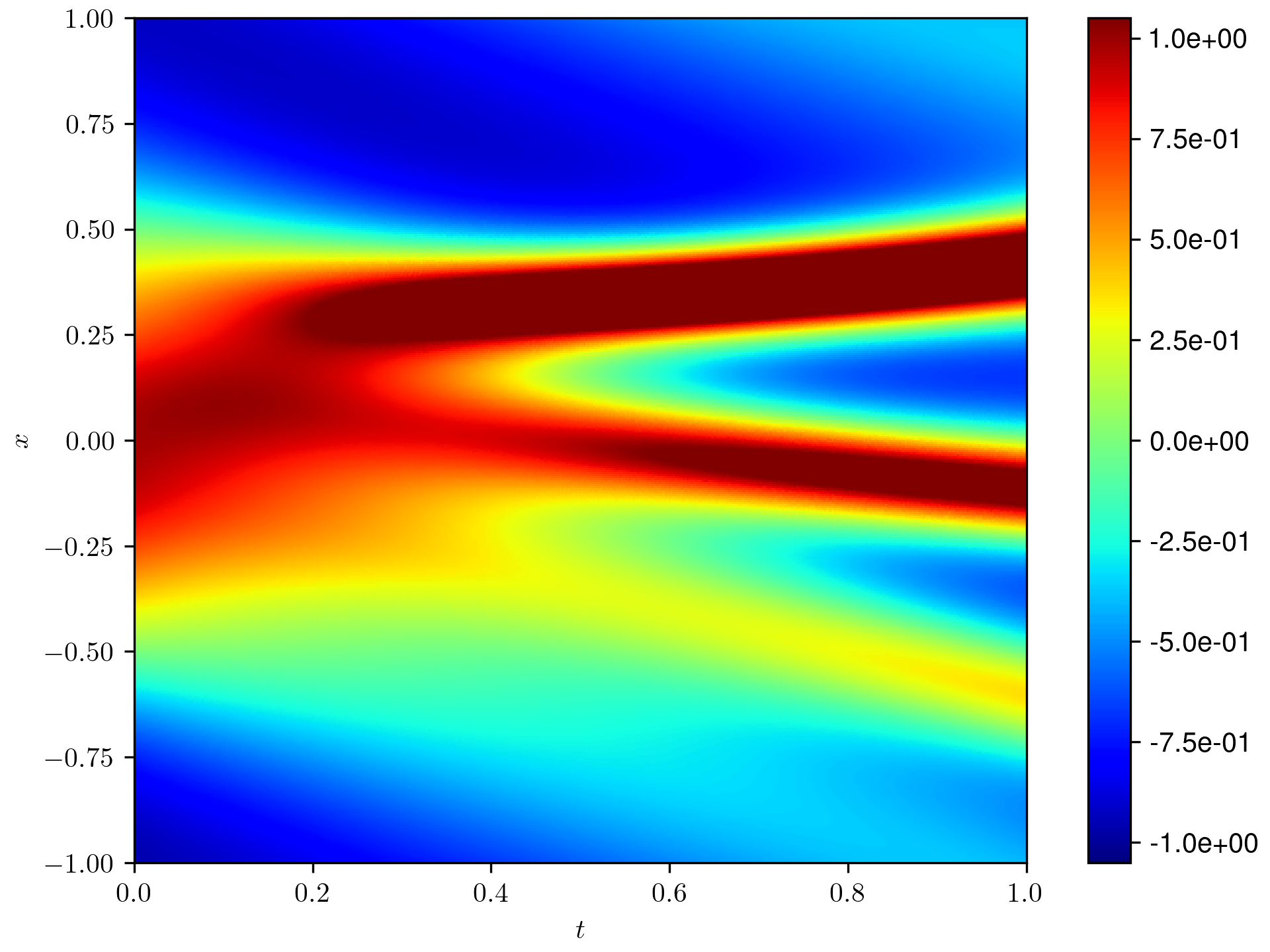}
    \end{subfigure}%
    \begin{subfigure}{.25\textwidth}
        \centering
        \includegraphics[height=0.75\textwidth,width=1.0\textwidth]{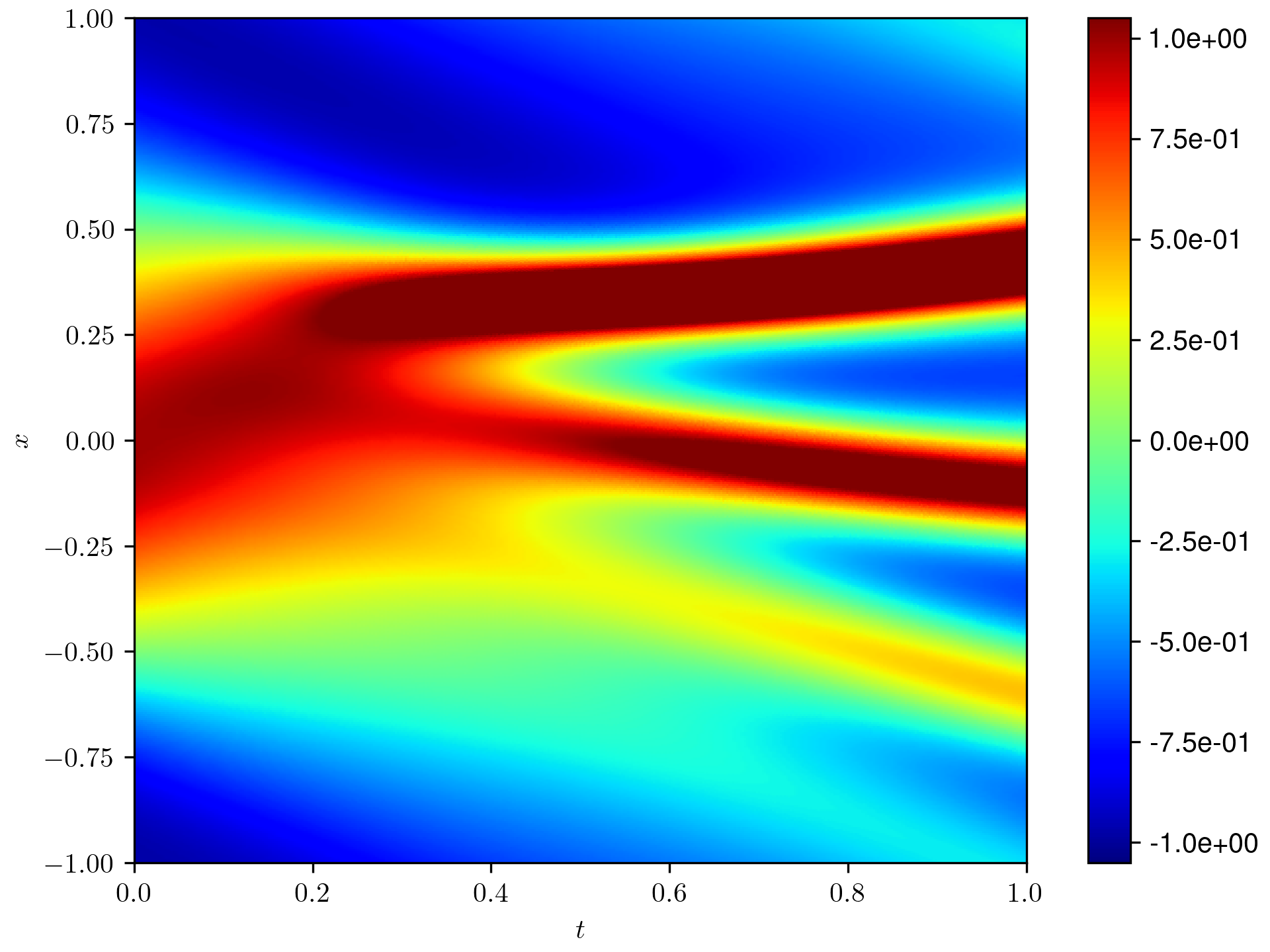}
    \end{subfigure}%
    \newline
    \raggedleft
    \begin{subfigure}{.25\textwidth}
        \centering
        \includegraphics[height=0.75\textwidth,width=1.0\textwidth]{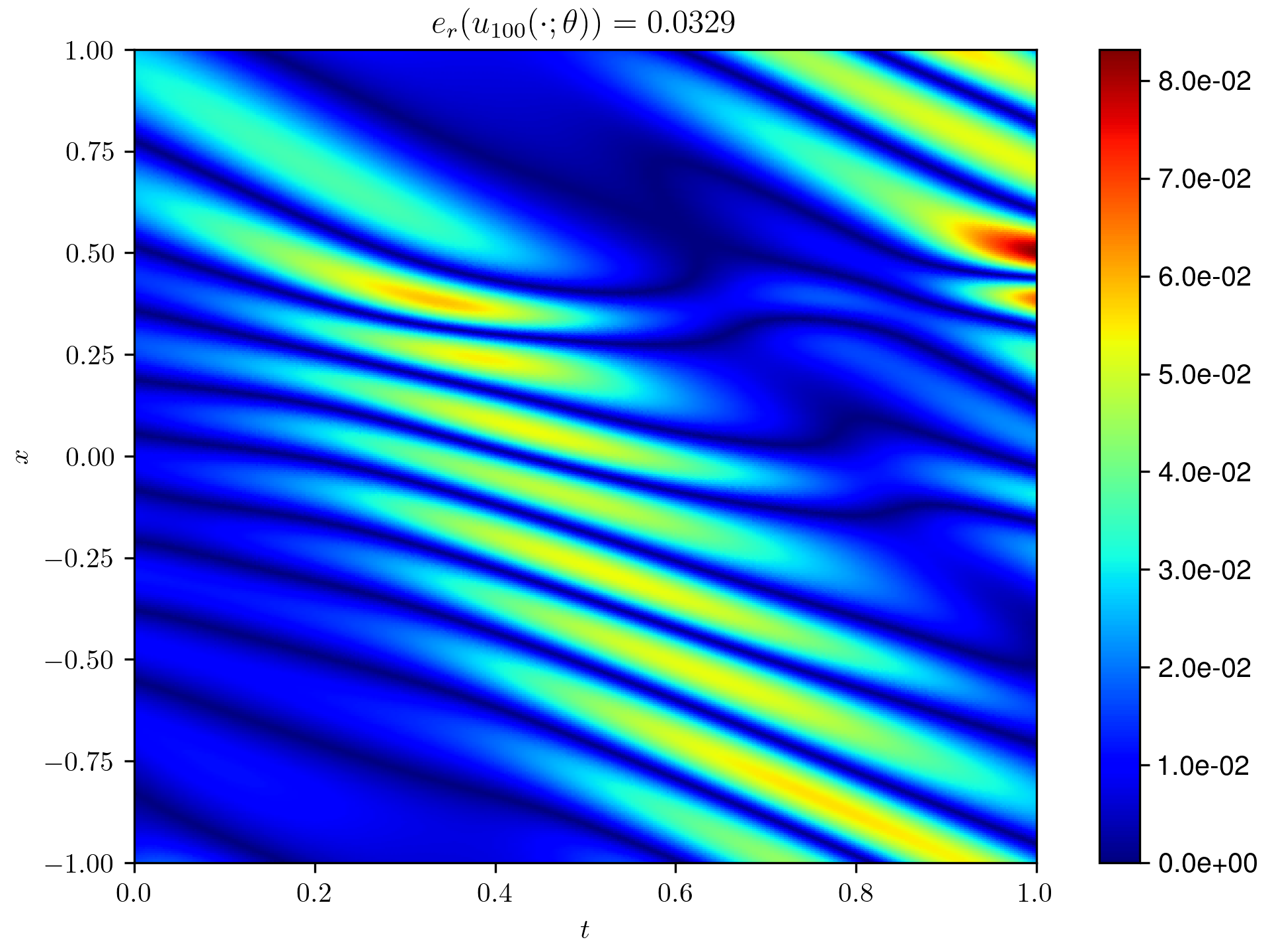}
        \caption{\textit{Uni}}
    \end{subfigure}%
    \begin{subfigure}{.25\textwidth}
        \centering
        \includegraphics[height=0.75\textwidth,width=1.0\textwidth]{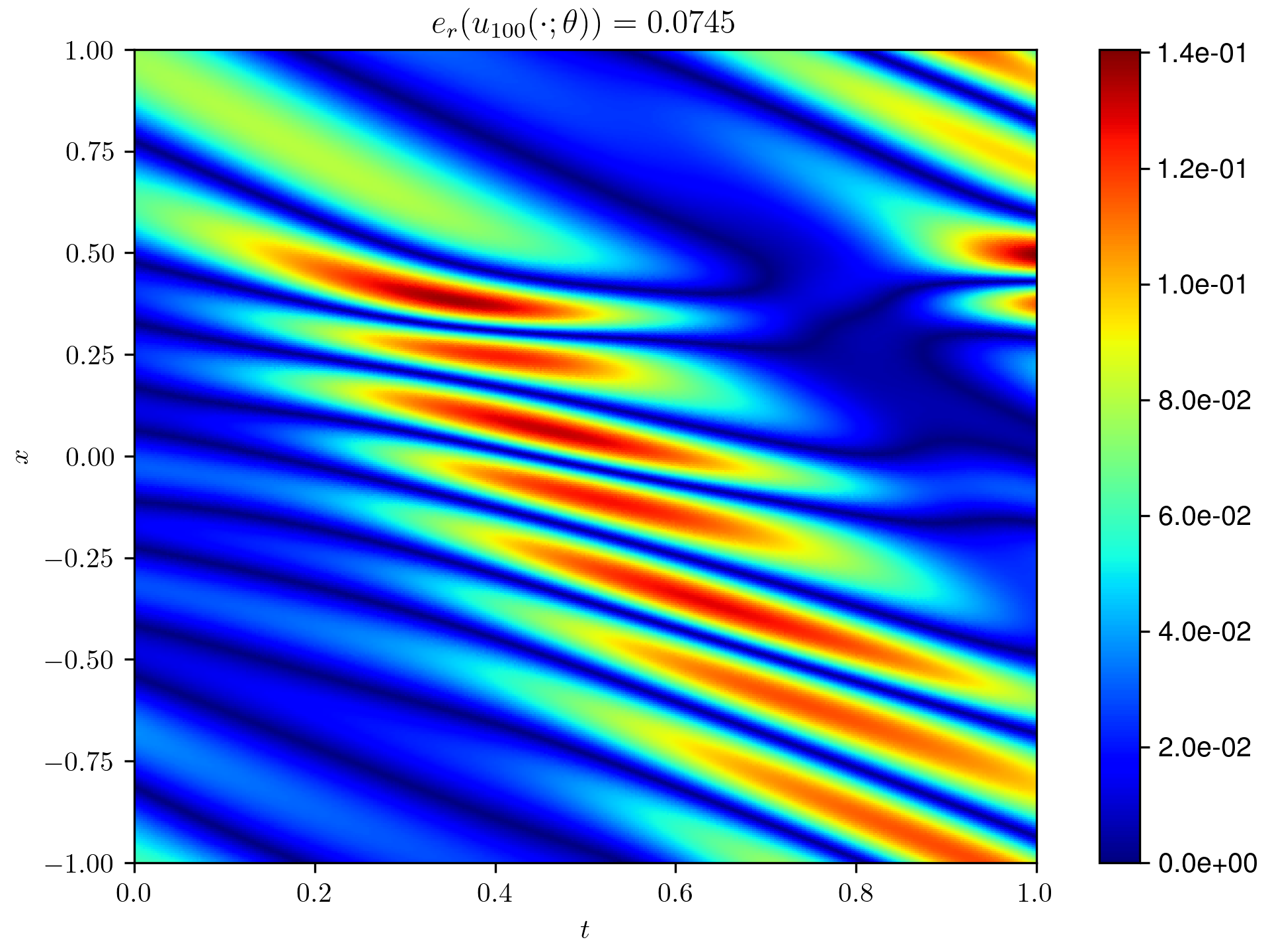}
        \caption{\textit{RAD}}
    \end{subfigure}%
    \begin{subfigure}{.25\textwidth}
        \centering
        \includegraphics[height=0.75\textwidth,width=1.0\textwidth]{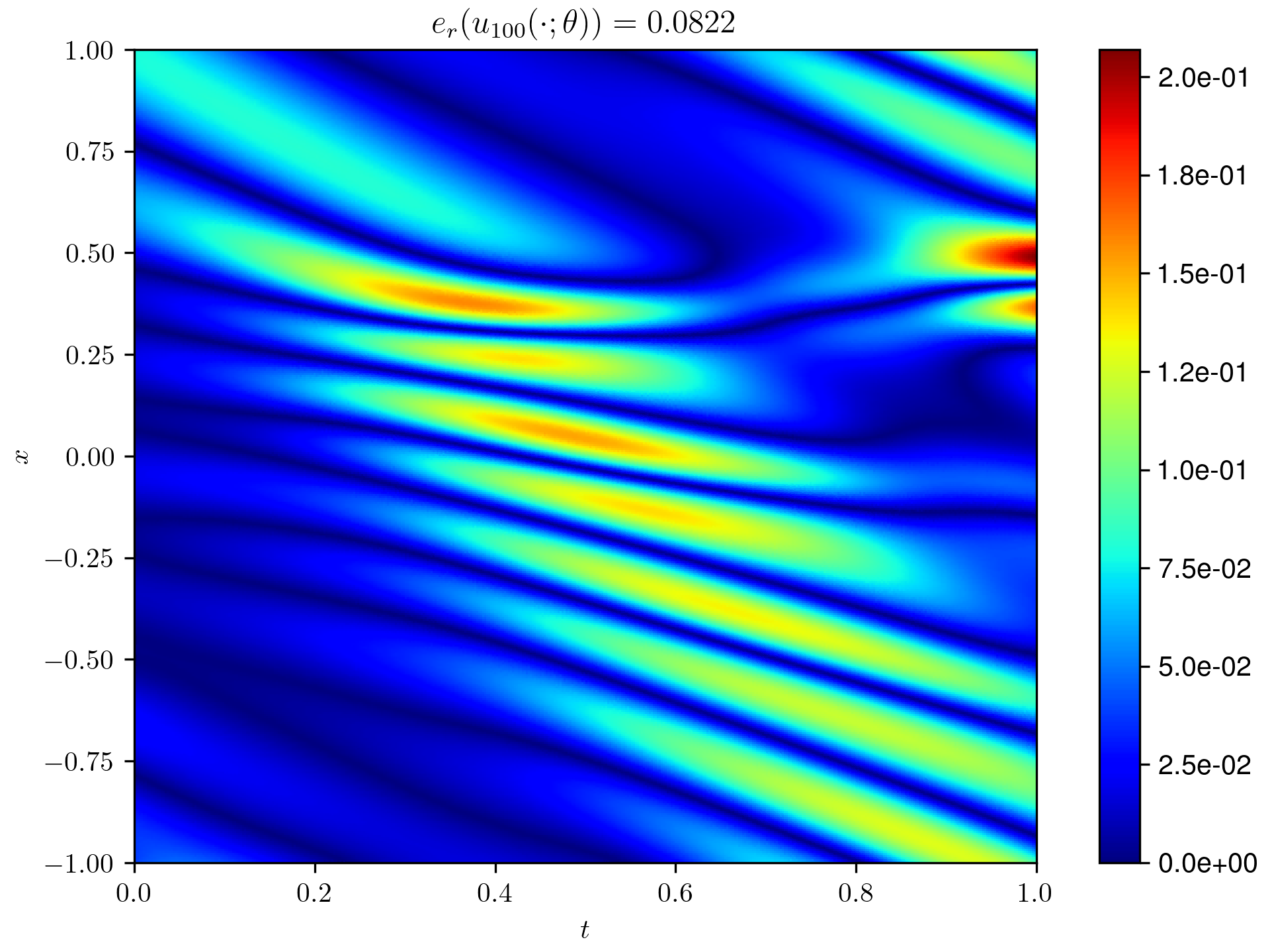}
        \caption{\textit{AAIS-g}}
    \end{subfigure}%
    \begin{subfigure}{.25\textwidth}
        \centering
        \includegraphics[height=0.75\textwidth,width=1.0\textwidth]{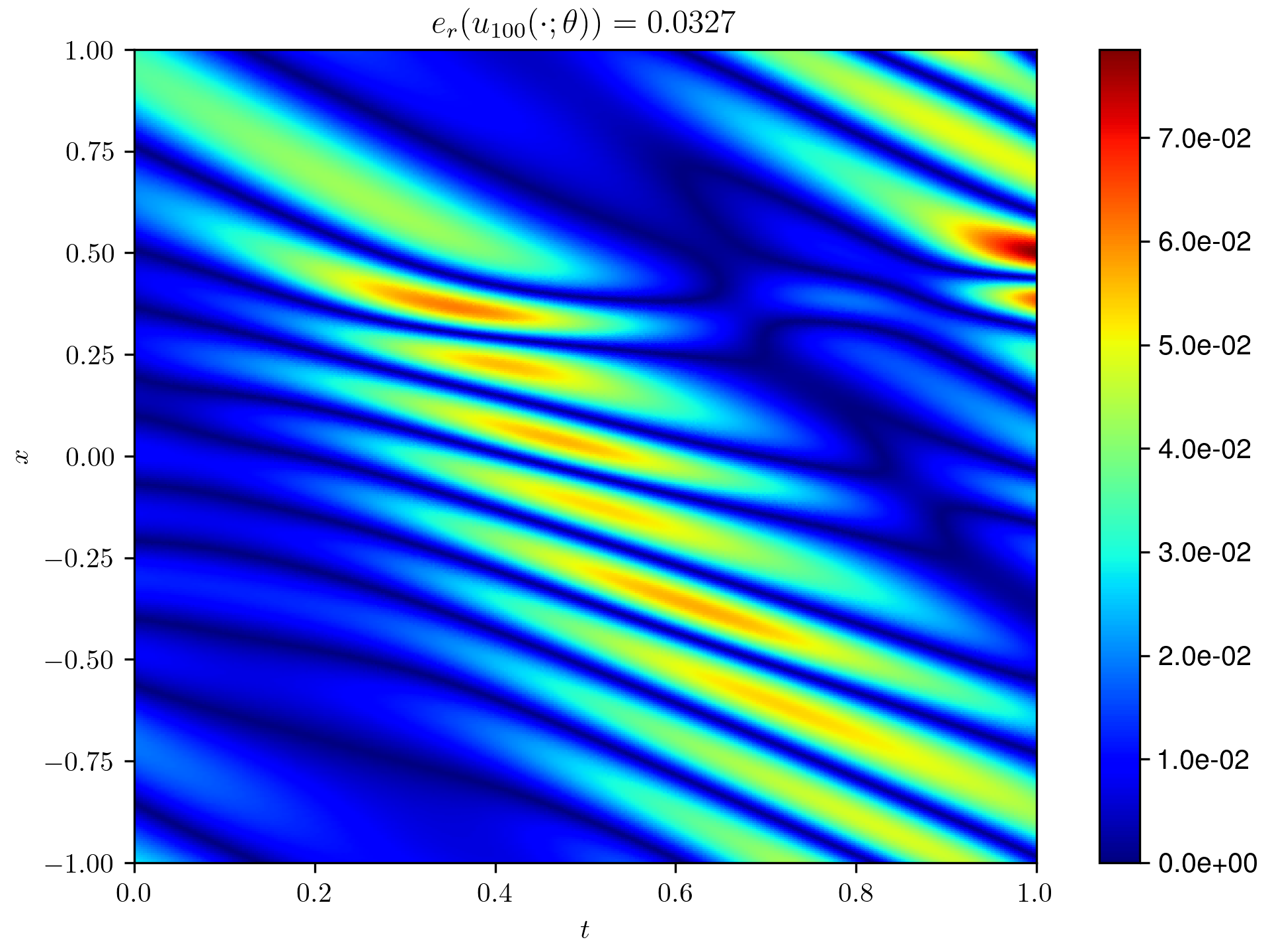}
        \caption{\textit{AAIS-t}}
    \end{subfigure}%
    \caption{Profiles of absolute error and neural network solutions for KdV equation after 100th training. }
    \label{fig:KdVErr5000e}
\end{figure}

Finally, we test the following Korteweg-de Vries(KdV) equation(reference solution is given in Figure \ref{fig:KdVexact}):
\begin{equation}
    \label{pde:KdV}
\begin{aligned}
&\partial_t u+u\partial_x u+00025\partial_{xxx}u =0,~~(t,x)\in(0,1)\times(-1,1),\\
&u(t,-1) = u(t,1)~~ u_t(t,-1)=u_t(t,1),~~t\in(0,1),\\
&u(0, x) = \cos(\pi x),~~x\in(0,1),
\end{aligned}
\end{equation}
with periodic boundary condition.
With the maximum iteration $M=100$, in the pre-train stage and each iteration there are 500 epochs for Adam of learning rate 0.0001 and 500 epochs for lbfgs of learning rate 1.0. $N_A$ for AAIS algorithms is set to be 4000. The loss and relative error after training at each iteration are listed in Figure \ref{fig:KdVLossErr500e}. Different like above experiments, the adaptive sampling methods would not help a lot during training, meaning that all methods give similar results. The residual and nodes are presented in Figure \ref{fig:KdVNode5000e}, we could see that although the \textit{Uni} method would generate a solution with smaller relative error, the residual would not be smoothed away and show some local behaviors. The profiles of solutions and absolute errors are given in Figure \ref{fig:KdVErr5000e}. Like Allen-Cahn equation, the higher absolute error  would not reflect in the residual may be the limitation of the adaptive sampling methods.

\section*{Data availability}
The code in this work is available from the GitHub repository \url{https://github.com/Zenki229/AAIS_PINN}.

\bibliographystyle{elsarticle-num} 
\bibliography{ref}

\end{document}